\documentclass{amsart}
\usepackage{amssymb,amsmath,amsfonts,amsthm,epsfig,amscd,hyperref}
\usepackage{stmaryrd,comment}
\usepackage[all,cmtip,poly]{xy}
\usepackage{svn, url}
\usepackage{mathrsfs}
\usepackage{xr}

\newtheorem{sublemma}{Sub-lemma}
\newtheorem{theorem}[subsection]{Theorem}
\newtheorem{thm}[subsubsection]{Theorem}
\newtheorem{lemma}[subsubsection]{Lemma}
\newtheorem{lem}[subsubsection]{Lemma}

\newtheorem{cor}[subsubsection]{Corollary}

\newtheorem{prop}[subsubsection]{Proposition}

\newtheorem{defn}[subsubsection]{Definition}

\theoremstyle{remark}
\newtheorem{remark}[subsubsection]{Remark}
\newtheorem{rem}[subsubsection]{Remark}

\newtheorem{assumption}[subsubsection]{Assumption}
\newtheorem{hyp}[subsubsection]{Hypothesis}

\numberwithin{equation}{subsection}

\newif\iffinalrun
\iffinalrun
\else
\fi

\iffinalrun
  \newcommand{\need}[1]{}
  \newcommand{\mar}[1]{}
\else
  \newcommand{\need}[1]{{\tiny *** #1}}
  \newcommand{\mar}[1]{\marginpar{\raggedright\tiny #1}}
\fi

\newcommand{\A}{\AA}
\newcommand{\C}{\CC}
\newcommand{\F}{\FF}

\newcommand{\Q}{\QQ}
\newcommand{\R}{\RR}
\newcommand{\Z}{\ZZ}

\renewcommand{\O}{\cO}

\newcommand{\m}{\frakm}
\newcommand{\p}{\frakp}

\newcommand{\sm}{\mathrm{sm}}

\renewcommand{\AA}{{\mathbb A}}

\newcommand{\CC}{{\mathbb C}}

\newcommand{\FF}{{\mathbb F}}

\newcommand{\LL}{{\mathbb L}}

\newcommand{\PP}{{\mathbb P}}
\newcommand{\QQ}{{\mathbb Q}}
\newcommand{\RR}{{\mathbb R}}

\newcommand{\TT}{{\mathbb T}}

\newcommand{\ZZ}{{\mathbb Z}}

\newcommand{\bA}{\ensuremath{\mathbf{A}}}

\newcommand{\bD}{\ensuremath{\mathbf{D}}}

\newcommand{\bL}{\ensuremath{\mathbf{L}}}

\newcommand{\bQ}{\ensuremath{\mathbf{Q}}}

\newcommand{\bZ}{\ensuremath{\mathbf{Z}}}

\renewcommand{\bf}{\ensuremath{\mathbf{f}}}

\newcommand{\bv}{\ensuremath{\mathbf{v}}}

\newcommand{\cA}{{\mathcal A}}

\newcommand{\cC}{{\mathcal C}}
\newcommand{\cD}{{\mathcal D}}

\newcommand{\cF}{{\mathcal F}}

\newcommand{\cH}{{\mathcal H}}
\newcommand{\cI}{{\mathcal I}}

\newcommand{\cK}{{\mathcal K}}
\newcommand{\cL}{{\mathcal L}}
\newcommand{\cM}{{\mathcal M}}

\newcommand{\cO}{{\mathcal O}}
\newcommand{\cP}{{\mathcal P}}
\newcommand{\cQ}{{\mathcal Q}}

\newcommand{\cS}{{\mathcal S}}

\newcommand{\cV}{{\mathcal V}}

\newcommand{\frakm}{\mathfrak{m}}
\newcommand{\frakn}{\mathfrak{n}}
\newcommand{\frakp}{\mathfrak{p}}

\newcommand{\frakX}{\mathfrak{X}}
\newcommand{\ainf}{\mathfrak{a}_\infty}

\newcommand{\Fbar}{\overline{\F}}
\newcommand{\Qbar}{\overline{\Q}}

\newcommand{\Fp}{\F_p}

\newcommand{\Fpbar}{\Fbar_p}

\newcommand{\Fpbarx}{\Fpbar^{\times}}

\newcommand{\Qp}{\Q_p}

\newcommand{\Qpbar}{\Qbar_p}

\newcommand{\Qpbarx}{\Qpbar^{\times}}

\DeclareMathOperator{\Aut}{Aut}

\DeclareMathOperator{\End}{End}

\DeclareMathOperator{\Gal}{Gal}
\DeclareMathOperator{\GL}{GL}

\DeclareMathOperator{\Hom}{Hom}

\DeclareMathOperator{\Ind}{Ind}

\DeclareMathOperator{\Mod}{Mod}

\DeclareMathOperator{\ord}{ord}
\DeclareMathOperator{\PGL}{PGL}
\DeclareMathOperator{\Pic}{Pic}

\DeclareMathOperator{\PSL}{PSL}

\DeclareMathOperator{\SL}{SL}

\DeclareMathOperator{\Spec}{Spec}

\DeclareMathOperator{\WD}{WD}
\DeclareMathOperator{\Supp}{Supp}
\DeclareMathOperator{\Tor}{Tor}

\newcommand{\ab}{\mathrm{ab}}
\newcommand{\cris}{\mathrm{cris}}
\newcommand{\dR}{\mathrm{dR}}

\newcommand{\Frob}{\mathrm{Frob}}

\newcommand{\HT}{\mathrm{HT}}

\newcommand{\st}{\mathrm{st}}

\newcommand{\ur}{\mathrm{ur}}
\newcommand{\BT}{\mathrm{BT}}
\newcommand{\ad}{\mathrm{ad}}
\newcommand{\loc}{\mathrm{loc}}
\newcommand{\Iw}{\mathrm{Iw}}
\newcommand{\rb}{\mathrm{b}}
\newcommand{\rs}{\mathrm{s}}
\newcommand{\rns}{\mathrm{ns}}

\newcommand{\rhobar}{\overline{\rho}}

\newcommand{\del}{\partial}

\newcommand{\llb}{\llbracket}
\newcommand{\rrb}{\rrbracket}

\newcommand{\onto}{\twoheadrightarrow}
\newcommand{\toisom}{\buildrel\sim\over\to}

\newcommand{\Art}{{\operatorname{Art}}}

\newcommand{\epsilonbar}{\overline{\epsilon}}

\newcommand{\rbar}{\overline{r}}

\newcommand{\CNL}{\operatorname{CNL}}
\newcommand{\Grm}{\mathrm{G}}
\newcommand{\Trm}{\mathrm{T}}
\newcommand{\Wrm}{\mathrm{W}}
\newcommand{\Brm}{\mathrm{B}}
\newcommand{\Mrm}{\mathrm{M}}

\newcommand{\Urm}{\mathrm{U}}
\newcommand{\Prm}{\mathrm{P}}

\newcommand{\Krm}{\mathrm{K}}
\newcommand{\tG}{\widetilde{G}}
\newcommand{\tB}{\widetilde{B}}
\newcommand{\tK}{\widetilde{K}}
\newcommand{\tDelta}{\widetilde{\Delta}}

\newcommand{\diag}{\mathrm{diag}}
\newcommand{\rec}{\mathrm{rec}}
\newcommand{\tX}{\widetilde{X}}
\newcommand{\PWP}{{^\Prm\Wrm^\Prm}}
\newcommand{\Sh}{\mathrm{Sh}}
\newcommand{\barx}{\overline{x}}
\newcommand{\tA}{\widetilde{A}}
\newcommand{\Gbar}{\overline{G}}
\newcommand{\avoid}{\mathrm{avoid}}
\renewcommand{\ss}{\mathrm{ss}}
\newcommand{\IQFassm}{Let $v \notin T$ be a finite place, with residue characteristic $\ell$. Then either 
	$T$ contains no $\ell$-adic places and $\ell$ is unramified in $F$, or there exists an imaginary quadratic subfield of $F$ in which $\ell$ splits.}

\begin{document}

\begin{abstract}

In this paper, we establish the modularity of every elliptic curve $E/F$, 
where $F$ runs over infinitely many imaginary quadratic fields, including 
$\Q(\sqrt{-d})$ for $d=1,2,3,5$.  More precisely, let $F$ be imaginary quadratic 
and assume that the modular curve
$X_0(15)$, which is an elliptic curve of rank $0$ over $\Q$, also has rank $0$ 
over $F$. Then we prove that all elliptic curves over $F$ 
are modular. More generally, when $F/\mathbb{Q}$ is an imaginary CM field
that does not contain a primitive $5$th root of unity, we prove the modularity 
of elliptic curves $E/F$ under a technical assumption on the 
image of the representation of $\Gal(\overline{F}/F)$ on $E[3]$ or $E[5]$. 

The key new technical ingredient we use is a local-global compatibility 
theorem for the $p$-adic Galois representations associated to torsion in the cohomology 
of the relevant locally symmetric spaces. 
We establish this result in 
the crystalline case, under some technical assumptions, 
but allowing arbitrary dimension, arbitrarily large regular Hodge--Tate weights, 
and allowing $p$ to be small and highly ramified in the imaginary CM field $F$. 

\end{abstract}

\title[Modularity of elliptic curves over imaginary quadratic fields]{On the modularity of elliptic curves over imaginary quadratic fields}
\author{Ana Caraiani} \email{caraiani.ana@gmail.com}
\address{Department of
  Mathematics, Imperial College London,
  London SW7 2AZ, UK}
  
\author{James Newton}\email{james.newton@maths.ox.ac.uk}
\address{Mathematical Institute,
University of Oxford,
Oxford OX2 6GG, UK}

\maketitle

\tableofcontents

\section{Introduction}

Let $F$ be a number field. We say that an elliptic curve 
$E/F$ is \emph{modular} if either $E$ has complex multiplication or if there 
exists a cuspidal automorphic representation $\pi$ of $\GL_2(\A_{F})$ 
of parallel weight $2$ whose associated $L$-function is the same as the 
$L$-function of $E$\footnote{The reason for the two cases is that if $E$ has CM by a field which embeds in $F$, then it cannot be associated to a cuspidal automorphic representation.}. 

In this paper, we establish 
the modularity of every elliptic curve $E/F$, 
where $F$ runs over infinitely many imaginary quadratic fields, including 
$\Q(\sqrt{-d})$ for $d=1,2,3,5$. 

Recall that 
the modular curve $X_0(15)$ is an elliptic curve of 
rank zero over $\mathbb{Q}$ 
-- it is the curve with Cremona label 15A1. 
We prove the following result. 

\begin{theorem}[Corollary~\ref{cor:modularity IQF improved}]\label{theorem:modularity IQ}
Let $F$ be an imaginary quadratic field such that the Mordell--Weil group 
$X_0(15)(F)$ is finite. Then every elliptic curve $E/F$ is modular.
\end{theorem}

\noindent We can compute the ranks of $X_0(15)$ over 
imaginary quadratic fields of small discriminant using Sage \cite{sagemath} or Magma \cite{magma} and check that 
the theorem applies to $F=\Q(\sqrt{-d})$ for the above values of $d$. By~\cite[Theorem 2.3]{mikic-najman}, 
the theorem applies to an infinite class of imaginary quadratic fields. 
Moreover, a celebrated conjecture of Goldfeld~\cite{goldfeld}, 
when coupled with the Birch--Swinnerton-Dyer 
conjecture, predicts that $X_0(15)$ should have rank $0$ over 
$50\%$ of quadratic fields, when these are ordered
by the absolute value of the discriminant. The conjecture predicts rank $0$ over slightly more than half of imaginary quadratic fields. 
More precisely, $X_0(15)$ is predicted to have rank $0$ over $100\%$ of those imaginary quadratic fields $\Q(\sqrt{-d})$ with $d$ positive square-free and $d$ mod $15 \in \{0, 1, 2, 3, 4, 5, 8,12\}$. This congruence condition corresponds to the global root number of $X_0(15)$ over $\Q(\sqrt{-d})$ being $+1$ (see for example \cite[Corollary 2]{Dokchitser-root-numbers}). Smith~\cite{smith-new} recently verified this prediction (since $X_0(15)$ has a rational cyclic degree $4$ isogeny, the previous results of~\cite{smith} exclude this case).  In fact, Smith shows that the $2^\infty$-Selmer corank is $0$ for 100\% of discriminants satisfying these congruence conditions, which implies that the rank is $0$ with no dependence on BSD.

The modularity of an elliptic curve $E$ over a number field $F$ implies that the 
associated $L$-function has analytic continuation to the entire complex plane. 
This is needed in order to formulate the Birch--Swinnerton-Dyer conjecture 
for $E$ unconditionally. Furthermore, modularity has historically played
a key role in progress 
on the BSD conjecture, going back to the use of Heegner points by Gross and Zagier
for (modular) elliptic curves over $\mathbb{Q}$. 
Recently, Loeffler and Zerbes made significant progress on the BSD 
conjecture for modular elliptic curves defined over imaginary quadratic fields~\cite{loeffler-zerbes}, 
making Theorem~\ref{theorem:modularity IQ} 
particularly timely. Other Diophantine applications of modularity, for example to the Fermat equation, require more information about mod $p$ Galois representations than comes from modularity (cf.~Serre's `$\epsilon$-conjecture' when $F=\Q$, and \cite{MR3811755, MR4053081} for conditional results over imaginary quadratic fields). This additional information appears to be beyond the reach of current techniques when $F$ is not totally real.

More generally, when $F/\mathbb{Q}$ is an imaginary CM field
that does not contain a primitive fifth root of unity, we prove the modularity 
of elliptic curves $E/F$ under a technical assumption on the 
image of the representation of $\Gal(\overline{F}/F)$ on $E[3]$ or $E[5]$.
As a consequence, we obtain the following result.  

\begin{theorem}[Corollary~\ref{cor:100 percent}]\label{theorem:modularity CM}
Let $F$ be an imaginary CM field that is Galois over $\mathbb{Q}$ and such that $\zeta_5\notin F$. 
Then 100\% of Weierstrass equations over $F$, ordered by their height, define a modular elliptic curve. 
\end{theorem}

The modularity of elliptic curves $E/\mathbb{Q}$ was pioneered
by Wiles and Taylor--Wiles in~\cite{wiles, MR1333036} and completed 
by Breuil--Conrad--Diamond--Taylor in~\cite{MR1839918}. The modularity 
of elliptic curves defined over \emph{real quadratic fields} was established, more recently, 
in~\cite{flhs}. Compared to the rational case, the real quadratic case 
relies on the improvements to the Taylor--Wiles method 
due to Kisin~\cite{kis04}, on supplementing the traditional 3-5 prime switch with an ingenious 3-7 switch, 
and on a sophisticated analysis of quadratic points on 
several modular curves of small level. Further results have been obtained
for more general \emph{totally real fields}, including cubic and quartic fields \cite{dns-cubic,box-quartic}. As another example, Thorne~\cite{thorne-jems} has proved
the modularity of every elliptic curve defined over the cyclotomic $\Z_p$-extension of $\Q$ for any prime $p$. 

The modularity of elliptic curves defined over \emph{imaginary CM fields}
has historically been more difficult to establish. This is because the systems
of Hecke eigenvalues that conjecturally match such elliptic 
curves contribute to the cohomology of locally symmetric spaces such as
Bianchi $3$-manifolds, which are not directly related to Shimura varieties. The situation has been extensively investigated numerically (for example, \cite{cremona-84,cremona-92,lmfdb}), and modularity of specific elliptic curves can be verified using the Faltings--Serre method \cite{gdp-faltings-serre}. Inspired by a programme outlined by Calegari--Geraghty in~\cite{CG},  
the potential modularity of such elliptic curves was 
established independently in~\cite{10author} and in~\cite{BCGP}. Since
then, Allen--Khare--Thorne proved many instances of actual 
modularity in~\cite{AKT}. More precisely, they established the modularity 
of a positive proportion of elliptic curves over imaginary CM 
fields together with strong residual modularity results modulo $3$ and modulo $5$. 

\begin{remark}
In fact,~\cite{BCGP} establish the potential modularity of 
elliptic curves defined over a general quadratic extension of a totally real field. A recent preprint~\cite{whitmore}
by Whitmore builds on their method and on the results of~\cite{AKT} 
to prove actual modularity for a positive proportion of such elliptic curves.
\end{remark}

To prove Theorem~\ref{theorem:modularity CM}, we combine the residual modularity results of~\cite{AKT}
with a modularity lifting theorem in the Barsotti--Tate case in the style of Kisin \cite{kis04}. The crucial
ingredient needed to prove our Barsotti--Tate modularity lifting theorem is 
a \emph{local-global compatibility} result 
for the Galois representations constructed by Scholze in~\cite{scholze-torsion}. This 
is a result of independent interest, which we now discuss. 

Let $K\subset \GL_n(\A_{F, f})$ be a neat compact open subgroup
and let $X_K$ be the corresponding locally symmetric space for $\GL_n/F$. A 
highest weight vector $\lambda$ for $\mathrm{Res}_{F/\Q}\GL_n$ 
determines a $\Z_p$-local system $\cV_{\lambda}$ on $X_K$ and
we are interested in understanding the systems of Hecke eigenvalues occurring 
in $H^*(X_K,\cV_{\lambda})$. 
Let $\mathbb{T}$ be the usual abstract spherical Hecke algebra acting on $H^*(X_K, \cV_{\lambda})$
by correspondences, let $\mathbb{T}(K,\lambda)$ be the maximal quotient of $\mathbb{T}$
through which this action is faithful, and let $\m \subset \mathbb{T}(K,\lambda)$ be a maximal ideal. 
When $\m$ is non-Eisenstein, Scholze constructed a continuous Galois representation
\[
\rho_{\m}: \Gal(\overline{F}/F)\to \GL_n(\mathbb{T}(K,\lambda)_{\m}/I),
\]
where $I\subset \mathbb{T}(K,\lambda)_{\m}$ is an ideal satisfying $I^4=0$, cf.~\cite{new-tho}. 
This $p$-adically interpolates the Galois representations attached to torsion
classes occurring in $H^*(X_K,\cV_{\lambda})_{\m}$ 
as well as those attached to characteristic $0$ automorphic forms, first constructed by~\cite{hltt}. 
We let $\bar{\rho}_{\m}$ denote the absolutely irreducible residual representation 
obtained by reducing $\rho_{\m}$ modulo $\m$. 

For applications to modularity,
it is extremely important to understand the properties of $\rho_{\m}$, cf.~\cite[Conjecture B]{CG}. 
One needs to know whether $\rho_{\m}$ satisfies some form of local-global 
compatibility: if $v\mid \ell$ is a prime of $F$ and $G_{F_v}:=\Gal(\overline{F}_v/F_v)$, how does the level $K_v$
at which $\m$ occurs (together with the weights $\lambda_v$ if $\ell = p$) determine 
the ramification of $\rho_{\m}|_{G_{F_v}}$? The case when $\ell = p$ is particularly subtle
because it is not (a priori) clear how to formulate integral $p$-adic Hodge theory 
conditions which should be satisfied by the Galois representations $\rho_{\m}$, %
and because the $\rho_{\m}$ are constructed in~\cite{scholze-torsion} 
via a $p$-adic interpolation
argument that loses track of the weight $\lambda$ and of the level $K_v$ for $v\mid p$. 

In~\cite{10author}, we established such a local-global 
compatibility result at $\ell = p$ 
in two restricted families of cases described by natural integral conditions: the ordinary case and certain Fontaine--Laffaille cases. In the 
present paper, we go much further than this and establish the desired result in 
the crystalline case, under some technical assumptions, 
but allowing arbitrary $n$, arbitrary weight $\lambda$, 
and allowing $p$ to be small and highly ramified in $F$. 
In this generality, the formulation via integral $p$-adic Hodge theory is still mysterious, but the local-global compatibility conjecture 
can be formulated as in~\cite[Conjecture 5.1.12]{gee-newton}, using the crystalline deformation
rings first constructed by Kisin~\cite{kisindefrings}. More precisely, 
we have a composition 
\begin{equation}\label{eq:galois to hecke map}
\xymatrix{R^{\square}_{\bar{\rho}_{\m}|_{G_{F_v}}}\ar[r]\ar@{->>}[dr] & R^{\square}_{\bar{\rho}_{\m}}\ar[r] & \mathbb{T}(K,\lambda)_{\m}/I 
\\ \  & R^{\square, \mathrm{crys}}_{\bar{\rho}_{\m}|_{G_{F_v}}}(\lambda_v)\ar@{-->}[ur]  &\ }
\end{equation}
where the first horizontal map is the usual map from the local deformation
ring of $\bar{\rho}_{\m}|_{G_{F_v}}$ to the global deformation ring 
of $\bar{\rho}_{\m}$ and the second horizontal map is induced by the existence of $\rho_{\m}$. 
When $K_v = \GL_n(\cO_{F_v})$ is a maximal compact subgroup, 
the natural conjecture is that the composition~\eqref{eq:galois to hecke map} 
factors through $R^{\square, \mathrm{crys}}_{\bar{\rho}_{\m}|_{G_{F_v}}}(\lambda_v)$, 
the crystalline deformation ring with Hodge--Tate weights determined by $\lambda_{v}$. We prove  
this conjecture in Theorem~\ref{thm:RtoT factors thru Kisin def ring} 
under some technical assumptions --- roughly, the statement is as follows. 

\begin{theorem}\label{theorem:LGC torsion} 
Let $F$ be an imaginary CM field that contains an imaginary quadratic field
$F_0$ and with maximal totally real subfield $F^+$\footnote{The field $F$ has to satisfy some
additional technical assumptions so that we can appeal to the unconditional
base change results of~\cite{shin-basechange}.}. 
Let $p$ be a rational prime that splits in $F_0$, let
$\bar{v}\mid p$ be a prime of $F^+$, and assume the following. 
\begin{enumerate}
\item Setting $\bar{v} = v\cdot v^{c}$, we have $K_v = \GL_n(\cO_{F_v})$ 
and $K_{v^c} = \GL_n(\cO_{F_{v^c}})$. 
\item There exists a prime $\bar{v}'\mid p$ of $F^+$ distinct from $\bar{v}$ such that
\[
\sum [F^+_{\bar{v}''}:\Q_p] \geq \frac{1}{2}[F^+:\Q],
\]
where the sum runs over primes $\bar{v}''\mid p$ of $F^+$ 
distinct from both $\bar{v}$ and $\bar{v}'$. 
\item $\m$ is a non-Eisenstein maximal ideal such that $\bar{\rho}_{\m}$ 
	is decomposed generic, cf. Definition~\ref{defn:generic}. 
\end{enumerate}

\noindent Then, up to possibly enlarging the nilpotent ideal $I$, 
the composition~\eqref{eq:galois to hecke map} factors through
$R^{\square, \mathrm{crys}}_{\bar{\rho}_{\m}|_{G_{F_v}}}(\lambda_v)$ as expected. 
\end{theorem}

\begin{remark}\leavevmode
\begin{enumerate}

\item In this paper we restrict attention to the crystalline case, hence the first assumption.  Building on the approach taken here, Bence Hevesi generalized Theorem~\ref{theorem:LGC torsion} to the potentially semi-stable case as part of his PhD thesis \cite{bence-thesis}. 

\item The second assumption is more serious than the first and excludes in particular 
the case where $F = F_0$. The third assumption is 
needed in order to appeal 
to the results of~\cite{caraiani-scholze-noncompact}
on unitary Shimura varieties, 
or alternatively to those of~\cite{koshikawa}. 

\item We also obtain in Theorem~\ref{thm:LGC in char 0} a local-global compatibility result 
for the characteristic $0$ Galois representations attached to 
regular algebraic cuspidal automorphic representations of 
$\GL_n(\mathbb{A}_F)$. 
In this setting, the local-global compatibility at $\ell\not = p$ is already 
known up to semi-simplification by work of Varma~\cite{ilavarma}.
More recently, A'Campo~\cite{ACampo} proved that these 
automorphic Galois representations are de Rham at all primes above $p$, and determined the Hodge--Tate weights of these representations, using Wang-Erickson's work on $p$-adic Hodge theoretic conditions for pseudorepresentations \cite{CWE-sst-pseudo}. These results still assume that the residual representation is absolutely irreducible and decomposed generic, but later work of A'Campo, Hevesi, Thorne and Whitmore \cite{ACHTW} proves local-global compatibility in characteristic 0 without these assumptions. They combine a version of the `$P$-ordinary' degree shifting introduced in this article with (amongst other innovations) a new idea involving control of the cohomology of unitary Shimura varieties outside middle degree. 

\item Motivated by our applications to elliptic curves, we prove a slightly more general result which includes semistable ordinary representations.

\end{enumerate}
\end{remark}

There are two key new ideas that allow us to prove much stronger
local-global compatibility results than in~\cite{10author}. The first idea 
is to work with $P$-ordinary parts at the prime $\bar{v}\mid p$ of $F^+$
where we want to prove local-global compatibility. The second idea,
which was suggested to us by Peter Scholze,  
is to increase the level at the auxiliary primes $\bar{v}''\mid p$ of $F^+$ 
in order to simplify the analysis of the boundary of the Borel--Serre compactification
in the relevant unitary Shimura varieties. Fortunately, these can be 
implemented simultaneously. 

To explain how the first new idea is useful, recall that the crystalline
deformation rings $R^{\square, \mathrm{crys}}_{\bar{\rho}_{\m}|_{G_{F_v}}}(\lambda_v)$
were defined by Kisin first after inverting $p$, and then integrally by taking 
Zariski closure from the generic fibre. On the other hand, the Galois 
representations $\rho_{\m}$ could be torsion. They are constructed 
by congruences using a subtle argument that involves $2n$-dimensional Galois 
representations. If we had a characteristic $0$ lift of $\rho_{\m}|_{G_{F_v}}$, which 
we knew was crystalline at $v$ with Hodge--Tate weights determined by $\lambda_v$, 
we would deduce that the diagram~\eqref{eq:galois to hecke map}
factors as desired. Conversely, if the diagram factored as desired, we would expect
the crystalline lift to exist by results of Tong Liu~\cite{liu-fourier}. 
It seems hard to guarantee that there is a characteristic $0$ 
crystalline lift of the global representation $\rho_{\m}$. However, by working with $P$-ordinary parts 
at $\bar{v}$ throughout, we construct for each $m\in \Z_{\geq 1}$
a $2n$-dimensional characteristic $0$
global representation $\rho_{\widetilde{\m}}$ such that 
\[
\rho_{\widetilde{\m}}|_{G_{F_v}} \simeq \left(\begin{smallmatrix} * & * 
\\ 0 & *\end{smallmatrix}\right) 
\]
with $n$-dimensional diagonal blocks and such that one of these blocks 
is congruent to the local representation $\rho_{\m}|_{G_{F_v}}\pmod{p^m}$. Moreover, we can ensure that each of these characteristic $0$ lifts is crystalline with the correct Hodge--Tate weights. We expect the global representation $\rho_{\widetilde{\m}}$ to be irreducible, and we do not produce characteristic $0$ lifts of the global representation $\rho_{\m} \pmod{p^m}$. %

The second new idea is useful for the ``degree-shifting'' argument 
needed to relate the cohomology groups $H^*(X_K, \cV_{\lambda})_{\m}$ 
to a middle degree boundary cohomology group 
$H^d(\partial \widetilde{X}_{\widetilde{K}},\cV_{\tilde{\lambda}})_{\widetilde{\m}}$
of some unitary Shimura variety $\widetilde{X}_{\widetilde{K}}$. 
We can control the latter using the main theorem of~\cite{caraiani-scholze-noncompact}. 
However, one only has a spectral sequence of Leray--Serre type from the former cohomology groups
to the latter -- controlling the behaviour of this spectral sequence seems to be a tricky 
problem in modular representation theory. In~\cite{10author}, 
we showed that the spectral sequence degenerates if $p$ is strictly greater than $n^2$ and 
is unramified in $F$. 
In the present paper, we increase the level at auxiliary primes 
$\bar{v}'' \mid p$ and, through a delicate induction argument, we keep track
of the terms in the spectral sequence modulo powers of $p$ without 
imposing the additional assumptions that $p > n^2$ and is unramified in $F$. 

To prove a modularity lifting theorem in the Barsotti--Tate case and 
deduce Theorem~\ref{theorem:modularity CM}, 
we apply Theorem~\ref{theorem:LGC torsion} in the case
when $n=2$ and $\lambda$ is trivial. For our applications, it is 
crucial to allow $p$
to be small and highly ramified in $F^+$. 
(We can then ensure
that the second condition of Theorem~\ref{theorem:LGC torsion}
is satisfied using an appropriate solvable base change.) 
This is why the local-global compatibility results of~\cite{10author} in
the Fontaine--Laffaille case were not strong enough and why~\cite{AKT} 
appealed instead to the results in the more restrictive ordinary case. 

We expect Theorem~\ref{theorem:LGC torsion} to have many more
applications to modularity over CM fields in the near future. Some examples have appeared since the first version of this article was posted: work of one of us (J.N.) with Boxer, Calegari, Gee and Thorne \cite{bcgnt} proving the Ramanujan and Sato--Tate conjectures for regular algebraic cuspidal automorphic representations of $\GL_2(\AA_F)$ with $F$ an imaginary quadratic field, and the work of Matsumoto \cite{matsumoto} proving the same results for arbitrary CM fields $F$. The case of weight $0$ was proved in \cite{10author}; to carry out a version of Harris's tensor product trick handling higher weights, it is crucial to use an improved modularity lifting theorem allowing ramification above $p$.

To deduce Theorem~\ref{theorem:modularity IQ}, we analyse
the imaginary quadratic points on 
several modular curves with small level at $3$ and $5$, classifying
elliptic curves for which both the $3$-torsion
and the $5$-torsion are exceptional. For a prime $p$, we let $\rb p \subset \GL_2(\F_p)$ denote
the upper-triangular Borel subgroup, $\rs p \subset \GL_2(\F_p)$ denote the normalizer of the standard split Cartan subgroup
and $\rns p \subset \GL_2(\F_p)$ denote the normalizer of the standard non-split Cartan subgroup $\rns p^\circ$.  
After some reductions using group theory,  
there turn out to 
be six modular curves of interest: 
\begin{enumerate}
\item $X(\rb3,\rb5)$ (also denoted by $X_0(15)$ above);
\item $X(\rs3, \rb5)$;
\item $X(\rns3^\circ, \rb5)$;
\item $X(\rb3, \rns5)$;
\item $X(\rs3, \rns5)$;
\item $X(\rns3^\circ, \rns5)$.
\end{enumerate}

\noindent The modular curves $X(\rb3,\rb5)$ and
$X(\rs3,\rb5)$ are isogenous elliptic curves of Mordell--Weil 
rank $0$ over $\Q$. They are the obstruction to extending 
Theorem~\ref{theorem:modularity IQ} to every imaginary quadratic 
field $F$, although we can at least understand how their torsion subgroup  
grows in imaginary quadratic extensions. 

The modular curve 
$X(\rns3^\circ, \rb5)$ is a genus $1$ curve without a rational point. 
This case does not occur in the real quadratic case because $\rns3^\circ$ 
does not contain an odd element, which should represent complex conjugation. 
The curve contains two infinite families of imaginary quadratic points, 
for which, miraculously, it is still possible to prove modularity! 
The elliptic curves in the first family turn out to all have 
rational $j$-invariant. The elliptic curves in the second family 
turn out to all be $\Q$-curves (isogenous to their conjugates over $\overline{\Q}$).

The remaining cases also do not occur in the real quadratic setting. The 
modular curve $X(\rb3, \rns5)$ is a genus $2$ hyperelliptic curve and
we study its imaginary quadratic points using similar methods to those of~\cite{flhs}.   
The modular curves $X(\rs3, \rns5)$ and $X(\rns3^\circ, \rns5)$ are bi-elliptic curves of genus $3$ whose
Jacobians have Mordell--Weil rank $1$. We analyse the imaginary quadratic points 
on these curves using the relative symmetric power Chabauty method developed by Siksek~\cite{siksek-chabauty} and Box~\cite{box-quadratic}. 

\begin{remark} 
It seems much more subtle to implement 
the 3-7 switch over an imaginary CM field than over a totally real field, 
as in~\cite[\S 3.2]{flhs}. The modular curve with full level structure at $7$ 
is isomorphic to the Klein quartic curve 
\[
x^3y + y^3z + z^3x = 0. 
\]
To implement the 3-7 switch, one
needs to produce points on a quadratic twist of the Klein quartic 
that are defined over solvable CM extensions of the original 
CM field. In the totally real case, this can be done with a clever 
application of Hilbert irreducibility, obtaining rational points over a 
degree $4$, thus solvable, totally real extension. This argument 
does not apply in the imaginary CM case: the direct argument gives points 
defined over a degree $4$ extension of the original field $F$, 
but this is not necessarily a CM field. 
By working with Weil restrictions of scalars to the maximal totally real subfield $F^+$,
the degree increases. One can obtain 
points defined over a CM extension of $F$ but it seems hard to guarantee 
that this extension is always solvable. 
\end{remark}

The organization of the paper is as follows. In Section 2, we collect 
preliminaries on locally symmetric spaces and develop $P$-ordinary
Hida theory in the setting of their Betti cohomology. In Section 3, we 
study the $P$-ordinary condition on the Galois side and record a
key argument with determinants that will be used for local-global compatibility. 
In Section 4, we prove Theorem~\ref{theorem:LGC torsion} and its characteristic
$0$ counterpart.  
In Section 5, we use this result together with the techniques developed in~\cite{10author} 
to prove a modularity lifting theorem over imaginary 
CM fields in the potentially Barsotti--Tate case. 
In Section 6, we combine this modularity lifting theorem with the results of~\cite{AKT}
to prove Theorem~\ref{theorem:modularity CM}. In Section 7, we analyse the 
imaginary quadratic points on several modular curves of small level and prove Theorem~\ref{theorem:modularity IQ}. 

\subsection*{Acknowledgements} Many thanks to the anonymous referee for a careful reading and helpful comments. We are very grateful to Peter Scholze 
for suggesting to us that increasing the level at auxiliary primes may help with degree-shifting 
and for sketching the proof of Lemma~\ref{lem:splitting cohomology}.
We thank Lambert A'Campo, Frank Calegari, John Cremona, Jessica Fintzen, Bence Hevesi, 
Samir Siksek, Matteo Tamiozzo and Jack Thorne for useful conversations. 
We thank Kiran Kedlaya and Steven Sivek for help running computations
in Magma and Sage. We thank Toby Gee, Florian Herzig, Alexander Smith and Matteo Tamiozzo 
for comments on an earlier version of this manuscript.

This project has received funding from the European Research Council (ERC) under the 
European Union’s Horizon 2020 research and innovation programme (grant agreement No. 804176). 
A.C. was supported in part by a Royal Society University Research Fellowship and by a Leverhulme 
Prize. J.N. was supported by a UKRI Future Leaders Fellowship, grant
MR/V021931/1. For the purpose of Open Access, the authors have applied a CC BY public copyright licence to any Author Accepted Manuscript version arising from this submission.

\subsection*{AI declaration} When revising this article, the authors prompted AI tools to check correctness. Some edits were made by the authors in response, with all content written by the authors. 

\subsection*{Notation} Our notation largely matches the one introduced in~\cite[\S 1.2]{10author}.  
If $F$ is a perfect field, we let $\overline{F}$ denote an algebraic 
closure of $F$ and $G_F$ denote the absolute Galois group $\Gal(\overline{F}/F)$. 

If $F$ is a number field, we let $S_p(F)$ be the set of 
places of $F$ above $p$. If $S$ is a finite set of finite places of a
number field $F$, we let 
$G_{F, S}$ denote the Galois group of the maximal extension
of $F$ that is unramified outside $S$. For a prime $\ell$, we let 
$\epsilon_{\ell}$ denote the $\ell$-adic cyclotomic character 
and $\bar{\epsilon}_{\ell}$ denote its reduction modulo $\ell$.

If $\pi$ is an irreducible admissible representation of $\GL_n(\A_F)$
and $\lambda\in (\Z^n)^{\Hom(F, \C)}$ is dominant for the standard 
upper triangular Borel subgroup, we say that $\pi$ is regular algebraic 
of weight $\lambda$ if the infinitesimal character of $\pi_\infty$ is the 
same as that of $V_{\lambda}^\vee$, where $V_{\lambda}$ is the 
algebraic representation of $\mathrm{Res}_{F/\Q}\GL_n$ of highest 
weight $\lambda$. See \S\ref{sec: unitary group} for a discussion 
of highest weight representations. 

If $K$ is a finite extension of $\Q_p$ for some prime $p$, 
we write $I_K$ for the inertia subgroup of $G_K$, $\Frob_K\in G_K/I_K$ 
for the geometric Frobenius and $W_K$ for the Weil group. 
We write $\Art_K: K^\times\toisom W_K^{\ab}$ for the Artin map of
local class field theory, normalized to take uniformizers to 
geometric Frobenius elements. We let $\mathrm{rec}_K$ 
denote the local Langlands correspondence of~\cite{ht}, 
which sends an irreducible smooth (admissible) representation $\pi$ 
of $\GL_n(K)$ over $\C$ to a Frobenius semi-simple Weil--Deligne representation
$\mathrm{rec}_K(\pi)$ of $W_K$, also over $\C$. We
also write $\mathrm{rec}_K^T$ for the arithmetic normalization
of the local Langlands correspondence, as defined for example
in~\cite[\S 2.1]{Clo14}; this normalization is defined for 
coefficients in any field which is abstractly isomorphic to 
$\C$, such as $\overline{\Q}_{\ell}$. We define labelled Hodge--Tate weights of $p$-adic representations of $G_K$ as in \cite[\S1.2]{10author}. In particular, $\epsilon_p$ has Hodge--Tate weight $-1$. 

If $G$ is a locally profinite group and $K$ is an open 
subgroup, we write $\cH(G,K)$ for the $\Z$-algebra of 
compactly supported bi-$K$-invariant functions $f:G\to \Z$,
cf.~\cite[Lemma 2.3]{new-tho}.

We let $E/\Q_p$ be a $p$-adic field
which will be our coefficient field, with ring of integers $\cO$, uniformizer
$\varpi$ and finite residue field $k:= \cO/\varpi$. We let $\CNL_{\cO}$
denote the category of complete, local, Noetherian $\cO$-algebras with 
residue field $k$.  

\section{The cohomology of locally symmetric spaces}

\subsection{Preliminaries}
In this section, we gather some preliminaries, and we
largely follow~\cite[\S 2]{10author} without giving complete details. 

\subsubsection{Locally symmetric spaces}\label{sec:locally symmetric spaces} 
Let $F$ be a number field
and $\mathrm{G}$ be a connected linear algebraic group over $F$,
with a model over $\cO_F$ that we will still denote by $\mathrm{G}$. 
We will denote by $X^{\mathrm{G}}$ the \emph{symmetric space} 
for $\mathrm{Res}_{F/\Q}\mathrm{G}$, which is 
a homogeneous space for $\mathrm{G}(F\otimes_{\Q}\R)$
as in~\cite[\S 2]{borel-serre} and \cite[Definition 3.1]{new-tho} (and which is determined by $\mathrm{G}$
up to isomorphism of homogeneous spaces).  

Let $K_{\mathrm{G}}\subset \mathrm{G}(\A_{F,f})$ be a \emph{good} compact
open subgroup in the sense of~\cite[\S 2.1]{10author}: namely 
it is \emph{neat} and of the form $\prod_{v}K_{\mathrm{G},v}$, where 
$v$ runs over the finite places of $F$. 
We consider the double quotient
\[
X^{\mathrm{G}}_{K_{\mathrm{G}}} 
:= \mathrm{G}(F) \backslash X^{\mathrm{G}}\times 
\mathrm{G}(\A_{F,f}) / K_{\mathrm{G}},
\]
which is a smooth, orientable Riemannian manifold. We 
also consider the partial Borel--Serre compactification 
$\overline{X}^{\mathrm{G}}$ of $X^{\mathrm{G}}$ as 
in~\cite[\S 7.1]{borel-serre} and form the double quotient 
\[
\overline{X}^{\mathrm{G}}_{K_{\mathrm{G}}} 
:= \mathrm{G}(F) \backslash \overline{X}^{\mathrm{G}}\times 
\mathrm{G}(\A_{F,f}) / K_{\mathrm{G}},
\]
which is a compact, smooth manifold with corners
with interior $X^{\mathrm{G}}_{K_{\mathrm{G}}}$. We note that the spaces $X^{\Grm}$ are always connected; when $\Grm(\R)$ is not connected, it is sometimes better to work with $X^{\Grm} \times \pi_0(\Grm(\R))$ (equivalently, replacing the isotropy subgroup in the definition of the symmetric space with its identity connected component). Since $\Grm(\R)$ will be connected in all the cases of interest to us, this will not concern us. Finally, we consider the boundaries $\partial X^{\mathrm{G}}: = 
\overline{X}^{\mathrm{G}}\setminus X^{\mathrm{G}}$ 
and $\partial X^{\mathrm{G}}_{K_{\mathrm{G}}}:= 
\overline{X}^{\mathrm{G}}_{K_{\mathrm{G}}} 
\setminus X^{\mathrm{G}}_{K_{\mathrm{G}}} $. 

We define $\mathfrak{X}_{\mathrm{G}} :=  
\varprojlim_{K_{\mathrm{G}}} X^{\mathrm{G}}_{K_{\mathrm{G}}}$,
endowed with the projective limit topology,
where $K_{\mathrm{G}}\subset \mathrm{G}(\A_{F,f})$
runs over good compact open subgroups.
We also consider the analogous spaces $\overline{\mathfrak{X}}_{\mathrm{G}}$ 
and $\partial\mathfrak{X}_{\mathrm{G}}$. All these spaces are
equipped with a continuous action of $\mathrm{G}(\A_{F,f})$, which
is equipped with the locally profinite topology. 
Note also that the spaces $\overline{\mathfrak{X}}_{\mathrm{G}}$ 
and $\partial\mathfrak{X}_{\mathrm{G}}$ are compact Hausdorff, 
being projective limits of compact Hausdorff spaces. We denote 
by $j:\mathfrak{X}_{\Grm}\hookrightarrow \overline{\mathfrak{X}}_{\Grm}$ the 
natural open immersion. 
As a consequence 
of~\cite[Lemma 6.2.1]{arizona} and~\cite[Lemma 2.31]{new-tho}, 
we see that the actions of any good subgroup $K_{\mathrm{G}}$ on
$\overline{\mathfrak{X}}_{\mathrm{G}}$ and $\partial\mathfrak{X}_{\mathrm{G}}$
are free in the sense of~\cite[Definition 2.23]{new-tho}. These limits have been considered previously by Rohlfs \cite{rohlfs}. It follows from the properness of the action of arithmetic groups on the symmetric space and its compactification (cf.~\cite[Proposition 1.9]{rohlfs}) that we have
\begin{align*}\mathfrak{X}_{\mathrm{G}} &= \mathrm{G}(F) \backslash {X}^{\mathrm{G}}\times 
\mathrm{G}(\A_{F,f}), ~\overline{\mathfrak{X}}_{\mathrm{G}} = \mathrm{G}(F) \backslash \overline{X}^{\mathrm{G}}\times 
\mathrm{G}(\A_{F,f}), \\ \text{and }\partial\mathfrak{X}_{\mathrm{G}} &= \mathrm{G}(F) \backslash \partial X^{\mathrm{G}}\times 
\mathrm{G}(\A_{F,f}),\end{align*} with topologies induced by the locally profinite topology on the adelic groups.  We prefer to work with these topological spaces, since they seem more natural than those used in \cite{new-tho,10author} which equip the adelic groups with the discrete topology. We compare these set-ups (`topological' and `discrete') in the next subsection.

\subsubsection{Hecke operators and coefficient systems}
\label{sec:Hecke formalism}

If $S$ is a finite set of finite places of $F$ 
we set $\Grm^S := \mathrm{G}(\A^{S}_{F,f})$ and $\Grm_S := 
\mathrm{G}(\A_{F, S})$, and similarly $K_\Grm^S = \prod_{v \not\in S}
 K_{\Grm, v}$ and $K_{\Grm, S} = \prod_{v \in S} K_{\Grm, v}$.
 We write $\cH(\Grm^S, K_\Grm^S)$ for the global Hecke algebra 
 over $\Z$ which is the restricted tensor product of the
 local Hecke algebras $\cH(\Grm(F_v), K_{\Grm,v})$ 
 for $v$ a finite place of $F$ not contained in $S$. 

Let $R$ be a commutative ring and let $\cV$ be a smooth 
$R[K_{\Grm,S}]$-module, which is finite free as an $R$-module. 
We now explain how to obtain from it a local system $\cV$ 
of $R$-modules on $X^{\Grm}_{K_{\Grm}}$ and how to equip
the usual and compactly supported cohomology groups 
$R\Gamma_{(c)}(X^{\Grm}_{K_{\Grm}},\cV)$ 
with an action of the Hecke algebra $\cH(\Grm^S, K_\Grm^S)\otimes_{\Z} R$,
by adapting the formalism of~\cite{new-tho} to our topological setting. 

Firstly, note that the $R[K_{\Grm, S}]$-module $\cV$ 
defines a $\Grm^S\times K_{\Grm,S}$-equivariant local system, 
which we denote by $\cV$ as well, 
on both $\mathfrak{X}_{\Grm}$ and 
$\overline{\mathfrak{X}}_{\Grm}$. Indeed, we first inflate $\cV$ to a smooth 
$R[\Grm^S\times K_{\Grm,S}]$-module, which is equivalent to a 
$\Grm^S\times K_{\Grm,S}$-equivariant sheaf on a point by~\cite[Lemma 2.26]{new-tho}, and then 
we pull back this sheaf to $\mathfrak{X}_{\Grm}$ and $\overline{\mathfrak{X}}_{\Grm}$,
respectively.  
From now on, we consider the $\Grm^S\times K_{\Grm,S}$-equivariant sheaves 
$\cV$ and $j_{!}\cV$ on $\overline{\mathfrak{X}}_{\Grm}$. 
By~\cite[\S 1, Lemma 1]{schneider}\footnote{It is 
assumed in \emph{loc. cit.} that the coefficients have characteristic $0$, but this is not
used in the proof.}, the category of $\Grm^S\times K_{\Grm,S}$-equivariant sheaves 
on $\overline{\mathfrak{X}}_{\Grm}$ has enough injectives. By~\cite[Lemma 2.25]{new-tho}, 
since $\overline{\mathfrak{X}}_{\Grm}$ is compact, 
the global sections of a $\Grm^S\times K_{\Grm,S}$-equivariant sheaf 
on $\overline{\mathfrak{X}}_{\Grm}$ form a smooth 
$R[\Grm^S\times K_{\Grm,S}]$-module. We therefore have a 
well-defined derived functor $R\Gamma(\overline{\mathfrak{X}}_{\Grm},-)$,
and we obtain 
\[
R\Gamma(\overline{\mathfrak{X}}_{\Grm}, \cV)\ \mathrm{and}\ 
R\Gamma(\overline{\mathfrak{X}}_{\Grm}, j_{!}\cV)
\] 
in the bounded below derived category of smooth 
$R[\Grm^S\times K_{\Grm,S}]$-modules. We apply the
functor $R\Gamma(K_{\Grm},-)$\footnote{This is the derived
functor of $K_{\Grm}$-invariants considered with its profinite topology, so
it computes the \emph{continuous} group cohomology of $K_{\Grm}$.}, 
which gives rise to objects
\[
R\Gamma(K_{\Grm}, R\Gamma(\overline{\mathfrak{X}}_{\Grm}, \cV))\ \mathrm{and}\ 
R\Gamma(K_{\Grm}, R\Gamma(\overline{\mathfrak{X}}_{\Grm}, j_{!}\cV))
\] 
 in the bounded below derived category of 
 $\cH(\Grm^S, K_\Grm^S)\otimes_{\Z} R$-modules. 
 
On the other hand, we can also view $\cV$ and $j_{!}\cV$ 
as $K_{\Grm}$-equivariant sheaves on $\overline{\mathfrak{X}}_{\Grm}$,
using the forgetful functor. %
Recall that the action of $K_{\Grm}$ on $\overline{\mathfrak{X}}_{\Grm}$ 
is free and that the quotient can be identified with $\overline{X}^{\Grm}_{K_{\Grm}}$; 
let $\pi:\overline{\mathfrak{X}}_{\Grm} \to \overline{X}^{\Grm}_{K_{\Grm}}$ denote
the projection map. 
The \emph{descent} functor $\cF\to \left(\pi_*\cF\right)^{K_{\Grm}}$ gives an equivalence between the category 
of $K_{\Grm}$-equivariant sheaves on $\overline{\mathfrak{X}}_{\Grm}$
and the category of sheaves on $\overline{X}^{\Grm}_{K_{\Grm}}$ by~\cite[Lemma 2.24]{new-tho}. 
We denote the corresponding sheaves on $\overline{X}^{\Grm}_{K_{\Grm}}$ by $\cV$ 
and $j_{!}\cV$ as well. 

\begin{prop}\label{prop:Hecke commutative diagram}
The following diagram of derived functors is commutative  
\[
 \xymatrix@C+1pc{ \mathrm{DSh}^+_{\mathrm{G}^S\times K_{\Grm,S}}(\overline{\mathfrak{X}}_\mathrm{G}) 
 \ar[r]^-{R\Gamma(\overline{\mathfrak{X}}_\mathrm{G},\ )} \ar[d]^-{\mathrm{forget}} & 
 \mathrm{D}^+_{\mathrm{sm}}(\mathrm{G}^S\times K_{\Grm, S}, R) \ar[r]^-{R\Gamma(K_{\Grm},\ )} & 
 \mathrm{D}^+(\cH(\mathrm{G}^S,K_{\Grm}^S)\otimes_{\Z}R) \ar[d]^{\mathrm{forget}} \\ 
 \mathrm{DSh}^+_{K_{\Grm}}(\overline{\mathfrak{X}}_\mathrm{G}) 
 \ar[r]^-{\mathrm{descent}} & \mathrm{DSh}^+(\overline{X}^{\mathrm{G}}_{K_{\Grm}}) 
 \ar[r]^-{R\Gamma\left(\overline{X}^{\mathrm{G}}_{K_{\Grm}},\ \right)} &  \mathrm{D}^+(R) }
\]
\end{prop}

\begin{proof} This is a topological version of~\cite[Prop. 2.18]{new-tho}. 
The corresponding diagram of underived functors commutes up to natural isomorphism. 
The forgetful functor from $\Grm^S\times K_{\Grm,S}$-equivariant
sheaves to $K_{\Grm}$-equivariant sheaves is exact and preserves injectives by~\cite[\S 3, Corollary 3]{schneider}. 
The descent functor is also exact and preserves injectives, since it is an equivalence of categories. 
The functor $\Gamma(\overline{\mathfrak{X}}_\mathrm{G},\ )$ preserves injectives by~\cite[Lemma 2.28]{new-tho}.
\end{proof}

\noindent Note that, since $\cV$ is locally free, we have a canonical isomorphism
\[
R\Gamma\left(\overline{X}^{\Grm}_{K_{\Grm}},\cV\right)
\toisom R\Gamma\left(X^{\Grm}_{K_{\Grm}},\cV\right)
\]
induced by the pullback map $j^*$, because $j$ is a homotopy equivalence, 
and that $R\Gamma(\overline{X}^{\Grm}_{K_{\Grm}},j_!\cV)$ 
precisely computes $R\Gamma_c(X^{\Grm}_{K_{\Grm}},\cV)$. This shows
how to construct morphisms 
\[
\cH(\Grm^S, K^S_{\Grm})\otimes_{\Z}R\to \End_{\mathrm{D}^+(R)}
\left(R\Gamma_{(c)}(X^{\Grm}_{K_{\Grm}},\cV)\right)
\footnote{We will only need this statement, which is slightly weaker than saying that these
	are objects in the bounded below derived category 
	of $\cH(\Grm^S, K^S_{\Grm})\otimes_{\Z}R$-modules.}.
\]
The same formalism
also applies to $R\Gamma(\partial X^{\Grm}_{K_{\Grm}},\cV)$. 

\begin{lemma}\label{lem:bounded_coh_dim} The functor $R\Gamma(\overline{\mathfrak{X}}_\mathrm{G},\ ): \mathrm{DSh}^+_{\mathrm{G}^S\times K_{\Grm,S}}(\overline{\mathfrak{X}}_\mathrm{G}) \to \mathrm{D}^+_{\mathrm{sm}}(\mathrm{G}^S\times K_{\Grm, S}, R)$ has bounded cohomological dimension.\end{lemma}
\begin{proof}
	We can check this after applying the forgetful functor to $K_{\Grm}$-equivariant sheaves. If $\cF \in \mathrm{Sh}_{ K_{\Grm}}(\overline{\mathfrak{X}}_\mathrm{G})$, then \cite[Lemma 2.35]{new-tho} implies that $R^i\Gamma(\overline{\mathfrak{X}}_\mathrm{G},\cF)$ vanishes for $i > \dim(X^{\Grm}_{K_{\Grm}})$.
\end{proof}

We now compare our set-up with that of \cite{new-tho}. We set  $\overline{\mathfrak{X}}_{\Grm}^{\mathrm{dis}} = \mathrm{G}(F) \backslash \overline{X}^{\mathrm{G}}\times 
\mathrm{G}(\A_{F,f})^{\mathrm{dis}}$, where the superscript indicates that we are considering $\mathrm{G}(\A_{F,f})$ with the discrete topology. If we have a good compact open subgroup $K_\Grm \subset \Grm(\AA_{F,f})$, then $(K_\Grm)^{\mathrm{dis}}$ acts freely on $\overline{\mathfrak{X}}_{\Grm}^{\mathrm{dis}}$ with quotient equal to $\overline{X}^\Grm_{K_\Grm}$. Using \cite[Lemma 2.19]{new-tho} to make the Hecke action on cohomology explicit, we see that, whether we use the topological or discrete set-up, we will obtain the same Hecke actions on the cohomology of $\overline{X}^{\mathrm{G}}_{K_{\Grm}}$. We will prove something a little stronger than this, to convince the reader that the two different set-ups really are naturally equivalent.

There is a natural map $\pi_{\mathrm{dis}}:\overline{\mathfrak{X}}_{\Grm}^{\mathrm{dis}} \to \overline{\frakX}_\Grm$ which induces an exact functor \[\pi_\mathrm{dis}^*: \mathrm{Sh}_{\mathrm{G}^S\times K_{\Grm,S}}(\overline{\mathfrak{X}}_\mathrm{G}) \to\mathrm{Sh}_{(\mathrm{G}^S\times K_{\Grm,S})^\mathrm{dis}}(\overline{\mathfrak{X}}_\mathrm{G}^{\mathrm{dis}}).\] Using descent to $\overline{X}^\Grm_{K_\Grm}$ for both the topological and discrete categories, we see that pullback by $\pi_{\mathrm{dis}}$ induces an equivalence
$\mathrm{Sh}_{K_{\Grm}}(\overline{\mathfrak{X}}_\mathrm{G}) = \mathrm{Sh}_{ K_{\Grm}^\mathrm{dis}}(\overline{\mathfrak{X}}_\mathrm{G}^{\mathrm{dis}})$.

\begin{lem}\label{lem:discretevstopHecke}
We have a natural isomorphism of functors from $\mathrm{DSh}^+_{\mathrm{G}^S\times K_{\Grm,S}}(\overline{\mathfrak{X}}_\mathrm{G})$ to $\mathrm{D}^+(\cH(\mathrm{G}^S,K_{\Grm}^S)\otimes_{\Z}R)$: \[R\Gamma(K_{\Grm},-)\circ R\Gamma(\overline{\mathfrak{X}}_\mathrm{G},-)  \cong R\Gamma(K_{\Grm}^\mathrm{dis},-)\circ R\Gamma(\overline{\mathfrak{X}}_\mathrm{G}^{\mathrm{dis}},-) \circ\pi^*_{\mathrm{dis}}.\]
\end{lem}
\begin{proof}
	It follows from \cite[Lemma 2.19]{new-tho} that we can identify the underived functors $\Gamma(K_{\Grm},-)\circ \Gamma(\overline{\mathfrak{X}}_\mathrm{G},-)  = \Gamma(K_{\Grm}^\mathrm{dis},-)\circ \Gamma(\overline{\mathfrak{X}}_\mathrm{G}^{\mathrm{dis}},-) \circ\pi^*_{\mathrm{dis}}.$ From this, we deduce that there is a natural transformation \[R\Gamma(K_{\Grm},-)\circ R\Gamma(\overline{\mathfrak{X}}_\mathrm{G},-)  \to R\Gamma(K_{\Grm}^\mathrm{dis},-)\circ R\Gamma(\overline{\mathfrak{X}}_\mathrm{G}^{\mathrm{dis}},-) \circ\pi^*_{\mathrm{dis}}.\] We can check that this is an isomorphism after composing with the forgetful map to $\mathrm{D}^+(R)$, and this follows from comparing Proposition \ref{prop:Hecke commutative diagram} and \cite[Proposition 2.18]{new-tho}.	
\end{proof}

We now recall an important finiteness result:
 \begin{lemma}\label{lem:perfect complex}
  Let $K_\Grm$ be a good subgroup, and let $K'_\Grm \subset K_\Grm$ 
 be a normal subgroup which is also good. 
 Let $R$ be a Noetherian ring, and let $\cV$ be a smooth 
 $R[K_\Grm]$-module, finite free as an $R$-module. 
 Then $R \Gamma_{(c)}( X^\Grm_{K_\Grm'}, \cV)$ are 
 perfect objects of 
 $\mathrm{D}^+(K_\Grm / K_\Grm', R)$; in other words, they 
 are isomorphic in this category to bounded complexes of projective 
 $R[K_\Grm / K_\Grm']$-modules.

 \end{lemma}
 
 \begin{proof} 
 The case of usual cohomology is essentially~\cite[Lemma 2.1.7]{10author}:
 we choose a finite triangulation of $\overline{X}^{\Grm}_{K_\Grm}$ and pull this back
 to a $K_{\Grm}$-invariant triangulation of 
 $\overline{\mathfrak{X}}_{\Grm}$, then consider the corresponding complex
 of simplicial chains $C_{\bullet}$. We notice that 
 $\Hom_{\Z[K'_{\Grm}]}\left(C_{\bullet}, \cV\right)$ is isomorphic 
 in $\mathrm{D}^+(K_\Grm / K_\Grm', R)$ to 
 $R \Gamma( \overline{X}^\Grm_{K_\Grm'}, \cV)$. 
 The case of cohomology with compact support can be
 done in a similar way, by choosing our triangulation of $\overline{X}^{\Grm}_{K_\Grm}$
 in such a way that a triangulation of $\partial X^{\Grm}_{K_\Grm}$ is a 
 simplicial subcomplex. Pulling back to $\overline{\mathfrak{X}}_{\Grm}$ we 
 obtain simplicial complexes $\partial C_{\bullet}\to C_{\bullet}$. Letting 
 $C^{\mathrm{BM}}_{\bullet}$ denote the cone of this map, we observe that 
 $\Hom_{\Z[K'_{\Grm}]}\left(C^\mathrm{BM}_{\bullet}, \cV\right)$ is isomorphic 
 in $\mathrm{D}^+(K_\Grm / K_\Grm', R)$ to 
 $R \Gamma_c(X^\Grm_{K_\Grm'}, \cV)$. 
 \end{proof}
 
Assume now that $R = \cO/\varpi^m$ for some $m\in \Z_{\geq 1}$. 
If $S\subseteq S_p(F)$ is a set of places of $F$ above $p$, and if
$\cV$ is a smooth $\cO/\varpi^m[K_{S_p\setminus S}]$-module, we define the 
\emph{completed cohomology} at $S$ of level $K^{S}_{\Grm}$ 
to be 
\[
R\Gamma(K^S_{\Grm}, R\Gamma(\overline{\mathfrak{X}}_{\Grm}, \cV))
\in \mathrm{D}^+_{\mathrm{sm}}(\Grm_{S},\cO/\varpi^m).
\]
Similarly, we define the \emph{completed cohomology with compact support} at $S$ 
of level $K^{S}_{\Grm}$ to be 
\[
R\Gamma(K^S_{\Grm}, R\Gamma(\overline{\mathfrak{X}}_{\Grm}, j_{!}\cV))
\in \mathrm{D}^+_{\mathrm{sm}}(\Grm_{S},\cO/\varpi^m).
\]
For a finite set of finite places $T\supseteq S_p(F)$ of $F$,
a variant of the above formalism equips these objects with 
actions of $\cH(\Grm^{T}, K^{T}_{\Grm})\otimes_{\Z}\cO/\varpi^m$.
The same formalism applies to 
\[
R\Gamma(K^S_{\Grm}, R\Gamma(\partial\mathfrak{X}_{\Grm}, \cV))
\in \mathrm{D}^+_{\mathrm{sm}}(\Grm_{S},\cO/\varpi^m).
\]
The following lemma offers a justification for the term \emph{completed cohomology}.  

\begin{lemma}\label{lem:recovering completed cohomology} 
For any $i\in \Z_{\geq 0}$, we 
have $\cH(\Grm^{T}, K^{T}_{\Grm})$-equivariant isomorphisms 
of admissible smooth $\cO/\varpi^m[\Grm_S]$-modules
\begin{equation}\label{eq:cc}
H^i\left( R\Gamma(K^S_{\Grm}, R\Gamma(\overline{\mathfrak{X}}_{\Grm}, \cV)) \right)
\toisom 
\varinjlim_{K_{\Grm,S}} H^i\left(X^{\Grm}_{K^S_{\Grm}K_{\Grm,S}}, \cV\right)
\end{equation}
and 
\begin{equation}\label{eq:ccc}
H^i\left( R\Gamma(K^S_{\Grm}, R\Gamma(\overline{\mathfrak{X}}_{\Grm}, j_{!}\cV)) \right)
\toisom 
\varinjlim_{K_{\Grm,S}} H^i_c\left(X^{\Grm}_{K^S_{\Grm}K_{\Grm,S}}, \cV\right). 
\end{equation}
\end{lemma}

\begin{proof} In the category
of compact Hausdorff spaces, we have 
\[
\overline{\mathfrak{X}}_{\Grm}/K^S_{\Grm} = 
\varprojlim_{K_{\Grm,S}} \overline{X}^{\Grm}_{K^S_{\Grm}K_{\Grm,S}}. 
\] 
Furthermore, $K^{S}_{\Grm}$ acts freely on $\overline{\mathfrak{X}}_{\Grm}$,
so we can functorially rewrite the LHS of~\eqref{eq:cc} and~\eqref{eq:ccc} 
in terms of the cohomology 
of either $\cV$ or $j_{!}\cV$ on the quotient $\overline{\mathfrak{X}}_{\Grm}/K^S_{\Grm}$. 
We can functorially rewrite the terms on the RHS in terms of the cohomology 
of either $\cV$ or $j_{!}\cV$ on $\overline{X}^{\Grm}_{K^S_{\Grm}K_{\Grm,S}}$. 
The result now follows from~\cite[Lemma 2.34]{new-tho}. Finally, admissibility 
of the cohomology groups follows from Lemma~\ref{lem:perfect complex}. 
\end{proof}

Our next result will imply an important property of completed 
cohomology: it is, in some sense, independent of the weight
$\cV$. It will be useful to work in a little more generality, so we 
assume that $\cV$ is a smooth $\cO/\varpi^m[K_{\Grm,S_p\backslash S}\times \Delta_S]$-module, 
flat over $\cO/\varpi^m$, for an open submonoid $\Delta_{S} \subset \Grm_{S}$ 
which contains an open subgroup $U_{S}$ of $\Grm_{S}$. 
As above, we associate to $\cV$ a $\Grm^T\times K_{\Grm, T\backslash S} \times 
U_S$-equivariant sheaf on $\overline{\mathfrak{X}}_{\Grm}$ by pulling back from a point. 

\begin{lemma}\label{lem:projection formula}
We have canonical isomorphisms 
\[
R\Gamma(\overline{\mathfrak{X}}_{\Grm}, \cO/\varpi^m)\otimes \cV \toisom R\Gamma(\overline{\mathfrak{X}}_{\Grm}, \cV)
\]
and
\[
R\Gamma(\partial\overline{\mathfrak{X}}_{\Grm}, \cO/\varpi^m)\otimes\cV\toisom R\Gamma(\partial \overline{\mathfrak{X}}_{\Grm}, \cV)
\]
in $\mathrm{D}^+_{\sm}(\Grm^T\times K_{\Grm, T\backslash S} \times U_S, \cO/\varpi^m)$. 
\end{lemma}

\begin{proof} We explain the case of $\overline{\mathfrak{X}}_{\Grm}$,
the case of $\partial\overline{\mathfrak{X}}_{\Grm}$ is the same. Let 
$f: \overline{\mathfrak{X}}_{\Grm}\to *$ be the $\Grm(\A_{F,f})$-equivariant 
projection to a point. We set $H:= G^T\times K_{\Grm, T\backslash S} \times U_S$. There is a pair of adjoint functors $(f^*, Rf_*)$ between
$D^+(\Sh_{H}(*), \cO/\varpi^m) \simeq D^+_{\sm}(H, \cO/\varpi^m)$ 
and $D^+(\Sh_{H}(\overline{\mathfrak{X}}_{\Grm}), \cO/\varpi^m)$.  
There is also a natural isomorphism 
\[
f^*(Rf_*(\cO/\varpi^m) \otimes \cV) \toisom f^*Rf_*(\cO/\varpi^m) \otimes f^*\cV
\] and hence by adjunction a natural map \[
f^*(Rf_*(\cO/\varpi^m) \otimes \cV) \to f^*\cV.
\]
We therefore have a morphism 
\begin{equation}\label{eq:projection morphism}
R\Gamma (\overline{\mathfrak{X}}_{\Grm}, \cO/\varpi^m)\otimes \cV
= Rf_*(\cO/\varpi^m) \otimes \cV\longrightarrow
Rf_* f^*\cV = R\Gamma (\overline{\mathfrak{X}}_{\Grm}, \cV) 
\end{equation}
in $\mathrm{D}^+_{\sm}(H, \cO/\varpi^m)$. It is enough to show that this
is an isomorphism after forgetting the equivariant structure. By~\cite[\S 1, Corollary 3]{schneider}, 
if we forget the equivariant structure for the sheaves on $\overline{\mathfrak{X}}_{\Grm}$,
 the resulting derived functors compute cohomology with compact support. 
Since $\overline{\mathfrak{X}}_{\Grm}$ is a compact Hausdorff space, 
we can apply~\cite[Prop. 2.6.6]{kashiwara-schapira}, which implies
that the morphism in~\eqref{eq:projection morphism} is an isomorphism (since $\cV$ is flat over $\cO/\varpi^m$, the assumption in \emph{loc.~cit.} that the coefficient ring has finite weak global dimension is not necessary). 
\end{proof}

We can use this lemma to \emph{define} an object \[R\Gamma(\overline{\mathfrak{X}}_{\Grm}, \cV):= R\Gamma(\overline{\mathfrak{X}}_{\Grm}, \cO/\varpi^m)\otimes \cV \in \mathrm{D}^+_{\sm}(\Grm^T\times K_{\Grm, T\backslash S} \times \Delta_S, \cO/\varpi^m).\] It is independent of the choice of $U_S$, by \cite[\S 1, Corollary 3]{schneider}. When $K_{\Grm,S}$ is a compact open subgroup of $\Delta_S$, we obtain $R\Gamma(K_{\Grm},R\Gamma(\overline{\mathfrak{X}}_{\Grm}, \cV)) \in \mathrm{D}^+(\cO/\varpi^m)$ with an action of $\cH(\Grm^T,K_\Grm^T)\otimes\cH(\Delta_{S},K_{\Grm,S})$.

In \S\ref{sec:boundary coh}, we will need a variant of this lemma with a coefficient system in a derived category. When we apply this lemma, we will just have group actions, not monoids, so we now assume $\cV \in \mathrm{D}^b_{\sm}(K_{\Grm,S_p\backslash S},\cO/\varpi^m)$. After inflation and pullback from a point, we get a corresponding object $\cV \in \mathrm{D}^b(\Sh_{\Grm^T\times K_{\Grm,T}}(\overline{\mathfrak{X}}_{\Grm},\cO/\varpi^m))$. To state the lemma, we need the derived tensor product functor $R\Gamma(\overline{\mathfrak{X}}_{\Grm}, \cO/\varpi^m)\otimes^{\LL}-$. Although $\Mod_{\sm}(\Grm^T\times K_{\Grm, T}, \cO/\varpi^m)$ does not have enough projectives, every object has a surjection from a $\cO/\varpi^m$-flat object.  Indeed, smoothness implies that there is a surjection from a direct sum of copies of compact inductions of trivial representations of compact open subgroups on $\cO/\varpi^m$\footnote{This is the usual proof of `enough projectives' over a characteristic $0$ field. The problem here is that the trivial representation with coefficients in $\cO/\varpi^m$ of a non-trivial compact $p$-adic group is not projective.}. This gives enough acyclic objects to compute derived tensor products on $\mathrm{D}^-_{\sm}(\Grm^T\times K_{\Grm, T}, \cO/\varpi^m)$. Since the functor $R\Gamma(\overline{\mathfrak{X}}_{\Grm},-)$ on $D^+(\Sh_{\Grm^T\times K_{\Grm,T}}(\overline{\mathfrak{X}}_{\Grm},\cO/\varpi^m))$ has bounded cohomological dimension (Lemma \ref{lem:bounded_coh_dim}), it takes bounded objects to bounded objects.

\begin{lem}\label{lem:projection formula derived}
	Let $\cV \in \mathrm{D}^b_{\sm}(K_{\Grm,S_p\backslash S},\cO/\varpi^m)$. We have a canonical isomorphism 
	\[R\Gamma(\overline{\mathfrak{X}}_{\Grm}, \cO/\varpi^m)\otimes^{\LL} \cV \toisom R\Gamma(\overline{\mathfrak{X}}_{\Grm}, \cV)
	\]
	in $\mathrm{D}^b_{\sm}(\Grm^T\times K_{\Grm, T}, \cO/\varpi^m)$. 
\end{lem}
\begin{proof}
	As in the proof of Lemma \ref{lem:projection formula}, we set $H = \Grm^T\times K_{\Grm,T}$ and consider $f:\overline{\frakX}_{\Grm}\to \ast$ the map to the point. By Lemma \ref{lem:bounded_coh_dim}, we have a pair of adjoint functors $(f^*,Rf_*)$ between unbounded derived categories  $\mathrm{D}_{\sm}(H,\cO/\varpi^m)$ and $\mathrm{D}(\Sh_{H}(\overline{\frakX}_{\Grm},\cO/\varpi^m))$. Since $f^*$ is exact, it is easy to see that we have a natural isomorphism \[
	f^*(Rf_*(\cO/\varpi^m) \otimes^{\LL} \cV) \toisom f^*Rf_*(\cO/\varpi^m) \otimes^{\LL} f^*\cV
	\] and we then obtain a map \[p_\cV:R\Gamma(\overline{\mathfrak{X}}_{\Grm}, \cO/\varpi^m)\otimes^{\LL} \cV \to R\Gamma(\overline{\mathfrak{X}}_{\Grm}, \cV)
	\] by adjunction. The fact that this is an isomorphism follows from the case where $\cV$ is a $\cO/\varpi^m$-flat module. More precisely, we can replace $\cV$ by a bounded above complex $\cF^\bullet$ of $\cO/\varpi^m$-flat objects in $\Mod_{\sm}(H, \cO/\varpi^m)$, and replace $\cO/\varpi^m$ by a bounded complex $\cI^\bullet$ of $Rf_*$-acyclic objects in $\Sh_{H}(\overline{\frakX}_{\Grm},\cO/\varpi^m)$. Using \cite[Prop. 2.6.6]{kashiwara-schapira} again, we see that each sheaf $\cI^i \otimes \cF^j$ is $Rf_*$-acyclic and the natural map \[f_*\cI^i\otimes\cF^j \to f_*(\cI^i \otimes \cF^j) \] is an isomorphism. The total complexes of the double complexes $f_*\cI^\bullet\otimes\cF^\bullet \to f_*(\cI^\bullet \otimes \cF^\bullet)$ respectively compute the source and target of $p_\cV$, so we see that $p_\cV$ is an isomorphism in $\mathrm{D}^b_{\sm}(\Grm^T\times K_{\Grm, T}, \cO/\varpi^m)$. 
\end{proof}
\begin{remark}
	With Lemma \ref{lem:projection formula derived} in hand, we can prove that $R\Gamma(\overline{\mathfrak{X}}_{\Grm}, \cO/\varpi^m)$ has bounded Tor-dimension, and then extend our projection formula to handle $\cV$ in the unbounded derived category. See for example \cite[Corollary 6.5.6]{Fu-etale} for the classical projection formula in $l$-adic cohomology.  
\end{remark}

Assume now that $R =\cO$, and that $\cV$ is an 
$\cO[K_{\Grm,S}]$-module, %
which is finite free as an $\cO$-module and such that $\cV/\varpi^m$ 
is a smooth $\cO/\varpi^m[K_{\Grm,S}]$-module for each $m\in \Z_{\geq 1}$. 
We then consider 
\[
R\Gamma(\overline{X}^{\Grm}_{K_{\Grm}}, \cV) := 
\varprojlim_{m} R\Gamma(\overline{X}^{\Grm}_{K_{\Grm}}, \cV/\varpi^m)
\]
in $\mathrm{D}^+(\cO)$, 
where the projective limit should be understood as a homotopy limit. 
We also consider the analogue with coefficient system $j_{!}\cV$. 
These limits can be endowed with an action of the Hecke algebra 
$\cH(\Grm^S,K_{\Grm}^S)\otimes_{\Z}\cO$\footnote{As
the proof of Lemma~\ref{lem:perfect complex} shows, we
have explicit perfect complexes that compute these derived functors, so 
we could simply take a projective limit on the level of complexes. To 
endow the projective limit with a Hecke action, we can instead 
consider adelic complexes that compute these derived
functors as in~\cite[\S 5.1]{borel}. Combining the fact that derived limits commute with cohomology \cite[\href{https://stacks.math.columbia.edu/tag/08U1}{Tag 08U1}]{stacks-project} and Lemma \ref{lem:discretevstopHecke} we can show that we will obtain the same Hecke actions as in \cite{new-tho}.}. 

Continue to assume $R=\cO$. If $S\subseteq S_p(F)$ is a set of places of $F$ above $p$, let
$\cV$ be an $\cO[K_{S_p\setminus S}]$-module which is finite 
free as an $\cO$-module and such that $\cV/\varpi^m$ is a smooth 
$\cO/\varpi^m[K_{S_p\setminus S}]$-module. We consider 
\[
R\Gamma(K^S_{\Grm}, R\Gamma(\overline{\mathfrak{X}}_{\Grm}, \cV)):=
\varprojlim_{m} R\Gamma(K^S_{\Grm}, R\Gamma(\overline{\mathfrak{X}}_{\Grm}, \cV/\varpi^m))
\]
in $\mathrm{D}^+(\Grm_{S}, \cO)$, where again the projective
limit should be understood as a homotopy limit. There is also 
the analogue with coefficient system $j_{!}\cV$. When $T\supseteq S_p(F)$
is a finite set of finite places of $F$, these limits can also 
be endowed with an action of the Hecke algebra 
$\cH(\Grm^T,K_{\Grm}^T)\otimes_{\Z}\cO$. 
 
We now assume that $\Grm$ is reductive, and let 
$\Prm = \mathrm{M} \mathrm{N}$ be a parabolic subgroup with 
Levi subgroup $\mathrm{M}$. Let $K_\Grm \subset \Grm(\A_{F,f})$ be 
a good subgroup. In this situation, we define $K_\Prm = K_\Grm \cap \Prm(\A_{F,f})$, 
$K_\mathrm{N} = K_\Grm \cap \mathrm{N}(\A_{F,f})$, and define $K_\mathrm{M}$ 
to be the image of $K_\Prm$ in $\mathrm{M}(\A_{F,f})$. 
We say that $K_\Grm$ is \emph{decomposed} with respect to $\Prm = \mathrm{M} 
\mathrm{N}$ if we have $K_\Prm = K_\mathrm{M} \ltimes K_{\mathrm{N}}$; 
equivalently, if $K_\mathrm{M} = K_\Grm \cap \mathrm{M}(\A_{F,f})$. 

Assume now that $K_\Grm$ is decomposed with respect to $\Prm = 
\mathrm{M} \mathrm{N}$, and let $S$ be a finite set of finite places of 
$F$ such that for all $v \not\in S$, $K_{\Grm, v}$ is a hyperspecial 
maximal compact subgroup of $\Grm(F_v)$. In this case, 
we can define homomorphisms 
\[ r_\Prm : \cH(\Grm^S, K_{\Grm}^S) \to \cH(\Prm^S, K_\Prm^S) 
\text{ and } r_{\mathrm{M}} : \cH(\Prm^S, K_\Prm^S) \to \cH(\mathrm{M}^S, K_\mathrm{M}^S), \]
given respectively by ``restriction to $\Prm$'' and 
``integration along $\mathrm{N}$''; see~\cite[\S 2.2.3]{new-tho}
and \cite[\S 2.2.4]{new-tho} respectively 
for the definitions of these maps, along with the proofs 
that they are indeed algebra homomorphisms,
and that $r_{\mathrm{M}}$ preserves integrality. We 
use $\cS:=r_{\mathrm{M}}\circ r_{\Prm}$ to denote the 
unnormalized Satake transform. 

Finally, we remark that the above formalism also 
applies to the case of the Hecke algebra of a monoid, 
see~\cite[\S 2.1.9]{10author}. 
 
\subsubsection{The general linear group and the quasi-split unitary group}\label{sec: unitary group}

From now on, we fix an integer $n\geq 2$ and we let 
$F$ be an imaginary CM field containing 
the maximal totally real subfield $F^+$. %
Let $c\in \Gal(F/F^+)$ denote complex conjugation. 
We set $\overline{S}_p:=S_p(F^+)$ and $S_p:=S_p(F)$. 
We let $\Psi_n$ be the matrix with 
$1$'s on the anti-diagonal and $0$'s elsewhere, and we let 
\[
J_n = \left(\begin{matrix} 0 & \Psi_n \\ -\Psi_n & 0 \end{matrix}\right). 
\]  
We let $\widetilde{G}/\cO_{F^+}$ be the group scheme 
defined by 
\[
\widetilde{G}(R) = \{g\in \GL_{2n}(R\otimes_{\cO_{F^+}}\cO_F)\mid
\ ^tgJ_ng^c = J_n\}
\]
for any $\cO_{F^+}$-algebra $R$. The generic fibre of 
$\widetilde{G}$ over $F^+$ is a quasi-split unitary group, 
which becomes isomorphic to $\GL_{2n}/F$ after base change 
from $F^+$ to $F$. In particular, if $\bar{v}$ is a place of
$F^+$ that splits in $F$, a choice of place $v\mid \bar{v}$ 
of $F$ determines a canonical isomorphism $\iota_v: G(F^+_{\bar{v}})
\toisom \GL_{2n}(F_v)$. 

We let $P\subset \widetilde{G}$ denote the Siegel parabolic 
consisting of block upper-triangular matrices with blocks 
of size $n\times n$. We let $P = U\rtimes G$ be a Levi 
decomposition such that we can identify 
$G$ with $\mathrm{Res}_{\cO_F/\cO_{F^+}}\GL_n$\footnote{
We use the same identification  
as in~\cite[\S 2.2.1]{10author}, namely $\left(\begin{smallmatrix}A & 0 
\\ 0 & D\end{smallmatrix}\right)\in G(R) \mapsto D\in \GL_n(R\otimes_{\cO_{F^+}}\cO_F)$
for an $\cO_{F^+}$-algebra $R$.}.
To simplify the notation, from now on we write 
$\widetilde{X}$ for $X^{\widetilde{G}}$ and 
$X$ for $X^G$. We also write $\widetilde{K}$ 
and $K$ for good subgroups of $\widetilde{G}(\A_{F^+,f})$ 
and of $G(\A_{F^+,f}) = \GL_n(\A_{F,f})$. Note that
the locally symmetric spaces $\tX_{\tK}$ are complex
manifolds of (complex) dimension $d:=n^2[F^+:\Q]$, whereas 
the locally symmetric spaces $X_{K}$ are real manifolds
of (real) dimension $d-1$. 

We now describe some explicit (integral and rational) 
coefficient systems for these symmetric spaces. These
will depend on a choice of a prime $p$ and on a choice
of a dominant weight for either $G$ or $\tG$. We fix a coefficient
field $E/\Q_p$ which is assumed to be sufficiently large, 
so that it contains the image of every embedding 
$\Hom(F,\overline{\Q}_p)$.  Let $T\subset \tB\subset \tG$
be the maximal torus of diagonal matrices and the upper triangular Borel
subgroup, respectively. Set $B :=\tB \cap G$, this can be identified
with the upper triangular Borel subgroup in $G$. 

We first treat the case of $G$. We identify the character group of 
$(\mathrm{Res}_{F^+/\Q}T)_E$ with $(\Z^n)^{\Hom(F, E)}$
in the usual way. A weight $(\lambda_{\tau,i})\in (\Z^n)^{\Hom(F,E)}$ 
with $\tau\in \Hom(F,E)$ and $i\in 1,\dots, n$ 
is {dominant} for $(\mathrm{Res}_{F^+/\Q}B)_E$ if it satisfies 
\[
\lambda_{\tau, 1}\geq \lambda_{\tau,2} \geq \dots \geq \lambda_{\tau,n}
\]
for each $\tau\in \Hom(F,E)$. We denote by $(\Z^n_+)^{\Hom(F,E)}$ the
subset of dominant weights. The expression `${\lambda}$ is a dominant weight for $G$' will indicate that a weight ${\lambda} \in X^\ast((\mathrm{Res}_{F^+/\Q}T)_E)$ is dominant for $(\mathrm{Res}_{F^+/\Q}B)_E$.

Assume now that $\lambda \in (\Z^n_+)^{\Hom(F,E)}$. 
We define the $G(\cO_{F^+,p}) = \prod_{{v}\in{S}_p}G(\cO_{F_{v}})$-representation $\cV_{\lambda}$ 
to be the integral dual Weyl module of highest weight $\lambda$ with coefficients in $\cO$, 
obtained from the Borel--Weil construction. 
More precisely, if
we let $\mathrm{B}_n\subset \GL_n$ denote the standard Borel 
consisting of upper-triangular matrices and $w_{0,n}$ denote
the longest element in the Weyl group of $\GL_n$, we consider 
the algebraic induction
\[
(\mathrm{Ind}_{\mathrm{B}_n}^{\GL_n} w_{0,n}\lambda_{\tau})_{/\cO}:= 
\{f\in \cO[\GL_n]\mid f(bg) = (w_{0,n}\lambda_{\tau})(b)f(g),
\]
\[
\forall\ \cO\to R, 
b\in \mathrm{B}_n(R), g\in \GL_n(R)\},
\] 
and we set $\cV_{\lambda_{\tau}}$ to be the finite free 
$\cO$-module obtained by evaluating this on $\cO$ and
$V_{\lambda_{\tau}}:=\cV_{\lambda_{\tau}}\otimes_{\cO}E$. When $\tau$ induces the place $v$ of $F$, these modules come with an action of $\GL_n(\cO_{F_v})$ and $\GL_n(F_v)$ respectively.

Finally, we set $\cV_{\lambda} := \otimes_{\tau,\cO} \cV_{\lambda_{\tau}}$
and $V_{\lambda}:=\cV_{\lambda}\otimes_{\cO}E$. Then 
$V_{\lambda}$ is the absolutely irreducible algebraic representation of $(\mathrm{Res}_{F^+/\Q}G)_{E}$
of highest weight $\lambda$ and it is finite-dimensional over $E$; the 
lattice $\cV_{\lambda}\subset V_{\lambda}$ is $G(\cO_{F^+,p})$-stable.  
For every $m\in \Z_{\geq 1}$, $\cV_{\lambda}/\varpi^m$
is a smooth $\cO/\varpi^m[G(\cO_{F^+,p})]$-module
that is finite free as an $\cO/\varpi^m$-module. Therefore,
the formalism of the previous section applies to $\cV_{\lambda}$. 

We now treat the case of $\widetilde{G}$. Assume that each 
place in $\overline{S}_p$ splits from $F^+$ to $F$ and that 
we have a partition of the form 
$S_p = \widetilde{S}_p\sqcup \widetilde{S}^c_p$, with $\tilde{v}\in \widetilde{S}_p$ 
the place lying above a place $\bar{v}\in \overline{S}_p$. This induces 
a partition on $\Hom(F,E)$, by choosing the embedding $\tilde{\tau}: F\hookrightarrow E$ 
above a given embedding $\tau: F^+\hookrightarrow E$ that induces a place
in $\widetilde{S}_p$. In turn, this induces an identification 
\[
(\mathrm{Res}_{F^+/\Q}\tG)_E = \prod_{\Hom(F^+,E)}\GL_{2n,E}
\]
and therefore an identification of the character group
of $(\mathrm{Res}_{F^+/\Q}T)_E$ with $(\Z^{2n})^{\Hom(F^+, E)}$. 
More precisely, this identifies a weight $\lambda = (\lambda_{\tilde{\tau},i})$
with a weight $\tilde{\lambda} = (\tilde{\lambda}_{\tau,i})$ where 
\begin{equation}\label{eq:identification of weights}
\tilde{\lambda}_{\tau} = (-\lambda_{\tilde{\tau} c,n},\dots , -\lambda_{\tilde{\tau} c,1},
\lambda_{\tilde{\tau},1},\dots,\lambda_{\tilde{\tau},n}). 
\end{equation}
The set of weights that are dominant for $(\mathrm{Res}_{F^+/\Q}\tB)_E$ 
are the ones in the subset $(\Z^{2n}_+)^{\Hom(F^+,E)}$. For such weights, 
we can therefore define the integral dual Weyl module of highest weight $\tilde{\lambda}$, $\cV_{\tilde{\lambda}} \subset V_{\tilde{\lambda}}$, a $\tG(\cO_{F^+,p})$-stable $\cO$-lattice in the highest weight $\tilde{\lambda}$ representation of $(\mathrm{Res}_{F^+/\Q}\tG)_E$. 
For every $m\in \Z_{\geq 1}$, $\cV_{\tilde{\lambda}}/\varpi^m$
is a smooth $\cO/\varpi^m[\prod_{\bar{v}\in \overline{S}_p} \tG(\cO_{F^+_{\bar{v}}})]$-module
that is finite free as an $\cO/\varpi^m$-module. Therefore,
the formalism of the previous section also 
applies to $\cV_{\tilde{\lambda}}$. We say `$\tilde{\lambda}$ is a dominant weight for $\tG$' to indicate that a weight $\tilde{\lambda} \in X^\ast((\mathrm{Res}_{F^+/\Q}T)_E)$ is dominant for $(\mathrm{Res}_{F^+/\Q}\tB)_E$.

We now define appropriate quotients of the Hecke algebras 
acting on the cohomology groups with these coefficient systems. Again,
we treat $G$ first. Let $S\supseteq S_p$ be a finite set of finite 
places of $F$ and let $K\subset \GL_n(\A_{F,f})$ be a good subgroup
such that $K_v = \GL_n(\cO_{F_v})$ for $v\not \in S$ and 
$K_v\subseteq \GL_n(\cO_{F_v})$ for $v\in S_p$. For any 
$\lambda\in (\Z^n_+)^{\Hom(F,E)}$, the complex $R\Gamma(X_{K},\cV_{\lambda})$ 
is well-defined as an object of $D^+(\cO)$ (up to unique isomorphism) 
and equipped with a Hecke action. We set 
$\mathbb{T}^S:= \cH(G^S, K^S)\otimes_{\Z}\cO$ 
and 
\[
\mathbb{T}^S(K,\lambda):= \mathrm{Im}\left(\mathbb{T}^S\to 
\mathrm{End}_{D^+(\cO)}(R\Gamma(X_{K},\cV_{\lambda}))\right). 
\]

In the case of $\tG$, let $S\supseteq S_p$ be a finite set of 
finite places of $F$ satisfying $S=S^c$. Let $\overline{S}$ denote
the set of finite places of $F^+$ below $S$. Let $\tK\subset \tG(\A_{F^+,f})$
be a good compact open subgroup such that $\tK_{\bar{v}} = \tG(\cO_{F^+_{\bar{v}}})$
for $\bar{v}\not\in\overline{S}$ and $\tK_{\bar{v}} \subseteq \tG(\cO_{F^+_{\bar{v}}})$ 
for $\bar{v}\in \overline{S}_p$. To simplify notation, we write $\tG^S = 
\tG^{\overline{S}}$ etc. For any $\tilde{\lambda}\in (\Z^{2n}_{+})^{\Hom(F^+,E)}$,
the complex $R\Gamma(X_{\tK},\cV_{\tilde{\lambda}})$ 
is well-defined as an object of $D^+(\cO)$ (up to unique isomorphism) 
and equipped with a Hecke action. We consider the abstract Hecke $\cO$-algebra
$\widetilde{\mathbb{T}}^S:= \cH(\tG^S, \tK^S)\otimes_{\Z}\cO$ and
its quotient $\mathbb{T}^S(\tK,\tilde{\lambda})$ acting faithfully on $R\Gamma(X_{\tK},\cV_{\tilde{\lambda}})$. 

As a consequence of Lemma~\ref{lem:perfect complex}, we 
see that both $\mathbb{T}^S(K,\lambda)$ and $\widetilde{\mathbb{T}}^S(\tK,\tilde{\lambda})$
are finite $\cO$-modules. 
There are obvious versions of all of this with $\cO/\varpi^m$-coefficients
and for compactly supported cohomology and for the cohomology of 
the boundary $\partial X_{\tK}$ of the Borel--Serre compactification 
of $X_{\tK}$. %

We will make use of particular elements of some Weyl groups, besides the longest element $w_{0,n}$ in the Weyl group of $\GL_n$ which we have already mentioned. For $\Grm = G$ or $\tG$, we will write $w_0^{\Grm}$ for the longest element in the Weyl group  $W((\mathrm{Res}_{F^+/\Q}\Grm)_E, (\mathrm{Res}_{F^+/\Q}T)_E)$. We set $w^P_0 = w_0^{G}w_0^{\tG}$. %
It is the longest element in the set $W^P$ of minimal length coset representatives for \[W((\mathrm{Res}_{F^+/\Q}\tG)_E, (\mathrm{Res}_{F^+/\Q}T)_E)/W((\mathrm{Res}_{F^+/\Q}G)_E, (\mathrm{Res}_{F^+/\Q}T)_E).\]

In our development of $P$-ordinary Hida theory, it will be important to compare coefficient systems for $\tG$ and $G$.  

Let $\mathrm{P}_{n,n}\subset \GL_{2n}$ be the parabolic subgroup 
of block-upper triangular matrices with Levi quotient $\GL_n\times \GL_n$.
By the transitivity of algebraic induction, $\cV_{\tilde{\lambda}_{\tau}}$ and $V_{\tilde{\lambda}_{\tau}}$ can be identified
with the evaluation on $\cO$ and $E$ respectively of the algebraic induction
\[
\left(\mathrm{Ind}_{\mathrm{P}_{n,n}}^{\GL_{2n}} \cV_{\lambda_{\tilde{\tau}}}\otimes
\cV_{-w_{0,n}\lambda_{\tilde{\tau}c}}\right)_{/\cO}. 
\] 

\begin{lemma}\label{lem:evsurjective}
The natural $\Prm_{n,n}(\cO)$-equivariant 
morphism 
\[
\cV_{\tilde{\lambda}_{\tau}}\to \cV_{\lambda_{\tilde{\tau}}} \otimes \cV_{-w_{0,n}\lambda_{\tilde{\tau}c}}
\]
given by evaluation of functions at the identity is surjective. 
\end{lemma}

\begin{proof} 
By transitivity of parabolic induction, we can identify
$\cV_{\tilde{\lambda}_{\tau}}$ with the evaluation on $\cO$ 
of 
\[
\left(\mathrm{Ind}_{\mathrm{B}_{2n}}^{\GL_{2n}}w_{0,2n}\tilde{\lambda}_{\tau}\right)_{/\cO}
\toisom
\left(\mathrm{Ind}_{\mathrm{P}_{n,n}}^{\GL_{2n}}\circ 
\mathrm{Ind}_{\mathrm{B}_n\times B_{n}}^{\GL_n\times \GL_n}w_{0,2n}\tilde{\lambda}_{\tau}\right)_{/\cO}, 
\]
where, by~\cite[\S I.3.5]{jantzen}, 
the map is given by $f\mapsto \tilde{f}(g)(h) = f(hg)$ for all $\cO$-algebras $R$, 
$h\in \GL_n(R)\times \GL_n(R)$, and $g\in \GL_{2n}(R)$. By Nakayama, 
it is enough to check surjectivity after base change to $\overline{\F}_p$, in which 
case the evaluation at identity map can be rewritten in geometric terms 
as the restriction map 
\[
H^0(X, \cL)\to H^0(X', \cL), 
\]
where $X = \mathrm{B}_{2n} \backslash \GL_{2n}$ is the full flag variety
for $\GL_{2n}$, $X'\subset X$ is the Schubert variety
for the longest Weyl group element in $\GL_n\times \GL_n$ and
$\cL$ is the line bundle on $X$ determined by $w_{0,2n}\tilde{\lambda}$. 
The result now follows from the main theorem of~\cite{andersen}
applied to $\mathrm{SL}_{2n}/\overline{\F}_p$.  
\end{proof}

\subsubsection{Explicit Hecke operators}\label{sec:explicit Hecke operators} 
Fix once and for all a choice $\varpi_{\bar{v}}$ of uniformizer of $F^+_{\bar{v}}$ for every finite place $\bar{v}$ of $F^+$. When $\bar{v}$ is unramified in $F$ we set $\varpi_v = \varpi_{\bar{v}}$ for $v|\bar{v}$.

We define some explicit Hecke operators at unramified primes first. If $v$ is a 
finite place of $F$ and $1 \leq i \leq n$ is an integer then 
we write $T_{v, i} \in \cH(\GL_n(F_v), \GL_n(\cO_{F_v}))$ 
for the double coset operator
\[
T_{v, i} = [ \GL_n(\cO_{F_v}) \diag(\varpi_v, \dots, \varpi_v, 1, \dots, 1) \GL_n(\cO_{F_v})], 
\]
where $\varpi_v$ appears $i$ times on the diagonal.
This is the same as the operator denoted by $T_{M, v, i}$ 
in ~\cite[Prop.-Def. 5.3]{new-tho}. 
We define a polynomial 
\begin{equation}\label{eqn:hecke_pol_for_GL_n} \begin{split} 
P_v(X) = X^n&-T_{v, 1}X^{n-1} + \dots + (-1)^iq_v^{i(i-1)/2}T_{v, i}X^{n-i}+\dots 
\\ & + (-1)^n q_v^{n(n-1)/2}T_{v,n} \in \cH(\GL_n(F_v), \GL_n(\cO_{F_v}))[X].
	\end{split}  
\end{equation}
It corresponds to the characteristic polynomial of a Frobenius 
element on $\rec^T_{F_v}(\pi_v)$, where $\pi_v$ is an 
unramified representation of $\GL_n(F_v)$. %

If $\overline{v}$ is a place of $F^+$ unramified in $F$, and 
$v$ is a place of $F$ above $\overline{v}$, and $1 \leq i \leq 2n$ 
is an integer, then we write $\widetilde{T}_{v, i} \in 
\cH( \widetilde{G}(F^+_{\overline{v}}),\widetilde{G}(\cO_{F^+_{\overline{v}}})) 
\otimes_\Z \Z[q_{\overline{v}}^{-1}]$ for the operator denoted 
$T_{G, v, i}$ in~\cite[Prop.-Def. 5.2]{new-tho}. We define a polynomial
\begin{equation}\label{eqn:hecke_pol_for_tilde_G} \begin{split} 
\widetilde{P}_v(X) = X^{2n} & - \widetilde{T}_{v, 1} X^{2n-1} +  
\dots + (-1)^j q_v^{j(j-1)/2} \widetilde{T}_{v, j}X^{2n-j} +  \dots \\&+ q_v^{n(2n-1)} 
\widetilde{T}_{v, 2n} \in \cH( \widetilde{G}(F^+_{\overline{v}}),\widetilde{G}
(\cO_{F^+_{\overline{v}}})) \otimes_\Z \Z[q_{\overline{v}}^{-1}][X]. \end{split} 
\end{equation}
It corresponds to the characteristic polynomial of a Frobenius element on 
$\rec^T_{F_v}(\pi_v)$, where $\pi_v$ is the base change of an 
unramified representation $\sigma_{\overline{v}}$ of 
the group $\widetilde{G}(F^+_{\overline{v}})$. 

We now describe the behaviour of these Hecke operators 
under the unnormalized Satake transform with respect to the Siegel parabolic. 
We use the following convention: if $f(X)$ is a polynomial of degree $d$,
with constant term a unit $a_0$, we set $f^\vee(X):=a_0^{-1}X^df(X^{-1})$. 

\begin{prop}\label{prop:satake_transform_unramified_case}
Let $v$ be a place of $F$, unramified over the place $\overline{v}$ of $F^+$. Let 
\[ 
\cS : \cH( \widetilde{G}(F^+_{\overline{v}}), \widetilde{G}(\cO_{F^+_{\overline{v}}})) \otimes_{\Z}\Z[q_{\bar{v}}^{-1}]
\to \cH( G(F^+_{\overline{v}}), G(\cO_{F^+_{\overline{v}}}))
\otimes_{\Z}\Z[q_{\bar{v}}^{-1}]
\]
denote the homomorphism defined at the end of \S \ref{sec:locally symmetric spaces}. 
Then we have
\[
\cS( \widetilde{P}_v(X)) = P_v(X) q_v^{n(2n-1)} P_{v^c}^\vee(q_v^{1-2n} X). 
\]
\end{prop}
\begin{proof}
	See \cite[\S 5.1]{new-tho}.
\end{proof}

We now discuss some Hecke operators at (possibly ramified) places
in $\overline{S}_p$. Assume that each prime $\bar{v}$ of $F^+$ above $p$ splits 
in $F$. 
Let $\bar{v}\in \overline{S}_p$, and recall that $\tilde{v}$ is a chosen prime
of $F$ above it. For integers $c\geq b\geq 0$, we define subgroups 
\[
\cP_{\bar{v}}(b,c)\subset 
\widetilde{G}(\cO_{F^+_{\bar{v}}}) = \GL_{2n}(\cO_{F_{\tilde{v}}})
\] 
which reduce to 
block upper-triangular matrices (with two $n\times n$ blocks) modulo
$\varpi_{\tilde{v}}^c$ and to block unipotent matrices modulo $\varpi_{\tilde{v}}^b$. 
We set $\cP_{\bar{v}} = \cP_{\bar{v}}(0,1)$, which is identified with the standard parahoric subgroup $\cP_{n,n}$ of $\GL_{2n}$. For each parabolic subgroup $Q_{\bar{v}}$ of $P_{\bar{v}}$ which contains $\widetilde{B}_{\bar{v}}$ we have an associated parahoric subgroup $\cQ_{\bar{v}} \subset \cP_{\bar{v}}$. We note that these subgroups all admit an Iwahori decomposition 
with respect to $P_{\bar{v}}$, and therefore the formalism of~\cite[\S 2.1.9]{10author} applies when we consider the Hecke algebras of monoids. 

Write 
$\tilde{u}_{\tilde{v},n}:=\diag(\varpi_{\tilde{v}}, 
\dots, \varpi_{\tilde{v}}, 1, \dots, 1)\in \GL_{2n}(F_{\tilde{v}})$,
where $\varpi_{\tilde{v}}$ appears exactly $n$ times on the diagonal. 
If $c\geq 1$, we write $\widetilde{U}_{\tilde{v},n}\in 
\cH\left(\widetilde{G}(\cO_{F^+_{\bar{v}}}), \cP_{\bar{v}}(b,c)\right)$ for the 
double coset operator
$\widetilde{U}_{\tilde{v},n} = [ \cP_{\bar{v}}(b, c) 
\iota_{\tilde{v}}^{-1}\tilde{u}_{\tilde{v},n} \cP_{\bar{v}}(b, c) ]$. Also write 
$\tilde{u}_{\tilde{v},2n}:= \diag(\varpi_{\tilde{v}},\dots, \varpi_{\tilde{v}})
\in \GL_{2n}(F_{\tilde{v}})$ 
and denote by $\widetilde{U}_{\tilde{v}, 2n}\in 
\cH\left(\widetilde{G}(\cO_{F^+_{\bar{v}}}), \cP_{\bar{v}}(b,c)\right)$ 
the corresponding double coset operator. Note
that these depend on both the choice of uniformizer $\varpi_{\tilde{v}}$ and the chosen level. We write $\widetilde{\Delta}_{\bar{v}}\subset \widetilde{G}(F^+_{\bar{v}})$
for the subset 
\[
\widetilde{\Delta}_{\bar{v}} := \iota_{\tilde{v}}^{-1}\left(\sqcup_{\mu_1\in \Z_+}\sqcup_{\mu_2\in \Z} \cP_{n,n} 
(\tilde{u}_{\tilde{v},n})^{\mu_1}(\tilde{u}_{\tilde{v},2n})^{\mu_2}\cP_{n,n}\right),
\] 
which is independent of the choice of $\tilde{v}\mid \bar{v}$. 

Considering cohomology at level $\cP_{\bar{v}}$ and the ordinary subspace for the Hecke operator $\widetilde{U}_{\tilde{v},n}$ will be most important for us. However, we will work a little more generally to allow us to keep track of additional Hecke operators at $\bar{v}$ and prove a local--global compatibility result for ordinary as well as crystalline representations.

So, more generally, we suppose we have a parabolic subgroup $\widetilde{B}_{\bar{v}} \subset Q_{\bar{v}} \subset P_{\bar{v}}$ corresponding to a subset $I \subset \Delta$ of the simple roots, and with Levi decomposition $Q_{\bar{v}} = M_{Q_{\bar{v}}}N_{Q_{\bar{v}}}$ compatible with the decomposition $P_{\bar{v}} = G_{\bar{v}}U_{\bar{v}}$. We consider the monoid of cocharacters \[X_{Q_{\bar{v}}} := \{\nu \in X_*(Z(M_{Q_{\bar{v}}})) : \langle \nu,\delta \rangle \ge 0 \text{ for all } \delta \in \Delta - I \}.\] In fact, this is simply the subset of $\widetilde{B}_{\bar{v}}$-dominant cocharacters in  $X_*(Z(M_{Q_{\bar{v}}}))$.  We then define a subset $\widetilde{\Delta}^{\cQ_{\bar{v}}}_{\bar{v}} \subset \widetilde{G}(F^+_{\bar{v}})$ containing the parahoric subgroup $\cQ_{\bar{v}} \subset \cP_{\bar{v}}$ by 
\[
\widetilde{\Delta}^{\cQ_{\bar{v}}}_{\bar{v}} := \coprod_{\nu \in X_{Q_{\bar{v}}}}\cQ_{\bar{v}} 
\nu(\varpi_{\bar{v}})\cQ_{\bar{v}}.
\]
We have $\widetilde{\Delta}^{\cP_{\bar{v}}}_{\bar{v}} =\widetilde{\Delta}_{\bar{v}}$.

We set $\Delta^{\cQ_{\bar{v}},+}_{\bar{v}} := \widetilde{\Delta}^{\cQ_{\bar{v}}}_{\bar{v}}\cap G_{F^+_{\bar{v}}}$ and $\Delta^{\cQ_{\bar{v}}}_{\bar{v}} := \Delta^{\cQ_{\bar{v}},+}_{\bar{v}}[\iota_{\tilde{v}}^{-1}(\tilde{u}_{\tilde{v},n}^{-1})]$ (the submonoid of $G(F^+_{\bar{v}})$ generated by $\Delta^{\cQ_{\bar{v}},+}_{\bar{v}}$ and its central element $\iota_{\tilde{v}}^{-1}(\tilde{u}_{\tilde{v},n}^{-1})$).

\begin{lemma}\label{lem:Q Hecke algebra}\begin{enumerate}
\item For $\nu \in X_{Q_{\bar{v}}}$, the element $\nu(\varpi_{\bar{v}})$ is \emph{$\cQ_{\bar{v}}$-positive}; i.e.~we have \begin{align*}\nu(\varpi_{\bar{v}})\left(N_{Q_{\bar{v}}} \cap \cQ_{\bar{v}}\right)\nu(\varpi_{\bar{v}})^{-1} & \subset N_{Q_{\bar{v}}} \cap \cQ_{\bar{v}} \\\text{ and } \nu(\varpi_{\bar{v}})^{-1}\left(\overline{N}_{Q_{\bar{v}}} \cap \cQ_{\bar{v}}\right)\nu(\varpi_{\bar{v}}) &\subset \overline{N}_{Q_{\bar{v}}} \cap \cQ_{\bar{v}}.\end{align*}
\item $\widetilde{\Delta}^{\cQ_{\bar{v}}}_{\bar{v}}$ is a monoid under multiplication.
\item The map $[(M_{Q_{\bar{v}}}\cap\cQ_{\bar{v}}) \nu(\varpi_{\bar{v}})(M_{Q_{\bar{v}}}\cap\cQ_{\bar{v}})] \mapsto [\cQ_{\bar{v}}\nu(\varpi_{\bar{v}})\cQ_{\bar{v}}]$ defines a ring isomorphism of Hecke algebras \[\cH(M_{Q_{\bar{v}}}\cap\widetilde{\Delta}^{\cQ_{\bar{v}}}_{\bar{v}},M_{Q_{\bar{v}}}\cap\cQ_{\bar{v}}) \toisom \cH(\widetilde{\Delta}^{\cQ_{\bar{v}}}_{\bar{v}},\cQ_{\bar{v}})\] which also factors through an isomorphism to $\cH(\Delta^{\cQ_{\bar{v}},+}_{\bar{v}},G(F^+_{\bar{v}})\cap\cQ_{\bar{v}})$.
\end{enumerate}
\end{lemma}
\begin{proof}
	The first part can be checked directly, or using root groups. The second part follows from the first, using the Iwahori decomposition of $\cQ_{\bar{v}}$. The third part is \cite[Corollary 6.12]{MR1643417}.
\end{proof}

\begin{remark}
	Our monoids are usually strictly contained in those defined in \cite[\S 2.1.9]{10author}. We only need to consider Hecke operators supported on double cosets of central elements in the Levi subgroup, which in particular implies (as shown in the preceding lemma) that the Hecke algebra $\cH(\widetilde{\Delta}^{\cQ_{\bar{v}}}_{\bar{v}},\cQ_{\bar{v}})$ is commutative.
\end{remark}

Fix $\bar{v}\in \bar{S}_p$
and $\tau\in \Hom(F^+_{\bar{v}}, E)$.  Let $\tilde{\lambda}\in \left(\Z^{2n}_+\right)^{\Hom(F^+,E)}$.
We define a character $\tilde{\alpha}^{\cQ_{\bar{v}}}_{\tilde{\lambda}_{\tau}}: \widetilde{\Delta}^{\cQ_{\bar{v}}}_{\bar{v}}\to E^\times$ 
by setting 
\[
\tilde{\alpha}^{\cQ_{\bar{v}}}_{\tilde{\lambda}_{\tau}}\left(\nu(\varpi_{\bar{v}})\right) = 
\tau(\varpi_{\bar{v}})^{\langle \nu, w_0^{\widetilde{G}}\tilde{\lambda}_\tau\rangle}
\] and setting the character to be trivial on $\cQ_{\bar{v}}$.

\begin{lemma}\label{lem:twisted action} Fix $\bar{v}\in \bar{S}_p$
	and $\tau\in \Hom(F^+_{\bar{v}}, E)$.  
	Define an action of $\cO[\widetilde{\Delta}^{\cQ_{\bar{v}}}_{\bar{v}}]$ on 
	$V_{\tilde{\lambda}_{\tau}}$ by 
	\begin{equation}\label{eq:twisted action}
	g\cdot_{\tilde{\lambda}_{\tau}}^{\cQ_{\bar{v}}} x:= \tilde{\alpha}^{\cQ_{\bar{v}}}_{\tilde{\lambda}_\tau}(g)^{-1} g \cdot x,
	\end{equation}
	where $g\cdot x$ is the usual action of $g\in \widetilde{\Delta}^{\cQ_{\bar{v}}}_{\bar{v}}\subset 
	\widetilde{G}(F^+_{\bar{v}})$ 
	on $x\in V_{\tilde{\lambda}_{\tau}}$. The lattice 
	$\cV_{\tilde{\lambda}_{\tau}}$ is
	stable under the $\cdot_{\tilde{\lambda}_{\tau}}^{\cQ_{\bar{v}}}$-action of $\cO[\widetilde{\Delta}^{\cQ_{\bar{v}}}_{\bar{v}}]$. 
\end{lemma}
\begin{proof} This follows from the fact that the rescaled action of $\widetilde{\Delta}^{\cQ_{\bar{v}}}_{\bar{v}}$ stabilizes each weight space in $\cV_{\tilde{\lambda}_{\tau}}$, which has lowest weight $w_{0,2n}\tilde{\lambda}_\tau$. Cf.~\cite[Definition 2.8]{geraghty}.
\end{proof}
Suppose we have a subset $\bar{S}\subseteq \bar{S}_p$ and standard parabolic subgroups $Q_{\bar{v}} \subset P_{\bar{v}}$ for each $\bar{v}\in \bar{S}$. Then we set $\widetilde{\Delta}^{\cQ_{\bar{S}}}_{\bar{S}} := \prod_{\bar{v}\in \bar{S}}\widetilde{\Delta}^{\cQ_{\bar{v}}}_{\bar{v}}$, with similar notation for $\Delta$. Let $\tilde{\lambda}\in \left(\Z^{2n}_+\right)^{\Hom(F^+,E)}$. If we omit the superscript $\cQ_{\bar{S}}$, we take $Q_{\bar{v}} = P_{\bar{v}}$ for all $\bar{v} \in \bar{S}$.

As a consequence of Lemma \ref{lem:twisted action}, we have constructed a twisted action of $\cO[\widetilde{\Delta}^{\cQ_{\bar{S}}}_{\bar{S}}]$ on $V_{\tilde{\lambda}}$ which stabilizes the lattice $\cV_{\tilde{\lambda}}$. This action is obtained from the usual action by rescaling with the inverse of the character
\[\tilde{\alpha}^{\cQ_{\bar{S}}}_{\tilde{\lambda}} := \prod_{\bar{v}\in\bar{S}}\prod_{\tau \in \Hom(F^+_{\bar{v}},E)}\tilde{\alpha}^{\cQ_{\bar{v}}}_{\tilde{\lambda}_\tau}.\]

We construct a similar rescaled action for the Levi subgroup $G$. Suppose $\lambda \in \left(\Z^{n}\right)^{\Hom(F,E)}$ is a dominant weight for $G$. Recall that we have identified $\lambda$ with a (not necessarily dominant) weight $\tilde{\lambda}$ of $\widetilde{G}$. We define a character $\alpha_{\lambda}^{\cQ_{\bar{v}}}:  \Delta^{\cQ_{\bar{v}}}_{\bar{v}} \to E^\times$ using the formula \[\alpha_{\lambda}^{\cQ_{\bar{v}}}(\nu(\varpi_{\bar{v}})) = \prod_{\tau\in\Hom(F^+_{\bar{v}},E)} \tau(\varpi_{\bar{v}})^{\langle \nu, w_0^G\tilde{\lambda}_{\tau}\rangle}\] and a rescaled action of $\Delta^{\cQ_{\bar{v}}}_{\bar{v}}$ on $\cV_{\lambda}$ by $g\cdot_{\lambda} x = \alpha_{\lambda}^{\cQ_{\bar{v}}}(g)^{-1} g\cdot x$. Note that the rescaling means that $\iota^{-1}_{\tilde{v}}(\tilde{u}_{\tilde{v},n})\cdot_{\lambda} x = x$.

Let $T\supseteq S_p$ be a finite set of finite places of $F$ that
satisfies $T = T^c$. The formalism 
of~\cite[\S 2.1.9]{10author} implies then that, for each good subgroup 
$\tK\subset \tG(\A_{F^+,f})$ such that $\tK_{\bar{v}} = \cP_{\bar{v}}(b,c)$ with $c\geq 1$
for each $\bar{v}\in \bar{S}\subseteq \bar{S}_p$, there is a canonical homomorphism 
\begin{equation}\label{eq:ord Hecke action Siegel}
\cH(\tG^T, \tK^T)\otimes_{\Z}\cH(\widetilde{\Delta}_{\bar{S}}, \tK_{\bar{S}})\to 
\mathrm{End}_{D^+(\cO)}\left(R\Gamma\left(\tX_{\tK}, \cV_{\tilde{\lambda}}\right)\right)
\end{equation}
and in particular all the Hecke operators $\widetilde{U}_{\tilde{v},n}$
and $\widetilde{U}_{\tilde{v},2n}$ for $\bar{v}\in \bar{S}$ 
act as endomorphisms of $R\Gamma\left(\tX_{\tK}, \cV_{\tilde{\lambda}}\right)$. 
In fact, the Hecke operators $\widetilde{U}_{\tilde{v}, 2n}$ for $\bar{v}\in \bar{S}$ 
act as automorphisms of $R\Gamma\left(\tX_{\tK}, \cV_{\tilde{\lambda}}\right)$
because the elements $\tilde{u}_{\tilde{v},2n}$ are central. 

We will also consider $\widetilde{K}$ with $\widetilde{K}_{\bar{v}} = \cQ_{\bar{v}}$ for each $\bar{v} \in \overline{S}$, and then we have a Hecke action 
\begin{equation}\label{eq:ord Hecke action}
\cH(\tG^T, \tK^T)\otimes_{\Z}\cH(\widetilde{\Delta}^{\cQ_{\bar{S}}}_{\bar{S}}, \tK_{\bar{S}})\to 
\mathrm{End}_{D^+(\cO)}\left(R\Gamma\left(\tX_{\tK}, \cV_{\tilde{\lambda}}\right)\right).
\end{equation} 

Similarly, for the Levi subgroup $G$, %
let $K \subset G(\AA_{F^+,f})$ be a good subgroup with $K_{\bar{v}} = \cQ_{\bar{v}} \cap G(F^+_{\bar{v}})$ for $\bar{v} \in \bar{S}$. 
Let $\{\lambda_{\tau}\}_{\tau\in \Hom(F^+_{\bar{v}},E)}$ be sets of dominant weights for $G$ at primes $\bar{v}\in \bar{S}$, giving rise to an 
$\cO[K_{\bar{S}}]$-module $\cV_{\lambda_{\bar{S}}}$. Let $\cV$ be an $\cO[K_{\bar{S}_p\setminus \bar{S}}]$-module,  
which is finite free as an $\cO$-module and such that $\cV/\varpi^m$ is a smooth $\cO/\varpi^m[K_{\bar{S}_p\setminus \bar{S}}]$-module 
for every $m\in \Z_{\geq 1}$. 
We get a Hecke action
\begin{equation}\label{eq:ord Hecke action Levi}
\cH(G^T, K^T)\otimes_{\Z}\cH({\Delta}^{\cQ_{\bar{S}}}_{\bar{S}}, K_{\bar{S}})\to 
\mathrm{End}_{D^+(\cO)}\left(R\Gamma\left(X_{K}, \cV\otimes_{\cO} \cV_{{\lambda_{\bar{S}}}}\right)\right).
\end{equation}
This Hecke action generalizes in the natural way 
to the case when $\cV$ is a complex of $\cO[K_{\bar{S}_p\setminus\bar{S}}]$-modules as above. 
For $v|\bar{v}$ we will be interested in the (invertible) Hecke operator $U_v$ corresponding to the central element $\iota_v^{-1}(\tilde{u}_{{v},n}^{-1}\tilde{u}_{v,2n})$. Under our identification of $G$ with $\mathrm{Res}_{F/F^+}\GL_n$, this element is $\diag(\varpi_v,\ldots,\varpi_v)$.

\subsubsection{Automorphic Galois representations and 
middle degree cohomology} We start by recalling
some well-known results about Galois representations 
associated to automorphic representations, and more generally
to systems of Hecke eigenvalues occurring in the cohomology 
of locally symmetric spaces with integral coefficients. 

\begin{thm}\label{thm:automorphic Galois reps unitary} Assume
that $F$ contains an imaginary quadratic field and that $\pi$ is
a cuspidal automorphic representation of $\tG(\A_{F^+})$ that
is $\xi$-cohomological\footnote{As in \cite{shin-basechange}, $\xi$-cohomological means that $\pi_\infty\otimes\xi$ has non-zero $(\mathfrak{g},K_\infty)$-cohomology.} for some irreducible algebraic representation
$\xi$ of $(\mathrm{Res}_{F^+/\Q}\tG)_{\C}$. For any isomorphism $\iota: \overline{\Q}_p\toisom \C$,
there exists a continuous, semisimple Galois representation 
\[
r_{\iota}(\pi):G_F\to \GL_{2n}(\overline{\Q}_p)
\] 
satisfying the following conditions:

\begin{enumerate} %
\item For each prime $\ell\not=p$ which is unramified in $F$ 
and above which $\pi$ is unramified, and for each prime $v\mid \ell$ of $F$, 
$r_{\iota}(\pi)|_{G_{F_v}}$ is unramified and the characteristic 
polynomial of $r_{\iota}(\pi)(\Frob_v)$ is equal to the image of $\widetilde{P}_v(X)$
in $\overline{\Q}_p[X]$ corresponding to the base change of $\iota^{-1}(\pi_v)$. 

\item For each prime $v\mid p$ of $F$, $r_{\iota}(\pi)$ is de Rham, and for each 
$\tau:F\hookrightarrow \overline{\Q}_p$, we have 
\[
\HT_{\tau}(r_{\iota}(\pi)) = 	\{\tilde{\lambda}_{\tau,1} + 2n-1,\tilde{\lambda}_{\tau,2} + 2n-2,
\dots,\tilde{\lambda}_{\tau, 2n}\},
\]
where $\tilde{\lambda}\in (\Z^{2n}_+)^{\Hom(F, \overline{\Q}_p)}$ is the highest weight 
of the representation $\iota^{-1}(\xi\otimes \xi)^\vee$ of $(\mathrm{Res}_{F/\Q}\GL_{2n})_{\overline{\Q}_p}$.\footnote{For each $\overline{\tau}:F^+\to \C$, $\xi$ gives a representation of $\tG_{\overline{\tau}}$ and hence for $\tau,\tau c$ extending $\overline{\tau}$ to $F$ we have a representation $\xi\otimes\xi$ of $\tG_{\overline{\tau}}\times \tG_{\overline{\tau}} = (\GL_{2n,F})_{\tau} \times (\GL_{2n,F})_{\tau c}$. Note that $r_{\iota}(\pi)^{\vee,c}(1-2n)\cong r_{\iota}(\pi)$, so $\HT_{\tau}(r_{\iota}(\pi))$ and $\HT_{\tau c}(r_{\iota}(\pi))$ can be read off from each other.}
\item If $F_0\subset F$ is an imaginary quadratic field and $\ell$ is a prime
which splits in $F_0$ (including possibly $\ell=p$), then for each prime $v\mid \ell$ of $F$ lying above 
a prime $\bar{v}$ of $F^+$, there is an isomorphism 
\[
\mathrm{WD}\left(r_{\iota}(\pi)|_{G_{F_v}}\right)^{F-ss} \simeq \mathrm{rec}^T_{F_v}
(\pi_{\bar{v}}\circ \iota_v). 
\]
\end{enumerate}
\end{thm}%

\begin{proof} This is~\cite[Theorem 2.3.3]{10author}. We mention that it 
relies on the base change result of~\cite{shin-basechange} and on the existence and properties
of the Galois representations associated to regular algebraic, conjugate self-dual 
cuspidal automorphic representations of $\GL_m$. 
\end{proof}

\begin{thm}\label{thm:torsion Galois reps residual} Let $\m\subset \mathbb{T}^T(K,\lambda)$ 
be a maximal ideal.
Suppose $F$ contains an imaginary quadratic field, the finite set of finite places $T$ of $F$ is stable under complex conjugation, and the following condition is satisfied:
\begin{itemize}\label{item:condition on S}
\item \IQFassm\footnote{This condition
can always be realized after enlarging $T$ and is used to ensure that the results of~\cite{scholze-torsion} 
that we appeal to are unconditional.} 
\end{itemize}
Then there exists a continuous, semi-simple Galois representation 
\[
\bar{\rho}_{\m}:G_{F,T}\to \GL_n\left(\mathbb{T}^T(K,\lambda)/\m\right)
\]
such that, for each finite place $v\not\in T$ of $F$, the characteristic polynomial
of $\bar{\rho}_{\m}(\Frob_v)$ is equal to the image of $P_v(X)$ in 
$\left(\mathbb{T}^T(K,\lambda)/\m\right)[X]$. 
\end{thm}

\begin{proof} This is~\cite[Theorem 2.3.5]{10author}: it essentially follows 
from~\cite[Corollary 5.4.3]{scholze-torsion}. 
\end{proof}

\begin{lemma}\label{lem:central character mod p} Let $\m \subset 
\mathbb{T}^T(K,\lambda)$ be a maximal ideal as in Theorem~\ref{thm:torsion Galois reps residual}.
Suppose $k = k(\m)$. Let $v$ be a $p$-adic place of $F$. The Hecke operator $U_v$ 
has a unique eigenvalue on $H^*(X_K, \cV_{\lambda}/\varpi)_{\m}$, equal to $\bar{\epsilon}_p^{\frac{n(n-1)}{2}}(\Art_{F_v}(\varpi_v)) \cdot \det \bar{\rho}_{\m}(\Art_{F_v}(\varpi_v))$. 
\end{lemma}

\begin{proof}
Our proof of this will be global, computing the action of $U_v$ in terms of a central Hecke operator at a suitable unramified prime. We write $Z_n$ for the centre of $\GL_{n,F}$. We have a right action of $Z_n(\AA_{F})$ on $X_K$ (by right multiplication of the finite adelic part on $\GL_n(\AA_{F,f})$ and of the archimedean part on $X^{\GL_{n,F}}$). We rescale the action of $Z_n(\prod_{w\in S_p}F_w)$, as in \S\ref{sec:explicit Hecke operators}, so that it stabilizes $\cV_\lambda$. Explicitly, this means that the elements $\diag(\varpi_w,\ldots,\varpi_w)$ act trivially (recall that we have fixed uniformizers at every finite place), whilst $Z_n(\cO_{F,p})$ acts via its natural action. This allows us to define an action of $Z_n(\AA_{F})$ on $\cV_{\lambda}$ which factors through its $p$-adic part and is compatible with its existing $K$ action. We obtain an action of $Z_n(\AA_{F})$ on $H^*(X_K, \cV_{\lambda}/\varpi)_{\m}$ which factors through the quotient $Z_n(\AA_{F})/F_\infty^\times (Z_n(\AA_{F,f})\cap K)$ for continuity reasons\footnote{In fact we do not even need continuity of the $F_\infty^\times$ action. The cohomology groups are finite, so the action of $F_\infty^\times$ gives a homomorphism from a product of copies of $\C^\times$ to a finite group. This is necessarily trivial, since $\C^\times$ has no proper finite index subgroups.}. 

We also have a continuous character $\psi_\m: \AA_F^\times = Z_n(\AA_F) \to k^\times$ determined by \[\psi_\m = \bar{\epsilon}_p^{\frac{n(n-1)}{2}}\det \bar{\rho}_{\m}\circ\Art_F.\]

Since $F_\infty^\times$ is connected, and $\bar{\rho}_{\m}$ is unramified away from $T$, $\psi_\m$ factors through the quotient $F^\times \backslash Z_n(\AA_F)/F_\infty^\times K_Z$ for a compact open subgroup $K_Z = \prod_w K_{Z,w}$ of $Z_n(\AA_{F,f})$ with $K_{Z,w} = \cO_{F_w}^\times$ for $w \notin T$. Shrinking $K_Z$ if necessary, we assume that $K_Z \subset Z_n(\AA_{F,f})\cap K$ and that $K_Z$ acts trivially on $\cV_\lambda/\varpi$. Note that for a finite place $w \notin T$, $\psi_\m(\Frob_w)$ is equal to $T_{w,n} \mod \m$. 

By Chebotarev density, we can find a place $w \notin T$ such that the uniformizer $\varpi_w$ and $\varpi_v$ map to the same element in the ray class group $F^\times \backslash Z_n(\AA_{F})/F_\infty^\times K_Z$. We write $\varpi_v = \alpha\varpi_w z_\infty k$, for $\alpha \in F^\times, z_\infty \in F_\infty^\times, k \in K_Z$. Since $z_\infty, k, \varpi_v, \varpi_w$ all act trivially on $\cV_\lambda/\varpi$, so does $\alpha$.

Now we can compare the actions of $z_\infty \varpi_w$ and $\varpi_v$ on $H^*(X_K, \cV_{\lambda}/\varpi)_{\m}$. Their ratio is equal to $\alpha k$. Since $k \in K$, it acts trivially. The principal id\`{e}le $\alpha$ acts trivially on both $\cV_\lambda/\varpi$ and the cover of $X_K$ trivializing $\cV_\lambda/\varpi$. So $\alpha$ also acts trivially on cohomology. Since the action of $F_\infty^\times$ on cohomology is trivial, we deduce that $\varpi_w$ and $\varpi_v$ have the same action on $H^*(X_K, \cV_{\lambda}/\varpi)_{\m}$. The unique eigenvalue of $\varpi_w$ (i.e.~of the Hecke operator $T_{w,n}$) on $H^*(X_K, \cV_{\lambda}/\varpi)_{\m}$ is $\psi_\m(\varpi_w)$. Our choice of $w$ means this is equal to $\psi_{\m}(\varpi_v)$, and we are done.
\end{proof}
\begin{rem}
	A more conceptual proof for Lemma \ref{lem:central character mod p} can be given by arguing with (mod $\varpi$) completed cohomology at level $K^p$, localized at $\m$. We can then assume $\lambda$ is trivial, in which case completed cohomology is equipped with a continuous action of $F^\times\backslash\AA_{F}^\times/F_\infty^\times(\AA_{F,f}\cap K^p)$ with a unique system of eigenvalues corresponding to the character $\bar{\epsilon}_p^{\frac{n(n-1)}{2}}\cdot \det \bar{\rho}_{\m}$.
\end{rem}

\begin{defn} We say that a maximal ideal $\m\subset \mathbb{T}^T(K,\lambda)$
is \emph{non-Eisenstein} if $\bar{\rho}_{\m}$ is absolutely irreducible. 
\end{defn}
Our convention is that when we ask for a maximal ideal $\m$ to be non-Eisenstein, we are implicitly imposing the assumptions of Theorem \ref{thm:torsion Galois reps residual}.

\begin{thm}\label{thm:torsion Galois reps Hecke} 
Let $\m\subset \mathbb{T}^T(K,\lambda)$ be a non-Eisenstein maximal ideal. There exist an integer $N\geq 1$, 
which depends only on $n$ and $[F:\Q]$, an ideal $I\subset \mathbb{T}^T(K,\lambda)_{\m}$
satisfying $I^N=0$, and a continuous homomorphism 
\[
\rho_{\m}: G_{F,T} \to \GL_n\left(\mathbb{T}^T(K,\lambda)_{\m}/I\right)
\]
such that, for each finite place $v\not\in T$ of $F$, the characteristic polynomial
of $\rho_{\m}(\Frob_v)$ is equal to the image of $P_v(X)$ in 
$\left(\mathbb{T}^T(K,\lambda)_{\m}/I\right)[X]$. 
\end{thm}

\begin{proof} This is~\cite[Corollary 5.4.4]{scholze-torsion}. 
\end{proof}

Let $\bar{S}\subseteq \bar{S}_p$ and let 
$\tK\subset \tG(\A_{F^+,f})$ be a compact open subgroup, which is 
decomposed with respect to the Levi decomposition $P = G\ltimes U$, 
and which satisfies $\tK_{\bar{v}} = \cP_{\bar{v}}(b,c)$ or $\cQ_{\bar{v}}$ for all 
$\bar{v}\in \bar{S}$. 
We set $K:=\tK\cap G(\A_{F^+,f})$ and $K_{U}:=\tK\cap U(\A_{F^+,f})$. 
Let $\m\subset \mathbb{T}^T(K,\lambda)$ be a non-Eisenstein maximal 
ideal and let $\widetilde{\m}\subset \widetilde{\mathbb{T}}^T$ denote its pullback 
under the unnormalized Satake transform $\widetilde{\mathbb{T}}^T \to \mathbb{T}^T$. 
Recall that the boundary $\partial \widetilde{X}_{\tK}$ of the Borel--Serre 
compactification of $\widetilde{X}_{\tK}$ has a $\tG(\A_{F^+,f})$-equivariant
stratification indexed by the rational parabolic subgroups of $\tG$ which contain $\tB$. 
See~\cite[\S 3.1.2]{new-tho}, especially~\cite[Lemma 3.10]{new-tho} for more details. 
For such a standard parabolic 
subgroup $Q$, we denote by $\widetilde{X}^Q_{\widetilde{K}}$ the stratum labelled
by $Q$. This stratum can be written as a double quotient: 
\[
\widetilde{X}^Q_{\widetilde{K}}= Q(F^+)\backslash X^Q \times \tG(\A_{F^+,f})/\widetilde{K}
\]
By applying the formalism in \S~\ref{sec:Hecke formalism}, 
there is, for any $\tilde{\lambda} \in (\Z^{2n}_+)^{\Hom(F^+,E)}$,
a homomorphism 
\[
\widetilde{\mathbb{T}}^T \to \mathrm{End}_{D^+(\cO)}
\left(R\Gamma\left(\tX^Q_{\tK}, \cV_{\tilde{\lambda}}\right)\right). 
\]
Therefore, we can define the localization 
$R\Gamma\left(\tX^Q_{\tK}, \cV_{\tilde{\lambda}}\right)_{\widetilde{\m}}$. 

\begin{thm}\label{thm:boundary only Siegel parabolic} Let $\m\subset
\mathbb{\mathbb{T}}^T(K,\lambda)$ be a non-Eisenstein maximal ideal and let 
$\widetilde{\m} := \cS^*(\m)\subset \widetilde{\mathbb{T}}^T$. Let 
$\tilde{\lambda}\in (\Z^{2n}_+)^{\Hom(F^+,E)}$. 
Pullback along the natural inclusion induces a 
$\widetilde{\mathbb{T}}^T$-equivariant isomorphism in $D^+(\cO)$: 
\[
R\Gamma\left(\partial \tX_{\tK},\cV_{\tilde{\lambda}}\right)_{\widetilde{\m}}
\toisom R\Gamma\left(\tX^P_{\tK}, \cV_{\tilde{\lambda}}\right)_{\widetilde{\m}}.
\] 
\end{thm}

\begin{proof} This is~\cite[Thm. 2.4.2]{10author}.
\end{proof}

\begin{thm}\label{thm:torsion Galois reps residual Shimura}
Let $\widetilde{\m}\subset 
\widetilde{\mathbb{T}}^T(\tK,\tilde{\lambda})$
be a maximal ideal. 
Suppose $F$ contains an imaginary quadratic field, 
the finite set of finite places $T$ of $F$ is stable under 
complex conjugation, and the following condition is satisfied:
\begin{itemize}
\item \IQFassm
\end{itemize}

\noindent Then there exists a continuous, semi-simple Galois representation 
\[
\bar{\rho}_{\widetilde{\m}}:G_{F,T}\to \GL_{2n}
\left(\widetilde{\mathbb{T}}^T(\tK,\tilde{\lambda})/\widetilde{\m}\right)
\]
such that, for each finite place $v\not\in T$ of $F$, the characteristic polynomial
of $\bar{\rho}_{\widetilde{\m}}(\Frob_v)$ is equal to the image of $\widetilde{P}_v(X)$ in 
$\left(\widetilde{\mathbb{T}}^T(\tK,\tilde{\lambda})/\widetilde{\m}\right)[X]$. 

\end{thm}

\begin{proof} We will use the Hecke algebra 
$\mathbb{T}^T_{cl}$ defined in~\cite[\S 6.5]{arizona}. 
By reducing to $\F$-coefficients and increasing the level 
to sufficiently small compact opens $\tK'_{\bar{v}}$ at all 
primes $\bar{v}\in \bar{S}_p$, we can trivialize the local system 
$\cV_{\tilde{\lambda}}\otimes_{\cO}\F$. The proof of~\cite[Theorem 6.5.3]{arizona}
shows that the map 
\[
\mathbb{T}^T_{cl}\to \mathrm{End}_{D^+\left(\prod_{\bar{v}\in \bar{S}_p}
\tG(\cO_{F^+_{\bar{v}}})/\tK'_{\bar{v}},\F\right)}\left(R\Gamma(\tX_{\tK'},\F)_{\widetilde{\m}}
\otimes_{\F} \left(\cV_{\tilde{\lambda}}\otimes_{\cO}\F\right) \right)
\]
is continuous for the discrete topology on the target. By the Hochschild--Serre spectral 
sequence, this implies
that the map 
\[
\mathbb{T}^T_{cl}\to \mathrm{End}_{D^+(\F)}\left(R\Gamma\left(\tX_{\tK},
\cV_{\tilde{\lambda}}\otimes_{\cO}\F\right)_{\widetilde{\m}}\right)
\]
is also continuous for the discrete topology on the target. 
The existence of a determinant 
valued in $\widetilde{\mathbb{T}}^T(\tK,\tilde{\lambda})/\widetilde{\m}$ now
follows 
from~\cite[Lemma 6.5.2]{arizona}, which is the version of~\cite[Corollary 5.1.11]{scholze-torsion}
for usual cohomology and which can be made unconditional by using 
Theorem~\ref{thm:automorphic Galois reps unitary} 
as an input. This 
determinant corresponds to a semi-simple $\overline{\F}_p$-valued representation by~\cite[Theorem A]{chenevier_det}. 
This representation can be realized over $\widetilde{\mathbb{T}}^T(\tK,\tilde{\lambda})/\widetilde{\m}$ 
by the same argument as in the proof of~\cite[Theorem 2.3.5]{10author}. 
\end{proof}

We now introduce the key technical condition that our residual
representations must satisfy in order to appeal to the main
result of~\cite{caraiani-scholze-noncompact}.  

\begin{defn}\label{defn:generic} A continuous representation $\bar{\rho}: G_{F}\to 
\GL_m(\F)$ is \emph{decomposed generic}\footnote{This is slightly weaker than the condition called
\emph{decomposed generic} in~\cite{caraiani-scholze-compact}. 
See~\cite[Remark 1.4]{caraiani-scholze-noncompact}
and~\cite[Corollary 5.1.3]{caraiani-scholze-noncompact} for an explanation.} if there exists a prime $\ell\not = p$
such that the following are satisfied: 
\begin{enumerate}
\item the prime $\ell$ splits completely in $F$; 
\item for every prime $v\mid \ell$ of $F$, the representation
$\bar{\rho}|_{G_{F_v}}$ is unramified and the eigenvalues
$\alpha_1,\dots,\alpha_m$ of $\bar{\rho}(\Frob_v)$ satisfy 
$\alpha_i/\alpha_j \not = \ell$ for $i\not =j$. 
\end{enumerate}
\end{defn}

\noindent We note that if a representation $\bar{\rho}$ is decomposed generic,
then by the Chebotarev density theorem there exist infinitely
many primes $\ell\not = p$ as in Definition~\ref{defn:generic}, 
cf.~\cite[Lemma 4.3.2]{10author}. 

\begin{thm}\label{thm:middle degree cohomology} 
Keep the same assumptions on $F$ and $T$ as in Theorem~\ref{thm:torsion Galois reps
residual Shimura}. 
Let $\widetilde{\m}\subset 
\widetilde{\mathbb{T}}^T(\tK,\tilde{\lambda})$
be a maximal ideal such that the associated Galois representation
$\bar{\rho}_{\widetilde{\m}}$ constructed in Theorem~\ref{thm:torsion Galois reps
residual Shimura} 
is decomposed generic, in the sense of Definition~\ref{defn:generic}. Recall 
that $d=\dim_{\C}\tX_{\tK}$. Then 
we have $\widetilde{\mathbb{T}}^T$-equivariant morphisms 
\[
H^d\left(\tX_{\tK}, \cV_{\tilde{\lambda}}[1/p]\right)_{\widetilde{\m}}
\hookleftarrow 
H^d\left(\tX_{\tK},\cV_{\tilde{\lambda}}\right)_{\widetilde{\m}}\twoheadrightarrow 
H^d\left(\partial \tX_{\tK},\cV_{\tilde{\lambda}}\right)_{\widetilde{\m}}. 
\]
\end{thm}%

\begin{proof} This follows from~\cite[Thm. 1.1]{caraiani-scholze-noncompact}
as in the proof of~\cite[Thm. 4.3.3]{10author}. Moreover, we can remove the technical hypotheses that $[F^+:\Q]>1$ and $\bar{\rho}_{\widetilde{\m}}$ has length at most two by appealing to Koshikawa's work \cite[Theorem 1.3]{koshikawa}. 
\end{proof}
\subsection{$P$-ordinary Hida theory}\label{sec:P-ord Hida}
In this section, we develop a 
$P$-ordinary version of Hida theory for the group $\tG$ and 
the Betti cohomology of the locally
symmetric spaces $\tX_{\tK}$, 
by extending the theory developed in~\cite{tilouine-urban} for
$\mathrm{GSp}_4$ and the Betti cohomology of 
Siegel modular threefolds. We relate
this construction to the $P$-ordinary part of completed cohomology. 

For this entire section, fix a subset $\bar{S}\subset \bar{S}_p$, 
where we will take the $P$-ordinary (or, slightly more generally, $Q$-ordinary) 
part of the cohomology of the $\tX_{\tK}$. 

\subsubsection{$P$-ordinary Hida theory at finite level} 
In this section, we will only consider good subgroups 
$\tK\subset \tG(\A_{F^+,f})$ such that the tame level $\tK^p$ is fixed 
and such that, for all $\bar{v}\in \bar{S}$, 
$\tK_{\bar{v}} =\cP_{\bar{v}}(b,c)$ for 
some integers $c\geq b\geq 0$. We denote
such a good subgroup by $\tK(b,c)$ and assume 
from now on that $c\geq 1$. Also set $\cP_{\bar{S}}(b,c):=\prod_{\bar{v}\in \bar{S}}
\cP_{\bar{v}}(b,c)$.  

Recall from~\eqref{eq:ord Hecke action Siegel} %
that we have well-defined actions of the Hecke algebras
$\cH(\widetilde{\Delta}_{\bar{S}}, \tK_{\bar{S}})$
on the complexes $R\Gamma(\tX_{\tK(b,c)}, \cV_{\tilde{\lambda}})$. 
These actions are compatible 
with the natural pullback maps as $b,c$ vary. 
We define the $P$-ordinary part $R\Gamma(\tX_{\tK(b,c)}, \cV_{\tilde{\lambda}})^{\ord}$ 
of the complex $R\Gamma(\tX_{\tK(b,c)}, \cV_{\tilde{\lambda}})$
to be the maximal direct summand %
on which all the $\widetilde{U}_{\tilde{v},n}$ act invertibly. This is a well-defined object 
of $D^+\left(\cP_{\bar{S}}(0,c)/\cP_{\bar{S}}(b,c), \cO\right)$ by
Lemma~\ref{lem:perfect complex} and by 
the theory of ordinary parts, cf.~\cite[\S 2.4]{KT}. Moreover, 
it inherits an action of the abstract Hecke algebra  
$\widetilde{\mathbb{T}}^T\otimes_{\Z} (\bigotimes_{\bar{v}\in \bar{S}}\cH(\widetilde{\Delta}_{\bar{v}}, \tK_{\bar{v}})[\widetilde{U}_{\tilde{v},n}^{-1}])$. 

We similarly have well-defined actions of the Hecke algebras $\cH(\widetilde{\Delta}_{\bar{S}}, \tK_{\bar{S}})$
on the complexes $R\Gamma(\partial \tX_{\tK(b,c)}, \cV_{\tilde{\lambda}})$. 
Therefore, we can also define 
the $P$-ordinary part $R\Gamma(\partial \tX_{\tK(b,c)}, \cV_{\tilde{\lambda}})^{\ord}$ 
of the complex $R\Gamma(\partial \tX_{\tK(b,c)}, \cV_{\tilde{\lambda}})$,
which is an object of $D^+\left(\cP_{\bar{S}}(0,c)/\cP_{\bar{S}}(b,c), \cO\right)$ equipped
with an action of 
$\widetilde{\mathbb{T}}^T\otimes_{\Z} (\bigotimes_{\bar{v}\in \bar{S}}\cH(\widetilde{\Delta}_{\bar{v}}, \tK_{\bar{v}})[\widetilde{U}_{\tilde{v},n}^{-1}])$. 

\subsubsection{The $P$-ordinary part of a smooth representation}\label{sec:ord part}
In this section, we will define various functors that 
will allow us to study $P$-ordinary Hida theory using
completed cohomology. These will be variants of 
the functors considered in~\cite[\S 5.2.1]{10author}, with essentially the same
properties, and we 
will appeal to the basic results in \emph{loc. cit.} throughout. We will introduce a variant with more general parahoric level in \S\ref{sec:ordQlevel} --- we find it clearer to introduce the simplest version of the theory first, which is already sufficient for our results on local-global compatibility in the crystalline case.

For an integer $b\geq 0$, we set 
\[
K_{\bar{v}}(b):=\ker \left(G(\cO_{F^+_{\bar{v}}})\to G(\cO_{F^+_{\bar{v}}}/\varpi^b_{\bar{v}})\right) , 
K_{\bar{S}}(b):= \prod_{\bar{v}}K_{\bar{v}}(b), K_{\bar{S}}:=K_{\bar{S}}(0). 
\]
We also set $U^0_{\bar{S}}:=\prod_{\bar{v}\in \bar{S}} U(\cO_{F^+_{\bar{v}}})$. 
Let $\Delta^+_{\bar{S}}\subset G_{\bar{S}}$ denote the monoid 
generated by $K_{\bar{S}}$ and by $\{\tilde{u}_{\tilde{v},n}\mid \bar{v}\in \bar{S}\}$
and $\Delta_{\bar{S}}\subset G_{\bar{S}}$ the subgroup obtained by adjoining the inverses of the elements $\tilde{u}_{\tilde{v},n}$. We set $\widetilde{\Delta}_{\bar{S},P} := \widetilde{\Delta}_{\bar{S}} \cap \prod_{\bar{v}\in\bar{S}}P(F^+_{\bar{v}})$ and note that we have $\widetilde{\Delta}_{\bar{S},P} = \Delta^+_{\bar{S}} \ltimes U^0_{\bar{S}}$. 

We let 
\[
\Gamma(U^0_{\bar{S}}, \ ):  \mathrm{Mod}_{\mathrm{sm}}(\widetilde{\Delta}_{\bar{S}},\cO/\varpi^m)\to
\mathrm{Mod}_{\mathrm{sm}}(\Delta^+_{\bar{S}},\cO/\varpi^m). 
\]
denote the functor of $U^0_{\bar{S}}$-invariants. 

For $V\in \mathrm{Mod}_{\mathrm{sm}}(\widetilde{\Delta}_{\bar{S}},\cO/\varpi^m)$, we define the action of
an element $g\in \Delta^+_{\bar{S}}$ 
on $v\in \Gamma(U^0_{\bar{S}}, V)$ by the formula 
\begin{equation}\label{eq:Hecke action}
v\mapsto g\cdot v := \sum_{n\in U^0_{\bar{S}}/gU^0_{\bar{S}}g^{-1}} ngv,
\end{equation}
cf.~\cite[\S 3]{emordone}. 
We obtain a derived functor  
\[
R\Gamma(U^0_{\bar{S}}, \ ): D^+_{\mathrm{sm}}(\widetilde{\Delta}_{\bar{S}},\cO/\varpi^m)\to
D^+_{\mathrm{sm}}(\Delta^+_{\bar{S}},\cO/\varpi^m).
\]

Since $U^0_{\bar{S}}$ is compact, an injective object in $\mathrm{Mod}_{\mathrm{sm}}(\widetilde{\Delta}_{\bar{S}},\cO/\varpi^m)$ remains  $\Gamma(U^0_{\bar{S}}, \ )$-acyclic on restriction to $\widetilde{\Delta}_{\bar{S},P}$. So this derived functor factors through the restriction functor to $D^+_{\mathrm{sm}}(\widetilde{\Delta}_{\bar{S},P},\cO/\varpi^m)$.

We also define a functor 
\[
\mathrm{ord}: \mathrm{Mod}_{\mathrm{sm}}(\Delta^+_{\bar{S}},\cO/\varpi^m)\to 
\mathrm{Mod}_{\mathrm{sm}}(\Delta_{\bar{S}},\cO/\varpi^m)
\]
that is the composition of the localization functors 
$\otimes_{\cO/\varpi^m[\tilde{u}_{\tilde{v},n}]}\cO/\varpi^m[(\tilde{u}_{\tilde{v},n})^{\pm 1}]$
for all $\bar{v}\in \bar{S}$. Note that
$\ord$ is an exact functor because localization is an
exact functor, and it preserves injectives by the same
argument as in~\cite[Lemma 5.2.7]{10author}. 

\begin{defn}\label{defn:P-ordinary part} We have a functor 
of \emph{$P$-ordinary parts} 
\[
D^+_{\mathrm{sm}}(\widetilde{\Delta}_{\bar{S}},\cO/\varpi^m) \to 
D^+_{\mathrm{sm}}(\Delta_{\bar{S}}, \cO/\varpi^m), \pi\mapsto 
R\Gamma(U^0_{\bar{S}},\pi)^{\ord}
\]
obtained by composing the functor $R\Gamma(U^0_{\bar{S}},\ )$ with the
functor $\mathrm{ord}$. 
\end{defn}

For $b\geq 0$, we also have a functor 
\[
\Gamma(U^0_{\bar{S}}\rtimes K_{\bar{S}}(b),\ ): 
\mathrm{Mod}_{\mathrm{sm}}(\widetilde{\Delta}_{\bar{S}},\cO/\varpi^m)\to 
\mathrm{Mod}(\Delta^+_{\bar{S}}/K_{\bar{S}}(b), \cO/\varpi^m),
\]
where the action of $g\in \Delta^+_{\bar{S}}/K_{\bar{S}}(b)$ is 
given by the same formula~\eqref{eq:Hecke action}. We 
denote the corresponding derived functor by 
$R\Gamma(U^0_{\bar{S}}\rtimes K_{\bar{S}}(b),\ )$. 
We also define the functor 
\[
\mathrm{ord}_b: \mathrm{Mod}(\Delta^+_{\bar{S}}/K_{\bar{S}}(b),\cO/\varpi^m)\to 
\mathrm{Mod}(\Delta_{\bar{S}}/K_{\bar{S}}(b),\cO/\varpi^m)
\]
by localization. Note that
$\mathrm{ord}_b$ is also an exact functor that preserves injectives. 
Finally, for $c\geq b\geq 0$ and $c\geq 1$, we also have a functor
\[
\Gamma(\cP_{\bar{S}}(b,c),\ ): 
\mathrm{Mod}_{\mathrm{sm}}(\widetilde{\Delta}_{\bar{S}},\cO/\varpi^m)\to 
\mathrm{Mod}(\Delta^+_{\bar{S}}/K_{\bar{S}}(b), \cO/\varpi^m),
\]
where the action of $g\in \Delta^+_{\bar{S}}/K_{\bar{S}}(b)$ is 
given by the same formula~\eqref{eq:Hecke action}. Note here that on the right hand side we are considering the natural action of the Hecke algebra \[\cH(\cP_{\bar{S}}(0,c)\Delta^+_{\bar{S}}\cP_{\bar{S}}(0,c),\cP_{\bar{S}}(b,c))=\cO[\Delta_{\bar{S}}^+/K_{\bar{S}}(b)].\]

\begin{lemma}\label{lem:ordinvariants} There is 
a natural isomorphism
\[
\ord_b\circ \Gamma(K_{\bar{S}}(b),\ ) \simeq 
\Gamma(K_{\bar{S}}(b),\ )\circ \ord
\]
of functors $\mathrm{Mod}_{\mathrm{sm}}(\Delta^+_{\bar{S}},\cO/\varpi^m)
\to \mathrm{Mod}(\Delta_{\bar{S}}/K_{\bar{S}}(b),\cO/\varpi^m)$, which extends to an isomorphism of derived functors
\[
\ord_b\circ R\Gamma(K_{\bar{S}}(b),\ ) \simeq 
R\Gamma(K_{\bar{S}}(b),\ )\circ \ord.
\]
\end{lemma}

\begin{proof}
The same argument as for~\cite[Lemma 5.2.6]{10author} works for the un-derived statement. Since $\mathrm{ord}$ is exact and preserves injectives and $\ord_b$ is exact, the statement for derived functors follows.
\end{proof}

\begin{lemma}\footnote{Compare with~\cite[Lemma 5.2.8]{10author}}\label{lem:level control 1} 
For all $c\geq b\geq 0$ with $c\geq 1$, 
there is a natural isomorphism 
\[
\ord_b \circ \Gamma(U^0_{\bar{S}}\rtimes K_{\bar{S}}(b),\ ) \simeq 
\ord_b \circ \Gamma(\cP_{\bar{S}}(b,c),\ )
\] 
of functors 
\[
\mathrm{Mod}_{\mathrm{sm}}(\widetilde{\Delta}_{\bar{S}},\cO/\varpi^m) \to 
\mathrm{Mod}(\Delta_{\bar{S}}/K_{\bar{S}}(b), \cO/\varpi^m). 
\]
\end{lemma}

\begin{proof} Let $V\in \mathrm{Mod}_{\mathrm{sm}}(\widetilde{\Delta}_{\bar{S}}, \cO/\varpi^m)$. 
We first claim that the natural inclusion $\Gamma(\cP_{\bar{S}}(b,c), V)\hookrightarrow 
\Gamma(U^0_{\bar{S}}\rtimes K_{\bar{S}}(b), V)$ 
is a morphism of $\cO/\varpi^m[\Delta^+_{\bar{S}}]$-modules. Indeed, 
by the formula~\eqref{eq:Hecke action}, it is enough to check that, for all 
$g\in \Delta^+_{\bar{S}}/K_{\bar{S}}(b)$, 
the map 
\[
U^0_{\bar{S}}/gU^0_{\bar{S}}g^{-1}\to \cP_{\bar{S}}(b,c)/\cP_{\bar{S}}(b,c)\cap g\cP_{\bar{S}}(b,c)g^{-1}
\]
is bijective, which holds by the Iwahori decomposition of $\cP_{\bar{S}}(b,c)$ with respect to $P$. 
By the exactness of $\ord_b$, we obtain an injection $\ord_b \Gamma(\cP_{\bar{S}}(b,c), V)\hookrightarrow 
\ord_b\Gamma(U^0_{\bar{S}}\rtimes K_{\bar{S}}(b), V)$. 

We are left to show that this injection is an equality. Let $\tilde{u}_{\bar{S}}:=\prod_{\bar{v}\in \bar{S}}\tilde{u}_{\tilde{v},n}$ 
and  $\widetilde{U}_{\bar{S}}:=\prod_{\bar{v}\in \bar{S}}\widetilde{U}_{\tilde{v},n}$. The result will follow if we show that, 
for any $v\in \Gamma(U^0_{\bar{S}}\rtimes K_{\bar{S}}(b), V)$, there exists $N\geq 0$ such that 
$(\tilde{u}_{\bar{S}})^N v \in V^{\cP_{\bar{S}}(b,c)}$. Since $V$ is smooth, there exists $c'\geq c$ 
such that $v\in V^{\cP_{\bar{S}}(b,c')}$. However, if $c'>c$, then $\widetilde{U}_{\bar{S}}v\in 
V^{\cP_{\bar{S}}(b,c'-1)}$ by~\cite[Lemma 3.3.2]{emordone}. We conclude
by induction. 
\end{proof}

\begin{lemma}\label{lem:level control 2} Let $\pi\in D^+_{\mathrm{sm}}(\widetilde{\Delta}_{\bar{S}}
,\cO/\varpi^m)$. Then for any $c\geq b\geq 0$ with $c\geq 1$ there is a natural 
isomorphism 
\[
R\Gamma(K_{\bar{S}}(b),\ord R\Gamma(U^0_{\bar{S}}, \pi))
\toisom \ord_b R\Gamma(\cP_{\bar{S}}(b,c), \pi)
\]
in $D^+(\Delta_{\bar{S}}/K_{\bar{S}}(b), \cO/\varpi^m)$. 
\end{lemma}

\begin{proof} This is proved in the same way as~\cite[Lemma 
5.2.9]{10author}, except we appeal to Lemma~\ref{lem:level control 1} above
instead of Lemma 5.2.8 of~\emph{op. cit.} 
\end{proof}

\subsubsection{Independence of level}

Recall that the finite free $\cO$-module 
$\cV_{\tilde{\lambda}}$ is equipped with an action
of $\widetilde{\Delta}_{\bar{S}}$. 
We consider the completed cohomology at $\bar{S}$ with coefficients
in $\cV_{\tilde{\lambda}}/\varpi^m$: 
\[
\pi(\tK^{\bar{S}},\tilde{\lambda}, m):= 
R\Gamma\left(\tK^{\bar{S}},R\Gamma\left(\mathfrak{\overline{X}}_{\tG}, 
\cV_{\tilde{\lambda}}/\varpi^m\right)\right). 
\]
We are using the discussion following Lemma \ref{lem:projection formula} to regard $R\Gamma\left(\mathfrak{\overline{X}}_{\tG}, 
\cV_{\tilde{\lambda}}/\varpi^m\right)$ as an object of $D^+_{\mathrm{sm}}(\tG^{T}\times \tK_{T\backslash\bar{S}}\times\widetilde{\Delta}_{\bar{S}}, \O/\varpi^m)$. Then $\pi(\tK^{\bar{S}},\tilde{\lambda}, m)$ is an object of $D^+_{\mathrm{sm}}(\widetilde{\Delta}_{\bar{S}}, \O/\varpi^m)$
equipped with an action of $\widetilde{\mathbb{T}}^T$. 
We note that, for any $c\geq b\geq 0$, $\cP_{\bar{S}}(b,c)\subset \widetilde{\Delta}_{\bar{S}}$
and we have a canonical $\widetilde{\mathbb{T}}^T$-equivariant 
isomorphism 
\[
R\Gamma\left(\cP_{\bar{S}}(b,c), \pi(\tK^{\bar{S}},\tilde{\lambda}, m) \right) \toisom  
R\Gamma\left(\tX_{\tK(b,c)}, \cV_{\tilde{\lambda}}/\varpi^m\right)
\]
in $D^+\left(\cP_{\bar{S}}(0,c)/\cP_{\bar{S}}(b,c), \cO\right/\varpi^m)$. 

We define $\pi^{\ord}(\tK^{\bar{S}}, \tilde{\lambda}, m)$ 
to be the $P$-ordinary part of $\pi(\tK^{\bar{S}},\tilde{\lambda}, m)$
as in Definition~\ref{defn:P-ordinary part}. 
This is an object in $D^+_{\mathrm{sm}}(\Delta_{\bar{S}}, \cO/\varpi^m)$
equipped with an action of $\widetilde{\mathbb{T}}^T$. 

\begin{prop}\label{prop:going to finite level} For 
any integers $m\geq 1$ and $c\geq b\geq 0$ with $c\geq 1$, 
there is a natural $\widetilde{\mathbb{T}}^T$-equivariant 
isomorphism 
\[
R\Gamma\left(K_{\bar{S}}(b), \pi^{\ord}(\tK^{\bar{S}}, \tilde{\lambda}, m)\right)\simeq 
R\Gamma\left(\tX_{\tK(b,c)}, \cV_{\tilde{\lambda}}/\varpi^m\right)^{\ord}
\]
in $D^+(K_{\bar{S}}/K_{\bar{S}}(b), \cO/\varpi^m)$. 
\end{prop}

\begin{proof} This is proved in the same way as~\cite[Prop. 5.2.15]{10author}, given 
Lemma~\ref{lem:level control 2} as an input. 
Again, the key point is that the two definitions of \emph{$P$-ordinary parts} 
(via Hida's idempotent as in~\cite[\S 2.4]{KT} or via localization as in \S~\ref{sec:ord part}) agree on
finite $\cO/\varpi^m$-modules, and that the cohomology groups
of $R\Gamma\left(\tX_{\tK(b,c)}, \cV_{\tilde{\lambda}}/\varpi^m\right)$ 
are finite $\cO/\varpi^m$-modules by Lemma~\ref{lem:perfect complex}. 
\end{proof}

\begin{cor}[Independence of level\footnote{Compare with~\cite[Cor. 5.2.16]{10author}}]
\label{cor:independence of level} For 
any integers $m\geq 1$ and $c\geq b\geq 0$ with $c\geq 1$, the natural 
$\widetilde{\mathbb{T}}^T$-equivariant 
morphism 
\[
R\Gamma\left(\tX_{\tK(b,\mathrm{max}\{1,b\})}, \cV_{\tilde{\lambda}}
/\varpi^m\right)^{\ord}\to 
R\Gamma\left(\tX_{\tK(b,c)}, \cV_{\tilde{\lambda}}/\varpi^m\right)^{\ord}
\]
is an
isomorphism in $D^+(K_{\bar{S}}/K_{\bar{S}}(b), \cO/\varpi^m)$. 
\end{cor}

\begin{proof} This follows immediately from Proposition~\ref{prop:going
to finite level}, upon noting that the LHS of the main isomorphism 
is independent of $c\geq \mathrm{max}\{1,b\}$. 
\end{proof}

\noindent The same results hold for the cohomology of 
the Borel--Serre boundary, with the same proof. 

\begin{prop}\label{prop:going to finite level boundary} For 
any integers $m\geq 1$ and $c\geq b\geq 0$ with $c\geq 1$, 
there is a natural $\widetilde{\mathbb{T}}^T$-equivariant 
isomorphism 
\[
R\Gamma\left(K_{\bar{S}}(b), \pi_{\partial}^{\ord}(\tK^{\bar{S}}, 
\tilde{\lambda}, m)\right)\simeq 
R\Gamma\left(\partial \tX_{\tK(b,c)}, \cV_{\tilde{\lambda}}/\varpi^m\right)^{\ord}
\]
in $D^+(K_{\bar{S}}/K_{\bar{S}}(b), \cO/\varpi^m)$. 
\end{prop}

\subsubsection{A variant at parahoric level}\label{sec:ordQlevel} %
We now suppose we have standard parabolic subgroups $Q_{\bar{v}} \subset P_{\bar{v}}$ for each $\bar{v} \in \bar{S}$, with corresponding parahoric subgroup $\cQ_{\bar{v}} \subset \widetilde{G}(\cO_{F^+_{\bar{v}}})$. We will still be interested in $P$-ordinary parts --- in other words, we only invert the Hecke operator for the element $\tilde{u}_{\overline{S}}$, as in the previous subsection. Our degree shifting arguments will apply to this $P$-ordinary part. Later on, we will take $Q$-ordinary parts after specializing to finite level $\cQ_{\overline{S}}$, but this means that we need to keep track of additional Hecke operators on our cohomology groups.

Recall that we have defined monoids $\widetilde{\Delta}_{\bar{S}}^{\cQ_{\bar{S}}}, {\Delta}_{\bar{S}}^{\cQ_{\bar{S}},+}, {\Delta}_{\bar{S}}^{\cQ_{\bar{S}}}$. We will also make use of $\widetilde{\Delta}_{\bar{S},P}^{\cQ_{\bar{S}}} :=  \widetilde{\Delta}_{\bar{S}}^{\cQ_{\bar{S}}} \cap \prod_{\bar{v}\in\bar{S}}P(F^+_{\bar{v}})$. Note that $\widetilde{\Delta}_{\bar{S},P}^{\cQ_{\bar{S}}} = {\Delta}_{\bar{S}}^{\cQ_{\bar{S}},+} \ltimes U^0_{\bar{S}}$. We set $K_{\bar{v}} := \cQ_{\bar{v}}\cap G(F^+_{\bar{v}})$ for each $\bar{v} \in \bar{S}$. 
	
Using the formula \eqref{eq:Hecke action}, we define a functor
\[R\Gamma(U^0_{\bar{S}},-): D^+_{\sm}(\widetilde{\Delta}_{\bar{S}}^{\cQ_{\bar{S}}},\cO/\varpi^m)\to  D^+_{\sm}({\Delta}_{\bar{S}}^{\cQ_{\bar{S}},+},\cO/\varpi^m)\] which factors through $D^+_{\sm}(\widetilde{\Delta}_{\bar{S},P}^{\cQ_{\bar{S}}},\cO/\varpi^m)$, and the localization, inverting $\tilde{u}_{\tilde{v},n}$ for all $\bar{v}\in\bar{S}$, \[\ord:D^+_{\sm}({\Delta}_{\bar{S}}^{\cQ_{\bar{S}},+},\cO/\varpi^m)\to D^+_{\sm}({\Delta}_{\bar{S}}^{\cQ_{\bar{S}}},\cO/\varpi^m)\] with composition denoted by $\pi \mapsto R\Gamma(U^0_{\bar{S}},\pi)^{\ord}$.  

Note that $K_{\bar{S}}$ is not in general normal in ${\Delta}_{\bar{S}}^{\cQ_{\bar{S}},+}$, but it is in the submonoid $K_{\bar{S}}[\tilde{u}_{\tilde{v},n}|\bar{v}\in\bar{S}]$. This means we can define a functor \[\ord_0R\Gamma\left(K_{\bar{S}},-\right): D^+_{\sm}({\Delta}_{\bar{S}}^{\cQ_{\bar{S}},+},\cO/\varpi^m)  \to  D^+_{\sm}(K_{\bar{S}}[\tilde{u}_{\tilde{v},n}^{\pm 1}|\bar{v}\in\bar{S}]/K_{\bar{S}},\cO/\varpi^m)\] as above, and equip an object in the image of this functor with an action of the Hecke algebra $\cH({\Delta}_{\bar{S}}^{\cQ_{\bar{S}}},K_{\bar{S}})$. We can also compute this functor as $R\Gamma(K_{\bar{S}},-)\circ\ord$, as in Lemma \ref{lem:ordinvariants}.
We have an analogue of Lemma \ref{lem:level control 2}:
\begin{lem}\label{lem:control level extra Hecke}
	Let $\pi \in D^+_{\sm}(\widetilde{\Delta}_{\bar{S}}^{\cQ_{\bar{S}}},\cO/\varpi^m)$. There is a natural isomorphism \[R\Gamma(K_{\bar{S}},\ord R\Gamma(U^0_{\bar{S}},\pi)) \toisom \ord_0R\Gamma\left(\cQ_{\bar{S}},\pi\right)\] in $D^+_{\sm}(K_{\bar{S}}[\tilde{u}_{\tilde{v},n}^{\pm 1}|\bar{v}\in\bar{S}]/K_{\bar{S}},\cO/\varpi^m)$ under which the action of $[K_{\bar{v}} \nu(\varpi_{\bar{v}}) K_{\bar{v}}] \in \cH({\Delta}_{\bar{S}}^{\cQ_{\bar{S}}},K_{\bar{S}})$ on the left hand side matches with the action of $[\cQ_{\bar{v}}\nu(\varpi_{\bar{v}})\cQ_{\bar{v}}]$ on the right hand side.
\end{lem}
\begin{proof}
	The proof is essentially the same as that of Lemma \ref{lem:level control 2}. The description of the Hecke action follows from the decomposition $\widetilde{\Delta}_{\bar{S},P}^{\cQ_{\bar{S}}} = {\Delta}_{\bar{S}}^{\cQ_{\bar{S}},+} \ltimes U^0_{\bar{S}}$ and the fact that the formula defining the double coset operator $[\cQ_{\bar{v}}\cap P(F^+_{\bar{v}})\nu(\varpi_{\bar{v}})\cQ_{\bar{v}}\cap P(F^+_{\bar{v}})]$ also defines $[\cQ_{\bar{v}}\nu(\varpi_{\bar{v}})\cQ_{\bar{v}}]$.
	
	More precisely, if $V \in \Mod_{\sm}(\widetilde{\Delta}_{\bar{S}}^{\cQ_{\bar{S}}},\cO/\varpi^m)$, we consider the natural inclusion
	\[\Gamma(\cQ_{\bar{S}},V)\hookrightarrow \Gamma( K_{\bar{S}}\ltimes U^0_{\bar{S}},V).\] For $\delta \in {\Delta}_{\bar{v}}^{\cQ_{\bar{v}},+}$, the action of $[\cQ_{\bar{v}}\delta\cQ_{\bar{v}}]$ on the left hand side is given by $\sum_{\gamma \in \cQ_{\bar{v}}/(\cQ_{\bar{v}}\cap \delta\cQ_{\bar{v}}\delta^{-1})} \gamma\delta$. The action on the right hand side is given by $\sum_{k \in K_{\bar{v}}/(K_{\bar{v}}\cap\delta K_{\bar{v}} \delta^{-1})} \sum_{n \in U_{\bar{v}}^0/k\delta U_{\bar{v}}^0 (k\delta)^{-1}} nk\delta$. The Iwahori decomposition of $\cQ_{\bar{v}}$ with respect to $P$ shows that the map $(k,n) \mapsto nk$ gives a well-defined bijection \[\{(k,n):k \in K_{\bar{v}}/(K_{\bar{v}}\cap \delta K_{\bar{v}} \delta^{-1}), n \in U_{\bar{v}}^0/k\delta U_{\bar{v}}^0 (k\delta)^{-1}\} \toisom \cQ_{\bar{v}}/(\cQ_{\bar{v}}\cap\delta\cQ_{\bar{v}}\delta^{-1}).\]
	
	It remains to show that the inclusion \[\ord_0\Gamma(\cQ_{\bar{S}},V)\hookrightarrow \ord_0\Gamma(K_{\bar{S}}\ltimes U^0_{\bar{S}},V)\] is bijective. As in the proof of Lemma \ref{lem:level control 1}, if $v$ is in the right hand side, there exists $c \ge 1$ such that $v \in V^{\cQ_{\bar{S}}\cap\cP_{\bar{S}}(0,c)}$. We can apply \cite[Lemma 3.3.2]{emordone} again to show that $\widetilde{U}_{\bar{S}}^N v \in V^{\cQ_{\bar{S}}}$ for some $N \ge 0$.
\end{proof}

\begin{remark}\label{rem:more general monoid}
Lemma~\ref{lem:control level extra Hecke} also holds, with the same proof, if, for each $\bar{v}\in \bar{S}$, we replace each  $\Delta^{\cQ_{\bar{v}},+}_{\bar{v}}$ with a more general monoid $\Sigma^{+}_{\bar{v}}\subset G(F^+_{\bar{v}})$ that contains $K_{\bar{v}}$ and $\tilde{u}_{\tilde{v},n}$ and that satisfies the positivity condition 
\[
\sigma U^0_{\bar{v}}\sigma^{-1}\subseteq U^{0}_{\bar{v}}, \forall \sigma\in \Sigma^{+}_{\bar{v}} 
\]
with respect to $U^0_{\bar{v}}$. For such a monoid $\Sigma^{+}_{\bar{v}}$, we define $\widetilde{\Sigma}_{\bar{v}}:=\cQ_{\bar{v}}\Sigma^+_{\bar{v}}\cQ_{\bar{v}}$, $\Sigma_{\bar{v},P}:=\widetilde{\Sigma}_{\bar{v}}\cap P(F^+_{\bar{v}})$ and $\Sigma_{\bar{v}}:=\Sigma_{\bar{v}}^+[\tilde{u}_{\tilde{v},n}^{-1}]$. We have $\Sigma_{\bar{v},P} = \Sigma^+_{\bar{v}}\ltimes U_{\bar{v}}^0$. The construction of the ordinary parts functors in the beginning of \S~\ref{sec:ordQlevel} goes through verbatim. The fact that the natural inclusion 
\[
\Gamma(\cQ_{\bar{v}}, V)\hookrightarrow\Gamma(K_{\bar{v}}\ltimes U_{\bar{v}}^0, V)
\]
is equivariant for the Hecke action of $\sigma\in \Sigma^+_{\bar{v}}$ follows from the positivity condition. We obtain that, for each $\pi\in D^+_{\mathrm{sm}}(\widetilde{\Sigma}_{\bar{v}},\cO/\varpi^m)$, there is a natural isomorphism 
\[
R\Gamma(K_{\bar{v}}, \mathrm{ord}R\Gamma(U^0_{\bar{v}},\pi))\stackrel{\sim}{\to}
\mathrm{ord}_0R\Gamma(\cQ_{\bar{v}}, \pi),
\]
in $D^+_{\mathrm{sm}}(K_{\bar{v}}[\tilde{u}_{\tilde{v},n}^{\pm 1}],\cO/\varpi^m)$, 
which matches the Hecke actions of $[K_{\bar{v}}\sigma K_{\bar{v}}]$ on the left hand side and $[\cQ_{\bar{v}}\sigma \cQ_{\bar{v}}]$ on the right hand side, for every $\sigma\in \Sigma_{\bar{v}}$.  
\end{remark}

Using Lemma~\ref{lem:control level extra Hecke}, we get natural analogues of our independence of level statements. In this context, we will have $\pi^{\ord}(\widetilde{K}^{\bar{S}},\tilde{\lambda},m) \in D^+_{\sm}(\widetilde{\Delta}^{\cQ_{\bar{S}}}_{\bar{S}},\cO/\varpi^m)$ and the same for boundary cohomology. 

\subsubsection{Independence of weight}

We retain our current set-up, with parabolic subgroups $Q_{\bar{v}}$ for $\bar{v} \in \bar{S}$ and a parahoric level subgroup $K_{\bar{S}}$. Assume that we have a dominant weight $\lambda$ for $G$. For a subset $\bar{S}\subseteq \bar{S}_p$, set 
\[
\cV_{\lambda_{\bar{S}}}:= \bigotimes_{{v}\in {S}}\bigotimes_{\tau\in \Hom(F_{{v}}, E), \cO} 
\cV_{\lambda_{{\tau}}}. 
\] 
This is, a priori, a finite free $\cO$-module with an action of $K_{\bar{S}}$. We have explained how to extend the $K_{\bar{S}}$-action 
to an action of $\Delta^{\cQ_{\bar{S}}}_{\bar{S}}$ (in particular, $\tilde{u}_{\tilde{v},n}$ acts trivially for each $\bar{v}\in \bar{S}$). 
\begin{lemma}
Using our identification of $K_{\bar{S}}$ with the block diagonal Levi subgroup of $\prod_{\bar{v}\in \bar{S}}\Prm_{n,n}(\cO_{F_{\tilde{v}}})$, we can identify
\[
\cV_{\lambda_{\bar{S}}}= \bigotimes_{{\bar{v}}\in {\bar{S}}}\bigotimes_{\tau\in \Hom(F^+_{\bar{v}}, E), \cO} 
\cV_{-w_{0,n}\lambda_{\tilde{\tau}c}}\otimes\cV_{\lambda_{\tilde{\tau}}}
\] with the action on both factors $\cV_{-w_{0,n}\lambda_{\tilde{\tau}c}}\otimes\cV_{\lambda_{\tilde{\tau}}}$ defined using the embedding $\tilde{\tau}$.\end{lemma}  %
\begin{proof}
Recall that for each $\bar{v}\in\overline{S}$ and place $\tilde{v}|\bar{v}$ of $F$ we have an isomorphism  $\iota_{\tilde{v}}:\widetilde{G}(F^+_{\bar{v}}) \cong \GL_{2n}(F_{\tilde{v}})$ identifying the Levi subgroup $G(F^+_{\bar{v}}) = \GL_n(F_{\tilde{v}}) \times \GL_n(F_{\tilde{v}^c})$ with block diagonal matrices in $\Prm_{n,n}(F_{\tilde{v}})$ via $(A_{\tilde{v}},A_{\tilde{v}^c}) \mapsto \begin{pmatrix}(\Psi_n\ ^tA_{\tilde{v}^c}^{-1}\Psi_n)^c & 0 \\ 0 & A_{\tilde{v}} \end{pmatrix}$. We write $\theta$ for the isomorphism $\GL_n(F_{\tilde{v}^c}) \cong \GL_n(F_{\tilde{v}})$ defined by \[\theta(A) = (\Psi_n\ ^tA^{-1}\Psi_n)^c,\]
which preserves our chosen Borel subgroup $\Brm$. 

For each $\tau \in \Hom(F^+_{\bar{v}}, E)$,  $K_{\bar{v}} = \GL_n(\cO_{F_{\tilde{v}}})\times\GL_n(\cO_{F_{\tilde{v}^c}})$ acts `factor-by-factor' on $\cV_{\lambda_{\tilde{\tau}}}\otimes \cV_{\lambda_{\tilde{\tau}c}}$. Describing this representation in terms of block diagonal matrices, we get $\theta^{-1}\cV_{\lambda_{\tilde{\tau}c}}\otimes \cV_{\lambda_{\tilde{\tau}}}$, where $\theta^{-1}\cV_{\lambda_{\tilde{\tau}c}}$ denotes the representation of $\GL_n(\cO_{F_{\tilde{v}}})$ given by pulling back the representation $\cV_{\lambda_{\tilde{\tau}c}}$ by $\theta^{-1}$. To finish the proof, we need to explain why $\theta^{-1}\cV_{\lambda_{\tilde{\tau}c}} \cong \cV_{-w_{0,n}\lambda_{\tilde{\tau}c}}$. To see this, consider the map

\begin{align*}(\Ind_{\Brm_n}^{\GL_n}w_{0,n}\lambda_{\tilde{\tau}c})_{/\cO} &\to (\Ind_{\Brm_n}^{\GL_n}-\lambda_{\tilde{\tau}c})_{/\cO}  \\ 
f &\mapsto \left(g \mapsto f(\Psi_n\ ^{t}g^{-1}\Psi_n)\right)\end{align*} and note that it gives the desired isomorphism.
\end{proof}

\begin{prop}\label{prop:independence of weight} Given a dominant weight 
$\tilde{\lambda}$ for $\widetilde{G}$ and a subset $\bar{S}\subseteq \bar{S}_p$, 
let $\tilde{\lambda}^{\bar{S}}$ be defined as follows: 
\begin{itemize}
\item if $\bar{v}\in \bar{S}$ and $\tau\in \Hom(F^+_{\bar{v}}, E)$, then $\tilde{\lambda}^{\bar{S}}_{\tau}:=(0,\dots,0)$. 
\item if $\bar{v}\in\bar{S}_p\setminus \bar{S}$ and 
$\tau\in \Hom(F^+_{\bar{v}}, E)$, then $\tilde{\lambda}^{\bar{S}}_{\tau}:=\tilde{\lambda}_{\tau}$. 
\end{itemize}
Identify $\tilde{\lambda}$ with a dominant weight $\lambda$ for $G$ as in \ref{eq:identification of weights}. For any integer $m\geq 1$ 
there is a natural $\widetilde{\mathbb{T}}^T$-equivariant isomorphism 
\[
\pi^{\ord}(\tK^{\bar{S}}, \tilde{\lambda}, m)\toisom 
\pi^{\ord}(\tK^{\bar{S}}, \tilde{\lambda}^{\bar{S}}, m) \otimes \cV_{w^P_0\lambda_{\bar{S}}}/\varpi^m
\]%
in $D^+_{\mathrm{sm}}(\Delta^{\cQ_{\bar{S}}}_{\bar{S}}, \cO/\varpi^m)$. 
\end{prop}

\begin{proof} Note first that $\pi^{\ord}(\tK^{\bar{S}}, \tilde{\lambda}, m)$
only depends on $\cV_{\tilde{\lambda}_{\bar{S}}}/\varpi^m$ as 
an object in $D^+_{\sm}(U^0_{\bar{S}}\rtimes \Delta^{\cQ_{\bar{S}},+}_{\bar{S}}, \cO/\varpi^m)$. 
As in Lemma~\ref{lem:evsurjective}, there is a $U^0_{\bar{S}}\rtimes \Delta^{\cQ_{\bar{S}},+}_{\bar{S}}$-equivariant morphism 
of finite free $\cO$-modules 
\[
\mathrm{ev}: \cV_{\tilde{\lambda}_{\bar{S}}}
\to \cV_{w^P_0\lambda_{\bar{S}}}
\]
given by evaluation of functions at the identity, where the action of 
$U^0_{\bar{S}}\rtimes \Delta^{\cQ_{\bar{S}},+}_{\bar{S}} \subset \widetilde{\Delta}^{\cQ_{\bar{S}}}_{\bar{S}}$ on the LHS is as in Lemma~\ref{lem:twisted action} 
and the action on the RHS factors through the 
action of $\Delta^{\cQ_{\bar{S}},+}_{\bar{S}}$. By Lemma~\ref{lem:evsurjective}, this morphism 
is surjective. Let $\cK_{\tilde{\lambda}_{\bar{S}}}:=\ker(\mathrm{ev})$, 
a finite free $\cO$-module with an action of $U^0_{\bar{S}}\rtimes \Delta^{\cQ_{\bar{S}},+}_{\bar{S}}$. 

For any $m\geq 1$, we have a short exact sequence of 
$\cO/\varpi^m[U^0_{\bar{S}}\rtimes \Delta^{\cQ_{\bar{S}},+}_{\bar{S}}]$-modules
\[
0\to \cK_{\tilde{\lambda}_{\bar{S}}}/\varpi^m \to \cV_{\tilde{\lambda}_{\bar{S}}}/\varpi^m
\to \cV_{w^P_0\lambda_{\bar{S}}}/\varpi^m\to  0. 
\]
We first claim that $\mathrm{ev}\pmod{\varpi^m}$ induces an isomorphism between 
\[
\pi^{\ord}(\tK^{\bar{S}}, \tilde{\lambda}, m) \stackrel{\mathrm{def}}{=} 
\ord R\Gamma\left(U^0_{\bar{S}}, R\Gamma\left(\tK^{\bar{S}}, R\Gamma(\overline{\mathfrak{X}}_{\widetilde{G}},
\cV_{\tilde{\lambda}}/\varpi^m) \right)\right)
\]
and 
\begin{equation}\label{eq:intermediate coh}
\ord  R\Gamma\left(U^0_{\bar{S}}, R\Gamma\left(\tK^{\bar{S}}, R\Gamma(\overline{\mathfrak{X}}_{\widetilde{G}},
\cV_{\tilde{\lambda}^{\bar{S}}}/\varpi^m)\otimes  \cV_{w^P_0\lambda_{\bar{S}}}/\varpi^m\right)\right). 
\end{equation}
By definition (cf.~Lemma \ref{lem:projection formula}), we have \[R\Gamma(\overline{\mathfrak{X}}_{\widetilde{G}},
\cV_{\tilde{\lambda}}/\varpi^m) = R\Gamma(\overline{\mathfrak{X}}_{\widetilde{G}},
\cV_{\tilde{\lambda}^{\bar{S}}}/\varpi^m)\otimes \cV_{\tilde{\lambda}_{\bar{S}}}/\varpi^m.\]  

To prove the claim, it is therefore enough to show that 
\[
\ord R\Gamma\left(U^0_{\bar{S}}, R\Gamma\left(\tK^{\bar{S}}, R\Gamma(\overline{\mathfrak{X}}_{\widetilde{G}},
\cV_{\tilde{\lambda}^{\bar{S}}}/\varpi^m)\otimes \cK_{\tilde{\lambda}_{\bar{S}}} /\varpi^m \right)\right)
\]
is trivial. This follows from Lemma~\ref{lem:ord kills kernel} below. 

Finally, we observe that there is a natural isomorphism 
in $D^+_{\mathrm{sm}}(\Delta^{\cQ_{\bar{S}}}_{\bar{S}}, \cO/\varpi^m)$
between~\eqref{eq:intermediate coh} 
and $\pi^{\ord}(\tK^{\bar{S}}, \tilde{\lambda}^{\bar{S}}, m) \otimes \cV_{w^P_0\lambda_{\bar{S}}}/\varpi^m$. 
Indeed, $\tilde{u}_{\tilde{v},n}$, $\tK^{\bar{S}}$ and $U^0_{\bar{S}}$ all act trivially on the finite free $\cO/\varpi^m$-module  
$\cV_{w^P_0\lambda_{\bar{S}}}/\varpi^m$, so we can pull this factor outside all the functors being applied in \eqref{eq:intermediate coh}.
\end{proof}

\begin{lemma}\label{lem:ord kills kernel} Let $\tau\in \Hom(F^+_{\bar{v}}, E)$ 
and let 
\[
\cK_{\tilde{\lambda}_{\tau}}:= \ker\left(\cV_{\tilde{\lambda}_{\tau}}\to
\cV_{\lambda_{\tilde{\tau}}}\otimes \cV_{-w_{0,n}\lambda_{\tilde{\tau}c}} \right)
\] 
be the kernel of the evaluation at identity map. For any $m\geq 1$, 
we have $(\tilde{u}_{\tilde{v},n})^m(\cK_{\tilde{\lambda}_{\tau}}/\varpi^m) = 0$.  
\end{lemma}

\begin{proof} 
Since $\cK_{\tilde{\lambda}_{\tau}}$ is a $\Trm_{2n}/\cO$-stable submodule of 
the evaluation at $\cO$ of an algebraic representation of $\GL_{2n}/\cO$, 
it has a decomposition into weight spaces for the diagonal torus $\Trm_{2n}$. The parabolic subgroup $\Prm_{n,n}$ corresponds to the set of simple roots $I=\Delta\backslash(e_n-e_{n+1})$. It follows from \cite[Proposition 4.1]{Cab84} that the weights which show up in $\cK_{\tilde{\lambda}_{\tau}}$ are the weights $\mu$ of $V_{\tilde{\lambda}_{\tau}}$ such that $\mu - w_{0,2n}\tilde{\lambda}_\tau$ contains a positive multiple of $(e_n-e_{n+1})$ in its decomposition into simple roots. This condition corresponds to the rescaled action of $\tilde{u}_{\tilde{v},n}$ on the $\mu$-weight space acting as multiplication by a positive power of $\tilde{\tau}(\varpi_{\tilde{v}})$. 
\end{proof}

The same result holds for the cohomology of 
the Borel--Serre boundary, with the same proof. 

\begin{prop}\label{prop:independence of weight boundary} 
Given a weight 
$\tilde{\lambda}$ and a subset $\bar{S}\subseteq \bar{S}_p$, 
let $\tilde{\lambda}^{\bar{S}}$ be defined as in Proposition~\ref{prop:independence of weight}. 
For any integer $m\geq 1$ 
there is a natural $\widetilde{\mathbb{T}}^T$-equivariant isomorphism 
\[
\pi^{\ord}_{\partial}(\tK^{\bar{S}}, \tilde{\lambda}, m)\toisom 
\pi^{\ord}_{\partial}(\tK^{\bar{S}}, \tilde{\lambda}^{\bar{S}}, m) \otimes \cV_{w^P_0\lambda_{\bar{S}}}/\varpi^m
\]
in $D^+_{\mathrm{sm}}(\Delta^{\cQ_{\bar{S}}}_{\bar{S}}, \cO/\varpi^m)$. 
\end{prop}

\subsubsection{A variant with dual coefficients}\label{sec:dual coeffs prelims}
We will also make use of a variant of $P$-ordinary Hida theory which can be applied with dual coefficient systems. It is formulated in terms of the ordinary parts $\ord^\vee$, $\ord_0^\vee$ defined using the Hecke action of $\tilde{u}_{\tilde{v},n}^{-1}$ on invariants under $\overline{U}^1_{\bar{v}}$  and $\cQ_{\bar{v}}$ respectively, where $\overline{U}^1_{\bar{v}}$ is the block-strictly-\emph{lower}-triangular part of the parahoric $\cP_{\bar{v}}$. We will fix $\bar{v} \in \bar{S}$ and start with a representation $\pi$ of the inverse monoid $(\widetilde{\Delta}^{\cQ_{\bar{v}}}_{\bar{v}})^{-1} = \coprod_{\nu \in X_{Q_{\bar{v}}}}\cQ_{\bar{v}}\nu(\varpi_{\bar{v}})^{-1}\cQ_{\bar{v}}$. Set $K_{\bar{v}}:=\cQ_{\bar{v}}\cap G(F^+_{\bar{v}})$.

\begin{lem}\label{lem:control level dual}
	Let $\pi \in D^+_{\sm}((\widetilde{\Delta}^{\cQ_{\bar{v}}}_{\bar{v}})^{-1},\cO/\varpi^m)$. There is a natural isomorphism \[R\Gamma(K_{\bar{v}},\ord^\vee R\Gamma(\overline{U}^1_{\bar{v}},\pi)) \toisom \ord_0^\vee R\Gamma\left(\cQ_{\bar{v}},\pi\right)\] in $D^+_{\sm}(K_{\bar{v}}[\tilde{u}_{\tilde{v},n}^{\pm1}]/K_{\bar{v}},\cO/\varpi^m)$ under which the action of $[K_{\bar{v}} \nu(\varpi_{\bar{v}})^{-1} K_{\bar{v}}] \in \cH((\Delta^{\cQ_{\bar{v}}}_{\bar{v}})^{-1},K_{\bar{v}})$ on the left hand side matches with the action of $[\cQ_{\bar{v}}\nu(\varpi_{\bar{v}})^{-1}\cQ_{\bar{v}}]$ on the right hand side.
\end{lem}
\begin{proof} We will reduce to Remark~\ref{rem:more general monoid} by conjugating by an appropriate group element. 
	Indeed, conjugation by $g:=\tilde{u}_{\tilde{v},n}^{-1}w^P_0$ sends $\overline{U}^1_{\bar{v}}$ to $U^0_{\bar{v}}$, $\cQ_{\bar{v}}$ to the parahoric $\overline{\cQ}^{w^P_0}_{\bar{v}}$ corresponding to the standard parabolic with Levi subgroup  $Q_{\bar{v}}^{w^P_0}\cap G(F^+_{\bar{v}})$, and it sends $K_{\bar{v}}$ to $K_{\bar{v}}^{w_0^P}$. Conjugation by $g$, moreover, sends $(\widetilde{\Delta}^{\cQ_{\bar{v}}}_{\bar{v}})^{-1}$ to $\widetilde{\Sigma}_{\bar{v}} = \overline{\cQ}^{w^P_0}_{\bar{v}}\Sigma^+_{\bar{v}}\overline{\cQ}^{w^P_0}_{\bar{v}}$, for 
    a monoid 
    \[
    K_{\bar{v}}^{w_0^P}\subset \Sigma^+_{\bar{v}}\subset G(F^+_{\bar{v}})
    \] 
    that satisfies the positivity condition in Remark~\ref{rem:more general monoid} with respect to $U^0_{\bar{v}}$. It sends $\tilde{u}_{\tilde{v},n}^{-1}$ to $\tilde{u}_{\tilde{v},n}\cdot \tilde{u}_{\tilde{v},2n}^{-1}$, where note that the central elements $\tilde{u}_{\tilde{v},2n}^{\pm 1}$ already belong to $\Sigma_{\bar{v}}^+$. 
	
	Conjugation by $g$ then identifies $R\Gamma(K_{\bar{v}},\ord^\vee R\Gamma(\overline{U}^1_{\bar{v}},\pi))$ and $\ord_0^\vee R\Gamma\left(\cQ_{\bar{v}},\pi\right)$ with $R\Gamma(K^{w^P_0}_{\bar{v}},\ord R\Gamma({U}^0_{\bar{v}},\pi^{g}))$ and $\ord_0 R\Gamma\left(\overline{\cQ}^{w^P_0}_{\bar{v}},\pi^{g}\right)$ respectively, where $\pi^g \in D^+_{\sm}(\widetilde{\Sigma}_{\bar{v}},\cO/\varpi^m)$ denotes $\pi$ with the action of $x \in \widetilde{\Sigma}_{\bar{v}}$ given by the action of $ g^{-1} xg \in (\widetilde{\Delta}^{\cQ_{\bar{v}}}_{\bar{v}})^{-1}$. Applying Remark \ref{rem:more general monoid} now gives the desired result.	
\end{proof}

An independence of weight statement for dual coefficients can be proved by following the proof of Proposition \ref{prop:independence of weight}, using the short exact sequence \[0 \to \cV_{w^P_0\lambda_{\bar{S}}}^\vee \xrightarrow{\mathrm{ev}^\vee} \cV_{\tilde{\lambda}_{\bar{S}}}^{\vee} \to \cK_{\tilde{\lambda}_{\bar{S}}}^\vee \to 0 \] and the topological nilpotence of $\tilde{u}_{\tilde{v},n}^{-1}$ on $\cK_{\tilde{\lambda}_{\bar{S}}}^\vee$.

\subsection{New ingredients for degree shifting}

\subsubsection{A computation of $P$-ordinary parts}\label{sec:local P-ord computation}

In this section, we compute the $P$-ordinary part of a 
parabolic induction from $G$ to $\widetilde{G}$, in the same 
spirit as the computation of ordinary parts 
in~\cite[\S 5.3]{10author}. Our calculations here 
are purely local. The global application is to the 
boundary of the Borel--Serre compactification of the locally
symmetric spaces for $\tG$ and it is 
carried out in \S~\ref{sec:LGC}. 

Fix a prime $\bar{v}$ of $F^+$ dividing $p$ and let $L:=F^+_{\bar{v}}$, 
a $p$-adic field (with ring of integers $\cO_L$ and uniformizer
$\varpi_{L}$). In this section, we let $\Grm/\cO_L$ be a 
split connected reductive group with split maximal torus $\Trm\subset \Grm$. 
Write $\Wrm:=W(\Grm,\Trm)$ for the Weyl group, and fix
a Borel subgroup $\Brm$ containing $\Trm$ and a parabolic subgroup 
$\Brm\subset \Prm \subset \Grm$ with Levi decomposition 
$\Prm = \Mrm \ltimes \Urm$. The Weyl group $\Wrm_{\Prm}$ 
of $\Mrm$ can be identified with a subgroup of $\Wrm$. We denote by 
$\Wrm^\Prm\subset \Wrm$ the subset of minimal length representatives 
of $\Wrm_P\backslash \Wrm$. We denote the length of an element $w\in \Wrm$
by $\ell(w)\in \Z_{\geq 0}$.  

Recall from~\cite[Cor. 5.20]{borel-tits} the (generalized) 
Bruhat decomposition 
\[
\Grm(L) = \bigsqcup_{w\in \Wrm^{\Prm}} \Prm(L)w\Brm(L). 
\]
Denote by $\PWP$ the intersection $\Wrm^{\Prm}\cap 
(\Wrm^\Prm)^{-1}$. This is a set of minimal length representatives
for the double cosets $\Wrm_{\Prm}\backslash \Wrm /\Wrm_{\Prm}$,
cf.~\cite[Lemma 3.2.2]{digne-michel}. 

\begin{lemma}\label{lem:gengenbruhat}
	We have a set-theoretic decomposition 
	\[
	\Grm(L) = \coprod_{w \in \PWP} \Prm(L)w\Prm(L).
	\] 
	The 
	closure relations (for the $p$-adic topology) are given by the Bruhat 
	ordering 
	\[
	\overline{\Prm(L)w\Prm(L)} = \coprod_{w' \le w \in \PWP} 
	\Prm(L)w'\Prm(L).  
	\] 
	Moreover, if $\Omega \subset \PWP$ is an upper subset\footnote{This means that, 
	if $w \in \Omega$ and $w' \in \PWP$ satisfies $w' \ge w$, then $w' \in \Omega$.},  
	then $\Prm(L)\Omega \Prm(L)$ is open in $\Grm(L)$.
\end{lemma}

\begin{proof}
	See \cite[Lemma 2.1.2]{hauseuxP}.
\end{proof}

We are interested in the parabolic induction functor 
\[
\mathrm{Ind}_{\Prm(L)}^{\Grm(L)} : D^+_{\sm}(\Prm(L), \cO/\varpi^m) \to 
D^+_{\sm}(\Grm(L), \cO/\varpi^m). 
\]
This functor is exact and preserves injectives. We define several 
functors related to it. 

For $w\in \PWP$, define $S_w: = \Prm(L)w\Prm(L)$ and 
$S_w^\circ: = \Prm(L)w \Mrm(L) \Urm^0$, where 
$\Urm^0 := \Urm(\cO_L)$. The subset $S_w^\circ \subset \Grm(L)$ is invariant under left multiplication by $\Prm(L)$ and right multiplication by inverses of elements in $\Mrm(L)^+\ltimes\Urm^0$, where $\Mrm(L)^+ = \{m \in \Mrm(L): m\Urm^0m^{-1}\subset \Urm^0\}$ (this means that functions with support in $S_w^\circ$ are stable under right translation by $\Mrm(L)^+\ltimes\Urm^0$).

For any $i\in \Z_{\geq 0}$, we define 
\[
\Grm_{\geq i}: =\bigsqcup_{\ell(w)\geq i} S_w, 
\]
which is an open subset of $\Grm(L)$ by Lemma~\ref{lem:gengenbruhat}, 
and which is invariant under left and right multiplication by 
$\Prm(L)$. For any $i\in \Z_{\geq 0}$, we define 
a functor 
\[
I_{\geq i} : \mathrm{Mod}_{\sm}(\Prm(L), \cO/\varpi^m)\to 
 \mathrm{Mod}_{\sm}(\Prm(L), \cO/\varpi^m) 
\]
by sending $\pi\in  \mathrm{Mod}_{\sm}(\Prm(L), \cO/\varpi^m)$ to 
\[ 
\begin{split} I_{\geq i}(\pi) = \{ f : \Grm_{\geq i} \to \pi \mid f  \text{ locally constant,}  
&\text{ of compact support modulo }\Prm(L),\\ & \forall p \in \Prm(L), 
g \in \Grm_{\geq i}, f(pg) = p f(g)  \}, \end{split} 
\]
where $\Prm(L)$ acts by right translation. For $w\in \PWP$, we define a functor
\[
I_{w} : \mathrm{Mod}_{\sm}(\Prm(L), \cO/\varpi^m)\to 
 \mathrm{Mod}_{\sm}(\Prm(L), \cO/\varpi^m) 
\]
by sending $\pi\in  \mathrm{Mod}_{\sm}(\Prm(L), \cO/\varpi^m)$ to 
\[ 
\begin{split} I_{w}(\pi) = \{ f : S_w \to \pi \mid f  \text{ locally constant,}  
&\text{ of compact support modulo }\Prm(L),\\ & \forall p \in \Prm(L), 
g \in S_w, f(pg) = p f(g)  \}, \end{split} 
\]
where again $\Prm(L)$ acts by right translation. 
Finally, for $w\in \PWP$, we also define a functor
\[
I^\circ_{w} : \mathrm{Mod}_{\sm}(\Prm(L), \cO/\varpi^m)\to 
\mathrm{Mod}_{\sm}(\Mrm(L)^+\ltimes \Urm^0, \cO/\varpi^m) 
\]
by defining $I^\circ_{w}(\pi)\subset I_w(\pi)$ to be the subset 
of functions with support in $S^\circ_w$. 

\begin{prop}\label{prop:ses}\leavevmode 
\begin{enumerate}
\item We have $I_{\geq 0} = \mathrm{Res}^{\Grm(L)}_{\Prm(L)} \circ \mathrm{Ind}_{\Prm(L)}^{\Grm(L)}$. 
\item Each functor $I_{\geq i}$, $I_w$ and $I^\circ_w$ is exact. 
\item For each $i\in \Z_{\geq 0}$ and $\pi \in  \mathrm{Mod}_{\sm}(\Prm(L), \cO/\varpi^m)$,
there is a functorial short exact sequence 
\[
0\to I_{\geq i+1}(\pi)\to I_{\geq i}(\pi)\to \bigoplus_{\ell(w) = i} I_{w}(\pi) \to 0. 
\]
\end{enumerate}
\end{prop}

\begin{proof} 
The first part follows from the definition of parabolic induction and 
from the fact that $ \Prm(L)\backslash \Grm(L)$ is compact. The second part 
follows from the fact that the natural map $\Grm(L)\twoheadrightarrow
\Prm(L)\backslash \Grm(L)$ admits a continuous section, 
which can be deduced from~\cite[Part II, \S 1.10]{jantzen} 
and~\cite[Lemma 2.1.1]{hauseux}.
The third part follows
from~\cite[Lemma 2.2.1]{hauseuxP}, noting that the length function $\ell$ 
is strictly monotonic for the Bruhat order. 
\end{proof}

We deduce that, for any $\pi\in D^+_{\sm}(\Prm(L), \cO/\varpi^m)$, 
there is a functorial distinguished triangle 
\begin{equation}\label{eq:distinguished triangle}
I_{\geq i+1}(\pi)\to I_{\geq i}(\pi)\to \bigoplus_{\ell(w) = i} I_{w}(\pi) \to I_{\geq i+1}(\pi)[1]
\end{equation}
in $D^+_{\sm}(\Prm(L), \cO/\varpi^m)$. 
 
\begin{prop}\label{prop:zeroconnecthom} Let $\pi\in D^+_{\sm}(\Prm(L), \cO/\varpi^m)$
and let $\cV$ be a finite free $\cO/\varpi^m$-module equipped with a smooth representation of an open submonoid $\Delta^+ \subset \Mrm(L)$ containing an open subgroup $\Krm \subset\Mrm(\cO_L)$. Assume that $g\Urm^0g^{-1}\subseteq \Urm^0$ for any $g\in \Delta^+$. 
For any $i\in \Z_{\geq 0}$ and $j\in \Z$, the sequence 
\[
0\to R^j\Gamma\left(\Krm\ltimes \Urm^0, \cV\otimes_{\cO/\varpi^m} I_{\geq i+1}(\pi) \right)
\to R^j\Gamma\left(\Krm\ltimes \Urm^0, \cV\otimes_{\cO/\varpi^m} I_{\geq i}(\pi) \right)
\]

\[
\to \bigoplus_{\ell(w) = i} R^j\Gamma\left(\Krm\ltimes \Urm^0, \cV\otimes_{\cO/\varpi^m} I_{w}(\pi) \right)\to 0
\]
associated to~\eqref{eq:distinguished triangle} is an exact sequence of $\cH(\Delta^+,\Krm)$-modules.  
\end{prop}

\begin{proof} The distinguished triangle \eqref{eq:distinguished triangle} gives a distinguished triangle \[
	\cV\otimes_{\cO/\varpi^m}I_{\geq i+1}(\pi)\to \cV\otimes_{\cO/\varpi^m}I_{\geq i}(\pi)\to \bigoplus_{\ell(w) = i} \cV\otimes_{\cO/\varpi^m}I_{w}(\pi) \to \cV\otimes_{\cO/\varpi^m}I_{\geq i+1}(\pi)[1]\]
	in $D^+_{\sm}(\Delta^+ \ltimes \Urm^0, \cO/\varpi^m)$, and taking cohomology gives us the desired sequence of $\cH(\Delta^+,\Krm)$-modules. It remains to check exactness, for which we can forget the Hecke action and just consider $\cV$ as a representation of $\Krm$. We consider 
decompositions $\Grm_{\geq i} = U_1\sqcup U_2$ 
into open and closed subsets that are $\Prm(L)$-invariant
on the left and $\Prm(\cO_L)$-invariant on the right, and
such that $U_1\subset \Grm_{\geq i+1}$. Any such 
decomposition induces a functorial decomposition
$I_{\geq i}(\pi) = I_{U_1}(\pi) \oplus I_{U_2}(\pi)$ 
in the category $\mathrm{Mod}_{\sm}(\Krm\ltimes \Urm^0,\cO/\varpi^m)$,
where $I_{U_1}$ denotes functions with support in $U_1$,
and similarly for $U_2$ (we can also tensor this decomposition with $\cV$). In particular, for any 
$\pi\in D^+_{\sm}(\Prm(L), \cO/\varpi^m)$, the associated
morphism 
\[
R^j\Gamma\left(\Krm\ltimes \Urm^0, \cV\otimes_{\cO/\varpi^m}I_{U_1}(\pi)\right) 
\to  R^j\Gamma\left(\Krm\ltimes \Urm^0, \cV\otimes_{\cO/\varpi^m} I_{\geq i}(\pi) \right)
\]
of $\cO/\varpi^m$-modules is injective. Lemma~\ref{lem:openclosedcover}
below implies that $I_{\geq i+1}$ can be written as a filtered direct limit
of functors of the form $I_{U_1}$. Since the tensor product $\cV\otimes_{\cO/\varpi^m}\ $ and
the functor $R^j\Gamma(\Krm\ltimes \Urm^0,\ )$ commute with filtered
direct limits\footnote{Since $\Krm\ltimes \Urm^0$ is compact, 
a filtered direct limit of injective objects in $\mathrm{Mod}_{\sm}(\Krm\ltimes \Urm^0,\cO/\varpi^m)$ 
is again an injective object in $\mathrm{Mod}_{\sm}(\Krm\ltimes \Urm^0,\cO/\varpi^m)$, 
cf.~\cite[Prop. 2.1.3]{emordtwo}.}, it follows that the morphism 
\[
R^j\Gamma\left(\Krm\ltimes \Urm^0, \cV\otimes_{\cO/\varpi^m}I_{\geq i+1}(\pi)\right) 
\to  R^j\Gamma\left(\Krm\ltimes \Urm^0, \cV\otimes_{\cO/\varpi^m} I_{\geq i}(\pi) \right)
\]
is injective. Since the injectivity applies for any $j\in \mathbb{Z}$,
the long exact sequence of cohomology groups attached 
to the distinguished triangle~\eqref{eq:distinguished triangle} (tensored with $\cV$) gives
the statement of the proposition. 
\end{proof}

\begin{lemma}\label{lem:openclosedcover}
For any $i \in \Z_{\geq 0}$, there exist decompositions 
$\Grm_{\geq i} = U^m_1\sqcup U^m_2$ into 
open and closed 
subsets, indexed by $m\in \Z_{\geq 1}$, that are $\Prm(L)$-invariant
on the left and $\Prm(\cO_L)$-invariant on the right, such that 
\[
\Grm_{\geq i+1} = \bigcup_{m\geq 1} U^m_1
\]
(in particular, each $U^m_1$ is a subset of $\Grm_{\geq i+1}$). 
\end{lemma}

\begin{proof} Let $\overline{S}_w:=\overline{\Prm(L)w\Prm(L)}$, 
a closed subset of $\Grm(L)$. We claim that it
is enough to find, for each 
$w\in \PWP$ with $\ell(w) = i$, 
decompositions $\Grm(L) = U_{w,1}^m \sqcup U_{w,2}^m$ into
open and closed subsets, indexed 
by $m\in \Z_{\geq 1}$, that are $\Prm(L)$-invariant
on the left and $\Prm(\cO_L)$-invariant on the right, 
such that 
\[
\overline{S}_w = \bigcap_{m\geq 1} U^m_{w,2}. 
\]
Indeed, once we have found such decompositions, we set 
\[
U_2^m:= \left(\bigcup_{\ell(w) = i} U^m_{w,2}\right)\cap \Grm_{\geq i},
\]
which is open and closed in $\Grm_{\geq i}$ 
because we have taken a finite union. The $U^m_2$ and their 
complements induce decompositions of $\Grm_{\geq i}$
with the desired properties. 

We now describe how to find
the decompositions for each $w$. 
Set $\mathscr{F}\ell:= \Prm \backslash \Grm$; 
this is a projective variety defined over $\cO_L$.  
The closed subset $\Prm(L)\backslash \overline{S}_w
\subset \Prm(L)\backslash \Grm(L)$ 
can be identified with the $\cO_L$-points of 
a closed Schubert subvariety $\overline{\cS}_w$ 
of $\mathscr{F}\ell$. 
For each $m\geq 1$, we consider the subset $\overline{V}_m
\subset \mathscr{F}\ell(\cO_L/\varpi^m_L)$ consisting of the 
points that satisfy modulo $\varpi^m_L$ 
the equations for $\overline{\cS}_w$. 
The preimage $V_m\subset \Prm(L)\backslash \Grm(L)$ of $\overline{V}_m$ 
is an open and closed subset that contains $\Prm(L)\backslash \overline{S}_w$, 
and is invariant under multiplication on
the right by the finite index normal subgroup of $\Prm(\cO_L)$ of elements that 
reduce to the identity modulo $\varpi_L^m$.
We intersect the translates of $V_m$ by finitely many coset
representatives in $\Prm(\cO_L)$ to obtain an open and closed
subset $W_m\subset \Prm(L)\backslash \Grm(L)$ 
that contains $\Prm(L)\backslash \overline{S}_w$, 
and is invariant under multiplication on
the right by $\Prm(\cO_L)$. Finally, we define
$U^m_{w,2}$ to be the preimage of $W_m$ in $\Grm(L)$. 
\end{proof}

Assume now that $\Grm =\GL_{2n}/L$ and that 
$\Prm$ is the standard upper-triangular parabolic with Levi $\GL_n\times \GL_n$. 
Write $\tilde{u}_L:=\mathrm{diag}(\varpi_{L},\dots,\varpi_{L}, 1,\dots,1)$, where the uniformizer
$\varpi_{L}$ of $\cO_L$ occurs $n$ times on the diagonal. Then 
$\Mrm(L)^+$ contains the monoid generated by $\Mrm(\cO_L)$ and $\tilde{u}_L$. %
This means that if $\pi \in D^+_{\sm}(\Mrm(L)^+\ltimes \Urm^0,\cO/\varpi^m)$, we can define $R\Gamma\left(\Urm^0,\pi\right) \in D^+_{\sm}(\Mrm(L)^+,\cO/\varpi^m)$ and the localization (inverting $\tilde{u}_L$) $\ord R\Gamma\left(\Urm^0,\pi\right) \in D^+_{\sm}(\Mrm(L)^+,\cO/\varpi^m)$.

Under these assumptions on $\mathrm{G}$ and $\mathrm{P}$, we now compute the complexes \[\ord R\Gamma\left(\Urm^0, I_{w}(\pi) \right)\]
in two special cases: 
for $w$ equal to either the longest element of $\PWP$ or to the identity element. 

\begin{lemma}\label{lem:longestelt1} Let $w^{\Prm}_0$ be the longest 
element of $\PWP$. Then:
\begin{enumerate}
\item $I^\circ_{w^{\Prm}_0}$ takes injectives to $\Gamma(\Urm^0,\ )$-acyclics. 
\item Let $\pi \in D^+_{\sm}(\Prm(L), \cO/\varpi^m)$. Then there is a natural isomorphism
\[
\ord R\Gamma(\Urm^0, I^\circ_{w^{\Prm}_0}(\pi))\toisom 
\ord R\Gamma(\Urm^0, I_{w^{\Prm}_0}(\pi)). 
\]
\end{enumerate}
\end{lemma}

\begin{proof} Since $w^{\Prm}_0$ is the longest element in $\PWP$, it normalizes
$\Mrm(L)$ (interchanging the two blocks) and therefore we have $S^\circ_{w^{\Prm}_0} = \Prm(L)w^{\Prm}_0\Urm^0$. 
The proof of~\cite[Lemma 5.3.4]{10author} applies verbatim;
for convenience, we reproduce it here. 

For the first part, let $\pi \in \mathrm{Mod}_{\sm}(\Prm(L), \cO/\varpi^m)$
and fix an $\cO/\varpi^m$-linear embedding $\pi\hookrightarrow I$, where $I$ 
is an injective $\cO/\varpi^m$-module. This gives rise to an embedding 
$\pi \hookrightarrow \mathrm{Ind}_{1}^{\Prm(L)} I$ of smooth 
$\cO/\varpi^m[\Prm(L)]$-modules. By~\cite[Lemma 2.1.10]{emordtwo}, 
it suffices to show that $I^\circ_{w^{\Prm}_0}\left(\mathrm{Ind}_{1}^{\Prm(L)} I\right)$ is 
an injective smooth $\cO/\varpi^m[\Urm^0]$-module. There is a natural 
$\Urm^0$-equivariant isomorphism 
\[
I^\circ_{w^{\Prm}_0}\left(\mathrm{Ind}_{1}^{\Prm(L)} I\right)\toisom 
\cC^\infty\left(\Prm(L)w^{\Prm}_0\Urm^0,I\right),
\]
where $\cC^\infty\left(\Prm(L)w^{\Prm}_0\Urm^0,I\right)$ denotes the set of functions $F:\Prm(L)w^{\Prm}_0\Urm^0\to I$
such that $u\mapsto (p\mapsto F(pw_0^{\Prm}u))$ defines a locally constant (equivalently, smooth) function $\Urm^0\to \mathrm{Ind}_1^{\Prm(L)}I$. 
This set is equipped with the usual action of $\Urm^0$ by right translation. 
The isomorphism sends a function $f\in I^\circ_{w^{\Prm}_0}\left(\mathrm{Ind}_{1}^{\Prm(L)} I\right)$ 
to $F\in \cC^\infty\left(\Prm(L)w^{\Prm}_0\Urm^0,I\right)$ given by $F(x):= f(x)(1)$.  
Since $\cC^\infty\left(\Prm(L)w^{\Prm}_0\Urm^0,I\right)$ is an injective
smooth $\cO/\varpi^m[\Urm^0]$-module, we conclude. 

For the second part, we first define an exact functor 
\[
J_{w^{\Prm}_0}: \mathrm{Mod}_{\sm}(\Prm(L), \cO/\varpi^m)\to 
\mathrm{Mod}_{\sm}(\Mrm(L)^+\ltimes \Urm^0, \cO/\varpi^m) 
\]
by the formula $J_{w^{\Prm}_0}(\pi) := I_{w^{\Prm}_0}(\pi)/I^\circ_{w^{\Prm}_0}(\pi)$. 
For each $\pi\in D^+_{\sm}(\Prm(L),\cO/\varpi^m)$, we have 
a distinguished triangle 
\[
\ord R\Gamma(\Urm^0, I^\circ_{w^{\Prm}_0}(\pi))\to \ord R\Gamma(\Urm^0, I_{w^{\Prm}_0}(\pi))
\to \ord R\Gamma(\Urm^0, J_{w^{\Prm}_0}(\pi))\xrightarrow{+1}. 
\]
It is enough to show that, for any $j\in \mathbb{Z}$, we have 
\[
\ord H^j(\Urm^0, J_{w^{\Prm}_0}(\pi)) = 0. 
\]
By direct computation, we see that $\tilde{u}_L$ acts locally nilpotently
on $J_{w^{\Prm}_0}(\pi)$: it shrinks the support of a function in $I_{w^{\Prm}_0}(\pi)$ 
towards $S^\circ_{w^{\Prm}_0}$. 
We conclude using the same argument as in~\cite[Lemma 3.3.1]{hauseux}. 
\end{proof}

\begin{lemma}\label{lem:longestelt2} Let $w^{\Prm}_0$ be the longest 
element of $\PWP$ and let $\pi \in D^+_{\sm}(\Prm(L),\cO/\varpi^m)$. 
Then there is a natural isomorphism 
\[
R\Gamma(\Urm^0, I^\circ_{w^{\Prm}_0}(\pi))\toisom \pi^{w^{\Prm}_0}
\]
in $D^+_{\sm}(\Mrm(L)^+, \cO/\varpi^m)$, where $m\in \Mrm(L)^+$ acts on 
$\pi^{w^{\Prm}_0}$ via the action of $w^{\Prm}_0m(w^{\Prm}_0)^{-1}$ on $\pi$. 
\end{lemma}

\begin{proof} By the first part of Lemma~\ref{lem:longestelt1}, it is
enough to show that there is a natural isomorphism 
\begin{equation}\label{eq:underived functors}
\Gamma(\Urm^0, I^\circ_{w^{\Prm}_0}(\pi))\toisom \pi^{w^{\Prm}_0}
\end{equation} 
of underived functors. The map sends an $\Urm^0$-invariant function 
$f: \Prm(L)w^{\Prm}_0\Urm^0\to \pi$ to the value $f(w^{\Prm}_0)\in \pi$. This 
is an isomorphism of $\cO/\varpi^m$-modules; it remains
to check that~\eqref{eq:underived functors} is $\Mrm(L)^+$-equivariant,
for the Hecke action of $\Mrm(L)^+$ on the LHS and the action twisted
by $w^{\Prm}_0$ on the RHS. In other words, we have to check that, 
for any $\Urm^0$-invariant function $f:  \Prm(L)w^{\Prm}_0\Urm^0\to \pi$ 
and any $m\in \Mrm(L)^+$, we have the equality 
\begin{equation}\label{eq:Hecke actions}
\sum_{\bar{u}\in \Urm^0/m\Urm^0m^{-1}} f(w^{\Prm}_0um) = w^{\Prm}_0m(w^{\Prm}_0)^{-1} f(w^{\Prm}_0). 
\end{equation}
This will hold if the only $\bar{u}\in \Urm^0/m\Urm^0m^{-1}$ 
that contributes to the LHS of~\eqref{eq:Hecke actions} is the identity element. 
We claim that $f(w^{\Prm}_0um) =0$ unless $u\in m\Urm^0m^{-1}$. Assume 
$f(w^{\Prm}_0um)\not = 0$; then $w^{\Prm}_0um \in \Prm(L)w^{\Prm}_0\Urm^0$. We write
$w^{\Prm}_0um = qw^{\Prm}_0u'$ with $q\in \Prm(L)$ and $u'\in \Urm^0$ and we obtain 
\[
u = (w^{\Prm}_0)^{-1} q w^{\Prm}_0 u'm^{-1} = ((w^{\Prm}_0)^{-1}qm' w^{\Prm}_0) (m u'm^{-1}), 
\]
with $m':= w^{\Prm}_0m^{-1}(w^{\Prm}_0)^{-1}\in \Mrm(L)$. But then  
$((w^{\Prm}_0)^{-1}qm' w^{\Prm}_0) \in (w^{\Prm}_0)^{-1}\Prm(L)w^{\Prm}_0\cap \Urm^0 = \{\mathrm{id}\}$
and we deduce $u\in m\Urm^0m^{-1}$.  
\end{proof}

Let $\chi: \Mrm(L)\to \cO^\times$ be the character
defined by the formula 
\[
\chi(m) = \frac{\mathrm{Nm}_{L/\Q_p}\det_{L}\left(\mathrm{Ad}(m)|_{\mathrm{Lie}\ U(L)}\right)^{-1}}
{|\mathrm{Nm}_{L/\Q_p}\det_{L}\left(\mathrm{Ad}(m)|_{\mathrm{Lie}\ U(L)}\right)|_p}.
\]

\begin{lemma}\label{lem:identityelt}
Let $\pi \in D^+_{\mathrm{sm}}(\Mrm(L),\cO/\varpi^m)$. 
Then there is a natural isomorphism 
\[
\ord R\Gamma\left(\Urm^0, \mathrm{Inf}_{\Mrm(L)^+}^{\Mrm(L)^+\ltimes \Urm^0}\pi \right)
\toisom \cO/\varpi^m(\chi) \otimes_{\cO/\varpi^m} \pi [-\mathrm{rk}_{\Z_p}\Urm^0]
\]
in $D^+_{\sm}(\Mrm(L)^+, \cO/\varpi^m)$. 

\end{lemma}

\begin{proof} Since $\Urm^0$ acts trivially on $\pi$, we have 
\[
\ord R\Gamma\left(\Urm^0, \mathrm{Inf}_{\Mrm(L)^+}^{\Mrm(L)^+\ltimes \Urm^0}\pi \right)
\toisom \pi\otimes_{\cO/\varpi^m} \ord R\Gamma\left(\Urm^0, \cO/\varpi^m \right). 
\]
By the proof of~\cite[Lemma 5.3.7]{10author}, the continuous cohomology groups
of $\Urm^0$ vanish below the top degree $\ell:=\mathrm{rk}_{\Z_p}\Urm^0$ after applying ordinary parts. 
It remains to show that 
\[
H^{\ell}(\Urm^0, \cO/\varpi^m) \simeq \cO/\varpi^m(\chi)
\]
as a representation of $\Mrm(L)^+$. This follows, just like~\cite[Prop. 3.1.8]{hauseux}, 
from the natural isomorphism in~\cite[Prop. 3.5.6]{emordtwo} and 
from the explicit description of the corestriction map in~\cite[Lemma 3.5.10]{emordtwo}. 
\end{proof}

\begin{remark}\label{rem:specialise to id} The functor $I_{\mathrm{id}}$ 
is the identity on $D^+_{\sm}(\Prm(L), \cO/\varpi^m)$
and the functor $I^\circ_{\mathrm{id}}$ is the natural restriction
$D^+_{\sm}(\Prm(L), \cO/\varpi^m)\to 
D^+_{\sm}(\Mrm(L)^+\ltimes \Urm^0, \cO/\varpi^m)$. 
Note also that in this case we have $S^\circ_{\mathrm{id}} = S_{\mathrm{id}}$, 
so that we have a natural identification
\[
R\Gamma(\Urm^0, I^\circ_{\mathrm{id}}(\ ))\toisom R\Gamma(\Urm^0, I_{\mathrm{id}}(\ ))
\]
in $D^+_{\sm}(\Mrm(L)^+, \cO/\varpi^m)$.  
Furthermore, if we start with $\pi \in D^+_{\mathrm{sm}}(\Mrm(L),\cO/\varpi^m)$, we can identify 
$\mathrm{Inf}_{\Mrm(L)^+}^{\Mrm(L)^+\ltimes \Urm^0}\pi$
with the functor $I^\circ_{\mathrm{id}}$ applied to $\mathrm{Inf}_{\Mrm(L)}^{\Prm(L)}\pi$. 
\end{remark}

\begin{cor}\label{cor:identityelt}
Let $\pi \in D^+_{\mathrm{sm}}(\Mrm(L),\cO/\varpi^m)$.
Then there is a natural isomorphism 
\[
\ord R\Gamma\left(\Urm^0, I_{\mathrm{id}}\left(\mathrm{Inf}_{\Mrm(L)}^{\Prm(L)}\pi\right) \right)
\toisom \cO/\varpi^m(\chi) \otimes_{\cO/\varpi^m} \pi [-\mathrm{rk}_{\Z_p}\Urm^0]
\]
in $D^+_{\sm}(\Mrm(L)^+, \cO/\varpi^m)$. 
\end{cor}

We now apply our results in the setting of Section \ref{sec:P-ord Hida}, so we adopt the notation of that section. We fix a place $\bar{v} \in \bar{S}$ and a standard parabolic $Q_{\bar{v}} \subset P_{\bar{v}}$. We apply the results of this section with $\Grm = \tG_{\cO_{F^+_{\bar{v}}}}$, $\Prm = P_{\bar{v}}$, $\Mrm = G_{\cO_{F^+_{\bar{v}}}}$ and $\Urm = U_{\bar{v}}$. Set $K_{\bar{v}} = \cQ_{\bar{v}}\cap G(F^+_{\bar{v}})$. Note that $\Delta_{\bar{v}}^{\cQ_{\bar{v}},+}$ is a submonoid of $G(F^+_{\bar{v}})^+$.  The following proposition summarizes the key result 
of this subsection. 

\begin{prop}\label{prop:ordsubquotients} 
Let $\pi\in D^+_{\sm}(G(F^+_{\bar{v}}),\cO/\varpi^m)$
and let $\cV$ be a finite free $\cO/\varpi^m$-module equipped with a smooth representation of $\Delta^{\cQ_{\bar{v}},+}_{\bar{v}}$ with $\tilde{u}_{\tilde{v},n}$ acting invertibly, which we inflate to an action of $\widetilde{\Delta}^{\cQ_{\bar{v}}}_{\bar{v},P} = \Delta^{\cQ_{\bar{v}},+}_{\bar{v}} \ltimes U^0_{\bar{v}}$.
Then 
\[
\ord_0 R^j\Gamma\left(K_{\bar{v}} \ltimes U^0_{\bar{v}}, \mathrm{Ind}_{P(F^+_{\bar{v}})}^{\tG(F^+_{\bar{v}})} \pi \otimes \cV\right)
\]
admits as $\cH(\Delta^{\cQ_{\bar{v}}}_{\bar{v}},K_{\bar{v}})$-module subquotients 
\[
R^j\Gamma(K_{\bar{v}}, \pi^{w^P_0}\otimes \cV)\ \text{and}\ 
R^{j-\mathrm{rk}_{\Z_p}U^0_{\bar{v}}}\Gamma(K_{\bar{v}}, \pi \otimes \cO/\varpi^m(\chi)\otimes \cV).
\]
\end{prop}

\begin{proof} By applying the exact functor $\mathrm{ord}_0$ to the 
short exact sequences in Proposition~\ref{prop:zeroconnecthom}, 
we see that 
$\ord_0 R^j\Gamma\left(K_{\bar{v}} \ltimes U^0_{\bar{v}}, \mathrm{Ind}_{P(F^+_{\bar{v}})}^{\tG(F^+_{\bar{v}})} \pi \otimes \cV\right)$ 
admits as subquotients 
$\mathrm{ord}_0R^j\Gamma(K_{\bar{v}} \ltimes U^0_{\bar{v}}, I_{w}(\pi)\otimes \cV)$
for $w = w^P_0$ and for $w=\mathrm{id}$. We have an
isomorphism of functors 
\[
\mathrm{ord}_0\circ R\Gamma(K_{\bar{v}}, \ ) \simeq R\Gamma(K_{\bar{v}}, \ )\circ \mathrm{ord},
\]
as in \S\ref{sec:ordQlevel}. By Lemma~\ref{lem:longestelt1} and
Lemma~\ref{lem:longestelt2}, we have an isomorphism 
\[
\ord R\Gamma(U_{\bar{v}}^0, I_{w^P_0}(\pi) \otimes \cV)\simeq 
\ord R\Gamma(U_{\bar{v}}^0, I_{w^P_0}(\pi))\otimes \cV \simeq \pi^{w^P_0}\otimes \cV. 
\]
By Corollary~\ref{cor:identityelt}, we have an isomorphism
\[
\ord R\Gamma(U_{\bar{v}}^0, I_{\mathrm{id}}(\pi)\otimes \cV) \simeq 
\ord R\Gamma(U_{\bar{v}}^0, I_{\mathrm{id}}(\pi))\otimes \cV \simeq 
\pi \otimes \cO/\varpi^m(\chi ) \otimes \cV [-\mathrm{rk}_{\Z_p}U_{\bar{v}}^0]. 
\] 
\end{proof}

\begin{remark}\label{rem:ordsubquotients}
Proposition~\ref{prop:ordsubquotients} still holds true, with the same proof, if we replace $\Delta^{\mathcal{Q}_{\bar{v}},+}_{\bar{v}}$ with a more general monoid $\Sigma^+_{\bar{v}}\subset G(F^+_{\bar{v}})$, as in Remark~\ref{rem:more general monoid}. 
\end{remark}

We have a corollary which will be applied with `dual' coefficients (cf.~\S\ref{sec:dual coeffs prelims}). This can be viewed as a computation of the ordinary part with respect to the opposite parabolic $\overline{P}$ applied to a parabolic induction from $P$.

\begin{cor}\label{cor:ordsubquotients dual}
Let $\pi\in D^+_{\sm}(G(F^+_{\bar{v}}),\cO/\varpi^m)$
and let $\cV$ be a finite free $\cO/\varpi^m$-module equipped with a smooth representation of $\left(\Delta^{\cQ_{\bar{v}},+}_{\bar{v}}\right)^{-1}$ with $\tilde{u}_{\tilde{v},n}^{-1}$ acting trivially, which we inflate to an action of $\left(\Delta^{\cQ_{\bar{v}},+}_{\bar{v}}\right)^{-1}\ltimes \bar{U}^1_{\bar{v}}$.
Then 
\[
\ord_0^\vee R^j\Gamma\left(K_{\bar{v}} \ltimes \overline{U}^1_{\bar{v}}, \mathrm{Ind}_{P(F^+_{\bar{v}})}^{\tG(F^+_{\bar{v}})} \pi \otimes \cV\right)
\]
admits $R^j\Gamma(K_{\bar{v}}, \pi \otimes \cV)$ as a $\cH(\left(\Delta^{\cQ_{\bar{v}}}_{\bar{v}}\right)^{-1},K_{\bar{v}})$-module subquotient.
\end{cor}
\begin{proof}
	We deduce the corollary from Remark~\ref{rem:ordsubquotients} by twisting. Set $g:=\tilde{u}_{\tilde{v},n}^{-1}w^P_0$ and $\Sigma^+_{\bar{v}}:=g(\Delta^{\cQ_{\bar{v}},+}_{\bar{v}})^{-1}g^{-1}$. Multiplication by $g$ induces an isomorphism \[R^j\Gamma\left(K_{\bar{v}} \ltimes \overline{U}^1_{\bar{v}}, \mathrm{Ind}_{P(F^+_{\bar{v}})}^{\tG(F^+_{\bar{v}})} \pi \otimes \cV\right) \toisom R^j\Gamma\left(K_{\bar{v}}^{w^P_0} \ltimes {U}^0_{\bar{v}}, \mathrm{Ind}_{P(F^+_{\bar{v}})}^{\tG(F^+_{\bar{v}})}\pi \otimes \cV^{w^P_0}\right)\] where the action of $\cH(\left(\Delta^{\cQ_{\bar{v}},+}_{\bar{v}}\right)^{-1},K_{\bar{v}})$ on the left is identified with the action of $\cH(\Sigma^+_{\bar{v}},K^{w_0^P}_{\bar{v}})$ on the right. 
    We recall from the proof of Lemma \ref{lem:control level dual} that $\overline{Q}^{w^P_0}_{\bar{v}}$ is the standard parabolic with Levi subgroup  $Q_{\bar{v}}^{w^P_0}\cap G(F^+_{\bar{v}})$. Now Proposition \ref{prop:ordsubquotients} tells us that we have a subquotient $R^j\Gamma(K_{\bar{v}}^{w^P_0}, \pi^{w^P_0}\otimes \cV^{w^P_0})$. After undoing the twisting by $g$, this can be identified with $R^j\Gamma(K_{\bar{v}}, \pi \otimes \cV)$ as a $\cH(\left(\Delta^{\cQ_{\bar{v}}}_{\bar{v}}\right)^{-1},K_{\bar{v}})$-module.\end{proof}

\subsubsection{Parabolic induction and cohomology} 
We return to the general situation of \S\ref{sec:local P-ord computation}, for a split reductive group $\Grm$ defined over (the ring of integers in) a local field $L$. In addition, we consider a compact Hausdorff
space $X$, equipped with a continuous action of the locally profinite group 
$\mathrm{P}(L)$. Set $\mathrm{K} = \Grm(\cO_L)$ and $\mathrm{K}_{\mathrm{P}}:= \mathrm{K} \cap \mathrm{P}(L)$. 
We assume that $X$ is a free $\Krm_{\Prm}$-space, 
in the sense of~\cite[Def. 2.23]{new-tho}. 

Denote by $X\times_{\mathrm{P}}\mathrm{G}$ the quotient
of $X\times \mathrm{G}(L)$ by the right action $(x,g)\cdot p = (xp,p^{-1}g)$ of $\mathrm{P}(L)$. Right multiplication by $\Grm(L)$ on itself gives $X\times_{\mathrm{P}}\mathrm{G}$ a right action of $\Grm(L)$.

The Iwasawa decomposition $\mathrm{G}(L) = \mathrm{K}\mathrm{P}(L)$  
implies
that there exists a $\Krm$-equivariant homeomorphism 
\[
X\times_{\mathrm{K}_\mathrm{P}} \mathrm{K} 
\toisom X\times_{\mathrm{P}}\mathrm{G}. 
\]
The LHS is visibly a compact Hausdorff space. 
Therefore, it makes sense to consider 
$R\Gamma\left(\mathrm{X}
\times_{\mathrm{P}} \mathrm{G}, \cO/\varpi^m \right)$ as
an element in $D^+_{\mathrm{sm}}(\mathrm{G}(L), \cO/\varpi^m)$. 

\begin{lemma}\label{lem:induction lemma} 
We have a natural isomorphism in 
$D^+_{\mathrm{sm}}(\mathrm{G}(L), \cO/\varpi^m)$ 
\begin{equation}\label{eq:induction formula}
R\Gamma\left(\mathrm{X}\times_{\mathrm{P}} \mathrm{G}, 
\cO/\varpi^m \right) \toisom \mathrm{Ind}_{\mathrm{P}(L)}^{\mathrm{G}(L)}\ 
R\Gamma(X,\cO/\varpi^m). 
\end{equation}
\end{lemma}

\begin{proof}

First, we explain how to construct a natural map 
from the LHS of~\eqref{eq:induction formula} 
to the RHS, then we observe that it is enough 
to prove the map is an isomorphism after restriction to 
$D^+_{\mathrm{sm}}(\Krm, \cO/\varpi^m)$, then we prove the latter. 

Consider the $\Prm(L)$-equivariant map $X\to \mathrm{X}
\times_{\mathrm{P}} \mathrm{G}$
given by $x\mapsto (x,1)$. We get an induced 
morphism 
\[
R\Gamma(X\times_{\Prm}\Grm, \cO/\varpi^m) 
\to R\Gamma(X, \cO/\varpi^m)
\]
in $D^+_{\mathrm{sm}}(\Prm(L), \cO/\varpi^m)$, 
which induces a morphism 
\[
R\Gamma(X\times_{\Prm}\Grm, \cO/\varpi^m) \to
\mathrm{Ind}_{\Prm(L)}^{\Grm(L)}\ R\Gamma(X, \cO/\varpi^m)
\]
in $D^+_{\mathrm{sm}}(\Grm(L), \cO/\varpi^m)$, 
by Frobenius reciprocity for smooth representations (also 
using that $\mathrm{Ind}_{\Prm(L)}^{\Grm(L)}$ is an exact functor). 

The same morphism can be constructed if we work with 
$\Krm_{\Prm}\subset \Krm$ replacing $\Prm\subset \Grm$. 
We now observe that we have a commutative 
diagram in $D^+_{\mathrm{sm}}(\Krm, \cO/\varpi^m)$ 
\begin{equation}\label{eq:commutative induction}
\xymatrix{\mathrm{Res}^{\Grm(L)}_{\Krm}\circ 
R\Gamma\left(X\times_{\Prm} \Grm, \cO/\varpi^m \right)\ar[r]\ar[d]^{\cong} & 
\mathrm{Res}^{\Grm(L)}_{\Krm}\circ \mathrm{Ind}_{\Prm(L)}^{\Grm(L)}\ 
R\Gamma (X, \cO/\varpi^m)\ar[d]^{\cong} \\
R\Gamma (X\times_{\Krm_{\Prm}} \Krm, \cO/\varpi^m)\ar[r] &
\mathrm{Ind}_{\Krm_{\Prm}}^{\Krm}\ R\Gamma (X, \cO/\varpi^m)}
\end{equation}
whose vertical arrows are isomorphisms by the Iwasawa
decomposition. 

We have a diagram with Cartesian outer square
\[
\xymatrix{\ & X\times \Krm\ar[dl]^{\varphi_2}\ar[dr]^{\varphi_1} & \ 
\\ X\times_{\Krm_{\Prm}} \Krm\ar[dr]^{\phi_2} & \ & X\ar[dl]^{\phi_1}\ar[ll]^{\phi} \\ \ & X/{\Krm}_{\Prm} &\ }
\]
The map $\phi$ is given by $x \mapsto (x,1)$; the top triangle is not commutative, 
but the bottom one is.  
Since $X\to X/\Krm_{\Prm}$ is a $\Krm_{\Prm}$-torsor, the same holds true 
for $X\times \Krm \to X\times_{\Krm_\Prm} \Krm$. This also implies that 
$X\times_{\Krm_{\Prm}}\Krm \to X/\Krm_{\Prm}$ is a $\Krm$-torsor. 
By~\cite[Lemma 2.24]{new-tho},
we have inverse equivalences of categories $(\varphi_2^*, \varphi_{2,*}^{\Krm_{\Prm}})$
between $\Sh_{\Krm}(X\times_{\Krm_{\Prm}}\Krm)$ and $\Sh_{\Krm_{\Prm}\times \Krm}(X\times \Krm)$. 
The analogous statements are true for $\varphi_1$, $\phi_1$ and
$\phi_2$. 

The bottom horizontal arrow in~\eqref{eq:commutative induction}
is induced by $\phi^*$. Let $\cO/\varpi^m \toisom \cI^{\bullet}$ be an injective 
resolution in $\mathrm{Sh}(X/\Krm_{\Prm})$. We claim that $\phi^*$ 
induces by Frobenius reciprocity a term-wise $\Krm$-equivariant 
isomorphism of complexes
\begin{equation}\label{eq:iso induction}
\Gamma (X\times_{\Krm_{\Prm}}\Krm, \phi_2^*\cI^\bullet) \toisom 
\mathrm{Ind}_{\Krm_{\Prm}}^{\Krm}\ \Gamma(X, \phi_1^*\cI^\bullet).
\end{equation}
To see the claim, note that $\phi_2^* = \varphi^{\Krm_{\Prm}}_{2,*} \varphi_2^* 
\phi_2^* = \varphi_{2,*}^{\Krm_{\Prm}}\varphi_1^*\phi_1^*$ 
and we can identify $\phi^*$ with the restriction to the identity of a function in
\[
\Gamma \left(X\times_{\Krm_{\Prm}}\Krm,  
\varphi_{2,*}^{\Krm_{\Prm}}\varphi_1^*(\phi_1^*\cI^\bullet)\right) 
\toisom
\mathrm{Ind}_{\Krm_{\Prm}}^{\Krm}\ \Gamma (X, \phi_1^*\cI^\bullet).
\]
The map in~\eqref{eq:iso induction} is precisely the horizontal map in the 
bottom row of~\eqref{eq:commutative induction}, which also implies
that the horizontal map in the top row is an isomorphism. 
\end{proof}

We will also make use of the following lemma about the interaction
between smooth induction and restriction. 

\begin{lemma}\label{lem:tensor identity} Let $\pi_1\in \mathrm{Rep}_{\sm}
(\mathrm{K}_{\mathrm{P}}, \cO/\varpi^m)$ and $\pi_2\in \mathrm{Rep}_{\sm}(\mathrm{K},\cO/\varpi^m)$ 
with the property that $\pi_2$ takes values in a finite free $\cO/\varpi^m$-module. There is a natural isomorphism 
\[
\left(\mathrm{Ind}_{\mathrm{K}_{\mathrm{P}}}^{\mathrm{K}}\pi_1\right)\otimes_{\cO/\varpi^m} \pi_2 \toisom 
\mathrm{Ind}_{\mathrm{K}_{\mathrm{P}}}^{\mathrm{K}} \left(\pi_1 \otimes_{\cO/\varpi^m} \mathrm{Res}_{\mathrm{K}_{\mathrm{P}}}^{\mathrm{K}} \pi_2\right)
\]
in $\mathrm{Rep}_{\sm}(\mathrm{K},\cO/\varpi^m)$. 
\end{lemma}

\begin{proof}

Since $\pi_2$ is a finite $\cO/\varpi^m$-module, there exists a compact open normal subgroup $U\subset \mathrm{K}$ that acts trivially on $\pi_2$.
One can define a natural map 
\[
\Phi: \left(\mathrm{Ind}_{\mathrm{K}_{\mathrm{P}}}^{\mathrm{K}}\pi_1\right)\otimes_{\cO/\varpi^m} \pi_2 \to
\mathrm{Ind}_{\mathrm{K}_{\mathrm{P}}}^{\mathrm{K}} \left(\pi_1 \otimes_{\cO/\varpi^m} \mathrm{Res}_{\mathrm{K}_{\mathrm{P}}}^{\mathrm{K}} \pi_2\right)
\]
by extending $\cO/\varpi^m$-linearly from the formula 
\[
f\otimes v \mapsto (g \mapsto f(g)\otimes \pi_2(g)v).
\]
One can check that $\Phi$ lands in $\mathrm{Ind}_{\mathrm{K}_{\mathrm{P}}}^{\mathrm{K}} \left(\pi_1 \otimes_{\cO/\varpi^m} \mathrm{Res}_{\mathrm{K}_{\mathrm{P}}}^{\mathrm{K}} \pi_2\right)$ and is a $\mathrm{K}$-equivariant homomorphism. To see that $\Phi$ does indeed take values in smooth functions, note that there is a compact open subgroup $U_f$ that stabilizes $f$. Then 
\[
f(gu)\otimes \pi_2(gu)v = f(g)\otimes \pi_2(g)v\ \forall u\in U_f\cap U\ \mathrm{and}\ g\in \mathrm{K}, 
\]
which shows that each $\Phi(f\otimes v)$ is a smooth function on $\mathrm{K}$.  

It remains to show that $\Phi$ is an isomorphism. For $m\in \Z_{\geq 1}$, we can find a decreasing sequence of compact open normal subgroups $U_m$ of $\mathrm{K}$ with $U_1 = U$ and that form a basis of neighbourhoods of identity in $\mathrm{K}$. It is enough to prove that $\Phi$ induces an isomorphism on the level of $U_m$-invariants for each $m$. Note that, for each $m$, the set of double cosets $\mathrm{K}_{\mathrm{P}}\backslash \mathrm{K}/U_m$ is finite, and for any smooth representation $\pi$ of $\mathrm{K}_{\mathrm{P}}$, we have an isomorphism 
\[
\left(\mathrm{Ind}_{\mathrm{K}_{\mathrm{P}}}^{\mathrm{K}}\pi\right)^{U_m} \toisom  \bigoplus_{\gamma\in \mathrm{K}_{\mathrm{P}}\backslash \mathrm{K}/U_m} 
\pi^{\mathrm{K}_{\mathrm{P}}\cap \gamma U_m\gamma^{-1}}, 
f\mapsto (f(\gamma))_{\gamma\in \mathrm{K}_{\mathrm{P}}\backslash \mathrm{K}/U_m}.
\]
Note also that each $\mathrm{K}_{\mathrm{P}}\cap \gamma U_m \gamma^{-1} = \mathrm{K}_{\mathrm{P}}\cap U_m$ acts trivially on $\pi_2$. 
We obtain a commutative diagram 
\[
\xymatrix{\left(\mathrm{Ind}_{\mathrm{K}_{\mathrm{P}}}^{\mathrm{K}}\pi_1\right)^{U_m}\otimes \pi_2 
\ar[r]_{\cong}\ar[d] & \bigoplus_{\gamma\in \mathrm{K}_{\mathrm{P}}\backslash \mathrm{K}/U_m} 
(\pi_1)^{\mathrm{K}_{\mathrm{P}}\cap \gamma U_m\gamma^{-1}}\otimes \pi_2 \ar[d] \\ 
\left(\mathrm{Ind}_{\mathrm{K}_{\mathrm{P}}}^{\mathrm{K}} \left(\pi_1 \otimes_{\cO/\varpi^m} 
\mathrm{Res}_{\mathrm{K}_{\mathrm{P}}}^{\mathrm{K}} \pi_2\right)\right)^{U_m}\ar[r]_{\cong}
 &  \bigoplus_{\gamma\in \mathrm{K}_{\mathrm{P}}\backslash \mathrm{K}/U_m} 
(\pi_1)^{\mathrm{K}_{\mathrm{P}}\cap \gamma U_m\gamma^{-1}}\otimes \pi_2}, 
\]
where the left vertical map is $\Phi$ and the 
right vertical map is $\mathrm{id} \otimes \gamma$ in the component corresponding to $\gamma$. This is visibly an isomorphism. 

\end{proof}

\subsubsection{A computation of group cohomology}

Keep the notation from the previous subsection, except we now write $\Krm$ for $\Mrm(\cO_L)$. Inside $\Krm$, consider the congruence 
subgroups $\Krm_{m}:= \{k\in \Krm\mid k\equiv \mathrm{id} \pmod{\varpi_L^m}\}$ 
indexed by $m\in \mathbb{Z}_{\geq 1}$. The group $\Urm^0$ is equipped 
with the adjoint action of $\Krm$, so that we can view the 
continuous cohomology $R\Gamma(\Urm^0, \ )$ 
as a functor $D^+_{\mathrm{sm}}(\Krm\ltimes \Urm^0, \cO/\varpi^m)\to 
D^+_{\mathrm{sm}}(\Krm, \cO/\varpi^m)$.  

\begin{lemma}\label{lem:splitting cohomology} Assume that $\Urm^0$ is a finite free $\mathbb{Z}_p$-module.\footnote{This assumption will be satisfied in the case we need, namely when $\mathrm{P}$ is the Siegel parabolic in $\widetilde{G}$.} For any $m\in \Z_{\geq 1}$, there 
exists $M = M(m)\geq m$ such that 
\[
R\Gamma(\Urm^0, \cO/\varpi^m) \simeq
\bigoplus^{\mathrm{rk}_{\Z_p}\Urm^0}_{i=0}  H^i(\Urm^0, \cO/\varpi^m)[-i]
\]
as an object in $D^+_{\mathrm{sm}}(\Krm_M, \cO/\varpi^m)$. Moreover, 
each $H^i(\Urm^0, \cO/\varpi^m)$ is non-zero and equipped with 
the trivial action of $\Krm_M$. 
\end{lemma}

\begin{proof} By~\cite[Prop. 2.1.11]{emordtwo}, an injective
smooth representation of $\Krm \ltimes \Urm^0$ is also injective
as a representation of $\Urm^0$. Therefore, if we apply the forgetful 
functor to $R\Gamma(\Urm^0, \cO/\varpi^m)$ in order to view it as an 
object of $D^+(\cO/\varpi^m)$, 
we obtain the continuous group cohomology of $\Urm^0$ with coefficients in $\cO/\varpi^m$. 
If we forget the $\Krm$-action, $\Urm^0$ is a free $\Z_p$-module. 
Therefore, we can compute its continuous group cohomology via
a Koszul complex, using~\cite[Lemma 7.3, part (ii)]{bms}, for example. We see
that all the differentials in the Koszul complex vanish, so we obtain
\[
R\Gamma(\Urm^0, \cO/\varpi^m) \simeq
\bigoplus_{i}  H^i(\Urm^0, \cO/\varpi^m)[-i]
\]
as an object in $D_{\mathrm{perf}}(\cO/\varpi^m)$, the full subcategory of 
$D^+(\cO/\varpi^m)$ consisting of perfect complexes. Moreover, we can identify 
each $H^i(\Urm^0, \cO/\varpi^m)$ with \[
\wedge^i_{\mathcal{O}/\varpi^m} \mathrm{Hom}_{\mathrm{cont}}(\Urm^0, \cO/\varpi^m),\]
which 
shows that the cohomology groups are non-zero precisely in the range 
$[0, \mathrm{rk}_{\Z_p}\Urm^0]$. 

Let $D_{\mathrm{perf}}(\Krm_{(M)}, \cO/\varpi^m)$ be the full subcategory 
of $D^+_{\mathrm{sm}}(\Krm_{(M)}, \cO/\varpi^m)$ whose essential image 
under the forgetful functor is $D_{\mathrm{perf}}(\cO/\varpi^m)$. It is enough
to show that, for every object $A$ of $D_{\mathrm{perf}}(\Krm, \cO/\varpi^m)$, 
there exists $M\geq m$ such that the restriction $A_M$ of $A$ 
to $D_{\mathrm{perf}}(\Krm_{M}, \cO/\varpi^m)$
is isomorphic to the constant object $B_M$ of $D_{\mathrm{perf}}(\Krm_{M}, \cO/\varpi^m)$ 
corresponding to the image $B$ of $A$ in $D_{\mathrm{perf}}(\cO/\varpi^m)$. 
In turn, this would follow from showing that the natural map 
\begin{equation}\label{eq:fully faithful}
\mathrm{colim}_{M} \left(\Hom_{D_{\mathrm{perf}}(\Krm_{M}, \cO/\varpi^m)}(B_M, A_M)\right)\to 
\Hom_{D_{\mathrm{perf}}(\cO/\varpi^m)}(B, A) 
\end{equation}
is a bijection. Indeed, on the RHS we have the identity morphism, which 
must correspond to a morphism on the LHS for some $M\geq m$. 
We can check that this morphism in $D_{\mathrm{perf}}(\Krm_{M}, \cO/\varpi^m)$
is an isomorphism after applying the forgetful functor to $D_{\mathrm{perf}}(\cO/\varpi^m)$. 

The objects of $D_{\mathrm{perf}}(\cO/\varpi^m)$ are precisely the dualizable %
objects of $D^+(\cO/\varpi^m)$ and the dual is given by applying
$\Hom_{\cO/\varpi^m}(\ , \cO/\varpi^m)$. We have the adjunction (cf.~\cite[\href{https://stacks.math.columbia.edu/tag/07VI}{Lemma 07VI}]{stacks-project})
\[
\Hom_{D_{\mathrm{perf}}(\cO/\varpi^m)}(B, A) = 
\Hom_{D_{\mathrm{perf}}(\cO/\varpi^m)}(\cO/\varpi^m, B^\vee \otimes^{\mathbb{L}} A) = H^0(B^\vee \otimes^{\mathbb{L}} A).  
\]

Now, we have an adjunction \[\mathrm{Hom}_{D_{\mathrm{perf}}(\Krm_M, \cO/\varpi^m)}\left(B_M, A_M\right) = \mathrm{Hom}_{D^+( \cO/\varpi^m)} \left(B, R\Gamma(\Krm_M,A_M)\right)\] and we can identify the right hand side with $H^0(B^\vee \otimes^{\mathbb{L}} R\Gamma(\Krm_M,A_M))$ (since $B \in D_{\mathrm{perf}}(\cO/\varpi^m)$). Fixing a strictly perfect complex representing $B$, and equipping it with the trivial $\Krm_M$ action to get a complex representing $B_M$, we identify $B^\vee \otimes^{\mathbb{L}} R\Gamma(\Krm_M,A_M) = R\Gamma(\Krm_M,B_M^\vee\otimes^{\mathbb{L}} A_M)$.

Therefore, the bijectivity in~\eqref{eq:fully faithful} 
reduces to the statement that 
\[
\mathrm{colim}_{M} H^0(\Krm_M, B_M^\vee \otimes^\mathbb{L} A_M)
\]
(degree $0$ hypercohomology)
is equal to the degree $0$ cohomology of the underlying complex of $B^\vee \otimes^\mathbb{L} A$. 
This follows from the fact that $B^\vee \otimes^\mathbb{L} A$ 
is a complex of smooth representations of $\Krm$.  
\end{proof}

We will apply the preceding lemma when proving the crucial Proposition \ref{prop:degree shifting}. In that proof, we will also need the following consequence of the Artin--Rees lemma:
\begin{lemma}\label{lem:artin-rees application} Let 
$N$ be a finite $\cO$-module and
$M$ be a subquotient of $N$. For any positive integer 
$m\in \Z_{\geq 1}$,
there exists an integer $m' \geq m$ 
such that $M/\varpi^mM$ is a subquotient of 
$N/\varpi^{m'}N$. 
\end{lemma}
\begin{proof} If $M$ is a quotient of $N$, then 
we may take $m' = m$. Therefore, it suffices to 
consider the case when $M \hookrightarrow N$ 
is a subobject. By the Artin--Rees lemma, cf.~\cite[Cor. 10.10]{atiyah-macdonald}, 
there exists $k\gg 0$ such that 
\[
\varpi^{m+k}N \cap M = \varpi^m (\varpi^kN \cap M)\subseteq \varpi^m M. 
\]
Set $m':=m+k$. We then have an inclusion 
$M/(\varpi^{m'}N \cap M)\hookrightarrow N/\varpi^{m'}N$
and a surjection $M/(\varpi^{m'}N \cap M) 
\twoheadrightarrow M/\varpi^mM$. 
\end{proof}

\section{Determinants, $P$-ordinary representations and deformation rings}\label{sec:determinants}

\subsection{The $P$-ordinary condition on the Galois side}\label{sec:P-ord}
We place ourselves in the setting of \S\ref{sec: unitary group}, so we have a CM field $F$, unitary group $\tG$, etc. We fix a standard parabolic $Q_{\bar{v}} \subset P_{F^+_{\bar{v}}}$. After fixing a place $v|\bar{v}$, this corresponds under $\iota_v$ to a standard parabolic $P_{n_1,\ldots,n_t} \subset \GL_{2n}$, with $(n_1,\ldots,n_t)$ a partition of $2n$ refining $(n,n)$.  Recall that $Q_{\bar{v}}=M_{\bar{v}}\ltimes N_{\bar{v}}$ is a Levi decomposition for $Q_{\bar{v}}$. We use the notation of \S\ref{sec:explicit Hecke operators}, so we have a parahoric subgroup $\cQ_{\bar{v}}$ associated with $Q_{\bar{v}}$.

 For $1 \le k \le t$, we define a cocharacter $\nu_k \in X_{Q_{\bar{v}}}$ by \[\nu_k(\varpi_{\bar{v}}):= \iota_v^{-1}\diag(\varpi_v,\ldots,\varpi_v,1,\ldots,1) \in \tDelta_{\bar{v}}^{\cQ_{\bar{v}}},\] where there are $n_1 + \cdots + n_k$ entries equal to $\varpi_v$, and denote the Hecke operator $[\cQ_{\bar{v}} \nu_k(\varpi_{\bar{v}}) \cQ_{\bar{v}}] \in \cH(\tDelta_{\bar{v}}^{\cQ_{\bar{v}}}, \cQ_{\bar{v}})$ by $\widetilde{U}_v^k$. It follows from Lemma \ref{lem:Q Hecke algebra} that \[\cH(\tDelta_{\bar{v}}^{\cQ_{\bar{v}}}, \cQ_{\bar{v}}) \cong \Z[\widetilde{U}_v^1,\ldots,\widetilde{U}_v^{t-1},(\widetilde{U}_v^{t})^{\pm 1}].\] 
If $\tilde{\lambda} \in (\ZZ_+^{2n})^{\Hom(F^+,\Qpbar)}$ is a dominant weight for $\tG$ and $\sigma$ is a smooth $\Qpbar$-representation of $\tG(F^+_{\bar{v}})$, we define the $\tilde{\lambda}$-\emph{rescaled} action of  $\cH(\tDelta_{\bar{v}}^{\cQ_{\bar{v}}}, \cQ_{\bar{v}})$ on $\sigma^{\cQ_{\bar{v}}}$ to be given by multiplying the usual double coset operator action of $[\cQ_{\bar{v}}g\cQ_{\bar{v}}]$ by $\tilde{\alpha}^{\cQ_{\bar{v}}}_{\bar{v}}(g)^{-1}$ (cf.~Lemma \ref{lem:twisted action}).

\begin{defn}
Let $\pi$ be a cuspidal automorphic representation of $\widetilde{G}(\A_{F^+})$ and fix an isomorphism $\iota: \overline{\Q}_p\toisom \C$. Let $\tilde{\lambda} \in (\ZZ_+^{2n})^{\Hom(F^+,\Qpbar)}$ be a dominant weight for $\tG$. We say that $\pi$ is $\iota$-$Q_{\bar{v}}$-ordinary of weight $\tilde{\lambda}$ if $\pi$ is $\iota V_{\tilde{\lambda}}^\vee$-cohomological and the $\tilde{\lambda}$-rescaled Hecke operators $\{\widetilde{U}_v^k : 1 \le k \le t\}$ have a simultaneous eigenvector with $p$-adic unit eigenvalues in $\iota^{-1}\pi^{\cQ_{\bar{v}}}_{\bar{v}}$.

If $\pi$ is $\iota$-$Q_{\bar{v}}$-ordinary of weight $\tilde{\lambda}$, we define the \emph{$\cQ_{\bar{v}}$-ordinary subspace} of $\iota^{-1}\pi_{\bar{v}}^{\cQ_{\bar{v}}}$ to be the largest $\cH(\tDelta_{\bar{v}}^{\cQ_{\bar{v}}}, \cQ_{\bar{v}})$-submodule on which the rescaled operators $\widetilde{U}_v^k$ have only $p$-adic unit eigenvalues for $1\le k\le t$.
\end{defn}

The goal of this section is to establish the following result, generalizing \cite[Corollary 2.33]{geraghty}, \cite[Theorem 2.4]{jackreducible}. 
\begin{thm}\label{thm:shape general}
	Suppose that $\pi$ is a cuspidal automorphic representation of $\widetilde{G}(\A_{F^+})$, $\iota$ is an isomorphism $\iota: \overline{\Q}_p\toisom \C$ and $\bar{v}$ is a $p$-adic place of $F^+$ 
	such that $\pi$ is $\iota$-$Q_{\bar{v}}$-ordinary of weight $\tilde{\lambda}$. Moreover, suppose that $F$ contains an imaginary quadratic field in which $p$ splits.
	
	Then we have the following conclusions: \begin{enumerate}\item The associated $p$-adic Galois representation 
	$r_{\iota}(\pi): G_F \to \GL_{2n}(\overline{\Q}_p)$ satisfies 
	\begin{equation}\label{eq:shape general}
	r_{\iota}(\pi)|_{G_{F_{v}}} \simeq \begin{pmatrix} r_1(\pi) & * & \cdots & *\\ 0 & r_2(\pi) & \cdots & *\\ 0 & 0 & \ddots & * \\ 0 & 0 & \cdots & r_t(\pi) \end{pmatrix},
	\end{equation}
	where $r_{j}(\pi): G_{F_v}\to \GL_{n_j}(\overline{\Q}_p)$ is a crystalline representation for each $j=1,2,\ldots,t$.

	\item The $\cQ_{\bar{v}}$-ordinary subspace of $\iota^{-1}\pi_{\bar{v}}^{\cQ_{\bar{v}}}$ is one-dimensional.
	
	\item For each embedding $\tau: F_v\hookrightarrow \overline{\Q}_p$, the 
	$\tau$-Hodge--Tate weights of the $r_j(\pi)$ are given by decomposing  $\tilde{\lambda}_{\tau, 2n}<
	\tilde{\lambda}_{\tau, 2n-1}+1 <\dots < \tilde{\lambda}_{\tau, 1} + 2n-1$ according to the partition $(n_1,\ldots,n_t)$. 
	
	\item The determinants $\det r_j(\pi)$ are given by the formulas:
	\begin{itemize}
		\item $\prod_{j=1}^{k}\det r_j(\pi)(\Art_{F_v}(u)) = \prod_{i=1}^{n_1+\cdots + n_k}\prod_{\tau: F_v \hookrightarrow \Qpbar}\tau(u)^{-\tilde{\lambda}_{\tau,2n-i+1}-i+1}$ for $u \in \cO_{F_v}^\times$.
		\item $\prod_{j=1}^{k}\det r_j(\pi)(\Art_{F_v}(\varpi_v))$ is equal to $\epsilon_p^{\sum_{i=1}^{n_1+\cdots + n_k}(1-i)}(\Art_{F_v}(\varpi_v))$ times the eigenvalue of $\widetilde{U}_v^k$ on the  $\cQ_{\bar{v}}$-ordinary subspace.
	\end{itemize}
\end{enumerate}
\end{thm}

\noindent Before giving the proof, we first establish a preliminary result. 

\begin{lemma}\label{lem:reducibility lemma} Assume that 
$r: G_{F_{v}}\to \GL_m(\overline{\Q}_p)$ is a semi-stable Galois representation.  
Let $v_1\leq v_2\leq  \dots \leq v_m$ denote the valuations of the eigenvalues 
of the geometric Frobenius acting on $\mathrm{WD}(r)$, and, 
for each $\tau: F_v\hookrightarrow \overline{\Q}_p$, let $h_{\tau,1}< h_{\tau, 2}< \dots 
< h_{\tau, m}$ denote the $\tau$-Hodge--Tate weights of $r$. Then 
\begin{equation}\label{eq:Newton above Hodge}
\sum_{i=1}^j v_i \geq \frac{1}{e_v} \sum_{i=1}^j \sum_{\tau} h_{\tau, i}
\end{equation}
for any $0\leq j \leq m$, where $e_v$ is the ramification degree of $F_v/\Q_p$. 

Furthermore, if we have an equality 
\[
\sum_{i=1}^j v_{\sigma(i)} = \frac{1}{e_v} \sum_{i=1}^j \sum_{\tau} h_{\tau, i}
\]
for some $1\leq j\leq m-1$ and permutation $\sigma\in S_m$, then 
$r \simeq  \begin{pmatrix} r_1 & * \\ 0 & r_2 \end{pmatrix}$
with $r_1: G_{F_v} \to \GL_j(\overline{\Q}_p)$
and $r_2: G_{F_v}\to \GL_{m-j}(\overline{\Q}_p)$. 
The representation $r_1$ has 
$\tau$-Hodge--Tate weights equal to $h_{\tau,1}<\dots < h_{\tau,j}$ for each $\tau$, 
the eigenvalues of the geometric 
Frobenius acting on $\mathrm{WD}(r_1)$ have valuations $v_1\leq v_2\leq \dots \leq v_j$ 
and we have $v_j<v_{j+1}$.   
\end{lemma}

\begin{proof} The first part follows from~\cite[Lemma 6.4.1]{hkv}. 
For the second part, if equality holds for some permutation $\sigma$, then 
by the first part and by the inequalities on the $v_i$, equality must also 
hold for the identity permutation. Furthermore,
we have $\sum_{i=1}^j v_i = \sum_{i=1}^j v_{\sigma(i)}$. 
If $\{\sigma(1),\dots, \sigma(j)\}\not = \{1, \dots, j\}$, then 
$\sigma(i)\geq j+1$ for some $i\leq j$,
and so we must have $v_j = v_{j+1}$. 
We show that $v_j = v_{j+1}$ is impossible. As in \emph{loc.~cit.}, using 
the inequality in~\eqref{eq:Newton above Hodge} 
for $j+1$ and $j-1$,
we deduce that the inequality for $j-1$ is actually an equality. 
Moreover, this implies that $h_{\tau, j} = h_{\tau, j+1}$ for 
all $\tau$, which is a contradiction. We deduce that 
$\{\sigma(1),\dots, \sigma(j)\}= \{1, \dots, j\}$
and that $v_{j}< v_{j+1}$. The last paragraph in 
the proof of~\cite[Lemma 6.4.1]{hkv} gives a $j$-dimensional
sub-$G_{F_v}$-representation $r_1$ of $r$ with 
$\tau$-Hodge--Tate weights equal to $h_{\tau,1}<\dots < h_{\tau,j}$
for each $\tau$ and such that the eigenvalues of the geometric
Frobenius acting on $\mathrm{WD}(r_1)$ are $v_1\leq v_2\leq \dots \leq v_j$. 
\end{proof}

\begin{proof}[Proof of Theorem \ref{thm:shape general}]
	We use $\iota_{{v}}$ to identify $\widetilde{G}(F_{\bar{v}}^+)$ with $\GL_{2n}(F_v)$. Since $\cQ_{\bar{v}}$ contains the Iwahori subgroup of $\widetilde{G}(F_{\bar{v}}^+)$, $\iota^{-1}\pi_{\bar{v}}$ is a subquotient of a normalized induction $\sigma = \mathrm{n-Ind}_{\Brm_{2n}(F_v)}^{\GL_{2n}(F_v)}(\chi_1\otimes\chi_2\otimes\cdots\otimes\chi_{2n})$ with each $\chi_i: F_v^\times\to \Qpbarx$ an unramified character. 
	
	We claim that, up to reordering the $\chi_i$, for each $1 \le k \le t$ the $p$-adic numbers $\prod_{i=1}^{n_1+\cdots n_k} \chi_i(\varpi_v)$ and $\delta_{\Brm_{2n}}(\nu_k(\varpi_{\bar{v}}))^{1/2}\alpha_{\tilde{\lambda}}^{\cQ_{\bar{v}}}(\nu_k(\varpi_{\bar{v}}))$ differ by a $p$-adic unit. 
	
	By the proof of \cite[Prop. 5.4]{jack}, we have an isomorphism $\sigma^{\cQ_{\bar{v}}} \cong (J_{Q_{\bar{v}}}(\sigma))^{\cQ_{\bar{v}}\cap M_{\bar{v}}(F^+_{\bar{v}})}$ (normalized Jacquet module) with the $\tilde{\lambda}$-rescaled action of the Hecke operator $\widetilde{U}_v^k$ on the left given by the action of $ \alpha_{\tilde{\lambda}}^{\cQ_{\bar{v}}}(\nu_k(\varpi_{\bar{v}}))^{-1}\delta_{Q_{\bar{v}}}^{-1/2}(\nu_k(\varpi_{\bar{v}}))\nu_k(\varpi_{\bar{v}})$ on the right. We can further compose with the injection $(J_{Q_{\bar{v}}}(\sigma))^{\cQ_{\bar{v}}\cap M_{\bar{v}}(F^+_{\bar{v}})} \hookrightarrow (J_{\Brm_{2n}}(\sigma))^{\cQ_{\bar{v}}\cap \Trm_{2n}(F_{{v}})}$, with $\widetilde{U}_v^k$ acting by $ \alpha_{\tilde{\lambda}}^{\cQ_{\bar{v}}}(\nu_k(\varpi_{\bar{v}}))^{-1}\delta_{\Brm_{2n}}^{-1/2}(\nu_k(\varpi_{\bar{v}}))\nu_k(\varpi_{\bar{v}})$ on the target. Note that since $\nu_k(\varpi_{\bar{v}})$ is central in $M_{\bar{v}}$, we have $\delta_{\Brm_{2n}}(\nu_k(\varpi_{\bar{v}})) = \delta_{Q_{\bar{v}}}(\nu_k(\varpi_{\bar{v}}))$. We have $J_{\Brm_{2n}}(\sigma)^{ss} = \bigoplus_{w \in W}\chi_{w(1)}\otimes\cdots\otimes\chi_{w(2n)}$. Choosing $w$ so that the unit-eigenvalue eigenvector for the $\widetilde{U}_v^k$ contributes to the summand indexed by $w$ gives the desired reordering of the $\chi_i$.
	
	By Theorem~\ref{thm:automorphic Galois reps unitary}, since $\tilde{\lambda}$ is dominant, 
	the $\tau$-Hodge--Tate weights of $r_{\iota}(\pi)|_{G_{F_{v}}}$ 
	are equal to 
	$\tilde{\lambda}_{\tau, 2n}< \tilde{\lambda}_{\tau, 2n-1} + 1<\dots 
	<\tilde{\lambda}_{\tau,1} + 2n-1$. In addition, by making part (3) 
	of Theorem~\ref{thm:automorphic Galois reps unitary} explicit at $v$, 
	the eigenvalues of the geometric Frobenius acting on $\mathrm{WD}(r_{\iota}(\pi)|_{G_{F_v}})$ 
	are $q_v^{\frac{2n-1}{2}} \chi_i(\varpi_v)$. Our claim about the $p$-adic valuations of the $\chi_i(\varpi_v)$ implies that the $p$-adic numbers 
	\[
	\prod_{i=1}^{n_1+\cdots+n_k} \chi_i(\varpi_v)
	\ \mathrm{and}
	\]
	\[
	q_v^{-\langle\nu_k,\rho_{\Brm_{2n}}\rangle}\prod_{\tau: F_v\hookrightarrow \overline{\Q}_p}
	\tau(\varpi_v)^{\langle\nu_k,w_0^{\widetilde{G}}\lambda_{\tau}\rangle}
	\]
	differ by a $p$-adic unit. The latter term has the same $p$-adic valuation as \[\prod_{i=1}^{n_1+\cdots+n_k}q_v^{\frac{1-2n}{2}}\prod_{\tau: F_v\hookrightarrow \overline{\Q}_p}
		\tau(\varpi_v)^{\lambda_{\tau, 2n-i+1}+ i-1}.\] In turn, this implies that the inequality in~\eqref{eq:Newton above Hodge}
	is an equality for $j=n_1+\cdots+n_k$ for each $k$. We conclude that $r_{\iota}(\pi)|_{G_{F_v}}$
	has the desired shape by Lemma~\ref{lem:reducibility lemma}. 

To show that the $\cQ_{\bar{v}}$-ordinary subspace of $\iota^{-1}\pi_{\bar{v}}^{\cQ_{\bar{v}}}$ is one-dimensional, it suffices to prove the analogous statement for $\sigma^{\cQ_{\bar{v}}} \cong (J_{Q_{\bar{v}}}(\sigma))^{\cQ_{\bar{v}}\cap M_{\bar{v}}(F^+_{\bar{v}})}$. By the geometric lemma \cite[Lemma 2.12]{BZ77ENS}, we have \begin{equation}\label{eq:jacquetdecomp}J_{Q_{\bar{v}}}(\sigma)^{ss} = \bigoplus_{w \in [W/W_{Q_{\bar{v}}}]} \left(\mathrm{n-Ind}^{M_{\bar{v}}}_{\Brm_{2n}\cap M_{\bar{v}}} \chi_{w(1)}\otimes\cdots\otimes\chi_{w(2n)} \right).\end{equation} The notation $w \in [W/W_{Q_{\bar{v}}}]$ means that we take the minimal length representative $w \in W$ for each coset. Explicitly, this is given by permutations $w \in S_{2n}$ which are order-preserving on each of the subsets $\{1,\ldots,n_1\}, \{n_1+1,\ldots,n_1+n_2\}, \ldots, \{n_1+\cdots+n_{t-1}+1,\ldots,2n\}$. Each summand in \eqref{eq:jacquetdecomp} has a one-dimensional space of invariants under the maximal compact subgroup $\cQ_{\bar{v}}\cap M_{\bar{v}}(F^+_{\bar{v}})$. On the invariants of the summand indexed by $w$, for each $1 \le k \le t$ the Hecke operator $\widetilde{U}_v^k$ acts by the scalar \[\alpha_{k,w} := \alpha_{\tilde{\lambda}}^{\cQ_{\bar{v}}}(\nu_k(\varpi_{\bar{v}}))^{-1}\delta_{\Brm_{2n}}^{-1/2}(\nu_k(\varpi_{\bar{v}}))\prod_{i=1}^{n_1+\cdots+n_k} \chi_{w(i)}(\varpi_v).\] We know that $\alpha_{k,1}$ is a $p$-adic unit for all $k$. The inequality in the conclusion of Lemma \ref{lem:reducibility lemma} shows that $v_p(\alpha_{1,w})> 0$ if $w$ does not preserve $\{1,\ldots,n_1\}$. For a fixed $w$, inducting on $k$ and repeating this argument shows that $v_p(\alpha_{k,w}) = 0$ for all $k$ if and only if $w = 1$.  This shows that the $\cQ_{\bar{v}}$-ordinary subspace of $\iota^{-1}\pi_{\bar{v}}^{\cQ_{\bar{v}}}$ is one-dimensional.

The statements about the Hodge--Tate weights and determinant of the $r_i(\pi)$ follow from Lemma \ref{lem:reducibility lemma} and the identifications above.  
	
	It remains to see that each $r_i(\pi)$ is crystalline. We write the argument for $r_1(\pi)$, it is similar for the other factors.  Let $w_i = v_p(\chi_i(\varpi_v))$ for $i=1,\dots, 2n$. 
	The eigenvalues of a geometric Frobenius on 
	$\mathrm{WD}(r_1(\pi))$ are $q_v^{\frac{2n-1}{2}}\chi_i(\varpi_v)$ for $i=1,\dots,n_1$. Lemma~\ref{lem:reducibility lemma} implies that, up to reordering the $\chi_i$ 
	for $i=1,\dots, n_1$, we may assume that $w_1\leq w_2\dots \leq w_{n_1} < w_{n_1+1}$.  
	If $r_1(\pi)$ is not crystalline, then by the Bernstein--Zelevinsky classification there exists some $i\in 2,\dots n_1$ 
	such that $\chi_{i-1} = \chi_{i}\cdot |\ |$ and
	$\pi_{\bar{v}}$ is a subquotient of the normalized parabolic induction
	\[
	\sigma' = \mathrm{n-Ind}_{Q'(F_v)}^{\GL_{2n}(F_v)}(
	\chi_1\otimes \dots \otimes \chi_{i-2}\otimes \mathrm{Sp}_2(\chi_i) \otimes \chi_{i+1}\otimes \dots \otimes \chi_{2n}),
	\]
	for an appropriate standard parabolic subgroup $Q' = M'N' \subset Q_{\bar{v}}$, where $\mathrm{Sp}_2(\chi_i)$ denotes a twist of the Steinberg representation. We let $\sigma'_0 = \chi_1\otimes \dots \otimes \chi_{i-2}\otimes \mathrm{Sp}_2(\chi_i) \otimes \chi_{i+1}\otimes \dots \otimes \chi_{2n}$.
	
	Applying the geometric lemma again, we have \[J_{Q_{\bar{v}}}(\sigma')^{ss} = \bigoplus_{w \in [W_{Q'}\backslash W/W_{Q_{\bar{v}}}]}  \left(\mathrm{n-Ind}^{M_{\bar{v}}}_{w^{-1}Q'w\cap M_{\bar{v}}}\left( J_{M'\cap wQ_{\bar{v}}w^{-1}}\sigma'_0 \right)^w\right)^{ss}\]
	
	As in the proof that the $\cQ_{\bar{v}}$-ordinary subspace is one-dimensional, it follows from considering valuations that the $\cQ_{\bar{v}}$-ordinary subspace of $\iota^{-1}\pi_{\bar{v}}^{\cQ_{\bar{v}}} \cong (J_{Q_{\bar{v}}}(\sigma'))^{\cQ_{\bar{v}}\cap M_{\bar{v}}(F^+_{\bar{v}})}$ can only contribute to the $w=1$ term in this decomposition.
We deduce that  \[\left(\mathrm{n-Ind}^{M_{\bar{v}}}_{Q'\cap M_{\bar{v}}}\left( \sigma'_0 \right)\right)^{\cQ_{\bar{v}}\cap M_{\bar{v}}(F^+_{\bar{v}})} \ne 0,\] but this is impossible because $\cQ_{\bar{v}}\cap M_{\bar{v}}(F^+_{\bar{v}})$ is maximal compact and $\sigma_0'$ has a Steinberg factor.
\end{proof}

\subsection{Determinants}\label{subsec:determinants}

In this section we prove a key proposition which will be combined with Theorem \ref{thm:shape general} to pass information about local--global compatibility for Galois representations with coefficients in $p$-torsion free Hecke algebras for $\tG$ to Galois representations with coefficients in torsion Hecke algebras for $G$. The initial set-up is as follows: recall our coefficient field $E \supset \cO 
\twoheadrightarrow k$ and assume we have an absolutely irreducible continuous
representation 
\[\rhobar_{\m}: G_F \to \GL_n(k)\] together with a continuous lift \[\rho_{\m}: G_F \to 
\GL_n(A)\] with 
coefficients in  $A \in \CNL_{\cO}$. 

We moreover have a finite flat 
$\cO$-algebra $\widetilde{A}\in \CNL_{\cO}$ with $\widetilde{A}[1/p] = \prod_{i = 1}^r K_i$ a 
product of fields, equipped with an $\cO$-algebra map 
$\widetilde{A}\rightarrow A$. Extending $E$ if necessary, we can assume 
that every $K_i$ has residue field $k$ (indeed, we can also assume that each $K_i$ is equal to $E$, but we won't need to do this). 
We suppose we have a continuous representation \[\widetilde{\rho}_{\m} = \prod 
\widetilde{\rho}_{i,\m} : G_F 
\to 
\prod \GL_{2n}(K_i) = \GL_{2n}(\widetilde{A}[1/p]) \]
such that the associated determinant of $G_F$ (in the sense of 
\cite{chenevier_det}) arises from a continuous $\widetilde{A}$-valued determinant 
$D_{\widetilde{\rho}_{\m}}: \widetilde{A}[G_F] \to \widetilde{A}$. This is equivalent to the characteristic polynomials of $\widetilde{\rho}_\m(g)$ having coefficients in $\widetilde{A}$ for all $g \in G_F$.

We recall from \cite{chenevier_det} the important notion of a \emph{Cayley--Hamilton} determinant: 
a determinant $D: R \to A$ is Cayley--Hamilton if the characteristic polynomial $\chi(r,t) = D(t-r) \in A[t]$ vanishes when evaluated at $t = r$ for all $r \in R$. 

The determinant $D_{\widetilde{\rho}_\m}$ factors through a 
Cayley--Hamilton determinant (which we also denote by 
$D_{\widetilde{\rho}_{\m}}$) of 
$\widetilde{B}:=\widetilde{\rho}_{\m}(\widetilde{A}[G_F]) \subset 
M_{2n}(\widetilde{A}[1/p])$. Since each element of $\widetilde{B}$ has  characteristic polynomial with coefficients in $\widetilde{A}$, $\widetilde{B}$ is 
integral over $\widetilde{A}$.
We also assume that we have a factorization of $A$-valued determinants of $A[G_{F}]$
\[D_{\widetilde{\rho}_{\m}}\otimes_{\widetilde{A}}A = 
D_{\rho_{\m}}D_{\rho_{\m}^{\vee,c}(1-2n)}.\]  

Fix a place $v|p$ of $F$. Extending $E$ if necessary, we may assume that
every irreducible constituent of the local representation 
$\rhobar_\m|_{G_{F_v}}$ is absolutely irreducible. 
We assume that, for each $i$, we have an $n$-dimensional 
$G_{F_v}$-sub-representation $(\widetilde{\rho}_{i,\m}^0 , V_i^0) \subset 
(\widetilde{\rho}_{i,\m}, K_i^{2n})$ with quotient $(\widetilde{\rho}_{i,\m}^1, 
V_i^1)$.
Finally, we assume that:

\begin{enumerate}
	\item\label{ass:det1} the isomorphism classes of the irreducible constituents of 
	$\rhobar_{\m}|_{G_{F_v}}$ are disjoint from those of 
	$\rhobar_{\m}^{\vee,c}(1-2n)|_{G_{F_v}}$;
	
	\item\label{ass:det2} for all $i$ the isomorphism classes of the irreducible 
	constituents of the residual representation 
	$\overline{(\widetilde{\rho}_{i,\m}^0)}$ coincide with those of 
	$\rhobar_{\m}|_{G_{F_v}}$.
\end{enumerate}

As in the proof of \cite[Thm.~2.22]{chenevier_det}, we are going to use a basic 
fact about idempotents, cf.~\cite[Chapter III, \S 4, Exercise 5 (a)]{Bourbaki-CA-1-7}. 

\begin{lemma}\label{lemma:basic idempotent}
	Let $S$ be a Henselian local ring, $R$ a (not necessarily commutative) 
	$S$-algebra which is integral over $S$. Let $I$ be a two-sided ideal of 
	$R$. Every idempotent in $R/I$ lifts to an idempotent in 
	$R$. 
\end{lemma}

\begin{proof} The statement can be reduced to the case
of a finite $S$-algebra $R$ generated by 
one element, in particular with $R$ commutative (take the subalgebra of $R$ generated by an arbitrary lift of the idempotent in $R/I$). 
In this case the statement is clear, since 
$R$ is then a product of finitely many Henselian local rings. 
\end{proof}

The following lemma constructs an element $\widetilde{e} \in M_{2n}(\widetilde{A}[1/p])$ which is an idempotent projection onto $\prod_i V_i^0$, lies in $\widetilde{\rho}_{\m}(\widetilde{A}[G_{F_v}])$ and moreover cuts out the first factor in the product decomposition $D_{\widetilde{\rho}_{\m}}\otimes_{\widetilde{A}}k = 
D_{\rhobar_{\m}}D_{\rhobar_{\m}^{\vee,c}(1-2n)}$ of residual determinants. Roughly speaking, this is possible because our assumptions on the irreducible constituents of $\rhobar_\m|_{G_{F_v}}$ mean the global factors $D_{\rhobar_{\m}}$, $D_{\rhobar_{\m}^{\vee,c}(1-2n)}$ can be distinguished locally at $v$. 

\begin{lem}\label{lem:idems}
	\begin{enumerate}
		\item The natural map $k[G_F] \xrightarrow{\rhobar_\m \times 
			\rhobar_{\m}^{\vee,c}(1-2n)} M_n(k) \times 
		M_n(k)$ induces an isomorphism $\iota: k[G_F]/\ker 
		(D_{\widetilde{\rho}_\m}\otimes_{\widetilde{A}}k) \cong M_n(k) \times 
		M_n(k)$. 
		
		\item The idempotent $e = \iota^{-1}(1_{M_n(k)}, 0)$ is contained 
		in the image of $k[G_{F_v}]$. 
		
		\item There is an idempotent lift $\widetilde{e} \in 
		\widetilde{B}$ of $e$ which is in the image $\widetilde{B}_v \subset 
		\widetilde{B}$ of 
		$\widetilde{A}[G_{F_v}]$. (Note that $D_{\widetilde{\rho}_\m}\otimes_{\widetilde{A}}k$ factors through $\widetilde{B}\otimes_{\widetilde{A}}k$, so $k[G_F]/\ker 
		(D_{\widetilde{\rho}_\m}\otimes_{\widetilde{A}}k)$ is a quotient of $\widetilde{B}$.)
		
		\item\label{lem:idems last part} Choose $\widetilde{e}$ as in the previous part. For each $i \in 
		\{1,\ldots,r\}$ let $\pi_i$ denote the projection 
		$\widetilde{A}[1/p]^{2n} 
		\to V_i^1$. Thinking of $\widetilde{e}$ as an element of 
		$\End(\widetilde{A}[1/p]^{2n})$, we have $\pi_i\circ  \widetilde{e} = 
		0$ and the image of $\widetilde{e}$ in $\End(V_i^0)$ is the identity. 
		In other words, $\widetilde{e}$ is an idempotent projection onto 
		$\prod_i V_i^0$.
	\end{enumerate}
\end{lem}
\begin{proof}
	The first claim follows from our assumption that $\rhobar_\m$ and 
	$\rhobar_{\m}^{\vee,c}(1-2n)$ are distinct and absolutely irreducible, by 
	\cite[Thm.~2.16]{chenevier_det}. 
	
	For the second claim we consider the subalgebra $B = \iota(k[G_{F_v}])$ of 
	$M_n(k) \times M_n(k)$. The semisimple quotient $B/\mathrm{Rad}(B)$ 
	contains the idempotent $\overline{e}_v$ which acts as the identity on each 
	irreducible 
	constituent of $\rhobar_{\m}|_{G_{F_v}}$ and as zero on each irreducible 
	constituent of $\rhobar_{\m}^{\vee,c}(1-2n)|_{G_{F_v}}$; here we are using 
	our 
	assumption that these two collections of irreducible constituents are 
	disjoint. Since $\mathrm{Rad}(B)$ is nilpotent, 
	we can lift $\overline{e}_v$ to an idempotent $e_v \in B$. By 
	considering 
	the composition series of $\rhobar_{\m}|_{G_{F_v}}$, we see that $e_v$ maps 
	to an idempotent unit (i.e.~the identity) under the first projection 
	to $M_n(k)$. Considering the composition series of 
	$\rhobar_{\m}^{\vee,c}(1-2n)$, $e_v$ maps to a nilpotent idempotent 
	(i.e.~zero) under the 
	second projection, so we have $e_v = e$.
	
	For the third claim, we know that $\widetilde{B}_v$ is 
	integral over $\widetilde{A}$ (since it is a subalgebra of $\widetilde{B}$). We apply Lemma \ref{lemma:basic idempotent}
	with $S=\widetilde{A}$, which is Henselian as it is a local ring and 
	finite over $\cO$, and $R = \widetilde{B}_v$. We deduce that the idempotent $e$ lifts to an idempotent $\widetilde{e}$ in $\widetilde{B}_v$.
	
	Now we come to the fourth claim. Fix an index $i$. We choose a $G_F$-stable 
	$\cO_{K_i}$-lattice $T \subset K_i^{2n}$. We have a short exact sequence \[ 
	0 \to V_i^0\cap T \to T \to \pi_i(T) \to 0\] of 
	$\cO_{K_i}[G_{F_v}]$-modules (in particular the submodule $V_i^0\cap T$ is 
	stable under $\widetilde{e}$). The image of $\widetilde{e}$ in 
	$\End_{k}(\pi_i(T)\otimes_{\cO_{K_i}}k)$ is equal to zero, since 
	$\widetilde{e}$ lifts the idempotent which acts as zero on each irreducible 
	constituent of the $G_{F_v}$ representation on 
	$\pi_i(T)\otimes_{\cO_{K_i}}k$ (by assumption 
	these coincide with the irreducible constituents of 
	$\rhobar_{\m}^{\vee,c}(1-2n)|_{G_{F_v}}$). We deduce that the image of 
	$\widetilde{e}$ in $\End_{\cO_{K_i}}(\pi_i(T))$ is equal to zero, since 
	it is an idempotent with image in $\varpi_{K_i} \pi_i(T)$. This shows that 
	$\pi_i \circ \widetilde{e} = 0$, as claimed. Similarly, the image of 
	$\widetilde{e}$ in $\End_k((V_i^0\cap T) \otimes_{\cO_{K_i}}k)$ is an 
	idempotent isomorphism, hence the identity, and so the image of  
	$\widetilde{e}$ in $\End_{\cO_{K_i}}(V_i^0\cap T)$ itself is an 
	idempotent isomorphism, hence equal to the identity.
\end{proof}

In the preceding lemma, we constructed an idempotent $\widetilde{e}$, compatible with both the local `$P$-ordinary' decomposition and the residual global decomposition $D_{\widetilde{\rho}_{\m}}\otimes_{\widetilde{A}}k = 
D_{\rhobar_{\m}}D_{\rhobar_{\m}^{\vee,c}(1-2n)}$. Following \cite[Theorem 1.4.4]{bellaiche_chenevier_pseudobook} (and its generalization \cite[Theorem 2.22]{chenevier_det}), we can now use the idempotents $e_1:=\widetilde{e}, e_2:=1-\widetilde{e}$ to equip $\widetilde{B}$ with a \emph{generalized matrix algebra} structure compatible with the determinant $D_{\widetilde{\rho}_\m}$. We now explain in detail what this means (see also \cite[\S2]{ANT}). By \cite[Lem.~2.4]{chenevier_det}, we have determinants of dimension $n$ \begin{align*}
D_{\widetilde{\rho}_\m,i}: 
e_i\widetilde{B}e_i &\to \widetilde{A} \\ x &\mapsto 
D_{\widetilde{\rho}_\m}(x+1-e_i)
\end{align*} for $i=1,2$, with $(D_{\widetilde{\rho}_\m,1}\otimes_{\widetilde{A}}k)(e_1 x e_1) =  D_{\rhobar_{\m}}(x)$ and $(D_{\widetilde{\rho}_\m,2}\otimes_{\widetilde{A}}k)(e_2 x e_2) =  D_{\rhobar_{\m}^{\vee,c}(1-2n)}(x)$. The $e_i\widetilde{B}e_i$ are equipped with $\widetilde{A}$-algebra isomorphisms \begin{equation}\label{eqn:displaypsi}\psi_i: e_i\widetilde{B}e_i \cong M_n(\widetilde{A})\end{equation} with $\det\circ\psi_i = D_{\widetilde{\rho}_\m,i}$. Moreover, the map $(\psi_1,\psi_2): e_1\widetilde{B}e_1 \oplus e_2\widetilde{B}e_2 \to M_n(\widetilde{A})\times M_n(\widetilde{A})$ lifts $\iota: k[G_F]/\ker 
(D_{\widetilde{\rho}_\m}\otimes_{\widetilde{A}}k) \cong M_n(k) \times 
M_n(k)$.

Define idempotents $E_i \in e_i\widetilde{B}e_i$ by asking for $\psi_i(E_i)$ to be the matrix with a one in the top left entry and zeroes elsewhere. We can now define $\widetilde{A}$-submodules $\cA_{i,j} := E_i \widetilde{B} E_j \subset \widetilde{B}$. For each $1 \le i,j,k \le 2$ we write $\cA_{i,j}\cA_{j,k}$ for the $\widetilde{A}$-module generated by products $xy$ for $x\in \cA_{i,j}, y \in \cA_{j,k}$. It is a submodule of $\cA_{i,k}$. The $\psi_i$ induce isomorphisms $\cA_{i,i} \cong \widetilde{A}$, so we have a multiplication $\cA_{i,j}\otimes_{\widetilde{A}}\cA_{j,i} \to \widetilde{A}$ and we identify $\cA_{i,j}\cA_{j,i}$ with an ideal of $\widetilde{A}$. If $i \ne j$, then $e_i\tB e_j\tB e_i$ maps to $\ker(D_{\widetilde{\rho}_{\m}}\otimes_{\tA}k)$, and therefore $\cA_{i,j}\cA_{j,i} \subset \m_{\tA}$. 

Putting everything together gives us isomorphisms \begin{align*}\widetilde{B}&\toisom \begin{pmatrix}
e_1\widetilde{B}e_1 & 
e_1\widetilde{B}e_2 \\ 
e_2\widetilde{B}e_1 & 
e_2\widetilde{B}e_2 
\end{pmatrix} \toisom \begin{pmatrix}
M_n(\widetilde{A}) & 
M_n(\cA_{1,2}) \\ 
M_n(\cA_{2,1}) & 
M_n(\widetilde{A})
\end{pmatrix}
\end{align*} with multiplication of the last matrix algebra given by using matrix multiplication and the maps $\cA_{i,j}\otimes_{\widetilde{A}}\cA_{j,i} \to \widetilde{A}$. The reducibility of the determinant $D_{\widetilde{\rho}_\m}$ is reflected in the GMA structure. More precisely, we apply  \cite[Proposition 2.5]{ANT} in our setting (a mild generalization of \cite[Proposition 1.5.1]{bellaiche_chenevier_pseudobook}) to deduce:

\begin{prop}\label{prop:redideal}
An ideal $J \subset \widetilde{A}$ contains $\cA_{1,2}\cA_{2,1}$ if and only if there are determinants $D_1, D_2 : \widetilde{B}\otimes_{\widetilde{A}}\tA/J \to \tA/J$ such that \begin{align*}D_{\widetilde{\rho}_\m}\otimes_{\widetilde{A}}\tA/J  &= D_1D_2, \\ D_1 \otimes_{\tA/J}k & = D_{\rhobar_{\m}},\\
\text{and }D_2\otimes_{\tA/J}k & = D_{\rhobar_{\m}^{\vee,c}(1-2n)}.\end{align*} If this property holds, $D_1$ and $D_2$ are uniquely determined and satisfy $\ker(D_{\widetilde{\rho}_\m}\otimes_{\widetilde{A}}\tA/J) \subset \ker(D_i)$ for $i=1,2$.
\end{prop}
\begin{proof}
	This follows immediately from \cite[Proposition 2.5]{ANT}. Note that the property `(COM)' in \cite[Lemma 1.3.5]{bellaiche_chenevier_pseudobook} (which follows from the fact that the trace of a determinant satisfies the identity $\mathrm{Tr}(xy) = \mathrm{Tr}(yx)$) means that the ideals $\cA_{1,2}\cA_{2,1}$ and $\cA_{2,1}\cA_{1,2}$ are equal.
\end{proof}

We can now state and prove the key proposition of this section, where we combine the GMA structure defined above with the fact that our idempotent $\widetilde{e}$ was chosen to have good local properties at $v$.

\begin{prop}\label{prop:char0lifts}
	Choose an idempotent $\widetilde{e} \in \widetilde{B_v}$ as in the third 
	part of Lemma \ref{lem:idems}. 
	\begin{enumerate}
		\item The map \begin{align*}A[G_{F}] &\to 
		\widetilde{e}\widetilde{B}\widetilde{e}\otimes_{\widetilde{A}}A\\ x 
		&\mapsto 
		\widetilde{e}x\widetilde{e}\otimes 1\end{align*} is a homomorphism and 
		it 
		induces the determinant $D_{\rho_\m}$ when we compose with 
		$\widetilde{e}\widetilde{B}\widetilde{e}\otimes_{\widetilde{A}}A \xrightarrow{\psi_1\otimes\mathrm{id}} 
		M_n(A)$ and the usual determinant (see \eqref{eqn:displaypsi} for $\psi_1$). 
		
		\item The map \begin{align*}\widetilde{A}[G_{F_v}] &\to 
		\widetilde{e}\widetilde{B}\widetilde{e}\\ x 
		&\mapsto 
		\widetilde{e}x\widetilde{e}\end{align*} is also a homomorphism and 
		it 
		induces the representation $\prod_{i=1}^r\tilde{\rho}_{i,\m}^0$ when we compose with the natural inclusion $\widetilde{e}\widetilde{B}\widetilde{e} \subset 
		\End_{\widetilde{A}[1/p]}(\prod_i V_i^0)$ (see the final part of Lemma \ref{lem:idems} for why we have this inclusion).
		
		\item There is an $\widetilde{A}$-valued lift of the 
		representation 
		$\rho_{\m}|_{G_{F_v}}$ which becomes isomorphic to $\prod_{i=1}^r 
		\widetilde{\rho}^0_{i, \m}$ when we invert $p$.
	\end{enumerate}
\end{prop}
\begin{proof}
	For the first part, we use Proposition \ref{prop:redideal}. Since the determinant 
	$D_{\widetilde{\rho}_\m}\otimes_{\widetilde{A}}A = 
	D_{\rho_{\m}}D_{\rho_{\m}^{\vee,c}(1-2n)}$ is reducible, this tells 
	us that the kernel 
	$J$ of $\widetilde{A} \to A$ contains the reducibility ideal 
	$\cA_{1,2}\cA_{2,1}$ (note that the reducibility ideal of the base-changed determinant $D_{\widetilde{\rho}_\m}\otimes_{\widetilde{A}}A$ is the $A$-ideal generated by the image of $\cA_{1,2}\cA_{2,1}$). It follows that the map 
	\begin{align*}\widetilde{A}/J[G_{F}] &\to 
	\widetilde{e}\widetilde{B}\widetilde{e}\otimes_{\widetilde{A}}\widetilde{A}/J\\ x 
	&\mapsto 
	\widetilde{e}x\widetilde{e}\otimes 1\end{align*} is a homomorphism, and so the same is true for its scalar extension from the subring $\widetilde{A}/J$ to $A$. The  determinant induced by $D_{\widetilde{\rho}_\m,1}\otimes_{\widetilde{A}}A$  is equal to $D_{\rho_{\m}}$ by the uniqueness part of Proposition \ref{prop:redideal} (applied over $A$). Note moreover that this determinant takes values in the subring $\widetilde{A}/J$ of $A$.
	
	For the second part, we can check that we have a homomorphism in $\End_{\widetilde{A}[1/p]}(\prod_i V_i^0)$, where it follows from the fact that $\prod_i V_i^0$ is $G_{F_{v}}$-stable. The identification of the representation with $\prod_{i=1}^r \tilde{\rho}_{i,\m}^0$ is now clear.
	
	For the third part, since $\rhobar_{\m}$ is absolutely irreducible as a $G_F$-representation, it follows from Skolem--Noether (see e.g.~\cite[Proposition IV.1.4]{milne-etcoh}) that, after conjugating $\rho_{\m}$ by an element of $\GL_n(A)$, we can assume that the representation $A[G_F] \to M_n(A)$ given by the first part is equal to $\rho_{\m}$. By the second part, $\psi_1$ also induces a \emph{local} representation $\widetilde{A}[G_{F_v}] \to M_n(\widetilde{A})$ which clearly lifts $\rho_{\m}|_{G_{F_v}}$. It also follows from the second part that after inverting $p$ this representation becomes isomorphic to $\prod_{i=1}^r \tilde{\rho}_{i,\m}^0$.
\end{proof}

\subsection{Local deformation rings}\label{sec:localdefrings}
In this section we fix a place $v \in S_p$ in $F$ and a residual local Galois representation
\[\rhobar_v: G_{F_v} \to \GL_n(k).\]

\begin{defn}
	Let $B$ be a finite $E$-algebra, and let $\lambda_v = (\lambda_{\tau,1}\ge \cdots \ge \lambda_{\tau,n})_{\tau \in \Hom(F_v, E)}$ be a dominant weight for $(\mathrm{Res}_{F_v/\Qp}\GL_n)_E$.
	\begin{enumerate}
		\item  A continuous representation $\rho: G_{F_v} \to \GL_n(B)$ is semistable-ordinary of weight $\lambda_v$ if it is conjugate to an upper triangular representation \[\begin{pmatrix} \chi_1 & * & \cdots & *\\ 0 & \chi_2 & \cdots & *\\ 0 & 0 & \ddots & * \\ 0 & 0 & \cdots & \chi_n \end{pmatrix}\] where for each $1\le j \le n$ and $\sigma \in I_{F_v}$ we have\[\chi_j(\sigma) = \prod_{\tau \in\Hom(F_v,E)}\tau(\Art_{F_v}^{-1}(\sigma))^{-\lambda_{\tau,n+1-j}-(j-1)} \]
		(cf.~\cite[Definition 3.8]{geraghty}.)
		\item We define a $p$-adic Hodge type (in the sense of \cite[\S 2.6]{kisindefrings}) $\bv_{\lambda_v}$ associated to $\lambda_v$ as in \cite[\S 3.3]{geraghty}. This is an $n$-dimensional $E$-vector space $D_E$ with a decreasing filtration on $D_E \otimes_{\Qp}F_v$ by $E\otimes_{\Qp}F_v$-submodules. More precisely $D_E\otimes_{\Qp}F_v$ is isomorphic as a filtered $E\otimes_{\Qp}F_v$-module to $D_{\dR}(\rho)$ where $\rho:G_{F_v} \to \GL_n(E)$ is a de Rham representation with labelled Hodge--Tate weights $(\lambda_{\tau,1}+n-1 > \lambda_{\tau,2}+n-2 > \cdots > \lambda_{\tau,n})_{\tau \in \Hom(F_v, E)}$.
		\item A semistable continuous representation $\rho: G_{F_v} \to \GL_n(B)$ has $p$-adic Hodge type $\bv_{\lambda_v}$ if for each $i$ there is an isomorphism of $B\otimes_{\Qp}F_v$-modules \[\mathrm{gr}^i(\rho\otimes_{\Qp}B_{\dR})^{G_{F_v}}\cong B\otimes_E (\mathrm{gr}^i D_{E}\otimes_{\Qp} F_v).\]
	\end{enumerate}
\end{defn} 
We recall (\cite[Lemma 3.9]{geraghty}) that a semistable-ordinary representation of weight $\lambda_v$ is semistable with $p$-adic Hodge type $\bv_{\lambda_v}$. 

\begin{lem}\label{lem:sst-ordinary deformations}
	Let $B$ be a finite local $E$-algebra and suppose $\rho_B: G_{F_v} \to \GL_n(B)$ is semistable of Hodge type $\bv_{\lambda_v}$ with $\rho:=\rho_B \otimes_B (B/\m_B)$ semistable-ordinary of weight $\lambda_v$. Then $\rho_B$ is semistable-ordinary of weight $\lambda_v$. 
\end{lem}
\begin{proof}
Extending $E$ if necessary, we assume $B$ has residue field $E$. Let $v_1 \le v_2 \le \cdots \le v_n$ denote the valuations of the eigenvalues of geometric Frobenius acting on $\WD(\rho)$ (with algebraic multiplicities). It follows from the proof of \cite[Lemma 6.4.1]{hkv} that in fact this sequence of slopes is strictly increasing and we have an equality \begin{equation}\label{eq:ordNPHP}v_i = \frac{1}{e_v}\sum_{\tau}(\lambda_{\tau,n+1-i}+i-1)\end{equation} for each $1 \le i \le n$. The ordinary filtration on $\rho$ corresponds to a filtration $\cF^\bullet$ of $D_{\st}(\rho)$ by admissible filtered $(\phi,N)$-modules which are free over $F_{v,0}\otimes_{\Qp} E$. For $i = 1,\ldots,n$, $\cF^i \subset D_{\st}(\rho)$ is the $E$-vector subspace spanned by the generalized eigenspaces of $\phi^{f_v}$ with eigenvalues of valuation $\le v_i$. 
The (admissible) filtered $(\phi,N)$-module $D_{\st}(\rho_B)$ is a successive extension of copies of $D_{\st}(\rho)$. We define $\cF^i_B$ to be the $B$-submodule spanned by the generalized eigenspaces of $\phi^{f_v}$ with eigenvalues of valuation $\le v_i$. This defines a filtration of $D_{\st}(\rho_B)$ by free $F_{v,0}\otimes_{\Qp}B$-submodules, stable under the actions of $\phi$ and $N$. Weak admissibility of $D_{\st}(\rho_B)$ and the equalities (\ref{eq:ordNPHP}) imply that each $\cF_B^i$ is weakly admissible and $\cF_B^i/\cF_B^{i-1} = D_{\st}(\chi_{i,B})$ for a crystalline character $\chi_{i,B}:G_{F_v}\to B^\times$ with labelled Hodge--Tate weights $(\lambda_{\tau,n+1-i}+i-1)_{\tau\in\Hom(F_v,E)}$.  
\end{proof}

We consider the lifting ring $R_{\rhobar_v}^{\square}$ with the universal lift $\rho_v^{\mathrm{univ}}: G_{F_v} \to \GL_n(R_{\rhobar_v}^{\square})$, and recall some results of Bellovin, Geraghty, Hartl--Hellmann and Kisin:

\begin{thm}\label{thm:localdefrings}
Let $\lambda_v = (\lambda_{\tau,1}\ge \cdots \ge \lambda_{\tau,n})_{\tau \in \Hom(F_v, E)}$ be a dominant weight for $(\mathrm{Res}_{F_v/\Qp}\GL_n)_E$.

\begin{enumerate}\item There is a unique $\cO$-flat quotient $R_{\rhobar_v}^{\st,\lambda_v}$ of $R_{\rhobar_v}^{\square}$ with the following property:
\begin{itemize}
	\item If $B$ is a finite $E$-algebra, an $\cO$-algebra map $\zeta: R_{\rhobar_v}^{\square} \to B$ factors through $R_{\rhobar_v}^{\st,\lambda_v}$ if and only if $\rho_v^{\mathrm{univ}}\otimes_{R_{\rhobar_v}^{\square},\zeta}B$ is semistable with $p$-adic Hodge type $\bv_{\lambda_v}$.
\end{itemize}
\item $R_{\rhobar_v}^{\st,\lambda_v}$ is reduced.
\item\label{part:crisdefring} There is a unique $\cO$-flat quotient $R_{\rhobar_v}^{\cris,\lambda_v}$ of $R_{\rhobar_v}^{\square}$ with the following property:
\begin{itemize}
	\item If $B$ is a finite $E$-algebra, an $\cO$-algebra map $\zeta: R_{\rhobar_v}^{\square} \to B$ factors through $R_{\rhobar_v}^{\cris,\lambda_v}$ if and only if $\rho_v^{\mathrm{univ}}\otimes_{R_{\rhobar_v}^{\square},\zeta}B$ is crystalline with $p$-adic Hodge type $\bv_{\lambda_v}$.
\end{itemize}
\item $R_{\rhobar_v}^{\cris,\lambda_v}[\frac{1}{p}]$ is regular (in particular, $R_{\rhobar_v}^{\cris,\lambda_v}$ is reduced).
\item  There is a unique $\cO$-flat quotient $R_{\rhobar_v}^{\triangle,\lambda_v}$ of $R_{\rhobar_v}^{\square}$ with the following property:
\begin{itemize}
	\item\label{part:orddefring} If $B$ is a finite $E$-algebra, an $\cO$-algebra map $\zeta: R_{\rhobar_v}^{\square} \to B$ factors through $R_{\rhobar_v}^{\triangle,\lambda_v}$ if and only if $\rho_v^{\mathrm{univ}}\otimes_{R_{\rhobar_v}^{\square},\zeta}B$ is semistable-ordinary of weight $\lambda_v$.
\end{itemize}
\item $\Spec(R_{\rhobar_v}^{\triangle,\lambda_v}[\frac{1}{p}])$ is an open and closed subspace of $\Spec(R_{\rhobar_v}^{\st,\lambda_v})[\frac{1}{p}]$. In particular, $R_{\rhobar_v}^{\triangle,\lambda_v}$ is reduced.
\end{enumerate}
\end{thm}
\begin{proof}
	The first and third parts are \cite[Theorem 2.7.6, Corollary 2.7.7]{kisindefrings}. The fourth part is a consequence of \cite[Theorem 3.3.8]{kisindefrings}. The second part is a consequence of (a very special case of) \cite[Theorem 3.3.3]{bellovin-gee} (see also \cite{bellovin,hartl-hellmann}). The fifth part is \cite[Lemma 3.10]{geraghty}. The sixth part follows from Lemma \ref{lem:sst-ordinary deformations}.
\end{proof}
\begin{remark}
We have emphasized the reducedness of our local deformation rings because in some parts of the literature the maximal reduced quotients of Kisin's deformation rings are introduced, characterized by their morphisms to finite field extensions of $E$ (e.g.~in \cite{BLGGT}). It won't make any difference to us in practice, because in the end we will be considering maps from local deformation rings to reduced finite flat $\cO$-algebras $\tA$ as in \S\ref{subsec:determinants}.
\end{remark}

\subsubsection{Fixed determinant deformation rings}\label{sec:fixeddetdefrings}
It is often useful to consider deformation rings with fixed determinant. Let $\lambda_v$ be a dominant weight for $(\mathrm{Res}_{F_v/\Qp}\GL_n)_E$ and suppose $\psi:G_{F_{v}}\to \cO^\times$ is a crystalline character with $\tau$-labelled Hodge--Tate weights $\sum_{i=1}^n \lambda_{\tau,i}+(n-i)$ for each $\tau:F_v\hookrightarrow E$. Suppose moreover that $\det\rhobar_{v}$ coincides with $\overline{\psi}:G_{F_v}\to k^\times$, the reduction of $\psi$. Then we have a quotient $R_{\rhobar_v}^{\square,\psi}$ of $R_{\rhobar_v}^{\square}$, classifying liftings with determinant $\psi$ (composed with the structure map from $\cO$ to the test $\cO$-algebra). 

We define quotients $R_{\rhobar_v}^{\cris,\lambda_v,\psi} = R_{\rhobar_v}^{\cris,\lambda_v}\otimes_{R_{\rhobar_v}^{\square}}R_{\rhobar_v}^{\square,\psi}$ and $R_{\rhobar_v}^{\triangle,\lambda_v,\psi} = R_{\rhobar_v}^{\triangle,\lambda_v}\otimes_{R_{\rhobar_v}^{\square}}R_{\rhobar_v}^{\square,\psi}$ of our local deformation rings.

\begin{lem}\label{lem:twistingdefrings} Suppose $p\nmid n$. Let $\lambda_v$ and $\psi$ be as above. Let $R = R_{\rhobar_v}^{\cris,\lambda_v}$ or $R_{\rhobar_v}^{\triangle,\lambda_v}$, and $R^\psi = R_{\rhobar_v}^{\cris,\lambda_v,\psi}$ or $R_{\rhobar_v}^{\triangle,\lambda_v,\psi}$ respectively. Then there is a section to the quotient map $R\to R^\psi$ extending to an isomorphism $R^\psi\llb X \rrb \cong R$. In particular, $R^\psi$ is $\cO$-flat and reduced and can be characterized using the properties of Theorem \ref{thm:localdefrings} parts \eqref{part:crisdefring} or \eqref{part:orddefring} respectively for maps from $R^{\square,\psi}_{\rhobar_v}$ to finite $E$-algebras.
\end{lem}
\begin{proof}
	This is essentially \cite[Lemma 4.3.1]{emertongeerefinedBM}. If we consider the universal lifting $\rho^{\mathrm{univ}}:G_{F_{v}}\to\GL_n(R)$, then the composition of the character $\psi^{-1}\det\rho^{\mathrm{univ}}$ with any map $R \to B$ to a finite $E$-algebra is crystalline with all labelled Hodge--Tate weights equal to zero. In other words these compositions are unramified. Since $R$ is $\cO$-flat and Noetherian, it follows that $\psi^{-1}\det\rho^{\mathrm{univ}}$ is itself unramified. This character is also residually trivial. Since $p\nmid n$, Hensel's lemma implies that there is an unramified character $\alpha:G_{F_v} \to 1+\m_R$ with $\alpha^n = \psi^{-1}\det\rho^{\mathrm{univ}}$. The representation $\alpha^{-1}\otimes\rho^{\mathrm{univ}}$ has determinant $\psi$ and defines a section of the quotient map $R\to R^\psi$ (we can now deduce that $R^\psi$ is $\O$-flat, since it injects into $R$). This extends to a map $R^\psi\llb X \rrb \to R$ sending $X$ to $\alpha(\Frob_v)-1$. We can identify $R^\psi\llb X \rrb$ with a quotient of $R^\square_{\rhobar_v}$ (in fact, of $R$) using the lifting $\rho^{\mathrm{univ},\psi}\otimes\ur(1+X)$, the twist of the universal lifting to $R^\psi$ with the unramified character taking $\Frob_v$ to $1+X$. Composing this lifting with our map $R^\psi\llb X \rrb \to R$ gives $\rho^{\mathrm{univ}}$, and it follows that our map is an isomorphism $R^\psi\llb X \rrb \cong R$.
\end{proof}

\section{Local-global compatibility in the crystalline case}\label{sec:LGC}

\subsection{A computation of boundary cohomology} \label{sec:boundary coh}
The goal of this section (Corollary \ref{cor:direct summand}) is to describe, in terms of the cohomology of the $G$-locally symmetric spaces, a particular direct summand of the completed cohomology for the boundary of the $\tG$-locally symmetric spaces. It will be one of the key ingredients allowing us to describe Hecke algebras acting on cohomology for $G$ in an arbitrary fixed cohomological degree in terms of the middle degree cohomology for $\tG$.

\subsubsection{Notation} Let $\bar{T}\supseteq \bar{S}_p$   
be a finite set of finite places of $F^+$ with preimage the finite set of 
finite places $T$ of $F$. 
Let $\bar{S}\subseteq \bar{S}_p$ be a set of 
primes of $F^+$ dividing $p$ with preimage 
$S\subseteq S_p$. For $\Grm = \tG, P, U$, or $G$, we set $\Grm^0_{\bar{S}} = \prod_{\bar{v}\in\bar{S}} \Grm(\cO_{F^+_{\bar{v}}})$.

Given
$\tilde{\lambda}\in (\mathbb{Z}^{2n}_+)^{\Hom(F^+, E)}$ 
a dominant weight for $\tG$, we define
\[
\cV_{\tilde{\lambda}_{\bar{S}}}:= \bigotimes_{\bar{v}\in \bar{S}} 
\bigotimes_{\tau\in \Hom(F^+_{\bar{v}}, E),\cO} \cV_{\tilde{\lambda}_{\tau}}. 
\]
We also define the object 
$\cV_{U}(\tilde{\lambda}_{\bar{S}}, m)$ in  
$D^b(\mathrm{Sh}_{G^{T}\times K_{T}}(\overline{\mathfrak{X}}_{G},\cO/\varpi^m))$
corresponding to the object 
$R\Gamma(U^0_{\bar{S}}, \cV_{\tilde{\lambda}_{\bar{S}}}/\varpi^m)$ of 
$D^b_{\sm}(K_{S},\cO/\varpi^m)$ 
(after inflation to $D^b_{\sm}(G^{T}\times K_{T},\cO/\varpi^m)$). 
We let $\cV^j_{U}(\tilde{\lambda}_{\bar{S}}, m)$ denote its 
cohomology sheaves; note that these are non-zero precisely when $j$ ranges from $0$ 
to $n^2\sum_{\bar{v}\in \bar{S}} [F^+_{\bar{v}}:\Q_p]$, by the K\"unneth formula
and a computation of group cohomology using a Koszul complex
(e.g.~using \cite[Lemma 7.3]{bms} as in the proof of Lemma~\ref{lem:splitting cohomology}).

If $K \subset G(\AA_{F^+,f})$ is a good subgroup, each $\cV_U(\tilde{\lambda}_{\bar{S}}, m)$ descends to an object in $D^b(\Sh(X_K,\cO/\varpi^m))$ with locally constant cohomology sheaves. Taking a homotopy limit gives $\cV_U(\tilde{\lambda}_{\bar{S}}) \in D^b(\Sh(X_K,\cO))$ again with locally constant cohomology sheaves (since $X_K$ is locally contractible we can find opens over which the $\cV^j_U(\tilde{\lambda}_{\bar{S}}, m)$ are constant for all $m \ge 1$). By passing to a limit over $m$, the cohomology  $R\Gamma(X_K,\cV_U(\tilde{\lambda}_{\bar{S}})\otimes \cV_{\lambda_{\bar{S}_p\backslash \bar{S}}})$ comes with whatever actions (e.g.~of $\TT^T$ or a monoid acting at $\bar{S}_p\backslash \bar{S}$) we have on $R\Gamma(X_K,\cV_U(\tilde{\lambda}_{\bar{S}},m)\otimes \cV_{\lambda_{\bar{S}_p\backslash \bar{S}}}/\varpi^m)$. 

We will often have an action of a commutative $\cO$-algebra $\TT$ (one of our Hecke algebras) on an $\cO$-module or, more generally, an object $C$ in $D(\cO)$. We will then write $\TT(C)$ for the image of $\TT$ in $\End_{D(\cO)}(C)$.

Given $\lambda\in (\mathbb{Z}^n_+)^{\Hom(F, E)}$ 
a dominant weight for $G$, we
define 
\[
\cV_{\lambda_{S}}:= \bigotimes_{v\in S} 
\bigotimes_{\tau\in \Hom(F_v, E),\cO} \cV_{\lambda_{\tau}}. 
\]

\subsubsection{Boundary cohomology}

\begin{thm}\label{thm:direct summand}

Assume that $\widetilde{K}\subset \widetilde{G}(\A_{F^+,f})$ is 
a good subgroup that is decomposed with respect to $P$, 
and with the property that, for each $\bar{v}\in \bar{S}_p$, 
$\widetilde{K}_{U,\bar{v}} = U^0_{\bar{v}}$. 
Let $\mathfrak{m}\subset \mathbb{T}^T$ be a non-Eisenstein
maximal ideal and let $\widetilde{\mathfrak{m}}:= \cS^*(\mathfrak{m})
\subset \widetilde{\mathbb{T}}^T$.  

Choose a partition 
\[
\bar{S}_p = \bar{S}_1\sqcup \bar{S}_2
\]
of the set $\bar{S}_p$ of primes of $F^+$ lying above $p$. 
Let $\tilde{\lambda}$ be a dominant weight 
for $\tG$ satisfying $\tilde{\lambda}_{\bar{v}} = 0$ for all $\bar{v}\in \bar{S}_2$. 
Then 
\[
\cS^*\circ \mathrm{Ind}_{P_{\bar{S}_2}}^{\tG_{\bar{S}_2}}R\Gamma\left(K^{\bar{S}_2}, 
R\Gamma\left(\overline{\mathfrak{X}}_{G}, \cV_{U}(\tilde{\lambda}_{\bar{S}_1}, m)\right)\right)_{\m}
\]
is a $\widetilde{\mathbb{T}}^T$-equivariant direct summand of 
\[
R\Gamma\left(\tK^{\bar{S}_2}, R\Gamma\left(\partial \mathfrak{X}_{\tG}, 
\cV_{\tilde{\lambda}}/\varpi^m\right)\right)_{\widetilde{\m}}
\]
in $D^+_{\sm}(\tG_{\bar{S}_2},\cO/\varpi^m)$. 

\end{thm}

\begin{proof} 
By combining Proposition~\ref{prop:borel-serre strata} and Lemma~\ref{lem:smooth induction},
we obtain a $\widetilde{\mathbb{T}}^T$-equivariant 
isomorphism in $D^+_{\sm}(\tG_{\bar{S}_2},\cO/\varpi^m)$
\begin{equation}\label{eq:induction}
R\Gamma\left(\tK^{\bar{S}_2}, R\Gamma\left(\partial \mathfrak{X}_{\tG}, 
\cV_{\tilde{\lambda}}/\varpi^m\right)\right)_{\widetilde{\m}}
\toisom R\Gamma\left(\tK^{\bar{S}_2}, \text{Ind}_{P^{\bar{S}_1}\times P^0_{\bar{S}_1}}
^{\tG^{\bar{S}_1}\times \tG^0_{\bar{S}_1}}
R\Gamma\left(\overline{\mathfrak{X}}_P, \cV_{\tilde{\lambda}}/\varpi^m\right)\right)_{\widetilde{\m}}. 
\end{equation}
Using the analogue for smooth representations of~\cite[Corollary 2.6]{new-tho}, 
we can rewrite the RHS of~\eqref{eq:induction}
as $r_P^*$ applied to 
\begin{equation}\label{eq:Iwasawa}
\bigoplus_{g}
\text{Ind}_{P_{\bar{S}_2}}^{\tG_{\bar{S}_2}} 
R\Gamma\left(g\tK^{\bar{S}_2}g^{-1}\cap P(\A_{F^+,f}), 
R\Gamma(\overline{\mathfrak{X}}_P, \cV_{\tilde{\lambda}}/\varpi^m)\right)_{r_G^*\m},  
\end{equation}
where $g$ runs over the double cosets 
\[
(P_{T\setminus \bar{S}_p}\times P^0_{\bar{S}_1})
\backslash (\tG_{T\setminus \bar{S}_p}\times \tG^0_{\bar{S}_1})
/(\tK_{T\setminus \bar{S}_p}\times \tK_{\bar{S}_1}),
\] 
and we view each $R\Gamma\left(g\tK^{\bar{S}_2}g^{-1}\cap P(\A_{F^+,f}), 
R\Gamma(\overline{\mathfrak{X}}_P, \cV_{\tilde{\lambda}}/\varpi^m)\right)$
as a $\mathbb{T}_P^{T}$-module in $D^+_{\sm}(P_{\bar{S}_2}, \cO/\varpi^m)$. We restrict to the direct summand corresponding to $g=1$ in~\eqref{eq:Iwasawa}. 
We now conclude by Proposition~\ref{prop:parabolic to levi}. 
\end{proof}

\begin{prop}\label{prop:borel-serre strata}\leavevmode
\begin{enumerate}
\item We have a $\tG(\A_{F^+,f})$-equivariant closed immersion 
\[
(\overline{\mathfrak{X}}_P\times \tG(\A_{F^+,f}))/
P(\A_{F^+,f})\hookrightarrow 
\partial \mathfrak{X}_{\tG}\footnote{Here, we consider the 
groups $\tG(\A_{F^+,f})$ and 
$P(\A_{F^+,f})$ 
as endowed with their natural profinite topologies.},
\]
whose complement is a disjoint union of locally closed subspaces of the form 
$(\mathfrak{X}_Q \times \tG(\A_{F^+,f}))/ Q(\A_{F^+,f})$ with 
$Q\subset \tG$ a standard $F^+$-rational 
parabolic such that $Q\not\subseteq P$. 
\item With all assumptions as in Theorem~\ref{thm:direct summand}, 
the natural pullback map induces a $\widetilde{\mathbb{T}}^T$-equivariant 
isomorphism in $D^+_{\sm}(\tG_{\bar{S}_2},\cO/\varpi^m)$
\[
R\Gamma\left(\tK^{\bar{S}_2}, R\Gamma\left(\partial \mathfrak{X}_{\tG}, 
\cV_{\tilde{\lambda}}/\varpi^m\right)\right)_{\widetilde{\m}}\toisom 
R\Gamma\left(\tK^{\bar{S}_2}, R\Gamma\left(
(\overline{\mathfrak{X}}_{P}\times \tG(\A_{F^+,f}))/P(\A_{F^+,f}), 
\cV_{\tilde{\lambda}}/\varpi^m\right)\right)_{\widetilde{\m}}.
\]

\end{enumerate}
\end{prop}

\begin{proof} The first part follows from the description of the boundary of
the Borel--Serre compactification in~\cite[\S 3.1.2]{new-tho}, see especially Lemma 3.10 of \emph{loc.~cit.} This reference uses the `discrete' versions $\mathfrak{X}_Q^{\mathrm{dis}}$ of the spaces $\frakX_Q$, but we can compare the two situations after taking quotients by compact open subgroups, and the space $(\mathfrak{X}_Q \times \tG(\A_{F^+,f}))/ Q(\A_{F^+,f})$ is equal to the inverse limit of its quotients by compact open subgroups $\tK$ of $\tG(\A_{F^+,f})$ (we can use the fact that $\tG(\A_{F^+,f}) \rightarrow Q(\A_{F^+,f})\backslash\tG(\A_{F^+,f})$ is a locally trivial $Q(\A_{F^+,f})$-torsor, as in the proof of Proposition \ref{prop:ses}). These quotients can in turn be computed as quotients of  $(\mathfrak{X}_Q^{\mathrm{dis}} \times \tG(\A_{F^+,f})^{\mathrm{dis}})/ Q(\A_{F^+,f})^{\mathrm{dis}}$, because $ Q(\A_{F^+,f})\backslash\tG(\A_{F^+,f})/\tK$ is a finite set for each compact open $\tK$.

For the second part, we 
first note that we can check whether a map is an isomorphism on
the level of cohomology groups. For a standard $F^+$-rational proper parabolic 
subgroup $Q\subset \tG$ and a good subgroup $\tK\subset \tG(\A_{F^+,f})$, set 
\[
\tX^Q_{\tK}:= (\mathfrak{X}_Q \times \tG(\A_{F^+,f}))/ Q(\A_{F^+,f})\tK, 
\]
which is a disjoint union of finitely many locally symmetric spaces for $Q$. 
Using excision and Lemma~\ref{lem:recovering completed cohomology}, 
we see that it is enough to show that, for $Q\subset \tG$ a standard $F^+$-rational 
parabolic with $Q\not \subseteq P$, we have 
\[
H^i_c (\tX^Q_{\tK}, \cV_{\tilde{\lambda}}/\varpi^m)_{\widetilde{\m}} = 0, 
\]
for any $i\in \Z_{\geq 0}$. This is standard by now, see for example
the proof of~\cite[Theorem 2.4.2]{10author}. 
\end{proof}

\begin{lemma}\label{lem:smooth induction} Keep the assumption on $\tilde{\lambda}$ from the statement
of Theorem~\ref{thm:direct summand}. There is a natural 
isomorphism in $D^+_{\sm}(\tG^{\bar{S}_1}\times \tG^0_{\bar{S}_1}, 
\cO/\varpi^m)$ 
\[
R\Gamma\left((\overline{\mathfrak{X}}_{P}\times \tG(\A_{F^+,f}))/P(\A_{F^+,f}), 
\cV_{\tilde{\lambda}}/\varpi^m\right) \toisom 
\mathrm{Ind}^{\tG^{\bar{S}_1}\times \tG^0_{\bar{S}_1}}_{P^{\bar{S}_1}\times P^0_{\bar{S}_1}} 
R\Gamma\left(\overline{\mathfrak{X}}_P, 
\cV_{\tilde{\lambda}}/\varpi^m\right).
\]
\end{lemma}

\begin{proof} 
The case of constant coefficients $\cO/\varpi^m$
follows from Lemma~\ref{lem:induction lemma} combined with 
the Iwasawa decomposition at primes in $\bar{S}_1$. For coefficients in a
local system $\cV_{\tilde{\lambda}}/\varpi^m$, where by assumption the action is non-trivial only at the 
primes in $\bar{S}_1$, 
we use the projection formula in Lemma~\ref{lem:projection formula}
and the tensor identity in Lemma~\ref{lem:tensor identity} (at primes in $\bar{S}_1$)
to reduce to the case of constant coefficients. 
\end{proof}

\begin{lemma}\label{lem:completed cohomology}
There is a natural isomorphism 
\[
\mathrm{Inf}_{G(\A_{F^+,f})}^{P(\A_{F^+,f})}
R\Gamma\left(\overline{\mathfrak{X}}_G, \cO/\varpi^m\right)
\toisom R\Gamma\left(\overline{\mathfrak{X}}_P, \cO/\varpi^m\right)
\]
in $D^+_{\sm}(P(\A_{F^+,f}), \cO/\varpi^m)$. 
\end{lemma}

\begin{proof} The map is given by
pullback along the $P(\A_{F^+,f})$-equivariant projection 
$\overline{\mathfrak{X}}_P\twoheadrightarrow \overline{\mathfrak{X}}_G$. 
To show that the map is an isomorphism, it is enough to check
that it induces an isomorphism after applying the forgetful functor to $D^+(\cO/\varpi^m)$,
and then it is enough to consider it on the level of cohomology groups. 
By Lemma~\ref{lem:recovering completed cohomology}, this 
reduces to establishing the isomorphism 
\[
\varinjlim_{K} H^i(X_K, \cO/\varpi^m)\toisom 
\varinjlim_{K_P} H^i(X^P_{K_P},\cO/\varpi^m),
\]
where $K_{P}\subset P(\A_{F^+,f})$ runs over compact open subgroups 
of the form $K \ltimes K_U$ with $K\subset G(\A_{F^+,f})$ and 
$K_{U}\subset U(\A_{F^+,f})$ compact open subgroups. 
We use the Leray--Serre spectral sequence for the fibration 
$X^U_{K_U}\to X^P_{K_P}\to X_K$. 
It is enough to check that 
\[
\varinjlim_{K_U} H^j(X^U_{K_U},\cO/\varpi^m) = 
\varinjlim_{\Gamma_U:= U(F^+) \cap K_U} H^j(\Gamma_U \backslash 
U(F^+\otimes_{\Q} \mathbb{R}),\cO/\varpi^m)  = 
\begin{cases}
\cO/\varpi^m & \text{for } j = 0 \\
0 & \text{otherwise.}
\end{cases}
\]
The first equality follows by strong approximation for $U$, 
the second by direct computation. 
We note that this computation should be thought of in the 
category of $\cO/\varpi^m[\![K]\!]$-modules, with a trivial action of $K$ on the RHS. 
\end{proof}

\begin{lemma}\label{lem:parabolic to levi}
With the notation as in the proof of Theorem~\ref{thm:direct summand},
there exists a natural isomorphism 
\[
\mathrm{Inf}_{G^{T\setminus \bar{S}_2}}^{P^{T\setminus \bar{S}_2}} 
R\Gamma\left (K_{T\setminus \bar{S}_2},  
R\Gamma\left(\overline{\mathfrak{X}}_{G}, \cV_{U}(\tilde{\lambda}_{\bar{S}_1}, m)\right)\right)
\toisom 
R\Gamma\left(K_{P, T\setminus \bar{S}_2}, R\Gamma (
\overline{\mathfrak{X}}_P, \cV_{\tilde{\lambda}}/\varpi^m )
\right)
\]
in $D^+_{\sm}(P^{T\setminus \bar{S}_2}, \cO/\varpi^m)$. 
\end{lemma}

\begin{proof} By Lemma~\ref{lem:projection formula derived}, 
we have isomorphisms  
\[
R\Gamma\left(\overline{\mathfrak{X}}_{G}, \cV_{U}(\tilde{\lambda}_{\bar{S}_1}, m)\right)
\toisom 
R\Gamma(\overline{\mathfrak{X}}_G, \cO/\varpi^m) \otimes^{\mathbb{L}}_
{\cO/\varpi^m} R\Gamma(K_{U, \bar{S}_1},\cV_{\tilde{\lambda}}/\varpi^m)
\]
\[
\toisom R\Gamma(\overline{\mathfrak{X}}_G, \cO/\varpi^m) 
\otimes^{\mathbb{L}}_{\cO/\varpi^m} 
R\Gamma(K_{U, T\setminus \bar{S}_2},\cV_{\tilde{\lambda}}/\varpi^m)
\]
in $D^+_{\sm}(G^{T\setminus \bar{S}_2}\times K_{T\setminus \bar{S}_2},\cO/\varpi^m)$.
We have a natural isomorphism 
\[
\mathrm{Inf}_{G^{T\setminus\bar{S}_2}}^{P^{T\setminus\bar{S}_2}}
\left(R\Gamma(\overline{\mathfrak{X}}_G, \cO/\varpi^m) 
\otimes^{\mathbb{L}}_{\cO/\varpi^m} 
R\Gamma(K_{U, T\setminus \bar{S}_2},\cV_{\tilde{\lambda}}/\varpi^m)\right)
\toisom 
\]
\[
R\Gamma \left(K_{U, T\setminus \bar{S}_2}, \mathrm{Inf}_{G(\A_{F^+, f})}^{P(\A_{F^+,f})} 
R\Gamma(\overline{\mathfrak{X}}_G, \cO/\varpi^m) 
\otimes^{\mathbb{L}}_{\cO/\varpi^m} \cV_{\tilde{\lambda}}/\varpi^m \right) 
\]
in $D^+_{\sm}(P^{T\setminus \bar{S}_2}\times K_{T\setminus \bar{S}_2}, \cO/\varpi^m)$
because an arbitrary colimit of injective representations of $K_{U, T\setminus \bar{S}_2}$ 
is injective. We conclude by Lemma~\ref{lem:completed cohomology}, by 
Hochschild--Serre applied to $K_{P, T\setminus \bar{S}_2} = 
K_{T\setminus \bar{S}_2}\ltimes K_{U, T\setminus \bar{S}_2}$, 
and by Lemma~\ref{lem:projection formula} applied to $\overline{\mathfrak{X}}_P$. 
\end{proof}

\begin{prop}\label{prop:parabolic to levi} 
With the notation as in the proof of Theorem~\ref{thm:direct summand}, 
we have a $\mathbb{T}^T_P$-equivariant isomorphism 
\[
R\Gamma\left(K_P^{\bar{S}_2}, 
R\Gamma (\overline{\mathfrak{X}}_P, \cV_{\tilde{\lambda}}/\varpi^m)\right)
\toisom 
r_G^* \circ \mathrm{Inf}_{G_{\bar{S}_2}}^{P_{\bar{S}_2}} R\Gamma\left(K^{\bar{S}_2},
R\Gamma\left(\overline{\mathfrak{X}}_{G}, \cV_{U}(\tilde{\lambda}_{\bar{S}_1}, m)\right)\right) 
\]
in $D^+_{\sm}(P_{\bar{S}_2}, \cO/\varpi^m)$. 
\end{prop}

\begin{proof} We separate the set of finite places of $F^+$ away from $\bar{S}_2$
into the union of the set of finite places away from $T$ and the set of finite places $T\setminus \bar{S}_2$, 
which contains $\bar{S}_1$. 
By the analogue for smooth representations 
of~\cite[Corollary 2.8]{new-tho}, 
there is a $\mathbb{T}^T_P$-equivariant 
natural transformation in $D^+_{\sm}(P_{\bar{S}_2},\cO/\varpi^m)$ between 
\[
r_G^* R\Gamma\left(K^T, \mathrm{Inf}_{G_{\bar{S}_2}}^{P_{\bar{S}_2}} 
R\Gamma\left(K_{T\setminus \bar{S}_2}, R\Gamma\left(\overline{\mathfrak{X}}_{G}, \cV_{U}(\tilde{\lambda}_{\bar{S}_1}, m)\right)\right)\right)
\]
and 
\[
R\Gamma\left(K_P^T, \mathrm{Inf}_{G^{T\setminus \bar{S}_2}}^{P^{T\setminus \bar{S}_2}} 
R\Gamma\left(K_{T\setminus \bar{S}_2}, R\Gamma\left(\overline{\mathfrak{X}}_{G}, \cV_{U}(\tilde{\lambda}_{\bar{S}_1}, m)\right)\right)\right).
\]
If we can show that this natural transformation is an isomorphism, we can 
conclude by Lemma~\ref{lem:parabolic to levi}. To prove that the natural transformation
is an isomorphism, we can forget the $\mathbb{T}^T_P$-action and work 
in $D^+(\cO/\varpi^m)$. There, it reduces to the statement that  
taking $K^T_U$-invariants is an exact functor, 
which holds true because $K^T_U$ is a profinite abelian group 
with trivial $p$-part. 
\end{proof}

The following is a result in the style of~\cite[Theorems 4.2.1 and 5.4.1]{10author}. 
This is the only part of this subsection that we will use in what follows. 

\begin{cor}\label{cor:direct summand}
Assume that $\widetilde{K}\subset \widetilde{G}(\A_{F^+,f})$ is 
a good subgroup that is decomposed with respect to $P$, 
and with the property that, for each $\bar{v}\in \bar{S}_p$, 
$\widetilde{K}_{U,\bar{v}} = U^0_{\bar{v}}$. 
Let $\mathfrak{m}\subset \mathbb{T}^T$ be a non-Eisenstein
maximal ideal and let $\widetilde{\mathfrak{m}}:= \cS^*(\mathfrak{m})
\subset \widetilde{\mathbb{T}}^T$.  

Choose a partition 
\[
\bar{S}_p = \bar{S}_1\sqcup \bar{S}_2\sqcup \bar{S}_3
\]
of the set $\bar{S}_p$ of primes of $F^+$ lying above $p$. 
Let $\tilde{\lambda}$ and 
$\lambda$ 
be dominant weights for $\widetilde{G}$ 
and $G$, respectively. Assume that the following conditions are satisfied: 

\begin{enumerate}
\item For each $\tau:F^+\hookrightarrow E$ 
inducing a place $\bar{v}\in \bar{S}_1$, $\tilde{\lambda}_{\tau} = (-w_{0,n}\lambda_{\tilde{\tau}c}, \lambda_{\tilde{\tau}})$ 
(identification as in~\eqref{eq:identification of weights}); 
\item For each $\tau:F^+\hookrightarrow E$ 
inducing a place $\bar{v}\in \bar{S}_2\sqcup \bar{S}_3$, $\tilde{\lambda}_{\tau} = 0$. 
\end{enumerate}

Then 
\[
\cS^*\circ \mathrm{Ind}_{P_{\bar{S}_3}}^{\tG_{\bar{S}_3}}
R\Gamma\left(K^{\bar{S}_3}, R\Gamma\left(\overline{\mathfrak{X}}_G, \cV_{\lambda_{\bar{S}_1}}/\varpi^m
\otimes \cV_{U}(\tilde{\lambda}_{\bar{S}_2}, m)\right)\right)_{\mathfrak{m}} 
\]
is a $\widetilde{\mathbb{T}}^T$-equivariant direct summand of 
\[
R\Gamma\left(\widetilde{K}^{\bar{S}_3}, R\Gamma\left(\partial \mathfrak{X}_{\tG}, \cV_{\tilde{\lambda}}/\varpi^m
\right)\right)_{\widetilde{\mathfrak{m}}}
\]
in $D^+_{\sm}(\tG_{\bar{S}_3},\cO/\varpi^m)$. 
\end{cor}

\begin{proof} This follows from Theorem~\ref{thm:direct summand} 
applied with $\bar{S}'_1:= \bar{S}_1\sqcup \bar{S}_2$ and 
$\bar{S}'_2:= \bar{S}_3$, 
as long as we can show that $\cV_{\lambda_{\bar{S}_1}}/\varpi^m$ 
is a $G^T\times K_T$-equivariant direct summand
of $\cV_{U}(\tilde{\lambda}_{\bar{S}_1}, m)$. This follows
from \cite[Corollary 2.11]{new-tho}, as in the proof of~\cite[Theorem 4.2.1]{10author} (and the proof of~\cite[Theorem 2.4.4]{10author}). 
\end{proof}

\subsection{The integral case}
In this section we prove our main results on local--global compatibility: Proposition \ref{prop:lgc torsion CTG} and Theorem \ref{thm:RtoT factors thru Kisin def ring}.
\subsubsection{Degree shifting}\label{subsec:degree shifting}

Let $\bar{S}\subset \bar{S}_p$. Recall that, for each $\bar{v}\in \bar{S}_p$, 
we have chosen a place $\tilde{v}\mid \bar{v}$ of $F$, with complex
conjugate $\tilde{v}^c$. The isomorphism $\iota_{\tilde{v}}:\widetilde{G}(F^+_{\bar{v}}) \cong \GL_{2n}(F_{\tilde{v}})$ identifies the Levi subgroup $G(F^+_{\bar{v}}) = \GL_n(F_{\tilde{v}}) \times \GL_n(F_{\tilde{v}^c})$ with block diagonal matrices in $\GL_{2n}(F_{\tilde{v}})$ via $(A_{\tilde{v}},A_{\tilde{v}^c}) \mapsto \begin{pmatrix}(\Psi_n\ ^tA_{\tilde{v}^c}^{-1}\Psi_n)^c & 0 \\ 0 & A_{\tilde{v}} \end{pmatrix}$.

We have standard parabolics $Q_{\bar{v}} \subset P_{\bar{v}}$ for each $\bar{v} \in \bar{S}$ and we set $K_{\bar{v}}:= \cQ_{\bar{v}} \cap G_{\bar{v}}(F^+_{\bar{v}})$, a parahoric subgroup of $G_{\bar{v}}(F^+_{\bar{v}})$. The conjugate parahoric $(K_{\bar{v}})^{w^P_0}$ can also be viewed as the intersection $\overline{\cQ}_{\bar{v}}^{w^P_0} \cap G_{\bar{v}}(F^+_{\bar{v}})$, where $\overline{Q}_{\bar{v}}^{w^P_0}$ is the standard parabolic with Levi subgroup $Q_{\bar{v}}^{w^P_0}\cap G_{\bar{v}}(F^+_{\bar{v}})$.

We introduce the following notation on the level of abstract 
Hecke algebras:
\begin{align*}\mathbb{T}^{\cQ_{\bar{S}},\bar{S}-\ord} &:= \mathbb{T}^T\otimes_{\ZZ}\left(\bigotimes_{\bar{v}\in\bar{S}}\cH(\Delta^{\cQ_{\bar{v}}}_{\bar{v}},K_{\bar{v}}) \right),
\mathbb{T}^{\cQ_{\bar{S}},\bar{S}-\ord}_{w^P_0}:= \mathbb{T}^T\otimes_{\ZZ}\left(\bigotimes_{\bar{v}\in\bar{S}}\cH((\Delta^{\cQ_{\bar{v}}}_{\bar{v}})^{w^P_0},K_{\bar{v}}^{w^P_0}) \right)\\ \mathrm{and}\ 
\widetilde{\mathbb{T}}^{\cQ_{\bar{S}},\bar{S}-\ord}&:= \widetilde{\mathbb{T}}^T\otimes_{\ZZ}\left(\bigotimes_{\bar{v}\in\bar{S}}\cH(\widetilde{\Delta}^{\cQ_{\bar{v}}}_{\bar{v}}, \cQ_{\bar{v}})[\widetilde{U}_{\tilde{v},n}^{-1}] \right). 
\end{align*}
The unnormalized Satake transform $\cS:\widetilde{\mathbb{T}}\to 
\mathbb{T}$ extends to a morphism 
denoted by 
\[
\cS^{w^P_0}:\widetilde{\mathbb{T}}^{\cQ_{\bar{S}},\bar{S}-\ord}
\to \mathbb{T}^{\cQ_{\bar{S}},\bar{S}-\ord}_{w^P_0},
\] 
given by $[\cQ_{\bar{v}}\nu(\varpi_{\bar{v}})\cQ_{\bar{v}}] \mapsto [K_{\bar{v}}^{w^P_0}\nu(\varpi_{\bar{v}})^{w^P_0}K_{\bar{v}}^{w^P_0}]$. In particular, $\widetilde{U}_{\tilde{v},n}$ maps to $U_{\tilde{v}} = [K_{\bar{v}}^{w^P_0}\diag(\varpi_{\tilde{v}},\cdots,\varpi_{\tilde{v}})K_{\bar{v}}^{w^P_0}]$ and $\widetilde{U}_{\tilde{v},2n}$ maps to 
\[
U_{\tilde{v}}U_{\tilde{v}^c}^{-1} = [K_{\bar{v}}^{w^P_0}\left(\diag(\varpi_{\tilde{v}},\cdots,\varpi_{\tilde{v}}),(\diag(\varpi_{\tilde{v}},\cdots,\varpi_{\tilde{v}})^{-1})^c\right)K_{\bar{v}}^{w^P_0}].
\]

To make use of Poincar\'{e} duality, we will need to introduce Hecke algebras twisted by the duality involutions $\iota, \tilde{\iota}$ which are defined by \[\iota[K^{w^P_0}gK^{w^P_0}] = [K^{w^P_0}g^{-1}K^{w^P_0}] \text{ and } \tilde{\iota}[\tK g \tK] = [\tK g^{-1} \tK]\] (see~\cite[\S 2.2.20]{10author}). They give isomorphisms \[\iota: \mathbb{T}^{\cQ_{\bar{S}},\bar{S}-\ord}_{w^P_0} \toisom \mathbb{T}^{\cQ_{\bar{S}},\bar{S}-\ord,\iota}_{w^P_0} := \mathbb{T}^T\otimes_{\ZZ}\left(\bigotimes_{\bar{v}\in\bar{S}}\cH(\left((\Delta^{\cQ_{\bar{v}}}_{\bar{v}})^{w^P_0}\right)^{-1},K_{\bar{v}}^{w^P_0}) \right)\] and \[\tilde{\iota}: \widetilde{\mathbb{T}}^{\cQ_{\bar{S}},\bar{S}-\ord} \toisom \widetilde{\mathbb{T}}^{\cQ_{\bar{S}},\bar{S}-\ord,\tilde{\iota}}:=\widetilde{\mathbb{T}}^T\otimes_{\ZZ}\left(\bigotimes_{\bar{v}\in\bar{S}}\cH(\left(\widetilde{\Delta}^{\cQ_{\bar{v}}}_{\bar{v}}\right)^{-1}, \cQ_{\bar{v}})[[\cQ_{\bar{v}}\tilde{u}_{\tilde{v},n}^{-1}\cQ_{\bar{v}}]^{-1}] \right).\] We will also make use of an untwisted Hecke algebra for the Levi, defined by $\mathbb{T}^{\cQ_{\bar{S}},\bar{S}-\ord,\iota} := \mathbb{T}^T\otimes_{\ZZ}\left(\bigotimes_{\bar{v}\in\bar{S}}\cH(\left(\Delta^{\cQ_{\bar{v}}}_{\bar{v}}\right)^{-1},K_{\bar{v}}) \right)$.

We define $\cS^{\iota}: \widetilde{\mathbb{T}}^{\cQ_{\bar{S}},\bar{S}-\ord,\tilde{\iota}}
\to \mathbb{T}^{\cQ_{\bar{S}},\bar{S}-\ord,\iota}$ extending $\cS$ according to the formula \[[\cQ_{\bar{v}}\nu(\varpi_{\bar{v}})^{-1}\cQ_{\bar{v}}] \mapsto [K_{\bar{v}}\nu(\varpi_{\bar{v}})^{-1}K_{\bar{v}}].\]

\begin{prop}\label{prop:initial subquotient_extra_Hecke} 
	Assume that $\widetilde{K}\subset \widetilde{G}(\A_{F^+,f})$ is 
	a good subgroup that is decomposed with respect to $P$ 
	and with the property that, for each $\bar{v}\in \bar{S}_p$, 
	$\widetilde{K}_{U,\bar{v}} = U^0_{\bar{v}}$. 
	Let $\mathfrak{m}\subset \mathbb{T}^T$ be a non-Eisenstein
	maximal ideal, let $\widetilde{\mathfrak{m}}:= \cS^*(\mathfrak{m})
	\subset \widetilde{\mathbb{T}}^T$,  
	and assume that $\bar{\rho}_{\widetilde{\m}}$ is decomposed generic
	in the sense of Definition~\ref{defn:generic}. 
	
	Choose a partition 
	\[
	\bar{S}_p = \bar{S}_1\sqcup \bar{S}_2\sqcup \bar{S}_3
	\]
	of the set $\bar{S}_p$ of primes of $F^+$ lying above $p$, together with standard parabolic subgroups $Q_{\bar{v}}\subset P_{\bar{v}}$ for each $\bar{v}\in \bar{S}_3$. 
	Let $\tilde{\lambda}$ and 
	$\lambda$ 
	be dominant weights for $\widetilde{G}$ 
	and $G$, respectively. Assume that the following conditions are satisfied: 
	
	\begin{enumerate}
		
		\item For each $\tau:F^+\hookrightarrow E$ 
		inducing a place $\bar{v}\in \bar{S}_1$, $\tilde{\lambda}_{\tau} = (-w_{0,n}\lambda_{\tilde{\tau}c}, \lambda_{\tilde{\tau}})$ 
		(identification as in~\eqref{eq:identification of weights}); 
		\item For each $\tau:F^+\hookrightarrow E$ 
		inducing a place $\bar{v}\in \bar{S}_2$, $\tilde{\lambda}_{\tau} = 0$. 
		\item For each $\bar{v}\in \bar{S}_3$, $\tK_{\bar{v}}=\cQ_{\bar{v}}$. For each
		$\tau:F^+\hookrightarrow E$ inducing such a place $\bar{v}$, we also have the standard identification
		$\tilde{\lambda}_{\tau} = (-w_{0,n}\lambda_{\tilde{\tau}c},\lambda_{\tilde{\tau}})$
	\end{enumerate}
	
	Then the unnormalized Satake transform 
	$\cS^{w^P_0}: \widetilde{\mathbb{T}}^{\cQ_{\bar{S}_3},\bar{S}_3-\ord} \to \mathbb{T}^{\cQ_{\bar{S}_3},\bar{S}_3-\ord}_{w^P_0}$ descends to a 
	homomorphism 
	\[
	\widetilde{\mathbb{T}}^{\cQ_{\bar{S}_3},\bar{S}_3-\ord} \left(H^d\left(\tX_{\tK}, 
	\cV_{\tilde{\lambda}}\right)_{\widetilde{\m}}^{\ord}\right) \to
	\mathbb{T}^{\cQ_{\bar{S}_3},\bar{S}_3-\ord}_{w^P_0}\left(\mathbb{H}^d\left( X_{K^{\bar{S}_3}K_{\bar{S}_3}^{w^P_0}}, \cV_{\lambda_{\bar{S}_1}}\otimes \cV_{U}(\tilde{\lambda}_{\bar{S}_2})\otimes 
	\cV_{\lambda_{\bar{S}_3}} \right)_{\m}\right),
	\]
	where $\mathbb{H}^d$ denotes the degree $d$ hypercohomology. The 
	Hecke action on the source is defined in~\eqref{eq:ord Hecke action} and
	the one on the target is defined in~\eqref{eq:ord Hecke action Levi}. The 
	Hecke operators $\widetilde{U}_{\tilde{v},n}$ for $\bar{v}\in \bar{S}_3$ are
	invertible on the source because we are considering the ordinary part of 
	cohomology at the primes in $\bar{S}_3$. 
	
	Moreover, we have an injection 
	\[
	\widetilde{\mathbb{T}}^{\cQ_{\bar{S}_3},\bar{S}_3-\ord}\left(H^d\left(\tX_{\tK}, 
	\cV_{\tilde{\lambda}}\right)_{\widetilde{\m}}^{\ord}\right) 
	\hookrightarrow 
	\widetilde{\mathbb{T}}^{\cQ_{\bar{S}_3},\bar{S}_3-\ord}\left(H^d\left(\tX_{\tK}, 
	\cV_{\tilde{\lambda}}[1/p]\right)_{\widetilde{\m}}^{\ord}\right). 
	\]%
\end{prop}

\begin{proof} We will apply Theorem~\ref{thm:middle degree cohomology}, noting that taking the ordinary part preserves injectivity and surjectivity of the maps in that statement.	This reduces us to proving that $\cS^{w^P_0}$ descends to a homomorphism 
	\[
	\widetilde{\mathbb{T}}^{\cQ_{\bar{S}_3},\bar{S}_3-\ord}\left(H^d\left(\partial\tX_{\tK}, 
	\cV_{\tilde{\lambda}}\right)_{\widetilde{\m}}^{\ord}\right) \to
	\mathbb{T}^{\cQ_{\bar{S}_3},\bar{S}_3-\ord}_{w^P_0}\left(\mathbb{H}^d\left( X_{K^{\bar{S}_3}K_{\bar{S}_3}^{w^P_0}}, \cV_{\lambda_{\bar{S}_1}}\otimes \cV_{U}(\tilde{\lambda}_{\bar{S}_2})\otimes 
	\cV_{\lambda_{\bar{S}_3}} \right)_{\m}\right).
	\]
	
	Combining Proposition \ref{prop:independence of weight boundary} and Lemma \ref{lem:control level extra Hecke} shows that $H^d(\partial \tX_{\tK}, \cV_{\tilde{\lambda}}/\varpi^m)^{\ord}_{\widetilde{\m}}$ is isomorphic as a $\widetilde{\mathbb{T}}^{\cQ_{\bar{S}_3},\bar{S}_3-\ord}$-module with \[R^d\Gamma\left(K_{\bar{S}_3}, \pi^{\ord}_{\del}(\widetilde{K}^{\bar{S}_3},\tilde{\lambda}^{\bar{S}_3},m)\otimes\cV_{\lambda_{\bar{S}_3}}^{w^P_0}/\varpi^m\right)%
	.\] For each $\bar{v}\in \bar{S}_3$, the action of $\cH(\widetilde{\Delta}^{\cQ_{\bar{v}}}_{\bar{v}}, \cQ_{\bar{v}})[\widetilde{U}_{\tilde{v},n}^{-1}]$ is via its isomorphism to $\cH({\Delta}^{\cQ_{\bar{v}}}_{\bar{v}}, K_{\bar{v}})$ (cf.~Lemma \ref{lem:Q Hecke algebra}).
	
	Using Corollary~\ref{cor:direct summand} and interchanging the $P$-ordinary part with the tensor product by $\cV^{w^P_0}_{\lambda_{\bar{S}_3}}$, we find a $\widetilde{\mathbb{T}}^{\cQ_{\bar{S}_3},\bar{S}_3-\ord}$-equivariant direct summand 	\begin{equation}\label{eq:ordinary part extra Hecke}
	\ord_0 \mathbb{H}^d\left(K_{\bar{S}_3} \ltimes U^0_{\bar{S}_3}, 
	\mathrm{Ind}_{P_{\bar{S}_3}}^{\tG_{\bar{S}_3}}
	R\Gamma\left(K^{\bar{S}_3}, R\Gamma\left(\overline{\mathfrak{X}}_G, \cV_{\lambda_{\bar{S}_1}}/\varpi^m
	\otimes \cV_{U}(\tilde{\lambda}_{\bar{S}_2}, m)\right)\right)
	\otimes \cV^{w^P_0}_{\lambda_{\bar{S}_3}}/\varpi^m \right)_{\widetilde{\m}}.
	\end{equation}
	
		We find a subquotient of the term~\eqref{eq:ordinary part extra Hecke} using the $w^P_0$-case 
	of Proposition~\ref{prop:ordsubquotients}. Note that, since $w_0^P$ indexes the open Bruhat cell, the subquotient in this case is actually a submodule. We claim that it is 
	$\widetilde{\mathbb{T}}^{\cQ_{\bar{S}_3},\bar{S}_3-\ord}$-equivariantly 
	isomorphic to 
	\begin{equation}\label{eq:twisted ordinary part extra Hecke}
	\mathbb{H}^d\left(X_{K^{\bar{S}_3}K_{\bar{S}_3}^{w^P_0}}, \cV_{\lambda_{\bar{S}_1}}/\varpi^m
	\otimes \cV_{U}(\tilde{\lambda}_{\bar{S}_2}, m)\otimes \cV_{\lambda_{\bar{S}_3}}/\varpi^m 
	\right)_{\widetilde{\m}},
	\end{equation}
	where the action of $\widetilde{\mathbb{T}}^{\cQ_{\bar{S}_3},\bar{S}_3-\ord}$ on~\eqref{eq:twisted ordinary part extra Hecke} is 
	via the extension $\cS^{w^P_0}$ of the unnormalized Satake transform defined above. This follows from the fact that the action of $G_{\bar{S}_3}$ on $R\Gamma\left(K^{\bar{S}_3},R\Gamma\left(\overline{\mathfrak{X}}_G, \cV_{\lambda_{\bar{S}_1}}/\varpi^m
	\otimes \cV_{U}(\tilde{\lambda}_{\bar{S}_2}, m)\right)\right)$ needs to be 
	pre-conjugated by $w^P_0$ when applying Proposition~\ref{prop:ordsubquotients}. Note that this is compatible with the rescaled Hecke actions, since the fact that $w^P_0$ commutes with $w_0^G$ implies that we have the equality \[\alpha^{\overline{\cQ}_{\bar{v}}^{w^P_0}}_{\lambda}(\nu(\varpi_{\bar{v}})^{w^P_0}) = \alpha^{\cQ_{\bar{v}}}_{\lambda^{w^P_0}}(\nu(\varpi_{\bar{v}})).\]
	 Taking stock, we have exhibited~\eqref{eq:twisted ordinary part extra Hecke} as a $\widetilde{\mathbb{T}}^{\cQ_{\bar{S}_3},\bar{S}_3-\ord}$-equivariant submodule of $H^d\left(\del\tX_{\tK}, 
\cV_{\tilde{\lambda}}
/\varpi^m\right)_{\widetilde{\m}}^{\ord}$, where the action on~\eqref{eq:twisted ordinary part extra Hecke} is via $\cS^{w^P_0}$. 

These presentations are compatible as $m$ varies. Because all the above cohomology groups are finitely generated $\cO$-modules, we can take inverse limits with respect to $m$ to deduce that
\begin{equation}\label{eq:twisted ordinary part extra Hecke integral}
    \mathbb{H}^d
	\left( X_{K^{\bar{S}_3}K_{\bar{S}_3}^{w^P_0}}, \cV_{\lambda_{\bar{S}_1}}\otimes \cV_{U}(\tilde{\lambda}_{\bar{S}_2})\otimes 
	\cV_{\lambda_{\bar{S}_3}}\right)_{\widetilde{\m}}
\end{equation}
is a $\widetilde{\mathbb{T}}^{\cQ_{\bar{S}_3},\bar{S}_3-\ord}$-equivariant submodule of $H^d\left(\del\tX_{\tK}, 
	\cV_{\tilde{\lambda}}\right)_{\widetilde{\m}}^{\ord}$, where the action of $\widetilde{\mathbb{T}}^{\cQ_{\bar{S}_3},\bar{S}_3-\ord}$
    on~\eqref{eq:twisted ordinary part extra Hecke integral} is via $\cS^{w_0^P}$. This gives a homomorphism 
    \[ 
	\widetilde{\mathbb{T}}^{\cQ_{\bar{S}_3},\bar{S}_3-\ord}\left(H^d\left(\del\tX_{\tK}, 
	\cV_{\tilde{\lambda}}\right)_{\widetilde{\m}}^{\ord}\right) \to
	\mathbb{T}^{\cQ_{\bar{S}_3},\bar{S}_3-\ord}_{w_0^P}\left(\mathbb{H}^d
	\left( X_{K^{\bar{S}_3}K_{\bar{S}_3}^{w^P_0}}, \cV_{\lambda_{\bar{S}_1}}\otimes \cV_{U}(\tilde{\lambda}_{\bar{S}_2})\otimes 
	\cV_{\lambda_{\bar{S}_3}}\right)_{\widetilde{\m}}\right).
	\]
	We conclude because $\cS^*({\m}) = \widetilde{\m}$, so $\mathbb{H}^d\left( X_{K^{\bar{S}_3}K_{\bar{S}_3}^{w^P_0}}, \cV_{\lambda_{\bar{S}_1}}
	\otimes \cV_{U}(\tilde{\lambda}_{\bar{S}_2})\otimes 
	\cV_{\lambda_{\bar{S}_3}} \right)_{\m}$
	is a Hecke-equivariant direct summand of
	$\mathbb{H}^d\left( X_{K^{\bar{S}_3}K_{\bar{S}_3}^{w^P_0}}, \cV_{\lambda_{\bar{S}_1}}\otimes \cV_{U}(\tilde{\lambda}_{\bar{S}_2})\otimes 
	\cV_{\lambda_{\bar{S}_3}} \right)_{\widetilde{\m}}$. 
\end{proof}

We have a similar statement to the above which will be useful after applying Poincar\'{e} duality. Most things go through very similarly. At places in $\bar{S}_2$, the weight $\tilde{\lambda}_{\bar{S}_2}$ is trivial so nothing changes when we take a dual. At places in $\bar{S}_3$ we use the results of \S\ref{sec:dual coeffs prelims} and Corollary \ref{cor:ordsubquotients dual}. However, we need to modify things at places in $\overline{S}_1$, because $R\Gamma(U^0_{\bar{S}},\cV_{\tilde{\lambda}_{\bar{S}}}^\vee/\varpi^m)$ may not admit $\cV_{{\lambda}_{\bar{S}}}^\vee/\varpi^m$ as a $K_{\bar{S}}$-equivariant direct summand. This causes a problem for the analogue of Corollary \ref{cor:direct summand} on the dual side. The following lemma will act as a replacement.

\begin{lem}
Let $\bar{S}\subset \bar{S}_p$, let $\tilde{\lambda}$ and $\lambda$ be dominant weights for $\widetilde{G}$ and $G$ respectively. Assume that the following condition is satisfied: \begin{enumerate}
	\item For each $\tau:F^+\hookrightarrow E$ inducing a place $\bar{v}\in\bar{S}$, $\tilde{\lambda}_{\tau} = (\lambda_{\tilde{\tau}},-w_{0,n}\lambda_{\tilde{\tau}c})$ (note that this is the $w^P_0$-conjugate of the standard identification).
\end{enumerate}

Let $m\in \ZZ_{\ge 1}$ be an integer. Then $R\Gamma\left(U^0_{\bar{S}},\cV_{\tilde{\lambda}_{\bar{S}}}^\vee/\varpi^m\right)$ admits $\cV_{\lambda_{\bar{S}}}^\vee/\varpi^m$ as a $K_{\bar{S}}$-equivariant direct summand.
\end{lem} 
\begin{proof}
Taking duals and applying the proof of \cite[Theorem 2.4.4]{10author}, it suffices to show that there is a surjective $P(\cO_{F^+,\bar{S}})$-equivariant map \[\cV_{\tilde{\lambda}_{\bar{S}}} \to \cV_{\lambda_{\bar{S}}}\] with a $K_{\bar{S}}$-equivariant splitting. It follows from \cite[Proposition 2.10]{new-tho} that there is a $K_{\bar{S}}$-equivariant decomposition $\cV_{\tilde{\lambda}_{\bar{S}}} = \cV_{\lambda_{\bar{S}}} \oplus W$, with $\cV_{\lambda_{\bar{S}}}$ invariant under the action of the unipotent subgroup $\overline{U}(\cO_{F^+,\bar{S}})$ in the parabolic opposite to $P$. Moreover, it follows from the main theorem of \cite{Cab84} that, after extending scalars to $\overline{E}$, $W\otimes_{\cO}\overline{E}$ is identified with the $P(F^+_{\bar{S}})$-stable subspace $(1-U(F^+_{\bar{S}}))\cV_{\tilde{\lambda}_{\bar{S},\overline{E}}} \subset \cV_{\tilde{\lambda}_{\bar{S},\overline{E}}}$ (thanks to Lambert A'Campo for pointing this out). Because $W$ is a direct summand of $\cV_{\tilde{\lambda}_{\bar{S}}}$ over $\cO$, it is a saturated $\cO$-submodule. This implies that $W = \cV_{\widetilde{\lambda}_{\bar{S}}}\cap (W\otimes_{\cO}\overline{E})$ is $P(\cO_{F^+,\bar{S}})$-stable, so quotienting out by $W$ gives the desired map $\cV_{\tilde{\lambda}_{\bar{S}}} \to \cV_{\lambda_{\bar{S}}}$.
\end{proof}

We can now state our version of Proposition \ref{prop:initial subquotient_extra_Hecke} with dual coefficients.
\begin{prop}
	\label{prop:initial subquotient_dual}
		Assume that $\widetilde{K}\subset \widetilde{G}(\A_{F^+,f})$ is 
	a good subgroup that is decomposed with respect to $P$, 
	and with the property that, for each $\bar{v}\in \bar{S}_p$, 
	$\widetilde{K}_{U,\bar{v}} = U^0_{\bar{v}}$. 
	Let $\mathfrak{m}\subset \mathbb{T}^T$ be a non-Eisenstein
	maximal ideal and assume that $\bar{\rho}_{\cS^*(\m^\vee)}$ is decomposed generic
	in the sense of Definition~\ref{defn:generic}. 
	
	Choose a partition 
	\[
	\bar{S}_p = \bar{S}_1\sqcup \bar{S}_2\sqcup \bar{S}_3
	\]
	of the set $\bar{S}_p$ of primes of $F^+$ lying above $p$, together with standard parabolic subgroups $Q_{\bar{v}}\subset P_{\bar{v}}$ for each $\bar{v}\in \bar{S}_3$. 
	Let $\tilde{\lambda}$ and 
	$\lambda$ 
	be dominant weights for $\widetilde{G}$ 
	and $G$, respectively. Assume that the following conditions are satisfied: 
	
	\begin{enumerate}
		
		\item For each $\tau:F^+\hookrightarrow E$ 
		inducing a place $\bar{v}\in \bar{S}_1$, $\tilde{\lambda}_{\tau} = (\lambda_{\tilde{\tau}}, -w_{0,n}\lambda_{\tilde{\tau}c})$ 
		(note the change compared to Proposition \ref{prop:initial subquotient_extra_Hecke}); 
		\item For each $\tau:F^+\hookrightarrow E$ 
		inducing a place $\bar{v}\in \bar{S}_2$, $\tilde{\lambda}_{\tau} = 0$. 
		\item For each $\bar{v}\in \bar{S}_3$, $\tK_{\bar{v}}=\cQ_{\bar{v}}$. For each
		$\tau:F^+\hookrightarrow E$ inducing such a place $\bar{v}$, we also have
		$\tilde{\lambda}_{\tau} = (\lambda_{\tilde{\tau}}, -w_{0,n}\lambda_{\tilde{\tau}c})$
	\end{enumerate}
	
	Then the unnormalized Satake transform 
	$\cS^{\iota}: \widetilde{\mathbb{T}}^{\cQ_{\bar{S}_3},\bar{S}_3-\ord,\tilde{\iota}} \to \mathbb{T}^{\cQ_{\bar{S}_3},\bar{S}_3-\ord,\iota}$ descends to a 
	homomorphism 
	\[
	\widetilde{\mathbb{T}}^{\cQ_{\bar{S}_3},\bar{S}_3-\ord,\tilde{\iota}} \left(H^d\left(\tX_{\tK}, 
	\cV_{\tilde{\lambda}}^\vee\right)_{\cS^*(\m^\vee)}^{\ord^\vee}\right) \to
	\mathbb{T}^{\cQ_{\bar{S}_3},\bar{S}_3-\ord,\iota}\left(\mathbb{H}^d\left( X_{K}, \cV_{\lambda_{\bar{S}_1}}^\vee\otimes \cV_{U}(\tilde{\lambda}_{\bar{S}_2})\otimes 
	\cV_{\lambda_{\bar{S}_3}}^\vee \right)_{\m^\vee}\right),
	\]
	where $\mathbb{H}^d$ denotes the degree $d$ hypercohomology. 
	Moreover, we have an injection 
	\[
	\widetilde{\mathbb{T}}^{\cQ_{\bar{S}_3},\bar{S}_3-\ord,\tilde{\iota}}\left(H^d\left(\tX_{\tK}, 
	\cV_{\tilde{\lambda}}^\vee\right)_{\cS^*(\m^\vee)}^{\ord^\vee}\right) 
	\hookrightarrow 
	\widetilde{\mathbb{T}}^{\cQ_{\bar{S}_3},\bar{S}_3-\ord,\tilde{\iota}}\left(H^d\left(\tX_{\tK}, 
	\cV_{\tilde{\lambda}}^\vee[1/p]\right)_{\cS^*(\m^\vee)}^{\ord^\vee}\right)\] and the target is isomorphic, via Poincar\'{e} duality, to $\widetilde{\mathbb{T}}^{\cQ_{\bar{S}_3},\bar{S}_3-\ord}\left(H^d_c\left(\tX_{\tK}, 
	\cV_{\tilde{\lambda}}[1/p]\right)_{\tilde{\iota}^*\cS^*(\m^\vee)}^{\ord}\right)$.
\end{prop}
\begin{proof}
	With the ingredients mentioned in the preamble, this is a straightforward modification of the proof of Proposition \ref{prop:initial subquotient_extra_Hecke}.
\end{proof}
Note that we have $\rhobar_{\tilde{\iota}^*\cS^*(\m^\vee)} = \rhobar_{\m}(-n)\oplus \rhobar_{\m}^{\vee,c}(1-n)$.

Assume we are given a non-Eisenstein maximal ideal $\m\subset \mathbb{T}$
with $\widetilde{\m}:=\cS^*(\m)$, and a subset $\bar{S}\subseteq \bar{S}_p$. 
We will use the following notation: 
\[
A(K,\lambda,q):= \mathbb{T}^{\overline{\cQ}_{\bar{S}}^{w^P_0},\bar{S}-\ord}_{w^P_0}\left(H^q(X_{K}, \cV_{\lambda})_{\m}\right),
\]
\[ 
A(K,\lambda,q,m):= \mathbb{T}^{\overline{\cQ}_{\bar{S}}^{w^P_0},\bar{S}-\ord}_{w^P_0}\left(H^q(X_{K}, \cV_{\lambda}/\varpi^m)_{\m}\right),\ \mathrm{and}
\]
\[
\tA(\tK,\tilde{\lambda}, \bar{S}):= \widetilde{\mathbb{T}}^{\overline{\cQ}_{\bar{S}}^{w^P_0},\bar{S}-\ord}
\left(H^d(\tX_{\tK}, \cV_{\tilde{\lambda}})^{\ord}_{\widetilde{\m}}\right).
\]
Given a neat compact open subgroup $K\subset \GL_n(\A_{F,f})$ and an integer 
$m\in \Z_{\geq 1}$, define the subgroup $K(m,\bar{S})\subset K$ by setting 
\[
K(m,\bar{S})_v: = K_v \cap \left\{ (1_n)\text{ mod }\varpi_v^m\right\}\subset \GL_n(\cO_{F_v})
\]
if $v$ is a $p$-adic place of $F$ which lies above a place in $\bar{S}$, and 
$K(m, \bar{S})_{v} := K_{v}$ otherwise. 
Also, given a good subgroup $\tK\subset \tG(\A_{F^+,f})$ and an integer $m\in \Z_{\geq 1}$, 
define the good subgroup $\tK(m,\bar{S})\subset \tK$ by setting 
\[
\tK(m, \bar{S})_{\bar{v}}: = \tK_{\bar{v}} \cap \left\{ \left( \begin{array}{cc} 1_n & * \\ 0 & 1_n 
\end{array}\right) \text{ mod }\varpi_{\tilde{v}}^{m}  \right\} \subset \GL_{2n}(\cO_{F_{\tilde{v}}}) = \widetilde{G}(\cO_{F^+_{\bar{v}}})
\]
if $\bar{v}$ is a $p$-adic place of $F^+$ contained in $\bar{S}$, and $\tK(m, \bar{S})_{\bar{v}}: = \tK_{\bar{v}}$
otherwise.

\begin{lemma}\label{lem:degree inequality}
Assume that the subset $\bar{S}\subset \bar{S}_p$ has the following property: 
$\sum_{\bar{v}\not\in \bar{S}} [F^+_{\bar{v}}:\Q_p] \geq \frac{1}{2}[F^+:\Q]$. 
Let $q\in \left[\left \lfloor{\frac{d}{2}}\right \rfloor, 
d-1\right]$. Then
$d-q \leq \sum_{\bar{v}\not\in \bar{S}} n^2[F^+_{\bar{v}}:\Q_p]$.
\end{lemma}

\begin{proof} If $d$ is odd, then so is $[F^+:\Q]$, and we have
\[
d-q \leq \frac{d+1}{2} \leq \frac{d+n^2}{2} \leq
\sum_{\bar{v}\not\in \bar{S}} n^2 [F^+_{\bar{v}}:\Q_p]. 
\]
If $d$ is even, then
$d-q \leq \frac{d}{2} \leq \sum_{\bar{v}\not\in \bar{S}} n^2 [F^+_{\bar{v}}:\Q_p]$. 
\end{proof}

\begin{prop}\label{prop:degree shifting} Let $\bar{v}$, $\bar{v}'$ be two distinct 
	places of $\bar{S}_p$. Let 
	$\bar{S}_1:=\{\bar{v}'\}$, $\bar{S}_3:=\{\bar{v}\}$ and $\bar{S}_2$ be their complement
	in $\bar{S}_p$. Let $\lambda\in \left(\Z^n_{+}\right)^{\Hom(F,E)}$
	be a highest weight for $G$. Let $m\in \Z_{\geq 1}$ be an integer and 
	$\tK \subset \tG(\A_{F^+,f})$ be a good subgroup, such that $\widetilde{K}^p$ is decomposed with respect to $P$. Assume 
	that the following conditions are satisfied. 
	
	\begin{enumerate}
		\item We have 
		\[
		\sum_{\substack{\bar{v}''\in \bar{S}_p \\ \bar{v}''\not = \bar{v},\bar{v}'}} [F^+_{\bar{v}''}:\Q_p] \geq \frac{1}{2}[F^+:\Q]. 
		\]
		\item For each $p$-adic place $\bar{v}''$ of $F^+$ not equal to $\bar{v}$ (including $\bar{v}'' = \bar{v}'$), 
		we have $U(\cO_{F^+_{\bar{v}''}}) \subset \tK_{\bar{v}''}$ and $\tK_{\bar{v}''} = \tK(m, \bar{S}_1\cup 
		\bar{S}_2)_{\bar{v}''}$; in other words $\tK_{\bar{v}''} \subset \left\{ \left( \begin{array}{cc} 1_n & * \\ 0 & 1_n 
		\end{array}\right) \text{ mod }\varpi_{\tilde{v}''}^{m}  \right\}$ for each of these places.
		Finally, we  
		have $\tK_{\bar{v}} = \overline{\cQ}_{\bar{v}}^{w^P_0}$ corresponding to the standard parabolic $\overline{Q}_{\bar{v}}^{w^P_0}\subset P_{F^+_{\bar{v}}}$ with Levi subgroup $Q_{\bar{v}}^{w^P_0}\cap G(F^+_{\bar{v}})$.   
		\item For each embedding $\tau:F\hookrightarrow E$ inducing the place $\bar{v}$ or $\bar{v}'$
		of $F^+$, we have $-\lambda_{\tau c, 1} - \lambda_{\tau, 1}\geq 0$. 
		\item $\m \subset \mathbb{T}$ is a non-Eisenstein maximal ideal such that $\bar{\rho}_{\widetilde{\m}}$ 
		is decomposed generic. 
	\end{enumerate}
	
	\noindent Define a weight $\tilde{\lambda} \in (\Z^{2n}_{+})^{\Hom(F^+, E)}$ as follows: if 
	$\tau: F^+\hookrightarrow E$ does not induce either $\bar{v}$ or $\bar{v}'$, set $\tilde{\lambda}_{\tau} = 0$. 
	If $\tau$ induces $\bar{v}$ or $\bar{v}'$, set $\tilde{\lambda}_{\tau} = (-\lambda_{\tilde{\tau}c}, \lambda_{\tilde{\tau}})$ 
	(identification as in~\eqref{eq:identification of weights})\footnote{The condition (3) in the statement
		of the Proposition guarantees that $\tilde{\lambda}$ will be dominant.}. 
	Set $K:=(\tK^{\bar{v}}\cap G(\A_{F^+, f}^{\bar{v}}))\cdot (\cQ_{\bar{v}}\cap G(F^+_{\bar{v}}))$, identified in the usual way with a neat subgroup of $\GL_n(\A_{F,f})$. 
	
	Let $q \in \left[\left \lfloor{\frac{d}{2}}\right \rfloor, 
	d-1\right]$. Then there exist an integer $m'\geq m$, an integer $N\geq 1$, 
	a nilpotent ideal $J\subset A(K, \lambda, q,m)$ satisfying $J^N = 0$, and a 
	commutative diagram 
	\[
	\xymatrix{\widetilde{\mathbb{T}}^{\overline{\cQ}^{w^P_0}_{\bar{v}},\{\bar{v}\}-\ord}\ar[d]^{\cS^{w^P_0}}\ar[r] & \tA\left(\tK(m',\bar{S}_1\cup \bar{S}_2), \tilde{\lambda}, \bar{v}\right)\ar[d]\\
		\mathbb{T}_{w^P_0}^{\overline{\cQ}^{w^P_0}_{\bar{v}},\{\bar{v}\}-\ord}\ar[r] & A(K,\lambda, q,m)/J.}
	\]
	Moreover, the integer $N$ can be chosen to only depend on $n$ and on $[F^+:\Q]$.  %
\end{prop}

\begin{proof} To simplify notation, we set $\TT:= 	\mathbb{T}_{w^P_0}^{\overline{\cQ}^{w^P_0}_{\bar{v}},\{\bar{v}\}-\ord}$, $\widetilde{\TT}:= \widetilde{\mathbb{T}}^{\overline{\cQ}^{w^P_0}_{\bar{v}},\{\bar{v}\}-\ord}$. We would like to show that there exist non-negative integers $m'\geq m$ 
	and $N$ such that 
	\[
	\cS^{w^P_0}\left(\mathrm{Ann}_{\widetilde{\mathbb{T}}}\ H^d(\tX_{\tK(m',\bar{S}_1\cup \bar{S}_2)}, 
	\cV_{\tilde{\lambda}})^{\ord}_{\widetilde{\m}}\right)^N \subseteq 
	\mathrm{Ann}_{\mathbb{T}}\ H^q(X_{K}, \cV_{\lambda}/\varpi^m)_{\m}. 
	\]
	This is similar to~\cite[Prop. 4.4.1]{10author};
	we will apply Proposition~\ref{prop:initial subquotient_extra_Hecke} repeatedly, which plays 
	the same role as Proposition 4.3.4 in \emph{loc. cit.}. 
	The argument is subtle, for two reasons. 
	\begin{itemize}
		\item We need to work with $\cO$-coefficients
		in order to access the Hecke algebras $\tA\left(\tK(m',\bar{S}_1\cup \bar{S}_2), \tilde{\lambda}, \bar{v}\right)$, 
		whilst the Hecke algebra $A(K,\lambda, q,m)$ acts on cohomology with torsion coefficients. 
		
		\item The spectral sequence computing the hypercohomology in Proposition~\ref{prop:initial subquotient_extra_Hecke} 
		is not known to degenerate. 
	\end{itemize}
	\noindent The second issue does not occur in~\cite{10author}. 
	To deal with both these issues, we will argue by descending induction
	on the degree $q$. The induction hypothesis is the following. 
	
	\begin{hyp} Let $q \in \left[\left \lfloor{\frac{d}{2}}\right \rfloor, d-1\right]$. 
		Then the Proposition holds for every cohomological degree $i\in [q+1, d-1]$ 
		and for every $m\in \Z_{\geq 1}$. Moreover, the integer $N$ can be chosen
		to depend only on $n$, $[F^+:\Q]$ and $q$. 
	\end{hyp}
	
	\noindent Note that the induction hypothesis is satisfied automatically for $q=d-1$. 
	Assume that the hypothesis is satisfied for some $q \in \left[\left \lfloor{\frac{d}{2}}\right \rfloor, d-1\right]$. 
	We will prove the Proposition for $q$, which will imply the induction hypothesis for $q-1$.
	
	Fix $m$, $\tK$ and $\lambda$ as in the statement. Let $M=M(m)\geq m$ be the integer 
	guaranteed by Lemma~\ref{lem:splitting cohomology}. We first increase the level, 
	going from $X_{K}$ to $X_{K(M, \bar{S}_1\cup \bar{S}_2)}$. This uses the same argument as in the 
	proof of~\cite[Prop. 4.4.1]{10author}, which we briefly recall here. Firstly, Poincar\'e 
	duality gives an equality
	\[
	\mathrm{Ann}_{\mathbb{T}}\ H^q(X_{K}, \cV_{\lambda}/\varpi^m)_{\m} = 
	\iota\left(\mathrm{Ann}_{\mathbb{T}^{\iota}}\ H^{d-1-q}(X_{K}, \cV_{\lambda}^\vee/\varpi^m)_{\m^\vee}\right),
	\]
	where $\m^\vee = \iota(\m) \subset \mathbb{T}^T$ and $\TT^{\iota}:=\mathbb{T}^{\cQ_{\bar{S}_3},\bar{S}_3-\ord,\iota}_{w^P_0}$.
	
	The Hochschild--Serre spectral 
	sequence gives an inclusion 
	\[
	\prod_{i=0}^{d-q-1} \mathrm{Ann}_{\mathbb{T}^\iota}\ H^{i}(X_{K(M, \bar{S}_1\cup\bar{S}_2)}, \cV_{\lambda}^\vee/\varpi^m)_{\m^\vee}
	\subset \mathrm{Ann}_{\mathbb{T}^\iota}\ H^{d-1-q}(X_{K}, \cV_{\lambda}^\vee/\varpi^m)_{\m^\vee},
	\]
	and Poincar\'e duality gives 
	\[
	\prod_{i=q}^{d-1} \mathrm{Ann}_{\mathbb{T}}\ H^{i}(X_{K(M, \bar{S}_1\cup\bar{S}_2)}, \cV_{\lambda}/\varpi^m)_{\m}
	\subset \mathrm{Ann}_{\mathbb{T}}\ H^{q}(X_{K}, \cV_{\lambda}/\varpi^m)_{\m}. 
	\] 
	We deal with the terms for $i\geq q+1$ using induction. 
	(See the last paragraph of the proof for more details on 
	how one applies the induction hypothesis.) Therefore, we are left 
	with the term for $i=q$ and we may assume that $K=K(M, \bar{S}_1\cup \bar{S}_2)$. 
	
	Now, note that the $\mathbb{T}$-algebra
	$A(K,\lambda,q,m)$ does not depend on $\lambda_{\bar{v}''}$ for $\bar{v}''\not=\bar{v}$,
	because the level $K_{\bar{v}''}$ is deep enough that the action on $\cV_{\lambda}/\varpi^m$ 
	is trivial. Therefore, we replace $\cV_{\lambda}/\varpi^m$ by 
	\[
	\cV_{\lambda_{\bar{v}}}\otimes \cV_{\lambda_{\bar{v}'}}\otimes \cV^{d-q}_{U}(\tilde{\lambda}_{\bar{S}_2},m). 
	\]
	This is non-zero by Lemmas~\ref{lem:splitting cohomology} and~\ref{lem:degree inequality}. 
	More generally, for any non-negative integer 
	$j\leq \sum_{\bar{v}''\in \bar{S}_2} n^2[F^+_{\bar{v}''}:\Q_p]$, 
	set 
	$\cV^{j}:= \cV_{\lambda_{\bar{v}}} \otimes \cV_{\lambda_{\bar{v}'}} 
	\otimes \cV^j_{U}(\tilde{\lambda}_{\bar{S}_2})$.
	We have a short exact sequence of $\mathbb{T}$-modules 
	\[
	0\to H^q(X_{K}, \cV^{d-q})_{\m} /\varpi^m \to H^q(X_{K}, \cV^{d-q}/\varpi^m)_{\m} 
	\to H^{q+1}(X_{K}, \cV^{d-q})_{\m}[\varpi^m] \to 0,
	\]
	where we are interested in understanding the Hecke algebra 
	$A(K,\lambda, q, m)$ acting on the term in the middle. We can understand the 
	Hecke algebra acting on the $\varpi^m$-torsion in $H^{q+1}(X_{K}, \cV^{d-q})_{\m}$
	using the induction hypothesis: the argument is identical to the argument used 
	in the proof of~\cite[Prop. 4.4.1]{10author}. (More precisely, the argument in \emph{loc. cit.} produces some suitable $m'>m$ and increases the level to $K(M(m'), \bar{S}_1\cup \bar{S}_2)$. The level-raising step is done as above, via Poincar\'e duality and the Hochschild--Serre spectral sequence.) Therefore, we are 
	left with understanding the faithful quotient of $\mathbb{T}$ acting on 
	$H^q(X_{K}, \cV^{d-q})_{\m} /\varpi^m$.  
	
	There is a $\mathbb{T}$-equivariant spectral sequence 
	\begin{equation}\label{eq:integral spectral sequence}
	E_{2}^{i,j} (\cO): = H^i( X_{K}, \cV^j )_{\m} \Rightarrow 
	\mathbb{H}^{i+j}( X_{K}, \cV_{\lambda_{\bar{v}}}\otimes 
	\cV_{\lambda_{\bar{v}'}}\otimes\cV_{U}(\tilde{\lambda}_{\bar{S}_2}))_{\m}. 
	\end{equation}
	If we knew that this spectral sequence
	degenerates on the $E_2$ page, we would deduce that 
	$H^q(X_{K}, \cV^{d-q})_{\m}$ 
	is a $\mathbb{T}$-equivariant 
	subquotient of $\mathbb{H}^{d}( X_{K}, \cV_{\lambda_{\bar{S}_1}}\otimes 
	\cV_{\lambda_{\bar{S}_3}}\otimes\cV_{U}(\tilde{\lambda}_{\bar{S}_2}))_{\m}$,
	and we would win by Proposition~\ref{prop:initial subquotient_extra_Hecke}.  
	However, it is not clear, in this generality, whether the spectral 
	sequence~\eqref{eq:integral spectral sequence} degenerates. 
	Instead, we compare it to the following $\mathbb{T}$-equivariant spectral sequence, 
	whose terms are $\cO/\varpi^m$-modules  
	\begin{equation}\label{eq:torsion spectral sequence}
	E_{2}^{i,j} (\cO/\varpi^m): = H^i( X_{K}, \cV^j/\varpi^m)_{\m} \Rightarrow 
	\mathbb{H}^{i+j}( X_{K}, \cV_{\lambda_{\bar{v}}}\otimes 
	\cV_{\lambda_{\bar{v}'}}\otimes\cV_{U}(\tilde{\lambda}_{\bar{S}_2}, m))_{\m}.
	\end{equation} 
	Since ${K}={K}(M, \bar{S}_1\cup\bar{S}_2)$, Lemma~\ref{lem:splitting cohomology} implies
	that all the differentials in~\eqref{eq:torsion spectral sequence} are zero. 
	Let 
	\[
	\phi_r^{i,j}: E_r^{i,j}(\cO)\to E_r^{i,j}(\cO/\varpi^m)
	\] 
	be the natural, $\mathbb{T}$-equivariant map between 
	the spectral sequences~\eqref{eq:integral spectral sequence} and~\eqref{eq:torsion spectral sequence}. 
	Let $F_r^{i,j}:=\mathrm{Im} (\phi_r^{i,j})$. When $r=2$, we have $F_2^{i,j} = E_2^{i,j}(\cO)/\varpi^m$.
	For $r\geq 3$, we at least have a surjection $E_r^{i,j}(\cO)/\varpi^m \twoheadrightarrow F_r^{i,j}$,
	because $F_r^{i,j}$ is an $\cO/\varpi^m$-module.  
	
	With this new notation, we are interested in relating $\mathrm{Ann}_{\mathbb{T}}(F_2^{q,d-q})$ to 
	$\mathrm{Ann}_{\mathbb{T}}(E_{\infty}^{d}(\cO))$. 
	For any $r\geq 2$, let $d_r^{-}: E_r^{q-r,d-q+r-1}\to E_r^{q,d-q}$ and $d_r^+: E_r^{q,d-q}\to E_r^{q+r, d-q+1-r}$
	denote the $r$th differentials. Because all the differentials in~\eqref{eq:torsion spectral sequence} are zero, 
	we have that $\mathrm{Im}(d_r^{-})\subseteq \mathrm{Ker}(\phi_r^{q,d-q})$. Since $\phi^{q,d-q}_{r}$  
	induces $\phi^{q,d-q}_{r+1}$, we deduce that we have an injection 
	\[
	F^{q,d-q}_{r+1} = \mathrm{Ker}(d_r^{+})/\left(\mathrm{Im}(d_r^{-}) + \mathrm{Ker}(d_r^{+})\cap
	\mathrm{Ker}(\phi_r^{q,d-q}) \right) \hookrightarrow F^{q,d-q}_r = E^{q,d-q}_r/\mathrm{Ker}(\phi_r^{q,d-q}). 
	\]
	Moreover, the cokernel of this injection becomes identified, under the map induced by $d^{+}_r$, 
	with $\mathrm{Im}(d^{+}_r)/ d^{+}_r(\mathrm{Ker}(\phi_r^{q,d-q}))$. Since $(\varpi^m)\subseteq \mathrm{Ker}(\phi_r^{q,d-q})$, 
	the latter is a quotient of $\mathrm{Im}(d^{+}_r)/\varpi^m$. By Lemma~\ref{lem:artin-rees application}, 
	there exists some $m'_r\geq m$ such that the latter is a subquotient of $E_2^{q+r, d-q+1-r}(\cO)/\varpi^{m'_r}$,
	or even of $E_2^{q+r, d-q+1-r}(\cO/\varpi^{m'_r})$. 
	We therefore have an inclusion 
	\[
	\mathrm{Ann}_{\mathbb{T}}\ E_{\infty}^d(\cO) \cdot \prod_{r=2}^{d-q-1} 
	\mathrm{Ann}_{\mathbb{T}}\ E_2^{q+r,d+1-q-r}(\cO/\varpi^{m'_r}) \subseteq 
	\mathrm{Ann}_{\mathbb{T}}\ F_2^{q,d-q}. 
	\]
	For the $E_{\infty}^d(\cO)$ term, Proposition~\ref{prop:initial subquotient_extra_Hecke} implies that 
	there is an inclusion 
	\[
	\cS^{w^P_0}\left(\mathrm{Ann}_{\widetilde{\mathbb{T}}}\ H^d(\tX_{\tK}, \cV_{\tilde{\lambda}})^{\ord}_{\widetilde{\m}}\right)
	\subseteq \mathrm{Ann}_{\mathbb{T}}\ E_{\infty}^d(\cO). 
	\]
	For each $E_2^{q+r,d+1-q-r}(\cO/\varpi^{m'_r}) = H^{q+r}(X_{K}, \cV^{d+1-q-r}/\varpi^{m'_r})_{\m}$, 
	we use the argument above via Poincar\'e duality and the Hochschild--Serre
	spectral sequence to increase the level to ${K}(M(m'_r), \bar{S}_1\cup \bar{S}_2)$. 
	We then apply the induction hypothesis. 
	
	Each time we apply the induction hypothesis, we find some integer
	$m'_i\geq m$ and some nilpotence degree $N_i$, which can be bounded in terms of $\dim(X_K)$, for $i$ running over some finite index set $I$ whose size can also be bounded in terms of $\dim(X_K)$. %
	To find a common $m'\geq m$, we let $m':=\mathrm{sup}_{i\in I} M(m'_i)$.
	For each $i$, we have 
	\[
	\mathrm{Ann}_{\widetilde{\mathbb{T}}}\ H^d(\tX_{\tK(m', \bar{S}_1\cup \bar{S}_2)}, \cV_{\tilde{\lambda}})^{\ord}_{\widetilde{\m}} \subseteq
	\mathrm{Ann}_{\widetilde{\mathbb{T}}}\ H^d(\tX_{\tK(m'_i, \bar{S}_1\cup \bar{S}_2)}, \cV_{\tilde{\lambda}})^{\ord}_{\widetilde{\m}}, 
	\]  
	because this is true rationally and the cohomology groups are torsion-free, by Theorem~\ref{thm:middle degree cohomology}. 
	We then let $J$ denote the image of the ideal 
	\[
    \cS^{w^P_0}\left(\mathrm{Ann}_{\widetilde{\mathbb{T}}}\ H^d(\tX_{\tK(m',\bar{S}_1\cup \bar{S}_2)}, 
	\cV_{\tilde{\lambda}})^{\ord}_{\widetilde{\m}}\right)
    \] 
	in $A(K, \lambda, q, m)$.  
	To find an appropriate nilpotence degree $N$, we set $N = 1+\sum_{i} N_i$. 
\end{proof}
Again we have a similar statement with dual coefficients. We introduce some more notation, depending on a decomposition $\bar{S}_p=\bar{S}_1\cup\bar{S}_2\cup\bar{S}_3$: \[
A^\vee(K,\lambda,q):= \mathbb{T}^{\cQ_{\bar{S}_3},\bar{S}_3-\ord,\iota}\left(H^q(X_{K}, \cV_{\lambda}^\vee)_{\m^\vee}\right),
\]
\[ 
A^\vee(K,\lambda,q,m):= \mathbb{T}^{\cQ_{\bar{S}_3},\bar{S}_3-\ord,\iota}\left(H^q(X_{K}, \cV_{\lambda}^\vee/\varpi^m)_{\m^\vee}\right),\ \mathrm{and}
\]
\[
\tA^\vee(\tK,\tilde{\lambda}, \bar{S}_3):= \widetilde{\mathbb{T}}^{\cQ_{\bar{S}_3},\bar{S}_3-\ord,\tilde{\iota}}
\left(H^d(\tX_{\tK}, \cV_{\tilde{\lambda}}^\vee)^{\ord^\vee}_{\cS^*{\m^\vee}}\right).
\]
\begin{prop}\label{prop:degree shifting dual} Let $\bar{v}$, $\bar{v}'$ be two distinct 
	places of $\bar{S}_p$. Let 
	$\bar{S}_1:=\{\bar{v}'\}$, $\bar{S}_3:=\{\bar{v}\}$ and $\bar{S}_2$ be their complement
	in $\bar{S}_p$. Let $\lambda\in \left(\Z^n_{+}\right)^{\Hom(F,E)}$
	be a highest weight for $G$. Let $m\in \Z_{\geq 1}$ be an integer and 
	$\tK \subset \tG(\A_{F^+,f})$ be a good subgroup such that $\tK^p$ is decomposed with respect to $P$. Assume 
	that the following conditions are satisfied. 
	
	\begin{enumerate}
		\item We have 
		\[
		\sum_{\substack{\bar{v}''\in \bar{S}_p \\ \bar{v}''\not = \bar{v},\bar{v}'}} [F^+_{\bar{v}''}:\Q_p] \geq \frac{1}{2}[F^+:\Q]. 
		\]
		\item For each $p$-adic place $\bar{v}''$ of $F^+$ not equal to $\bar{v}$ (including $\bar{v}'' = \bar{v}'$), 
		we have $U(\cO_{F^+_{\bar{v}''}}) \subset \tK_{\bar{v}''}$ and $\tK_{\bar{v}''} = \tK(m, \bar{S}_1\cup 
		\bar{S}_2)_{\bar{v}''}$.
		We  
		have $\tK_{\bar{v}} = {\cQ}_{\bar{v}}$ corresponding to the standard parabolic ${Q}_{\bar{v}}\subset P_{F^+_{\bar{v}}}$.   
		\item For each embedding $\tau:F\hookrightarrow E$ inducing the place $\bar{v}$ or $\bar{v}'$
		of $F^+$, we have $\lambda_{\tau c, n} + \lambda_{\tau, n}\geq 0$. 
		\item $\m \subset \mathbb{T}$ is a non-Eisenstein maximal ideal such that $\bar{\rho}_{\cS^*(\m^\vee)}$ 
		is decomposed generic. 
	\end{enumerate}
	
	\noindent Define a weight $\tilde{\lambda} \in (\Z^{2n}_{+})^{\Hom(F^+, E)}$ as follows: if 
	$\tau: F^+\hookrightarrow E$ does not induce either $\bar{v}$ or $\bar{v}'$, set $\tilde{\lambda}_{\tau} = 0$. 
	If $\tau$ induces $\bar{v}$ or $\bar{v}'$, set $\tilde{\lambda}_{\tau} = (\lambda_{\tilde{\tau}},-w_{0,n}\lambda_{\tilde{\tau}c})$. 
	Set $K:=\tK\cap G(\A_{F^+, f})$, identified in the usual way with a neat subgroup of $\GL_n(\A_{F,f})$. 
	
	Let $q \in \left[\left \lfloor{\frac{d}{2}}\right \rfloor, 
	d-1\right]$. Then there exist an integer $m'\geq m$, an integer $N\geq 1$, 
	a nilpotent ideal $J\subset A^\vee(K, \lambda, q,m)$ satisfying $J^N = 0$, and a 
	commutative diagram 
	\[
	\xymatrix{\widetilde{\mathbb{T}}^{{\cQ}_{\bar{v}},\{\bar{v}\}-\ord,\tilde{\iota}}\ar[d]^{\cS^{\iota}}\ar[r] & \tA^\vee\left(\tK(m',\bar{S}_1\cup \bar{S}_2), \tilde{\lambda}, \bar{v}\right)\ar[d]\\
		\mathbb{T}^{{\cQ}_{\bar{v}},\{\bar{v}\}-\ord,\iota}\ar[r] & A^\vee(K,\lambda, q,m)/J.}
	\]
	Moreover, the integer $N$ can be chosen to only depend on $n$ and on $[F^+:\Q]$.  %
\end{prop}

We will now be able to reduce questions about Galois representations with coefficients in the torsion Hecke algebras $A(K,\lambda,q,m)$ to understanding the properties of the Galois representations with coefficients in the $p$-torsion free Hecke algebras $\widetilde{A}(\cdots)$ and $\widetilde{A}^\vee(\cdots)$. To do this, we need some results about automorphic Galois representations in characteristic 0. 

For the statement of the next proposition, recall that we have introduced Hecke operators $\widetilde{U}_v^k$ for a place $v \in S_p$ in \S\ref{sec:P-ord}. 

\begin{prop}\label{prop:det mod p} 
	Assume that $p$ splits in an imaginary quadratic subfield of $F$. Let $\bar{v}$ be a $p$-adic place of $F^+$,
	let $\m\subset \mathbb{T}_{w^P_0}^{\overline{\cQ}_{\bar{v}}^{w^P_0},\{\bar{v}\}-\ord}$ be a non-Eisenstein maximal ideal in the 
	support of some $H^*(X_K, \cV_{\lambda})$, 
	and set $\widetilde{\m}:=(\cS^{w^P_0})^*(\m)$, a maximal 
	ideal of $\widetilde{\mathbb{T}}^{\overline{\cQ}_{\bar{v}}^{w^P_0},\{\bar{v}\}-\ord}$. Fix $v|\bar{v}$ in $F$ and suppose that $\widetilde{U}_v^k \notin \widetilde{\m}$ for $1 \le k \le t$. 
	
	Assume that $\pi$ is a cuspidal automorphic representation of $\widetilde{G}(\A_{F^+})$, $\iota$ is an isomorphism $\iota:\overline{\Q}_p\toisom \C$ and $\pi$ is $\iota$-$\overline{Q}_{\bar{v}}^{w^P_0}$-ordinary of weight $\tilde{\lambda}$. Suppose moreover that the Hecke eigenvalues on $(\iota^{-1}\pi^{\infty})^{\tK,\overline{\cQ}_{\bar{v}}^{w^P_0}-\ord}$ come from a map \[f: \widetilde{\mathbb{T}}^{\overline{\cQ}_{\bar{v}}^{w^P_0}-\ord}_{\widetilde{\m}} \to \Qpbar,\] where the superscript `$\overline{\cQ}_{\bar{v}}^{w^P_0}-\ord$' denotes that we replace the local factor $\pi_{\bar{v}}^{\overline{\cQ}_{\bar{v}}^{w^P_0}}$ with its one-dimensional $\overline{\cQ}_{\bar{v}}^{w^P_0}$-ordinary subspace.
	
	We have the associated $p$-adic Galois representation 
	$r_{\iota}(\pi): G_F \to \GL_{2n}(\overline{\Q}_p)$ (the existence of $f$ implies that we have an isomorphism of semi-simplified reductions $\rbar_{\iota}(\pi)\cong \rhobar_{\widetilde{\m}}$). Consider the $(n,n)$-block decomposition 
	\[
	r_{\iota}(\pi)|_{G_{F_{\tilde{v}}}} \simeq \begin{pmatrix} r_1(\pi) & * \\ 0 & r_2(\pi) \end{pmatrix}
	\]
	guaranteed by Theorem~\ref{thm:shape general} (noting that $r_1(\pi), r_2(\pi)$ may be further decomposed according to the shape of $\overline{\cQ}_{\bar{v}}^{w^P_0}$).  
	For $i=1,2$, assume that $E$ is large enough that $r_i(\pi)$ 
	can be defined over it, via the embedding $E \hookrightarrow\Qpbar$ coming from $f$, and let $\overline{r_i(\pi)}$ be the semi-simplification of the reduction modulo
	$\varpi$ of $r_i(\pi)$. Then 
	\[
	\det \overline{r_1(\pi)} \left(\Art_{F_{{v}}}(\varpi_{{v}})\right) 
	= \det \bar{\rho}_{\m}(\Art_{F_{{v}}}(\varpi_{{v}}))\ \mathrm{and} 
	\]
	\[
	\det \overline{r_2(\pi)} \left(\Art_{F_{{v}}}(\varpi_{{v}})\right) 
	= \det \left(\bar{\rho}^{\vee, c}_{\m}(1-2n)\right)(\Art_{F_{{v}}}(\varpi_{{v}})). 
	\]
\end{prop}

\begin{proof} By Theorem~\ref{thm:shape general}, 
	$\det r_1(\pi)\left(\Art_{F_{\tilde{v}}}(\varpi_{\tilde{v}})\right)\in \cO^\times$ 
	is equal to $\epsilon_p^{\frac{n(1-n)}{2}}(\mathrm{Art}_{F_v}(\varpi_v))$
	times the eigenvalue of $\widetilde{U}_{{v},n}$ acting on 
	the $\overline{\cQ}_{\bar{v}}^{w^P_0}$-ordinary subspace of  $\iota^{-1}\pi^{\overline{\cQ}_{\bar{v}}^{w^P_0}}$. 
	By the description of the map $\cS^{w^P_0}$, the reduction 
	of this eigenvalue modulo $\varpi$ is equal to 
	the image of $U_{{v}}$ in $\mathbb{T}_{w^P_0}^{\overline{\cQ}_{\bar{v}}^{w^P_0},\{\bar{v}\}-\ord}/\m$. 
	By Lemma~\ref{lem:central character mod p}, 
	$\det \bar{\rho}_{\m}(\Art_{F_{{v}}}(\varpi_{{v}}))$
	is equal to $\bar{\epsilon}_p^{\frac{n(1-n)}{2}}(\mathrm{Art}_{F_v}(\varpi_v))$
	times the eigenvalue of $U_{{v}}$ acting
	on $H^*(X_K, \cV_{\lambda}/\varpi)_{\m}$, 
	so we obtain the first equation.
	For the second equation, let $\overline{r_{\iota}(\pi)}$ be the 
	semi-simplification of the reduction modulo
	$\varpi$ of $r_{\iota}(\pi)$. 
	The same line of reasoning 
	implies that $\det \overline{r_{\iota}(\pi)} \left(\Art_{F_{{v}}}
	(\varpi_{{v}})\right)$ is equal to 
	$\bar{\epsilon}_p^{n(1-2n)}(\mathrm{Art}_{F_v}(\varpi_v))$ times
	the image of $U_{{v}}\cdot U_{{v}^c}^{-1}$ in $\mathbb{T}^{\{\bar{v}\}-\ord}/\m$. 
	We conclude by the first equation and by 
	another application of Lemma~\ref{lem:central character mod p}.
\end{proof}

To proceed, we recall the notion of a CTG weight from~\cite[Def. 4.3.5]{10author}. 

\begin{defn}\label{defn:CTG} A weight $\tilde{\lambda}\in (\Z^{2n}_+)^{\Hom(F^+, E)}$ 
is CTG (``cohomologically trivial for $G$'') if it satisfies the following condition
\begin{itemize}
\item Given $w\in W^P$, define $\lambda_w = w(\tilde{\lambda}+\rho) - \rho$, 
viewed as an element of $(\Z^n_+)^{\Hom(F, E)}$ in the usual way. For each
$w\in W^P$ and $i_0 \in \Z$, there exists $\tau\in \Hom(F,E)$ such that $\lambda_{w,\tau} - 
\lambda_{w,\tau c}^{\vee}\not = (i_0,i_0,\dots, i_0)$. 
\end{itemize}
\end{defn}

An important application of the CTG assumption is the following variant of \cite[Theorem 2.4.11]{10author}:
\begin{prop}\label{prop:CTG automorphic coh}
Assume that p splits in an imaginary quadratic subfield of F. Let $\mathfrak{m}\subset \mathbb{T}^T$ be a non-Eisenstein
maximal ideal. Fix a place $\bar{v} \in \bar{S}_p$ and a standard parabolic $Q_{\bar{v}}\subset P_{\bar{v}}$ and suppose $\widetilde{\mathfrak{m}}$ is a maximal ideal of $\widetilde{\TT}^{\cQ_{\bar{v}},\{\bar{v}\}-\ord}$ which extends $\cS^*(\mathfrak{m})$. Let $\tK \subset \tG(\AA_{F^+,f})$ be a good subgroup such that $\widetilde{\mathfrak{m}}$ is in the support of $H^*(\tX_{\tK},\cV_{\tilde{\lambda}})^{\ord}$ for a CTG weight $\tilde{\lambda}$. Suppose that $\widetilde{U}_v^k \notin \widetilde{\m}$ for $1 \le k \le t$ (in other words, these Hecke operators act with unit eigenvalues on $H^*(\tX_{\tK},\cV_{\tilde{\lambda}})^{\ord}_{\widetilde{\mathfrak{m}}}$). Let $d = \frac{1}{2}\dim_{\R}X^{\tG} = n^2[F^+:\Q]$. 

Then $H^d(\tX_{\tK},\cV_{\tilde{\lambda}})^{\ord}_{\widetilde{\mathfrak{m}}}[1/p]$ is a semisimple $\widetilde{\TT}^{\cQ_{\bar{v}},\{\bar{v}\}-\ord}[1/p]$-module, and for every homomorphism \[f: \widetilde{\TT}^{\cQ_{\bar{v}},\{\bar{v}\}-\ord}(H^d(\tX_{\tK},\cV_{\tilde{\lambda}})^{\ord}_{\widetilde{\mathfrak{m}}})\to \Qpbar,\] and every isomorphism $\iota:\Qpbar\to\C$ there exists a cuspidal automorphic representation $\pi$ of $\tG(\AA_{F^+})$ which is $\iota$-$Q_{\bar{v}}$-ordinary of weight $\tilde{\lambda}$ such that $f$ is associated to the Hecke eigenvalues on $(\iota^{-1}\pi^{\infty})^{\tK,\cQ_{\bar{v}}-\ord}$, where `$\cQ_{\bar{v}}-\ord$' indicates that we replace the local factor $\pi_{\bar{v}}^{\cQ_{\bar{v}}}$ with its one-dimensional $\cQ_{\bar{v}}$-ordinary subspace.
\end{prop}
\begin{proof}
	This follows from combining (the proof of) \cite[Theorem 2.4.11]{10author} with the fact that the $\cQ_{\bar{v}}$-ordinary subspace of $\pi_{\bar{v}}^{\cQ_{\bar{v}}}$ (which is all that contributes to cohomology localized at $\widetilde{\mathfrak{m}}$) is one-dimensional, which is part of Theorem \ref{thm:shape general}. 
\end{proof}

\begin{remark}\label{rem:CTG}\leavevmode
\begin{enumerate}
\item If the weight $\tilde{\lambda}$ is CTG, and if $\widetilde{\m}=\mathcal{S}^*(\m)$ for a non-Eisenstein $\m$, then we have a canonical, Hecke-equivariant isomorphism 
\[
 H^d_c(\widetilde{X}_{\widetilde{K}}, \cV_{\tilde{\lambda}}[1/p])_{\widetilde{\m}}\stackrel{\sim}{\to} H^d(\widetilde{X}_{\widetilde{K}}, \cV_{\tilde{\lambda}}[1/p])_{\widetilde{\m}},
\]
by~\cite[Theorem 2.4.11]{10author}. This implies that, in the CTG setting, the term $\widetilde{\mathbb{T}}^{\cQ_{\bar{S}_3},\bar{S}_3-\ord}\left(H^d_c\left(\tX_{\tK}, 
	\cV_{\tilde{\lambda}}[1/p]\right)_{\tilde{\iota}^*\cS^*(\m^\vee)}^{\ord}\right)$ of Proposition~\ref{prop:initial subquotient_dual} can be replaced with $\widetilde{\mathbb{T}}^{\cQ_{\bar{S}_3},\bar{S}_3-\ord}\left(H^d\left(\tX_{\tK}, 
	\cV_{\tilde{\lambda}}[1/p]\right)_{\tilde{\iota}^*\cS^*(\m^\vee)}^{\ord}\right)$. 
    
\item In Proposition~\ref{prop:degree shifting}, we may assume
that the weight $\tilde{\lambda}$ is CTG, without changing the Hecke algebra 
$A(K,\lambda, q, m)$. This is because~\cite[Lemma 4.3.6]{10author}
shows that a weight can be ensured to be CTG by modifying it at only one embedding
$\tau: F^+\hookrightarrow E$. Choose a $\tau$ which induces the place $\bar{v}'$ 
of $F^+$. Because the level at $\bar{v}'$ is assumed to be deep enough in Proposition~\ref{prop:degree shifting},
we may modify $\tilde{\lambda}_{\tau} = (-\lambda_{\tilde{\tau}c}, \lambda_{\tilde{\tau}})$ 
without changing the Hecke algebra $A(K,\lambda, q,m)$.  
\end{enumerate}
\end{remark}

We can now apply the results of \S\ref{sec:determinants} to obtain the main result of this section.

\begin{prop}\label{prop:lgc torsion CTG}

Assume that $p$ splits in an imaginary quadratic subfield of $F$. Let $K \subset G(\AA_{F^+,f}) = \GL_n(\AA_{F,f})$ be a good subgroup and fix distinct places $\bar{v},\bar{v'}\in S_p$. Let $\lambda$ be a dominant weight for $G$.
	
	Let $m \in \Z_{\ge 1}$ be an integer. Fix a standard parabolic $Q_{\bar{v}} \subset P_{\bar{v}}$, suppose that $K_{\bar{v}}  = \cQ_{\bar{v}} \cap G(F^+_{\bar{v}})$ and  let $\frakm \subset \TT^{\cQ_{\bar{v}},\{\bar{v}\}-\ord}$ be a maximal ideal in the support of $H^*(X_K,\cV_\lambda/\varpi^m)$. 
	
	Fix $\tilde{v}|\bar{v}$ in $F$. Using $\iota_{\tilde{v}}$, we identify $Q_{\bar{v}}$ with a standard block-upper-triangular parabolic subgroup of $\GL_{2n}$ corresponding to a decomposition $(n_1,\ldots,n_t)$ of $2n$. Suppose that $n = n_1 + \cdots+ n_r$.
	
		Assume that:
	\begin{enumerate} 
		\item We have \label{item:condition on size of local fields}
\[
\sum_{\substack{\bar{v}''\in \bar{S}_p \\ \bar{v}''\not = \bar{v},\bar{v}'}} [F^+_{\bar{v}''}:\Q_p] \geq \frac{1}{2}[F^+:\Q]. 
\]
\item $\m$ is a non-Eisenstein maximal ideal such that $\bar{\rho}_{\m}$ 
is decomposed generic. 
\item\label{assm:basechange} Let $T\supseteq S_p$ be a finite set of finite places of $F$ that
satisfies $T = T^c$ and that satisfies the following additional condition. Let $v \notin T$ be a finite place, with residue characteristic $l$. Then either 
$T$ contains no $l$-adic places and $l$ is unramified in $F$, or there is an 
imaginary quadratic subfield of $F$ in which $l$ splits. Assume that, whenever $v\not\in T$, $K_v = \mathrm{GL}_n(\cO_{F_{v}})$. 
\item For all $\nu \in X_{Q_{\bar{v}}}$, the Hecke operator $[K_{\bar{v}}\nu(\varpi_{\bar{v}})K_{\bar{v}}]$ is not contained in $\m$.
	\end{enumerate}
	
	Then for each $q \in [0,d-1]$ there exist an integer $N \ge 1$, depending only on $n$ and $[F^+:\Q]$, a nilpotent ideal $J$ of $\TT^{\cQ_{\bar{v}},\{\bar{v}\}-\ord}(H^q(X_K,\cV_\lambda/\varpi^m)_{\m})$ with $J^N = 0$ and a continuous $n$-dimensional representation
	\[\rho_{\m}: G_{F,T} \to \GL_n(\TT^{\cQ_{\bar{v}},\{\bar{v}\}-\ord}(H^q(X_K,\cV_\lambda/\varpi^m)_{\m})/J) \] such that the following conditions are satisfied:
	\begin{enumerate}
		\item For each place $v \notin T$ of $F$, the characteristic polynomial of $\rho_{\m}(\Frob_v)$ is equal to the image of $P_v(X)$. 
		\item For $v|\bar{v}$, the representation $\rho_{\m}|_{G_{F_{{v}}}}$ has a lift to $\tilde{\rho}_{{v}}: G_{F_{{v}}} \to \GL_n(\widetilde{A})$, where $\widetilde{A}$ is a finite flat local $\cO$-algebra equipped with a morphism \[f:\widetilde{A} \to A(K,\lambda, q, m)/J:=\TT^{\cQ_{\bar{v}},\{\bar{v}\}-\ord}(H^q(X_K,\cV_\lambda/\varpi^m)_{\m})/J.\]
		\item Inverting $p$, the lift $\tilde{\rho}_{{v}}[1/p]$ is semistable with labelled Hodge--Tate weights $(\lambda_{\tau,n} < \cdots < \lambda_{\tau,1}+n-1)_{\tau:F_{{v}}\hookrightarrow E}$. 
		\item Furthermore, these semistable lifts satisfy \[\tilde{\rho}_{\tilde{v}}[1/p] \simeq \begin{pmatrix} \tilde{\rho}_{\tilde{v},r+1} & * & \cdots & *\\ 0 & \tilde{\rho}_{\tilde{v},r+2} & \cdots & *\\ 0 & 0 & \ddots & * \\ 0 & 0 & \cdots & \tilde{\rho}_{\tilde{v},t} \end{pmatrix}\text{ and } \tilde{\rho}_{\tilde{v}^c}[1/p] \simeq \begin{pmatrix} \tilde{\rho}_{\tilde{v}^c,r} & * & \cdots & *\\ 0 & \tilde{\rho}_{\tilde{v}^c,r-1} & \cdots & *\\ 0 & 0 & \ddots & * \\ 0 & 0 & \cdots & \tilde{\rho}_{\tilde{v}^c,1}\end{pmatrix}\] where the representations $\tilde{\rho}_{v,j}:G_{F_v} \to \GL_{n_j}(\widetilde{A}[1/p])$ are crystalline with labelled Hodge--Tate weights determined by the requirement that they are increasing from top left to bottom right.
		\item For $j = r+1,\ldots, t$, the characters $\det\tilde{\rho}_{\tilde{v},j}$ take values in $\widetilde{A}$ and their image under $f$ is given by characters $\psi_j$ determined by:
			\begin{itemize}
			\item $\prod_{j=r+1}^{k} \psi_j(\Art_{F_{\tilde{v}}}(u)) = \prod_{i=1}^{n_{r+1}+\cdots + n_k}\prod_{\tau: F_{\tilde{v}}\hookrightarrow \Qpbar}\tau(u)^{-{\lambda}_{\tau,n-i+1}-i+1}$ for $u \in \cO_{F_{\tilde{v}}}^\times$.
			\item $\prod_{j=r+1}^{k} \psi_j(\Art_{F_{\tilde{v}}}(\varpi_{\tilde{v}}))$ is equal to $\epsilon_p^{\sum_{i=1}^{n_{r+1}+\cdots + n_k}(1-i)}(\Art_{F_{\tilde{v}}}(\varpi_{\tilde{v}}))\cS^{w_0^P}(\widetilde{U}_v^{k-r})$, where, by an abuse of notation, we identify the abstract Hecke operator $\cS^{w_0^P}(\widetilde{U}_v^{k-r})$ with its image in $A(K,\lambda,q,m)/J$. 
		\end{itemize}
	\end{enumerate}
\end{prop}
\begin{proof}
We already know the existence of $\rho_{\m}$ satisfying the first condition (local--global compatibility at unramified places), so we are free to enlarge $T$. As explained in the proof of \cite[Corollary 4.4.8]{10author}, we may assume (applying a twisting argument) that  $\rhobar_{\widetilde{\m}}=\rhobar_{\m}\oplus\rhobar_{\m}^{\vee,c}(1-2n)$ is decomposed generic, not just that $\bar{\rho}_{\m}$ is decomposed generic. We will use a similar twisting argument later in this proof. We can also use Hochschild--Serre to reduce to the case when $K_{v''} \subset \left\{ g: g = 1_n \text{ mod }\varpi_{{v}''}^{m}  \right\}$ for each $v''|\bar{v}'' \in \overline{S}_p - \{\bar{v}\}$. 
This means we can moreover assume that $\lambda_{\bar{v}''} = 0$ if $\bar{v}'' \in \bar{S}_p - \{\bar{v}\}$. 

Now we let $\widetilde{K}\subset \widetilde{G}(\AA_{F^+,f})$ be a good subgroup satisfying:\begin{itemize}
	\item $\widetilde{K}^{\bar{v}}\cap G(\AA^{\bar{v}}_{F^+,f}) = K^{\bar{v}}$.
	\item $\widetilde{K}^T = \widetilde{G}(\widehat{\cO}^T_{F^+})$.
	\item For each $\bar{v}'' \in \overline{S}_p - \{\bar{v}\}$, 
	we have $U(\cO_{F^+_{\bar{v}''}}) \subset \tK_{\bar{v}''}$ and  \[\tK_{\bar{v}''} \subset \left\{ \left( \begin{array}{cc} 1_n & * \\ 0 & 1_n 
	\end{array}\right) \text{ mod }\varpi_{\tilde{v}''}^{m}  \right\}.\]
	\item $\widetilde{K}_{\bar{v}} = \overline{\cQ}_{\bar{v}}^{w^P_0}$, so the corresponding standard parabolic has block sizes \[(n_{r+1},\ldots,n_t,n_1,\ldots,n_r).\]
\end{itemize}

Next, we use a twisting argument to reduce to the case when $-\lambda_{\tilde{\tau}c,1}-\lambda_{\tau,1} \ge 0$. Indeed, twisting by $\epsilon_p^\mu$ moves us from the weight $\lambda$ to the weight $\lambda':=(\lambda_{\tau,1}-\mu,\ldots,\lambda_{\tau,n}-\mu)_{\tau: F \hookrightarrow E}$ (cf.~\cite[Proposition 2.2.23]{10author}), and we satisfy the desired condition if we take $\mu$ to be sufficiently positive. 

At this point we assume that $q \ge \lfloor \frac{d}{2} \rfloor$. We will handle small $q$ at the end of the proof using Poincar\'{e} duality. We are now in a position to apply Proposition \ref{prop:degree shifting}. Following Remark \ref{rem:CTG}, we are free to modify $\lambda_{\bar{v}'}$ so that the weight $\tilde{\lambda}$ is CTG. We have $A(K,\lambda,q,m) = \TT^{\cQ_{\bar{v}},\{\bar{v}\}-\ord}(H^q(X_K,\cV_\lambda/\varpi^m)_{\m})$. 

Suppose we have a continuous character $\overline{\psi}: G_F \to k^\times$ (perhaps after extending $\cO$), which is unramified at $S_p$, and let $\psi:G_F\to\cO^\times$ denote the Teichm\"{u}ller lift of $\overline{\psi}$. Choose a finite set ${T}'\supset {T}$ (closed under complex conjugation) which also contains all the places where $\psi$ is ramified, and a good normal subgroup $K' \subset K$ satisfying:
\begin{itemize}
	\item $(K')^{T'-T} = K^{T'-T}$.
	\item $K/K'$ is abelian of order prime to $p$.
	\item For each place $v$ of $F$, the restriction of $\psi|_{G_{F_v}}\circ\Art_{F_v}$ to $\det(K'_v)$ is trivial.
	\item $T'$ satisfies assumption (\ref{assm:basechange}) from the Proposition.
\end{itemize}
We will then consider the Hecke algebras for the twist \[A(K',\lambda,q,m,\psi) := \TT^{\cQ_{\bar{v}},\{\bar{v}\}-\ord}(H^q(X_{K'},\cV_\lambda/\varpi^m)_{\m(\psi)}),\] (see \cite[\S2.2.20]{10author} for the definition of $\m(\psi)$ and note that $\rhobar_{\m(\psi)} = \rhobar_{\m}\otimes \overline{\psi}$). Establishing the proposition for any of these twists will imply it for the Hecke algebra $A(K,\lambda,q,m)$. We will always assume that $\overline{\psi}$ is chosen so that $\rhobar_{\widetilde{\m(\psi)}}$ remains decomposed generic.

We can, twisting by a suitable $\psi$ if necessary, assume that the isomorphism classes of the irreducible constituents of $\rhobar_\m|_{G_{F_{\tilde{v}}}}$ are disjoint from those of 
$\rhobar_{\m}^{\vee,c}(1-2n)|_{G_{F_{\tilde{v}}}}$. 

Applying Proposition \ref{prop:degree shifting}, for each $\psi$ as above we have a finite flat $\cO$-algebra $\widetilde{A}(\psi)$ and a nilpotent ideal $J_{\psi}$ with a map $f_{\psi}:\widetilde{A}(\psi) \to A(K',\lambda,q,m,\psi)/J_\psi$. Using Propositions \ref{prop:det mod p} and \ref{prop:CTG automorphic coh}, we deduce that there is a $\overline{Q}^{w^P_0}_{\bar{v}}$-ordinary Galois representation $\tilde{\rho}_{\m(\psi)} = \prod_{i=1}^r \tilde{\rho}_{\m(\psi)}^i: G_F \to \GL_{2n}(\widetilde{A}(\psi)[1/p]) = \prod_{i=1}^{r}E$ such that, for each $i$, the factor $\tilde{\rho}_{\m(\psi)}^i$ comes with a $(n,n)$-block decomposition \[\tilde{\rho}_{\m(\psi)}^i |_{G_{F_{\tilde{v}}}} \simeq \begin{pmatrix} r_{1,\psi}^i & * \\ 0 & r_{2,\psi}^i \end{pmatrix}
\]
We can take a semisimplified reduction to the residue field and then, by Proposition \ref{prop:det mod p}, we have
\begin{equation}\label{eqn:detordsub}
\det \overline{r^i_{1,\psi}} \left(\Art_{F_{{\tilde{v}}}}(\varpi_{\tilde{v}})\right) 
= \det (\bar{\rho}_{\m}\otimes\overline{\psi})(\Art_{F_{\tilde{v}}}(\varpi_{\tilde{v}})).
\end{equation}

\begin{sublemma}
	Possibly after enlarging $\cO$, there is a continuous character $\overline{\psi}: G_F \to k^\times$, unramified at $S_p$ and with $\rhobar_{\widetilde{\m(\psi)}}$ decomposed generic, such 
	that \begin{enumerate}
		\item the isomorphism classes of the irreducible constituents of 
		$\rhobar_{\m(\psi)}|_{G_{F_{\tilde{v}}}}$ are disjoint from those of 
		$\rhobar_{\m(\psi)}^{\vee,c}(1-2n)|_{G_{F_{\tilde{v}}}}$;
		
		\item for all $1 \le i \le r$ the isomorphism classes of the irreducible 
		constituents of the residual representation 
		$\overline{r^i_{1,\psi}}$ coincide with those of 
		$\rhobar_{\m(\psi)}|_{G_{F_{\tilde{v}}}}$.
	\end{enumerate}	
\end{sublemma}
\begin{proof}
		We denote the irreducible constituents (with multiplicity) of 
	$\rhobar_{\m}|_{G_{F_{\tilde{v}}}}$ over $\overline{\F}_p$ by $S_1, \ldots , S_m$, and let $d_j = \dim S_j$ 
	and $\delta_j = \det S_j$. 
	We denote the irreducible constituents of $\rhobar_{\m}^{\vee,c}(1-2n)|_{G_{F_{\tilde{v}}}}$ over $\overline{\F}_p$
	by $T_1,\ldots,T_{\ell}$
    and let $e_j = \dim T_j$ and $\eta_j = \det T_j$. We enlarge $\mathcal{O}$ and hence its residue field $k$ so that all these constituents are defined over $k$. 
    
    For each $\psi$, the 
	irreducible constituents of
	$\overline{r_{1,\psi}^i}$ are given by a multiset $\{S_j\otimes\bar{\psi}: 
	j \in I_i \}\coprod \{T_j\otimes\bar{\psi}^{\vee,c}: j \in J_i\}$ for two subsets $I_i \subset 
	\{1,\ldots,m\}$ and $J_i\subset \{1,\dots, \ell\}$ with $\sum_{j\in I_i}d_j + \sum_{j \in J_i}e_j = n$. 
	Comparing what this entails for $\det \overline{r^i_{1,\psi}}(\Art_{F_{\tilde{v}}}(\varpi_{\tilde{v}}))$ with the formula (\ref{eqn:detordsub}), we get:

    \begin{equation}\label{eq:det after twisting}
	\left(\bar{\psi}^n\prod_{j=1}^m\delta_j\right)(\Art_{F_{\tilde{v}}}(\varpi_{\tilde{v}})) = \left(\prod_{j \in I_i} \bar{\psi}^{d_{j}}\delta_{j} 
	\prod_{j \in J_i} 
	(\overline{\psi}^{\vee,c})^{e_{j}}\eta_{j}\right)(\Art_{F_{\tilde{v}}}(\varpi_{\tilde{v}})).
    \end{equation}

            We claim that we 
	can choose $\overline{\psi}$ 
	so that any equation of this form forces $\sum_{j \in J_i}e_{j} = 0$ and 
	$\sum_{j \in 
		I_i}d_{j} = n$ (whilst also preserving the first condition of the 
	sub-lemma). Indeed, let $\mathcal{L}> n$ be a prime such that $\mathcal{L}\nmid \#(k^\times)\cdot p$. We enlarge the coefficients so that the new residue field $k'$ contains an element $\alpha\in (k')^\times$ of order $\mathcal{L}$. Since Grunwald--Wang allows us to find $\overline{\psi}$ with specified behaviour at any finite set of places, we can choose $\overline{\psi}: G_F\to (k')^{\times}$ 
    as follows: 
    \begin{itemize}
    \item $\bar{\psi}(\mathrm{Art}_{F_{\tilde{v}}}(\varpi_{\tilde{v}})) = \alpha$; 
    \item $\bar{\psi}(\mathrm{Art}_{F_{\tilde{v}^c}}(\varpi_{\tilde{v}^c})) = 1$; 
    \item $\bar{\psi}$ is locally trivial at a prime which is decomposed generic for $\bar{\rho}_{\widetilde{\m}}$. 
    \end{itemize}
From equation~\eqref{eq:det after twisting}, we obtain 
\[
\alpha^{n}\prod_{j=1}^m \delta_j(\mathrm{Art}_{F_{\tilde{v}}}(\varpi_{\tilde{v}}))
= \alpha^{\sum_{j\in I_i}d_j}\prod_{j\in I_i} \delta_j(\mathrm{Art}_{F_{\tilde{v}}}(\varpi_{\tilde{v}})) \prod_{j\in J_i} \eta_j(\mathrm{Art}_{F_{\tilde{v}}}(\varpi_{\tilde{v}})),  
\]
which is equivalent to 
\[
\alpha^{\sum_{j\in J_i}e_j} = \prod_{j\not\in I_i} \delta^{-1}_j(\mathrm{Art}_{F_{\tilde{v}}}(\varpi_{\tilde{v}})) \prod_{j\in J_i}\eta_j(\mathrm{Art}_{F_{\tilde{v}}}(\varpi_{\tilde{v}})) \in k^\times. 
\]
As the order of $\alpha$ is $\mathcal{L}>n\geq \sum_{j\in J_i}e_j\geq 0$ and as $\mathcal{L}\nmid \#(k^\times)$, the only possibility is that $\sum_{j\in J_i}e_j = 0$.

	We conclude that $J_i = \emptyset$ and $I_i = 
	\{1,\ldots,m\}$ for every $i$, so in other words the isomorphism classes of 
	the irreducible 
	constituents of the residual representation 
	$\overline{r_{1,\psi}^i}$ coincide with those of 
	$\rhobar_{\m(\psi)}|_{G_{F_{\tilde{v}}}}$.
\end{proof}

We may now assume that the isomorphism classes of the irreducible constituents of 
$\rhobar_{\m}|_{G_{F_{\tilde{v}}}}$ are disjoint from those of 
$\rhobar_{\m}^{\vee,c}(1-2n)|_{G_{F_{\tilde{v}}}}$ and that for all $1 \le i \le r$ the isomorphism classes of the irreducible 
constituents of the residual representation 
$\overline{r^i_{1}}$ coincide with those of 
$\rhobar_{\m}|_{G_{F_{\tilde{v}}}}$. Applying Proposition \ref{prop:char0lifts} and Theorem \ref{thm:shape general} for the two places $v = \tilde{v}, \tilde{v}^c$, we deduce the statement of the Proposition for the Hecke algebra $A(K,\lambda,q,m)$. 

It remains to handle the case $q < \lfloor \frac{d}{2}\rfloor$. Poincar\'{e} duality gives an isomorphism:
\[\iota: A(K,\lambda,q,m) \cong A^\vee(K,\lambda,d-1-q, m).\] We can now run the same argument as above, using Proposition \ref{prop:degree shifting dual}. Since the Satake map in the dual degree shifting is untwisted, in this case we will have $2n$-dimensional $Q_{\bar{v}}$-ordinary representations $\tilde{\rho}_{\m^\vee}$ lifting \[\rhobar_{\tilde{\iota}^*\cS^*(\m^\vee)} = \rhobar_{\m}(-n)\oplus \rhobar_{\m}^{\vee,c}(1-n)\] with a decomposition \[\tilde{\rho}_{\m^\vee} |_{G_{F_{\tilde{v}}}} \simeq \begin{pmatrix} r_1 & * \\ 0 & r_{2} \end{pmatrix}
\] such that the \emph{lower right} block $r_2$ lifts $\rhobar_{\m}(-n)$. Proposition \ref{prop:char0lifts} and Theorem \ref{thm:shape general} then apply to $\rho_{\m}^{\vee,c}$, giving the desired properties of $\rho_{\m}$ by taking duals. This is compatible with the following: in the dual case, we twist by a sufficiently \emph{negative} power of the cyclotomic character to arrange that $\tilde{\lambda} = w^P_0\lambda$ is dominant for $\widetilde{G}$ (using our standard identification of weights for $G$ and $\widetilde{G}$). Then $\tilde{\rho}_{\m^\vee}$ is cohomological of weight $w^P_0\lambda$, so the $\tau$-labelled Hodge--Tate weights of $\tilde{\rho}_{\m^\vee} |_{G_{F_{\tilde{v}}}}$ are given by \[-\lambda_{\tau c,1}< \cdots < -\lambda_{\tau c, n}+n-1 < \lambda_{\tau,n} + n < \cdots < \lambda_{\tau,1}+2n-1.\] In particular, the Hodge--Tate weights of $r_2(n)$ are as expected. This completes the proof for all values of $q$.
\end{proof}
\subsubsection{Local--global--compatibility using deformation rings}

\noindent We formulate a consequence of Proposition \ref{prop:lgc torsion CTG} in terms of Galois deformation rings. The local deformation rings we need were defined in \S\ref{sec:localdefrings}.

\begin{thm}\label{thm:RtoT factors thru Kisin def ring}
Suppose that $F$ is an imaginary CM field that contains an imaginary quadratic field. Let $p$ 
be a prime which splits in an imaginary quadratic subfield of $F$. Let $T$ be a finite set of finite places of $F$, which 
contains $S_p$ and which 
is stable under complex conjugation, and such that the following condition is satisfied:
\begin{itemize}
\item \IQFassm
\end{itemize}

\noindent Let $K \subset G(\AA_{F^+,f}) = \GL_n(\AA_{F,f})$ be a good subgroup with $K_{v}=\GL_n(\cO_{F_v})$ $\forall v\not\in T$. 
Fix distinct places $\bar{v},\bar{v'}\in S_p$. Let $\lambda$ be a dominant weight for $G$.

We fix a standard parabolic $Q_{\bar{v}} \subset P_{\bar{v}}$ and suppose we are in one of three cases:

\begin{itemize}
	\item[(cr-ord)] $\iota_{\tilde{v}}(Q_{\bar{v}})$ is the standard parabolic given by the partition $(n,1,\ldots,1)$ of $2n$.
	\item[(ord)]  $\iota_{\tilde{v}}(Q_{\bar{v}}) = \Brm_{2n}$
	\item[(cr)] $Q_{\bar{v}} = P_{\bar{v}}$.
\end{itemize}

Suppose that $K_{\bar{v}}  = \cQ_{\bar{v}} \cap G(F^+_{\bar{v}})$ and let $\frakm \subset \TT^{\cQ_{\bar{v}},\{\bar{v}\}-\ord}$ be a maximal ideal in the support of $H^*(X_K,\cV_\lambda)$. 

Assume that:
\begin{enumerate}
	\item We have 
	\[
	\sum_{\substack{\bar{v}''\in \bar{S}_p \\ \bar{v}''\not = \bar{v},\bar{v}'}} [F^+_{\bar{v}''}:\Q_p] \geq \frac{1}{2}[F^+:\Q]. 
	\]
	\item $\m$ is a non-Eisenstein maximal ideal such that $\bar{\rho}_{\m}$ 
	is decomposed generic. 
	\item For all $\nu \in X_{Q_{\bar{v}}}$, the Hecke operator $[K_{\bar{v}}\nu(\varpi_{\bar{v}})K_{\bar{v}}]$ is not contained in $\m$.
\end{enumerate}

Then there exist an integer $N \ge 1$, depending only on $n$ and $[F^+:\Q]$, a nilpotent ideal $J$ of $\TT^{\cQ_{\bar{v}},\{\bar{v}\}-\ord}(R\Gamma(X_K,\cV_\lambda)_{\m})$ with $J^N = 0$ and a continuous $n$-dimensional representation
\[\rho_{\m}: G_{F,T} \to \GL_n(\TT^{\cQ_{\bar{v}},\{\bar{v}\}-\ord}(R\Gamma(X_K,\cV_\lambda)_{\m})/J) \] satisfying 
\begin{itemize}
\item[(lgc-ur)] For each place $v \notin T$ of $F$, $\rho_{\m}$ is unramified and 
the characteristic polynomial of $\rho_{\m}(\Frob_v)$ is equal to the image of $P_v(X)$.
\end{itemize}  
Moreover, the induced map $t_{\rho_{\m}}:R^\square_{\rhobar_{\m}} \to \TT^{\cQ_{\bar{v}},\{\bar{v}\}-\ord}(R\Gamma(X_K,\cV_\lambda)_{\m})/J$ has the following property in each of our three cases:
	\begin{itemize}
		\item[(cr-ord)] The restriction of $t_{\rho_{\m}}$ to $R_{\rhobar_{\m}|_{G_{F_{\tilde{v}}}}}^\square$ factors through $R_{\rhobar_{\m}|_{G_{F_{\tilde{v}}}}}^{\triangle,\lambda_{\tilde{v}}}$ and the restriction to $R_{\rhobar_{\m}|_{G_{F_{\tilde{v}^c}}}}^\square$ factors through $R_{\rhobar_{\m}|_{G_{F_{\tilde{v}^c}}}}^{\cris,\lambda_{\tilde{v}^c}}$.
		
		\item[(ord)] For $v|\bar{v}$, restriction of $t_{\rho_{\m}}$ to $R_{\rhobar_{\m}|_{G_{F_{{v}}}}}^\square$ factors through $R_{\rhobar_{\m}|_{G_{F_{{v}}}}}^{\triangle,\lambda_{{v}}}$.
		
		\item[(cr)] For $v|\bar{v}$, restriction of $t_{\rho_{\m}}$ to $R_{\rhobar_{\m}|_{G_{F_{{v}}}}}^\square$ factors through $R_{\rhobar_{\m}|_{G_{F_{{v}}}}}^{\cris,\lambda_{{v}}}$.
	\end{itemize}
\end{thm}
\begin{proof}
	The first point is that it is enough to prove our statement for Galois representations with coefficients in $\TT^{\cQ_{\bar{v}},\{\bar{v}\}-\ord}(R\Gamma(X_K,\cV_\lambda/\varpi^m)_{\m})$ for integers $m \in \ZZ_{\ge 1}$. This is because we have an isomorphism (by \cite[Lemma 3.11]{new-tho}) \[\TT^{\cQ_{\bar{v}},\{\bar{v}\}-\ord}(R\Gamma(X_K,\cV_\lambda)_{\m}) \toisom \varprojlim_m \TT^{\cQ_{\bar{v}},\{\bar{v}\}-\ord}(R\Gamma(X_K,\cV_\lambda/\varpi^m)_{\m}).\]
    We can pass to the inverse limit in $m$ because the degree of nilpotence of the ideal $J$ is bounded only in terms of $n$ and $[F^+:\Q]$, not in terms of $m$. 
	In fact, we can prove the statement one cohomological degree at a time (cf.~the proof of \cite[Theorem 4.5.1]{10author}). Indeed, the kernel of the map \[\TT^{\cQ_{\bar{v}},\{\bar{v}\}-\ord}(R\Gamma(X_K,\cV_\lambda/\varpi^m)_{\m}) \to \prod_{q}\TT^{\cQ_{\bar{v}},\{\bar{v}\}-\ord}(H^q(X_K,\cV_\lambda/\varpi^m)_{\m})\] is a nilpotent ideal with vanishing $d$th power, and a Galois representation with coefficients in a quotient of $\prod_{q}\TT^{\cQ_{\bar{v}},\{\bar{v}\}-\ord}(H^q(X_K,\cV_\lambda/\varpi^m)_{\m})$ which satisfies condition (lgc-ur) can be conjugated to take values in the image of  $\TT^{\cQ_{\bar{v}},\{\bar{v}\}-\ord}(R\Gamma(X_K,\cV_\lambda/\varpi^m)_{\m})$ (by Carayol's lemma, cf.~the proof of \cite[Corollary 4.4.8]{10author}). 
	
	Now our statement follows from Proposition \ref{prop:lgc torsion CTG}. This produces a lift of the map $R_{\rhobar_{\m}|_{G_{F_{{v}}}}}^\square \to \TT^{\cQ_{\bar{v}},\{\bar{v}\}-\ord}(H^q(X_K,\cV_\lambda/\varpi^m)_{\m})/J$ to a map with target a finite flat local $\cO$-algebra $R_{\rhobar_{\m}|_{G_{F_{{v}}}}}^\square \to \widetilde{A}$. This map factors through the appropriate (crystalline or ordinary) quotient by the characterizing property of this quotient. 
\end{proof}

We also need a small refinement which will help us in a `fixed determinant' setting:

\begin{cor}\label{cor:RtoT factors thru Kisin def ring fixed det}
In the setting of Theorem \ref{thm:RtoT factors thru Kisin def ring}, assume moreover that $p\nmid n$ and we have a quotient map $f: \TT^{\cQ_{\bar{v}},\{\bar{v}\}-\ord}(R\Gamma(X_K,\cV_\lambda)_{\m})/J \to A$ such that $\det(f_\ast(\rho_{\m})) = \psi$ for a character $\psi: G_{F,T} \to \cO^\times$ which is crystalline at all places in $S_p$ with $\tau$-labelled Hodge--Tate weights $\sum_{i=1}^n \lambda_{\tau,i}+(n-i)$ for each $\tau:F\hookrightarrow E$. 

Then for $v|\bar{v}$ the induced map $R_{\rhobar_{\m}|_{G_{F_{v}}}}^{\triangle,\lambda_{v}} \to A$ or $R_{\rhobar_{\m}|_{G_{F_{v}}}}^{\cris,\lambda_{v}} \to A$ factors through the appropriate fixed determinant $\psi$ lifting ring (cf.~\S\ref{sec:fixeddetdefrings}).
\end{cor}
\begin{proof}
This follows from Lemma \ref{lem:twistingdefrings}.
\end{proof}

\subsection{The characteristic $0$ case}
For simplicity, we restrict to the crystalline case here. %
For this subsection, we drop our running assumption that $F$ contains an imaginary quadratic field.
\begin{thm}\label{thm:LGC in char 0}
	Suppose $\pi$ is a cuspidal automorphic representation of $\GL_n(\AA_F)$, regular algebraic of weight $\lambda$, with $F$ a totally real or CM field. Let $v$ be a place of $F$ dividing $p$ and suppose $\pi^{\GL_n(\cO_{F_v})}$ and $\pi^{\GL_n(\cO_{F_{v^c}})}$ are both non-zero. (We allow the possibility that $v = v^c$, even in the CM case.) Let $\iota:\Qpbar \to \C$ be an isomorphism, and consider the continuous semisimple representation $r_{\iota}(\pi):G_F\to\GL_n(\Qpbar)$ constructed in \cite{hltt}. Assume that:
	\begin{enumerate}
		\item $\overline{r_{\iota}(\pi)}$ is irreducible and decomposed generic.
	\end{enumerate}

Then $r_{\iota}(\pi)|_{G_{F_v}}$ and $r_{\iota}(\pi)|_{G_{F_{v^c}}}$ are crystalline with $\tau$-labelled Hodge--Tate weights $\lambda_{\iota\tau,n} < \cdots < \lambda_{\iota\tau,1}+n-1$ for $\tau: F \hookrightarrow \Qpbar$ inducing $v$ or $v^c$ respectively.
\end{thm}
\begin{proof}
Fix a prime $\ell$ such that $\overline{r_{\iota}(\pi)}$ satisfies the decomposed generic condition at $\ell$. Using cyclic base change and case (cr) of Theorem \ref{thm:RtoT factors thru Kisin def ring}, it suffices to find a solvable CM extension $F'/F$ with the following properties:
\begin{enumerate}
	\item $F'$ is linearly disjoint from $\overline{F}^{\ker \overline{r_{\iota}(\pi)}}$ over $F$.
	\item $F'$ contains an imaginary quadratic field in which $p$ splits. 
	\item $\ell$ splits completely in $F'$.
	\item The places $v, v^c$ split completely in $F'$.
	\item There is a place $\bar{w}$ of $(F')^+$, lying over the place $\bar{v}|v$ of $F^+$, and another $p$-adic place $\bar{w}'$ of $(F')^+$ such that \[\sum_{\substack{\bar{w}''|p \text{ in }(F')^+\\\bar{w}''\ne \bar{w},\bar{w}'}}[(F')^+_{\bar{w}''}:\Qp]\ge \frac{1}{2}[(F')^+:\Q].\]
\end{enumerate}
We can achieve the final property by choosing $F'$ with $[(F')^+:F^+] \ge 4$ and with $\bar{v}$ split completely in $(F')^+$, and then choosing $\bar{w}, \bar{w}'$ to be two distinct places of $(F')^+$ lying over $\bar{v}$. We conclude that it is possible to find such an extension $F'$. %
\end{proof}

\section{Automorphy lifting}
\subsection{A potentially Barsotti--Tate modularity lifting theorem}
We begin by stating the main theorem in this section. In order to get an optimal result for applications to modularity of elliptic curves, we only consider Galois representations with inverse-cyclotomic determinant. 

\begin{theorem}\label{thm:pBT_lifting}
Let $F$ be an imaginary CM field and let $p$ be an odd prime. %
Suppose given a continuous representation $\rho : G_F \to \GL_2(\overline{\bQ}_p)$ satisfying the following conditions:
\begin{enumerate}
\item $\rho$ is unramified almost everywhere and $\det\rho = \epsilon_p^{-1}$.
\item\label{assm:pBT} For each place $v | p$ of $F$, the representation $\rho|_{G_{F_v}}$ is potentially semistable with all labelled Hodge--Tate weights equal to $(0,1)$.  
\item $\overline{\rho}$ is decomposed generic (Definition \ref{defn:generic}) and $\overline{\rho}|_{G_{F(\zeta_p)}}$ is irreducible. 
\item\label{assumption:p5} If $p = 5$ and the projective image of $\rhobar(G_{F(\zeta_5)})$ is conjugate to $\PSL_2(\F_5)$, we assume further that the extension of $F$ cut out by the projective image of $\rhobar$ does not contain $\zeta_5$.
\item\label{assumption:A4} If $\det\rhobar$ has order $3$ (hence $3|(p-1)$), then the projective image of $\rhobar$ is not isomorphic to $A_4$. 
\item There exists a cuspidal automorphic representation $\pi$ of $\PGL_2(\bA_F)$ and an isomorphism $\iota:\Qpbar\cong \C$ satisfying the following conditions:
\begin{enumerate}
\item $\pi$ is regular algebraic of weight $0$. %
\item\label{assm:samecomp} For each place $v |  p$ where $\rho|_{G_{F_v}}$ is potentially crystalline, $r_{\iota}(\pi)|_{G_{F_v}}$ is potentially ordinary of weight $0$ (in the sense of \cite[\S5.2]{geraghty}) if and only if $\rho|_{G_{F_v}}$ is potentially ordinary of weight $0$. We moreover assume that $\rec_{F_v}(\pi_v)$ has monodromy operator $0$. 
\item\label{assm:sameordcomp} For each place $v | p$ where $\rho|_{G_{F_v}}$ is not potentially crystalline, $\pi$ is $\iota$-ordinary of weight $0$ at $v$ and $r_{\iota}(\pi)|_{G_{F_v}}$ is not potentially crystalline. 
\item\label{assm:resmod} There is an isomorphism $\rhobar \cong \rbar_{\iota}(\pi)$.
\end{enumerate}
\end{enumerate}
Then $\rho$ is automorphic: there exists a cuspidal automorphic representation $\Pi$ of $\PGL_2(\bA_F)$, regular algebraic of weight $0$, such that $\rho \cong r_\iota(\Pi)$. 
\end{theorem}
\begin{rem}Using Theorem \ref{thm:RtoT factors thru Kisin def ring}, we can replace the assumption that $r_{\iota}(\pi)|_{G_{F_v}}$ is potentially ordinary of weight $0$ at certain places $v$ with an assumption on $\pi_v$.
\end{rem}
\begin{rem}
We restrict to Galois representations with inverse-cyclotomic determinant, as in \cite{AKT}, so that we can handle the case where $p=3$ and the image of $\rhobar(G_{F(\zeta_3)})$ is $\SL_2(\F_3)$. Assumption \eqref{assumption:A4} comes from the exceptional case in Lemma \ref{lem:DDT4.11} --- when this assumption fails, it will often be possible to apply the Theorem after restriction to $G_{F(\zeta_p)}$. It needs to be checked that the decomposed generic assumption is still valid after restriction, but this will often be the case (cf.~Lemma \ref{lem:quadratic irred implies generic}, for example).
\end{rem}

\begin{rem} We could consider also the Fontaine--Laffaille case, for arbitrary dimension $n$ Galois representations -- namely, assume that 
$p$ is unramified in $F$ and that the weight $\lambda$ satisfies the condition:
\begin{equation}\label{eq:fl up to twist}
\lambda_{\tau, 1} - \lambda_{\tau, n}< p-n \ \mathrm{for}\ \mathrm{all}\ \tau\in \Hom(F, E). 
\end{equation}
Using the local-global compatibility result given by Theorem~\ref{thm:RtoT factors thru Kisin def ring} instead of Theorem~\cite[Theorem 4.5.1]{10author}, one could also prove an automorphy lifting theorem in this case that strengthens~\cite[Theorem 6.1.1]{10author}. The restrictions on the Hodge--Tate weights in the automorphy lifting theorem in \emph{loc.~cit.} are stronger than the condition in~\eqref{eq:fl up to twist} only because of restrictions in the corresponding result on local-global compatibility. The rest of the argument would go through verbatim. 
\end{rem}

\subsection{Galois deformation theory} %
To prove our automorphy lifting theorem we apply the patching method in \cite[\S6]{10author}, making modifications  as in \cite{AKT} to avoid the assumption that there is a $\sigma \in G_F - G_{F(\zeta_p)}$ such that $\rhobar(\sigma)$ is scalar and to include cases with $p =3$ or $5$ where $\rhobar|_{G_{F(\zeta_p)}}$ does not have enormous image. In addition, we work with local lifting rings at $p$ that have two irreducible components  (either ordinary/non-ordinary or crystalline ordinary/non-crystalline ordinary), which is why we need assumption  (\ref{assm:samecomp}) in this theorem. The fact that these local lifting rings have generically reduced special fibre, which we will recall shortly, is important for implementing (derived) Ihara avoidance in a situation where local lifting rings have more than one component. 

We adopt all the terminology and notation of \cite[\S6.2.1]{10author}, although the coefficient ring `$\Lambda$' appearing there will always be $\cO$ for us, and our Galois representations will all be two-dimensional. So, we fix a continuous and absolutely irreducible $\rhobar: G_F \to \GL_2(k)$, and let $S$ be a finite set of finite places of $F$ containing $S_p$ and all the places where $\rhobar$ is ramified. 

For each $v \in S$, we have a local lifting ring $R_v^\square$ for $\rhobar_v:=\rhobar|_{G_{F_v}}$ and the notion of a local deformation problem $\cD_v$: a set valued functor on $\CNL_{\cO}$ satisfying some conditions which in particular imply that it is represented by a quotient $R_v$ of $R_v^\square$.

A global deformation problem is a tuple \[(\rhobar,S,\{R_v\}_{v \in S}),\] where for each $v \in S$, $R_v$ is a quotient of $R_v^\square$ representing a local deformation problem for $\rhobar_v$.

A global deformation problem with fixed determinant is a tuple \[(\rhobar,\psi,S,\{R_v\}_{v \in S}),\] where $\psi: G_{F,S} \to \cO^\times$ is a character which lifts $\det(\rhobar)$, $(\rhobar,S,\{R_v\}_{v \in S})$ is a global deformation problem, and the lifts parameterized by each $R_v$ have determinant $\psi|_{G_{F_v}}$. 

If $\cS$ is a global deformation problem (with or without fixed determinant), and $T$ is a subset of $S$, we have the functor $\cD_{\cS}^T$ of $T$-framed deformations of type $\cS$. When $T = \emptyset$, we denote the functor (of deformations of type $\cS$) by $\cD_{\cS}$.

The functors $\cD_{\cS}$ and $\cD_{\cS}^T$ are represented by $\CNL_{\cO}$-algebras $R_\cS$ and $R_{\cS}^T$, respectively.

For a global deformation problem $\cS= (\rhobar,S,\{R_v\}_{v \in S})$ (or one with fixed determinant, $\cS = (\rhobar,\psi,S,\{R_v\}_{v \in S})$) and $T \subset S$, we define $R_{\cS}^{T,\loc} = \widehat{\otimes}_{v \in T}R_v$. There is a natural local $\cO$-algebra map $R_{\cS}^{T,\loc} \to R_{\cS}^T$. 

We are assuming that $p$ is odd. We will make use of the following local deformation problems (which we identify in terms of their representing ring), where $\psi$ always denotes a fixed determinant character.

\begin{itemize}
\item Fixed determinant lifting rings $R_v^{\psi}$, parameterizing lifts of $\rhobar|_{G_{F_v}}$ with determinant $\psi|_{G_{F_v}}$. We will make use of these rings when $v \notin S_p$ and $H^2(F_v,\ad^0\rhobar) = 0$, in which case $R_v^{\psi}$ is formally smooth over $\cO$ of relative dimension $3$.
\item Fixed determinant `level raising' lifting rings $R_v^{\psi,\chi}$, for $v$ with $q_v \equiv 1 \mod p$ and $\rhobar|_{G_{F_v}}$ trivial, and a character $\chi:\cO_{F_v}^\times \to \cO^\times$ which is trivial modulo $\varpi$ (cf.~\cite[{\S}A.1.2]{AKT}). They classify lifts $\rho$ with determinant $\psi$ and characteristic polynomial \[\mathrm{char}_{\rho(\sigma)}(X) = (X-\chi(\Art_{F_v}^{-1}(\sigma)))(X-\chi^{-1}\psi(\Art_{F_v}^{-1}(\sigma)))\] for all $\sigma \in I_{F_v}$.

\item Barsotti--Tate lifting rings $R_v^{\psi,\BT}$ for $v|p$. These are the rings $R_{\rhobar_v}^{\cris,(0,0)_{\tau \in \Hom(F_v,E)},\psi}$ in the notation of \S\ref{sec:localdefrings}. Note that we assume that $\psi|_{G_{F_v}}$ is crystalline with all labelled Hodge--Tate weights equal to $1$ and that it lifts $\det\rhobar_v$ for this ring to be defined. In practice, we will take $\psi$ to be the inverse of the cyclotomic character.
\end{itemize}

\subsubsection{Ordinary deformation rings} We will also use, when $v\mid p$, the ordinary lifting ring $R_v^{\triangle,(0,0)_{\tau \in \Hom(F_v,E)},\psi}$ with fixed determinant and Hodge--Tate weights $(0,1)$. For simplicity, we only consider the case $\psi = \epsilon_p^{-1}$, with $\rhobar|_{G_{F_v}}$ trivial and $\epsilonbar_p$ trivial on $G_{F_v}$, and set 
\[R_v^{\triangle} :=  R_v^{\triangle,(0,0)_{\tau \in \Hom(F_v,E)},\psi}\] in the notation of \S\ref{sec:localdefrings}.

We recall some properties of $R_v^{\triangle}$, following \cite{snowden-ord}. The reader should refer to the published version, since this contains some updates compared to the arXiv preprint of \cite{snowden-ord}. Note that twisting by the cyclotomic character and using the reducedness of $R_v^{\triangle}$ (Theorem \ref{thm:localdefrings}) shows that our lifting ring can indeed be identified with the ring denoted by $R$ in \cite[Proposition 4.3.2]{snowden-ord}.

\begin{prop}\label{prop:ordcomps}\begin{enumerate}
\item $\Spec(R_v^{\triangle})$ is equidimensional of dimension $[F_v:\Qp]+4$, with two irreducible components $X^{cr}=\Spec(R_v^{\triangle,cr}), X^{st}=\Spec(R_v^{\triangle,st})$ characterized by their points valued in finite extensions $E'/E$:
\begin{itemize}
\item $x: R_v^{\triangle} \to E'$ factors through $R_v^{\triangle,cr}$ if and only if $\rho_x$ is crystalline.
\item $x: R_v^{\triangle} \to E'$ factors through $R_v^{\triangle,st}$ if and only if $\rho_x$ is conjugate to a representation of the form $\begin{pmatrix}1 & *\\ 0 & \epsilon_p^{-1} \end{pmatrix}.$
\end{itemize}
\item Each generic point of $\Spec(R_v^\triangle/\varpi)$ is the specialization of a unique generic point of $\Spec(R_v^\triangle)$. 
\end{enumerate}
\end{prop}
\begin{proof}
The first part follows from \cite[Proposition 4.3.2]{snowden-ord}. The second part also follows from Snowden's results, as we now explain. It suffices to show that the dimension of $X^{cr}\cap X^{st} \cap \Spec(R_v^\triangle/\varpi)$ is $< [F_v:\Qp]+3$, as this shows that there is no point of large enough dimension to be a generic point of $\Spec(R_v^\triangle/\varpi)$ generalizing to both generic points of $\Spec(R_v^\triangle)$.

Snowden defines another ring $\widetilde{R}_v^\triangle$ \cite[Proposition 4.4.3]{snowden-ord}. The ring $\widetilde{R}_v^\triangle$ is a quotient ($\cO$-flat and reduced) of the finite $R_v^{\square,\epsilon_p^{-1}}$-algebra given by adjoining a root of the characteristic polynomial of a lift of Frobenius under the universal lifting of $\rhobar|_{G_{F_v}}$. It comes with a finite morphism $\pi: \Spec(\widetilde{R}_v^\triangle) \to \Spec(R_v^\triangle)$ which is an isomorphism after inverting $p$. In particular, $\pi$ is surjective and induces a bijection between irreducible components. We denote the irreducible components of $\Spec(\widetilde{R}_v^\triangle)$ lying over $X^{cr}$ and $X^{st}$ by $\widetilde{X}^{cr}$ and $\widetilde{X}^{st}$ respectively. The map $\pi$ induces a finite surjective map $\widetilde{X}^{st}\cap \widetilde{X}^{cr}  \cap \Spec(\widetilde{R}_v^\triangle/\varpi) \to X^{cr}\cap X^{st} \cap \Spec(R_v^\triangle/\varpi)$. Indeed, it follows from \cite[Theorem 4.6.1, Lemma 4.6.4]{snowden-ord} that a point of $\widetilde{X}^{cr}  \cap \Spec(\widetilde{R}_v^\triangle/\varpi)$ is contained in $\widetilde{X}^{st}$ if and only if an element lifting Frobenius has unipotent image under the corresponding Galois representation. This condition holds for every element of $X^{st} \cap \Spec(R_v^\triangle/\varpi)$. We can now bound the dimension of the target of this map by bounding the dimension of the source. Snowden explicitly describes the mod $\varpi$ fibre $\Spec(\widetilde{R}_v^\triangle/\varpi)$ \cite[Lemma 4.6.4]{snowden-ord}, and its two irreducible components intersect in a proper closed subset (of dimension $[F_v:\Qp]+2$).
\end{proof}

\begin{lemma}\label{lem:cegs}
Each generic point of $\Spec(R_v^{\epsilon_p^{-1},\BT}/\varpi)$ is the specialization
of a unique generic point of $\Spec(R_v^{\epsilon_p^{-1},\BT})$.  
\end{lemma}
\begin{proof}
We let $R = R_v^{\epsilon_p^{-1},\BT}$. The lemma follows from generic reducedness of $\Spec(R/\varpi)$ and the fact that every generic point of $\Spec(R)$ has characteristic $0$. The generic reducedness follows from \cite[Theorem 1.3]{CEGSA}. The cited result applies to the lifting ring $R_v^{\BT} = R_v^{(0,0)_{\tau \in \Hom(F_v,\Qpbar)}}$ without fixed determinant, but since $p$ is odd $R_v^{\BT}$ is formally smooth over $R$. So we deduce that $R/\varpi$ is also generically reduced.  Every generic point of $\Spec(R)$ has characteristic $0$ because $R$ is $\cO$-flat. To deduce the claim about generic points of $\Spec(R/\varpi)$, let $\p \in \Spec(R)$ be the image of a generic point of $\Spec(R/\varpi)$. Since $\Spec(R/\varpi)$ is generically reduced, $R_\p/\varpi R_\p = (R/\varpi)_\p$ is a field. Now we know that $R_\p$ is a Noetherian local ring with a principal maximal ideal; it is therefore a local principal ideal ring with non-zero ideals generated by powers of $\varpi$. Since $\p$ is not a generic point of $\Spec(R)$, $R_\p$ has dimension $1$ and $(0)$ is its unique minimal prime (in particular, $R_\p$ is a DVR). 
\end{proof}

\begin{lemma}\label{lem:BTcomp}Let $v$ be a place of $F$ with $v|p$. Suppose that $\rhobar|_{G_{F_v}}$ is trivial, the residue field $k_v$ is not equal to $\F_p$, and $R_v^{\epsilon_p^{-1},\BT}$ is non-zero. Then $\Spec(R_v^{\epsilon_p^{-1},\BT})$ has exactly two irreducible components, one whose points correspond to ordinary Galois representations and one whose points correspond to non-ordinary Galois representations.
\end{lemma}
\begin{proof}
This follows from results of Kisin \cite[Corollary 2.5.16]{kis04} and Gee \cite[Proposition 2.3]{MR2280776}.
\end{proof}

\begin{defn}
A Taylor–Wiles datum for a global deformation problem $\cS$ is a tuple $(Q, N, (\alpha_{v,1},\alpha_{v,2})_{v \in Q})$ consisting
of:
\begin{itemize}
\item A finite set of finite places $Q$ of $F$, disjoint from $S$, and a positive integer $N$ such that $q_v \equiv 1 \mod p^N$ for each $v \in Q$.
\item For each $v \in Q$ and $i \in \{1,2\}$, distinct unramified $k$-valued characters $\alpha_{v,1}, \alpha_{v,2}$ such that $\rhobar|_{G_{F_v}} \cong \bigoplus_{i=1}^2 \alpha_{v,i}$.
\end{itemize}
We call $N$ the level of the Taylor--Wiles datum. 

If $(Q, N, (\alpha_{v,1}, \alpha_{v,2})_{v \in Q})$ is a Taylor--Wiles datum for $\cS$, then we define a new global deformation problem
\[\cS_Q= (\rhobar,S\cup Q, \{R_v\}_{v \in S}\cup\{R_v^{\square}\}_{v \in Q}) \]
(respectively, $\cS_Q = (\rhobar,\psi,S\cup Q, \{R_v\}_{v \in S}\cup\{R_v^{\psi}\}_{v \in Q} )$ if $\cS$ has fixed determinant).
\end{defn}

\subsection{Patching}\label{ssec:patching}
Our set-up is very close to that of \cite[\S6]{10author}. First we give an axiomatic description of the kinds of objects which will be the output of the patching method, and deduce a formal modularity lifting result.

We assume given the following objects:
\begin{enumerate}
\item A power series ring $S_\infty=\cO[[X_1,\cdots,X_r]]$ with augmentation ideal $\mathfrak{a}_\infty = (X_1, \dots, X_r)$. 
\item\label{ihara_complex_iso} Perfect complexes $C_\infty, 
C'_\infty$ of $S_\infty$-modules, and a fixed isomorphism 
\[ C_\infty \otimes^\bL_{S_\infty} S_{\infty} / \varpi \cong C'_\infty \otimes^\bL_{S_\infty} S_{\infty} / \varpi \]
in $\mathbf{D}(S_\infty / \varpi)$.
\item Two $S_\infty$-subalgebras
\[ T_\infty\subset 
\End_{\bD(S_\infty)}(C_\infty) \]
and
\[ T'_\infty \subset 
\End_{\bD(S_\infty)}(C'_\infty), \]
which have the same image in
\[ \End_{\bD(S_\infty/\varpi)}(C_\infty \otimes^\bL_{S_\infty} S_{\infty} / \varpi ) = \End_{\bD(S_\infty/\varpi)}(C'_\infty \otimes^\bL_{S_\infty} S_{\infty} / \varpi ), \]
where these endomorphism algebras are identified using the fixed isomorphism in (\ref{ihara_complex_iso}). Call this common image 
$\overline{T}_\infty$. Note that $T_\infty$ and $T'_\infty$ are finite 
$S_\infty$-algebras.
\item Two Noetherian complete local $S_\infty$-algebras $R_\infty$ and $R'_\infty$ and 
surjections $R_\infty\onto T_\infty/I_\infty$, $R'_\infty\onto 
T'_\infty/I'_\infty$, where $I_\infty$ and $I'_\infty$ are nilpotent ideals. 
We write $\overline{I}_\infty$ and $\overline{I}'_\infty$ for the images of these 
ideals in $\overline{T}_\infty$. Note that it then makes sense to 
talk about the support of $H^*(C_\infty)$ and $H^*(C'_\infty)$ over $R_\infty$, 
$R'_\infty$, even though they are not genuine modules over these rings. These 
supports actually belong to the closed subsets of $\Spec R_\infty$, $\Spec 
R'_\infty$ given by $\Spec T_\infty$, $\Spec T'_\infty$, and hence are finite 
over $\Spec S_\infty$.
\item An isomorphism $R_\infty/\varpi\cong R'_\infty/\varpi$ 
compatible with the $S_\infty$-algebra structure and the actions 
(induced from $T_\infty$ and $T'_\infty$) on 
\[ H^*( C_\infty \otimes^\bL_{S_\infty} S_{\infty} / \varpi)/(\overline{I}_\infty+\overline{I}'_\infty) = H^*( C'_\infty \otimes^\bL_{S_\infty} S_{\infty} / \varpi)/(\overline{I}_\infty+\overline{I}'_\infty), \]
where these cohomology groups are identified using the fixed isomorphism. 
\item Integers $q_0 \in \bZ$ and $l_0 \in \bZ_{\geq 0}$.
\end{enumerate}

\begin{assumption}    \label{assumptionsetup} Our set-up is assumed to satisfy the following: 
\begin{enumerate}
\item $\dim R_\infty=\dim R'_\infty=\dim S_\infty -l_0$, and $\dim R_\infty/\varpi=\dim R'_\infty/\varpi=\dim S_\infty-l_0-1$. 
\item (Behaviour of components) Assume that each generic point of $\Spec 
R_\infty/\varpi$ is the specialization of unique generic points of $\Spec 
R_\infty$ and $\Spec R'_\infty$. Moreover, we assume that $\Spec R_\infty$ and $\Spec R'_\infty$ are equidimensional and that all of their generic points have characteristic $0$. This implies that $\Spec R_\infty/\varpi$, $\Spec R'_\infty/\varpi$ are equidimensional of dimension one less (by the principal ideal theorem). 
\item (Generic concentration) \label{part:genericconcentration} 
We have
\[ H^\ast(C_\infty\otimes_{S_\infty}^{\bL} S_\infty/\mathfrak{a}_\infty )[\frac{1}{p}]\neq 0, \]
and these groups are non-zero only for degrees in the interval $[q_0,q_0+l_0]$.
\item (Automorphic point) \label{part:autpoint} We fix a characteristic $0$ point $x \in \Spec (T_\infty/\mathfrak{a}_{\infty}T_\infty)$.
\end{enumerate}
\end{assumption}

Note that $\Supp_{R_\infty}(H^*(C_\infty))=\Spec T_\infty$ and $\Supp_{R_\infty'}(H^*(C_\infty'))=\Spec T_\infty'$. (This is because the kernel of $T_\infty \to \End_{S_\infty}(H^*(C_\infty))$ is nilpotent, and the same holds for $T_\infty'$ and $C_\infty'$.)

\begin{prop}\label{prop:automorphic_component} Consider the automorphic subset of $\Spec R_\infty$: \[ \Supp_{R_\infty}(H^*(C_\infty)) = \Spec T_\infty \subset \Spec R_\infty.\]
	\begin{enumerate}
		\item There exists an irreducible component $C_a \subset \Spec R_\infty$, containing the automorphic point $x$, with $C_a \subset \Spec T_\infty$.
		\item\label{part:automorphic_component} Let $C_a \subset \Spec T_\infty$ be an irreducible component of $\Spec R_\infty$ which contains $x$. Suppose $C \subset \Spec R_\infty$ is an irreducible component such that the subsets $C\cap \Spec(R_\infty/\varpi)$ and $C_a\cap \Spec(R_\infty/\varpi)$ of $\Spec(R_\infty/\varpi)$ contain generic points $\barx_C$, $\barx_a$ respectively,  which generalize to the same generic point $x'$ of $\Spec R'_\infty$. Then $C \subset \Spec T_\infty$.
	\end{enumerate}
\end{prop}

\begin{proof} First we note that the pullback of $x$ to $S_\infty$ is $\mathfrak{a}_\infty$. We write $\tilde{x} \in \Spec(T_\infty\otimes_{S_\infty}S_{\infty,\ainf})$ for the prime ideal extending $x$. 

It follows from our assumptions and \cite[Lemma 6.2]{CG}, just as in the proof of \cite[Proposition 6.3.8]{10author}, that $H^*(C_{\infty,\ainf})$ is non-zero exactly in degree $q_0+l_0$ and that $M_\infty:= H^{q_0+l_0}(C_{\infty,\ainf})$ is a Cohen--Macaulay $S_{\infty,\ainf}$-module with depth and dimension equal to $\dim S_{\infty,\ainf} - l_0 = \dim S_{\infty} - l_0-1$. 

Since the image of $\ainf S_{\infty,\ainf}$ in $\End(H^{q_0+l_0}(C_{\infty,\ainf}))$ is contained in the image of $\tilde{x}$, we deduce that $\mathrm{depth}(\tilde{x},M_\infty) \ge \dim S_{\infty,\ainf} - l_0$. An $M_\infty$-regular sequence in $\tilde{x}$ remains regular on $M_{\infty,\tilde{x}}$, and the localization $(T_\infty\otimes_{S_\infty}S_{\infty,\ainf})_{\tilde{x}}$ is equal to $T_{\infty,x}$. So $\mathrm{depth}_{T_{\infty,x}}(M_{\infty,\tilde{x}}) \ge \dim S_{\infty,\ainf} - l_0$. In particular, $\dim T_{\infty,x} \ge \dim S_{\infty}- l_0-1$, so the one-dimensional prime $x$ is contained in an irreducible component of $\Spec T_{\infty}$ of dimension at least $\dim S_{\infty} - l_0$. By dimension considerations, this irreducible component can be identified with an irreducible component of $\Spec R_\infty$. This shows the existence of an irreducible component $C_a$ as in the first part.

For the second part, we let $x_a$ be the generic point of $C_a$. Since the pullback of $x_a$ to $S_\infty$ is contained in $\ainf$, the `localization' $C_{\infty,x_a}$ defined following \cite[Lemma 6.3.3]{10author} is quasi-isomorphic to the complex with $M_{\infty,x_a}$ in degree $q_0+l_0$ and zero elsewhere. In particular, with the length function on complexes defined in \cite[\S6.3.1]{10author}, we have $\lg_{T_{\infty,x_a}}(C_{\infty,x_a}) \neq 0$. The generic point $\barx_a$ of $\Spec R_\infty/(x_a,\varpi)$ given to us in the statement of the proposition has dimension $\dim S_\infty-l_0-1$ and is a generic point of $\Spec R_\infty/(\varpi)$ which lies in $\Spec T_\infty$. Let $\barx_a'$ denote the corresponding point of $\Spec R_\infty'/(\varpi)$. It has a unique generalization $x' \in \Spec R_\infty'$.

Now let $x_C$ be the generic point of $C$. We wish to show that it lies in $\Spec T_\infty$. We are given a generic point $\barx_C$ of $\Spec R_\infty/(x_C,\varpi)$, which must have dimension $\dim S_\infty-l_0-1$ and be a generic point of $\Spec R_\infty/(\varpi)$. Let $\barx_C'$ denote the corresponding point of $\Spec R_\infty'/(\varpi)$, which also generalizes to $x'$ by hypothesis.

We now repeatedly use \cite[Lemma 6.3.7]{10author}.
As $\lg_{T_{\infty,x_a}}(C_{\infty,x_a}) \neq (0)$, we deduce that  $\lg_{T_{\infty,\barx_a}}((C_{\infty} \otimes^{\bL}_{S_\infty} S_\infty/(\varpi))_{\barx_a}) \neq 0$. Hence $\barx_a' \in \Spec T_\infty'$ and $\lg_{T_{\infty,\barx_a'}'}((C_{\infty}' \otimes^{\bL}_{S_\infty} S_\infty/(\varpi))_{\barx_a'}) \neq 0$, from which we deduce that $x'\in \Spec T_\infty'$ and $\lg_{T_{\infty,x'}'}(C'_{\infty,x'}) \neq 0$. We further deduce that $\barx_C' \in \Spec T_\infty'$ and $\lg_{T_{\infty,\barx_C'}'}((C_{\infty}' \otimes^{\bL}_{S_\infty} S_\infty/(\varpi))_{\barx_C'}) \neq 0$. Hence $\barx_C \in \Spec T_\infty$ and $\lg_{T_{\infty,\barx_C}}((C_{\infty} \otimes^{\bL}_{S_\infty} S_\infty/(\varpi))_{\barx_C}) \neq 0$, from which we finally deduce that $x_C \in \Spec T_\infty$ (and $\lg_{T_{\infty,x_C}}(C_{\infty,x_C}) \neq (0)$). 
\end{proof}

\begin{cor}\label{cor:char0automorphy} 
Let $C$ be an irreducible component of $\Spec R_\infty$ satisfying the assumption of Proposition \ref{prop:automorphic_component}(\ref{part:automorphic_component}) for some `automorphic' component $C_a$. Let $x$ be a point of $C$, and let $y$ be the contraction of $x$ in 
$S_\infty$. Then the support of $H^\ast(C_\infty\otimes^{\bL}_{S_\infty} 
S_\infty/y)_y$ over $\Spec R_\infty$ contains 
$x$. If $y$ is one-dimensional of characteristic $0$ this says that $x$ is in the support of $H^\ast(C_\infty\otimes^{\bL}_{S_\infty} 
S_\infty/y)[1/p]$.
\end{cor}
\begin{proof}  
This follows from Proposition~\ref{prop:automorphic_component}(\ref{part:automorphic_component}) by considering the $\Tor$ spectral sequence computing the cohomology of $C_{\infty, y} \otimes^{\bL}_{S_{\infty,y}} S_{\infty,y}/y$, as in the proof of \cite[Corollary 6.3.9]{10author}. %
\end{proof}
\subsection{Hecke algebras and cohomology of locally symmetric spaces for $\PGL_2$}
We now go back to the constructions of \S\ref{sec:Hecke formalism}, which we apply to $\Grm = \Gbar =\PGL_{2,F}$ for an imaginary CM field $F$. We need to drop the assumption that $K = \prod_v K_v \subset \Gbar(\A_{F,f})$ is neat. We assume for convenience that $K \subset \PGL_2(\widehat{\cO}_F)$. Thanks to the results of \cite[\S5]{AKT}, all  the properties we need for cohomology of locally symmetric spaces for $\PGL_2$ can be deduced from the case of $\GL_2$ with neat level.

We fix a finite set $S$ of finite places of $F$ such that $S_p \subset S$ and $K_v = \PGL_2(\cO_{F_v})$ for $v \notin S$. We assume that $R = \cO$ or $\cO/\varpi^m$ for some $m \in \Z_{\ge 1}$. Let $\cV$ be an $R[K_S]$-module, finite free as an $R$-module, and such that $\cV/\varpi^r$ is a smooth $K_S$-module for each $r \ge 1$. 

We will make use of the Hecke algebra $\cH(\Gbar^S,K^S)$. For each finite place $v \notin S$ and $1 \le i \le 2$, we write $T_{v,i}$ for the image of $T_{v,i} \in \cH(\GL_2(F_v),\GL_2(\cO_{F_v}))$ in $\cH(\PGL_2(F_v),\PGL_2(\cO_{F_v}))$. In fact $T_{v,2} = 1$ in $\cH(\Gbar^S,K^S)$. We write $P_v(X)$ for the image of the polynomial \eqref{eqn:hecke_pol_for_GL_n} in $\cH(\PGL_2(F_v),\PGL_2(\cO_{F_v}))[X]$.

We consider the object
\[C^\bullet(K,\cV) := \varprojlim_r R\Gamma(K,R\Gamma(\overline{\frakX}_{\Gbar},\cV/\varpi^r))\] of $D^+(R)$ (we do not need to take a limit if $R$ is finite), which comes equipped with an action of $\cH(\Gbar^S,K^S)$.

More generally, if $K' = \prod K'_v \subset K$ is an open normal subgroup with $(K')^S = K^S$, then we consider the object
\[C^\bullet(K/K', \cV) = \varprojlim_r R\Gamma(K',R\Gamma(\overline{\frakX}_{\Gbar},\cV/\varpi^r))\] of $D^+(R[K/K'])$, which again comes with an action of $\cH(\Gbar^S,K^S)$.

We can construct $C^\bullet(K/K', \cV)$, with its Hecke action, by taking a derived limit of the complexes $\Hom_{\Z[K']}(C_\bullet,\cV/\varpi^r)$, where $C_\bullet$ denotes the complex of singular chains with $\Z$-coefficients on  $\overline{\frakX}_{\Gbar}$. Commuting the derived limit with cohomology of $K/K'$, using \cite[\href{https://stacks.math.columbia.edu/tag/08U1}{Tag 08U1}]{stacks-project}, we see that $R\Gamma(K/K',C^\bullet(K/K', \cV)) = C^\bullet(K,\cV)$. 

\begin{lem}
There are natural Hecke equivariant quasi-isomorphisms $A(K/K',\cV) \cong C^\bullet(K/K',\cV)$, where $A(K/K',\cV)$ are the complexes constructed in \cite[\S5.1]{AKT} using the singular chains of $\overline{\frakX}_{\Gbar}^\mathrm{dis}$.
\end{lem}
\begin{proof}We can reduce to the case where $\cV= \cV/\varpi^r$ and $K'$ is neat. Pullback by the continuous map $\overline{\frakX}_{\Gbar}^\mathrm{dis} \to \overline{\frakX}_{\Gbar}$ induces a Hecke equivariant map $A(K/K',\cV) \to C^\bullet(K/K',\cV)$ inducing the identity on the cohomology groups, which are identified with $H^*(\overline{X}^{\Gbar}_{K'},\cV)$ on both sides. 
\end{proof}

Since $K$ is not assumed to be neat, $C^\bullet(K,\cV)$ is not necessarily a perfect complex. However, using the Hochschild--Serre spectral sequence to compute its cohomology in terms of $C^\bullet(K',\cV)$ for $K' \subset K$ a neat open normal subgroup, we see that its cohomology groups are finitely generated $R$-modules \cite[Lemma 5.1]{AKT}.

A similar argument shows that the $\cO$-algebras \[\TT^S_{\Gbar}(C^\bullet(K/K',\cV)) \subset \End_{D^+(R[K/K'])}(C^\bullet(K/K',\cV))\] generated by the image of $\cH(\Gbar^S,K^S)$ are $\cO$-finite \cite[Lemma 5.2]{AKT}.

As a consequence, for a maximal ideal $\m \subset \TT^S_{\Gbar}(C^\bullet(K/K',\cV))$, we have a direct summand $C^\bullet(K/K',\cV)_{\m}$ cut out by an idempotent $e_\m$ in the Hecke algebra \cite[Prop.~3.6]{AKT}

\begin{prop}\label{prop:modpGal}
Suppose that $p$ is odd, and that $S$ satisfies the following conditions:
\begin{enumerate}
\item $S$ is stable under complex conjugation
\item $F$ contains an imaginary quadratic field. Let $v \notin S$ be a finite place, with residue characteristic $l$. Then either $S$ contains no $l$-adic places and $l$ is unramified in $F$, or there exists an imaginary quadratic subfield of $F$ in which $l$ splits.
\end{enumerate}
Then for any maximal ideal $\m \subset \TT^S_{\Gbar}(C^\bullet(K/K',\cV))$, there exists a semisimple continuous representation
\[\rhobar_\m: G_{F,S} \to \GL_2(\TT^S_{\Gbar}(C^\bullet(K/K',\cV))/\m)\]
 such that for each $v \notin S$,
\[\det(X - \rhobar_\m(\Frob_v)) = P_v(X) \mod \m.\]

Suppose moreover that $\rhobar_{\m}$ is absolutely irreducible. Then there exists an integer $N\ge 1$, depending only on $[F:\Q]$, an ideal $J \subset \TT^S_{\Gbar}(C^\bullet(K/K',\cV))$ with $J^N = 0$, and a continuous representation \[\rho_\m: G_{F,S}\to \GL_2(\TT^S_{\Gbar}(C^\bullet(K/K',\cV))_{\m}/J)\]
such that for each $v \notin S$,
\[\det(X - \rho_\m(\Frob_v)) = P_v(X) \mod J.\]

In particular, we have $\det \rho_{\m} = \epsilon_p^{-1}$.
\end{prop}
\begin{proof}
This is essentially \cite[Cor.~5.7, Cor.~5.9]{AKT} (the assumption there that $K/K'$ is abelian is not necessary). %
We let $K_G$ and $K'_G$ be the pre-images in $\GL_2(\widehat{\cO}_F)$ of $K$ and $K'$ respectively. Then \cite[Corollary 5.5]{AKT} identifies $\TT^S_{\Gbar}(H^\ast(C^\bullet(K/K',\cV)))$ as a quotient of a $\GL_2$-Hecke algebra $\TT^S_{G}(H^\ast(C^\bullet(K_G/K'_G,\cV)))$. We then reduce to the case where $K_G = K'_G$ is a neat level subgroup in $\GL_2$, as in the proofs of \cite[Theorems 5.6, 5.8]{AKT}, and finally appeal to \cite[Theorems 2.3.5, 2.3.7]{10author}. The determinant of $\rho_{\m}$ is $\epsilon_p^{-1}$ by Chebotarev density, considering the constant terms of the polynomials $P_v(X)$. 
\end{proof}

\begin{prop}\label{prop:aktperfect}
Let $\m$ be a maximal ideal of $\TT^S_{\Gbar}(C^\bullet(K,\cV))$ with residue field $k$. We make the following assumptions:
\begin{enumerate} \item $\cV \otimes k \cong k$ (with trivial action of $K_S$). \item $p$ is odd and $\rhobar_\m$ is absolutely irreducible.
\item $\zeta_p \in F$. 
\end{enumerate}
Then the cohomology groups $H^i(C^\bullet(K,\cV))_{\m}$ vanish for $i > \dim_{\R}X^{\Gbar}$. In particular, $C^\bullet(K,\cV)_{\m}$ is a perfect complex of $R$-modules.
\end{prop}
\begin{proof}
This follows from \cite[Thm.~5.11]{AKT}.
\end{proof}

\subsection{The proof of Theorem \ref{thm:pBT_lifting}}
We first prove a version of Theorem \ref{thm:pBT_lifting} with some additional assumptions. The general case will follow using solvable base change, as in \cite[\S6.5.12]{10author} and the proof of \cite[Theorem A.14]{AKT}.

We fix the following data:
\begin{enumerate}
\item\label{mlthyp1} An imaginary CM field $F$, an odd prime $p$ and an isomorphism $\iota: \Qpbar \cong \C$.
\item A finite set $S$ of finite places of $F$, including the places above $p$.
\item A (possibly empty) subset $R \subset S$ of places prime to $p$ and a decomposition $S_p = S_p^{cr}\coprod S_p^{st}$.
\item A cuspidal automorphic representation $\pi$ of $\PGL_2(\bA_F)$ which is regular algebraic of weight $0$. We may identify $\pi$ with a cuspidal automorphic representation of $\GL_2(\bA_F)$ with trivial central character.
\end{enumerate}
We assume the following conditions are satisfied:
\begin{enumerate}
\addtocounter{enumi}{4}
\item If $l$ is a prime lying below an element of $S$, or which is ramified in $F$, then $F$ contains an imaginary quadratic field in which $l$ splits. In particular, each place of $S$ is split over $F^+$ and the extension $F/F^+$ is everywhere unramified.
\item For each $v \in S_p$, let $\overline{v}$ denote the place of $F^+$ lying below $v$. Then there exists a place $\overline{v}' \neq \overline{v}$ of $F^+$ such that $\overline{v}' | p$ and 
\[ \sum_{\overline{v}'' \neq \overline{v}, \overline{v}'} [ F^+_{\overline{v}''} : \bQ_p ] > \frac{1}{2} [ F^+ : \bQ ]. \] Moreover, we assume that the residue field of $v$ is strictly bigger than $\F_p$.
\item $\pi_v$ is unramified for $v \notin R\cup S_p^{st}$.
\item If $v \in R \cup S_p^{st}$, then $\pi_v^{\mathrm{Iw}_v} \ne 0$.
\item\label{5.2assm-ord} If $v \in S_p^{st}$, then $\pi$ is $\iota$-ordinary of weight $0$ at $v$ %
and $r_\iota(\pi)|_{G_{F_v}}$ is non-crystalline ordinary.
\item If $S = S_p \cup R$, then $\zeta_p \in F$.
\item If $S \ne S_p\cup R$, then $S - (S_p \cup R)$ contains at least two places with distinct residue characteristics.
\item If $v \in S - (R\cup S_p)$, then $v \notin R^c$ and  $H^2(F_v,\ad^0\rbar_{\pi,\iota}) = 0$. 
\item $\rbar_{\pi,\iota}$ is decomposed generic and ${\rbar_{\pi,\iota}}|_{G_{F(\zeta_p)}}$ is irreducible. 
\item $\rbar_{\pi,\iota}|_{G_{F_v}}$ is the trivial representation for $v \in S_p \cup R$. (In particular, by considering determinants, $q_v \equiv 1 \mod p$ for $v \in R$.)
\item\label{mltlasthyp} If $p = 5$ and the projective image of $\rbar_{\pi,\iota}(G_{F(\zeta_5)})$ is conjugate to $\PSL_2(\F_5)$, we assume further that the extension of $F$ cut out by the projective image of $\rbar_{\pi,\iota}$ does not contain $\zeta_5$.
\end{enumerate}

\begin{prop}\label{prop:pBT_lifting_after_BC}
With notation and assumptions as in (\ref{mlthyp1})--(\ref{mltlasthyp}), suppose given a continuous representation $\rho : G_F \to \GL_2(\overline{\bQ}_p)$ satisfying the following conditions:
\begin{enumerate}
\item We have $\rhobar \cong \rbar_{\pi,\iota}$ and $\det(\rho) = \epsilon_p^{-1}$.
\item For each place $v \in S_p^{cr}$, $\rho|_{G_{F_v}}$ is Barsotti--Tate. 
\item\label{samecomp_after_BC} For each place $v \in S_p^{cr}$, $r_{\iota}(\pi)|_{G_{F_v}}$ is ordinary if and only if $\rho|_{G_{F_v}}$ is ordinary. 
\item\label{sameordcomp_after_BC} For each place $v \in S_p^{st}$, $\rho|_{G_{F_v}}$ is a non-crystalline extension of $\epsilon_p^{-1}$ by the trivial character.
\item For each finite place $v \notin S$ of $F$, $\rho|_{G_{F_v}}$ is unramified.
\item For each place $v \in R$, $\rho|_{G_{F_v}}$ is unipotently ramified.

\end{enumerate}

Then $\rho$ is automorphic: there exists a cuspidal automorphic representation $\Pi$ of $\PGL_2(\bA_F)$ of weight $0$ such that $\rho \cong r_\iota(\Pi)$. %
\end{prop}
\begin{proof}
We define a compact open subgroup $K = \prod_{v}K_v$ of $\PGL_2(\widehat{\cO}_F)$ as follows:
\begin{itemize}
\item If $v \notin S$ or $v \in S_p^{cr}$, then $K_v = \PGL_2(\cO_{F_v})$.
\item If $v \in R \cup S_p^{st}$, then $K_v = \Iw_v$.
\item If $v \in S - (S_p \cup R)$, then $K_v = \Iw_{v,1}$ is the pro-$v$ Iwahori subgroup of $\PGL_2(\cO_{F_v})$.
\end{itemize}
If $S - (S_p \cup R)$ is non-empty, $K$ is neat, by the same argument as \cite[Lemma 6.5.2]{10author}. We set $T = S\cup S^c$. 

Recalling the unitary group $\widetilde{G}$ from \S\ref{sec: unitary group} (for $n =2$), we need to define standard parabolic subgroups $Q_{\bar{v}}\subset \widetilde{G}_{F^+_{\bar{v}}}$ for each $\bar{v} \in \overline{S}_p$. For each $\bar{v}$ we choose a place of $F$, $\tilde{v}|\bar{v}$. We make this choice so that if at least one $v | \bar{v}$ is in $S_p^{st}$, then $\tilde{v}$ is in $S_p^{st}$. Then for each $\bar{v} \in \overline{S}_p$ we consider the following three cases:

\begin{itemize}
\item[(cr-ord)] $\tilde{v} \in S_p^{st}$ and $\tilde{v}^c \in S_p^{cr}$. Then $\iota_{\tilde{v}}(Q_{\bar{v}})$ is the standard parabolic given by the partition $(2,1,1)$.
\item[(ord)] $\tilde{v} \in S_p^{st}$ and $\tilde{v}^c \in S_p^{st}$. Then $\iota_{\tilde{v}}(Q_{\bar{v}}) = \Brm_{4}$, the Borel subgroup.
\item[(cr)] $\tilde{v} \in S_p^{cr}$ and $\tilde{v}^c \in S_p^{cr}$. Then $Q_{\bar{v}} = P_{\bar{v}}$, the Siegel parabolic.
\end{itemize}
We have a Hecke algebra $\TT^{\cQ_{\bar{S}_p},\bar{S}_p-\ord}_{\Gbar}(\cC^\bullet(K,\cO))$, defined as in \S\ref{subsec:degree shifting} by adding Hecke operators at places $v|p$ to $\TT^T(\cC^\bullet(K,\cO))$. 

Then we can find a coefficient field $E \subset \Qpbar$ and a maximal ideal \[\m \subset \TT^{\cQ_{\bar{S}_p},\bar{S}_p-\ord}_{\Gbar}(\cC^\bullet(K,\cO))\] such that $\rhobar_\m: G_{F,T} \to \GL_2(\Fpbar)$ satisfies $\rhobar_\m \cong \rbar_{\iota}(\pi)$ (cf.~\cite[Theorem 5.10]{AKT}). 

Moreover, for $v \in S_p^{st}$, since $\pi$ is $\iota$-ordinary at $v$, the Hecke operators $U_v:= [K_v \iota_v^{-1}\begin{pmatrix}\varpi_v & 0\\ 0 & 1 \end{pmatrix}K_v]$ are not in $\m$. 

Enlarging $E$ if necessary, we assume that the residue field of $\m$ is equal to $k$ and that $k$ contains all eigenvalues of the elements of $\rhobar_{\m}(G_F)$. 

We now describe the global deformation problems we will be working with. They will depend on a choice of character $\chi = \prod_{v \in R}\chi_v : \prod_{v \in R}\cO_{F_v}^\times \to \cO^\times$ which is trivial modulo $\varpi$.

For each $\chi$, we have the global deformation problem with fixed determinant \[\cS_{\chi} = (\rhobar,\epsilon_p^{-1}, S, \{R_v^{\epsilon_p^{-1},\BT}\}_{v \in S_p^{cr}}\cup\{R_v^{\triangle}\}_{v \in S_p^{st}} \cup \{R_v^{\epsilon_p^{-1},\chi_v}\}_{v \in R} \cup \{R_v^{\epsilon_p^{-1}}\}_{v \in S - (S_p\cup R)}).\]

The character $\chi_v: \cO_{F_v}^\times \to \cO^\times$ (which we note factors through $k_v^\times$) determines a character \begin{align*}\chi_v: \Iw_v & \to \cO^\times\\ \begin{pmatrix}
a & b\\ c & d 
\end{pmatrix} &\mapsto \chi_v(a/d).\end{align*}

We have an $\cO[K_S]$-module $\cO(\chi^{-1})$, where $K_S$ acts by the projection $K_S \to K_R = \prod_{v \in R}\Iw_v \xrightarrow{\prod_v \chi_v}\cO^\times$.

For each $\chi$, there is a canonical, surjective, $\cO$-algebra map
\[\TT^{\cQ_{\bar{S}_p},\bar{S}_p-\ord}_{\Gbar}(\cC^\bullet(K,\cO(\chi^{-1}))) \to \TT^{\cQ_{\bar{S}_p},\bar{S}_p-\ord}_{\Gbar}(\cC^\bullet(K,k))\] inducing a bijection on maximal ideals. So the maximal ideal $\m$ corresponds to a maximal ideal of $\TT^{\cQ_{\bar{S}_p},\bar{S}_p-\ord}_{\Gbar}(\cC^\bullet(K,\cO(\chi^{-1})))$ for each $\chi$. We abusively denote all these ideals by $\m$. The localization $\cC^\bullet(K,\cO(\chi^{-1}))_\m$ is a perfect complex of $\cO$-modules (using Proposition \ref{prop:aktperfect} when $K$ is not neat). We will also consider the Hecke algebra $\TT^{S,\cQ_{\bar{S}_p},\bar{S}_p-\ord}_{\Gbar}(\cC^\bullet(K,\cO(\chi^{-1}))_\m)$ obtained by adding in the spherical Hecke operators at places in $S^c - S$ to $\TT^{\cQ_{\bar{S}_p},\bar{S}_p-\ord}_{\Gbar}(\cC^\bullet(K,\cO(\chi^{-1}))_\m)$.

\begin{prop}\label{prop:lgcnoQ}
There exist an integer $N \ge 1$, depending only on $[F:\Q]$, an ideal $J \subset \TT^{S,\cQ_{\bar{S}_p},\bar{S}_p-\ord}_{\Gbar}(\cC^\bullet(K,\cO(\chi^{-1}))_\m)$ such that $J^N = 0$, and a continuous surjective homomorphism
\[f_{\cS_\chi}: R_{\cS_\chi} \to \TT^{S,\cQ_{\bar{S}_p},\bar{S}_p-\ord}_{\Gbar}(\cC^\bullet(K,\cO(\chi^{-1}))_\m)/J\] such that for each finite place $v \notin S$ of $F$, the characteristic polynomial of $f_{\cS_\chi}\circ \rho_{\cS_{\chi}}^{\mathrm{univ}}$ equals the image of $P_v(X)$ in $\TT^{S,\cQ_{\bar{S}_p},\bar{S}_p-\ord}_{\Gbar}(\cC^\bullet(K,\cO(\chi^{-1}))_\m)/J$.
\end{prop}
\begin{proof}
 Proposition \ref{prop:modpGal} already gives us a representation of $G_{F,S\cup S^c}$ with the right local properties at $v \notin S\cup S^c$, so it remains to check each prime $v \in S\cup S^c$. As in the proof of Proposition \ref{prop:modpGal}, we reduce to proving a similar local-global compatibility statement for a neat level in $\GL_2$. For $v \in S_p$, we apply Theorem \ref{thm:RtoT factors thru Kisin def ring}. For the remaining $v$, we proceed as in the proof of \cite[Proposition 6.5.3]{10author}, using \cite[Theorem 3.1.1]{10author}.
\end{proof}

To complete the proof of Proposition \ref{prop:pBT_lifting_after_BC} we need to show that the point of $\Spec(R_{\cS_1})$ given by $\rho$ is in the support of $H^*(C^\bullet(K,\cO)_\m)$. Then \cite[Theorem 5.10]{AKT} implies that $\rho$ is automorphic of weight $0$.

We need a local--global compatibility statement allowing ramification at Taylor--Wiles primes. So we suppose we have a Taylor--Wiles datum $(Q, N, (\alpha_{v,1},\alpha_{v,2})_{v \in Q})$ for $\cS_1$ (which is then also a Taylor--Wiles datum for every $\cS_{\chi}$). We assume that each place of $Q$ has residue characteristic split in an imaginary quadratic subfield of $F$. Now we define deformation problems
\[\cS_{\chi,Q} = (\rhobar,\epsilon_p^{-1}, S\cup Q, \{R_v^{\epsilon_p^{-1},\BT}\}_{v \in S_p^{cr}}\cup\{R_v^{\triangle}\}_{v \in S_p^{st}}\cup \{R_v^{\epsilon_p^{-1},\chi_v}\}_{v \in R} \cup \{R_v^{\epsilon_p^{-1}}\}_{v \in S\cup Q - (S_p\cup R)}).\]

For $v \in Q$, let $\Delta_v = k_v^\times(p)$, the maximal $p$-power quotient of $k_v^\times$. As in \cite[\S A.1.4]{AKT}, the local lifting ring $R_v^{\epsilon_p^{-1}}$ is equipped with the structure of an $\cO[\Delta_v]$-algebra. Setting $\Delta_Q = \prod_{v \in Q}\Delta_v$, we obtain an $\cO[\Delta_Q]$-algebra structure on $R_{\cS_{\chi,Q}}$. 

We define subgroups $K_1(Q) \subset K_0(Q) \subset K$, with $K_0(Q)/K_1(Q) \cong \Delta_Q$ as in \cite[\S6.5]{10author} (taking the image in $\PGL_2$ of the subgroups defined there). From this point, we follow loc.~cit. very closely, so we just explain the key points of the argument.  

There is a direct summand $C^\bullet(K_0(Q)/K_1(Q),\cO(\chi^{-1}))_{\frakn_1^Q}$ of $C^\bullet(K_0(Q)/K_1(Q),\cO(\chi^{-1}))$ in $\mathbf{D}(\cO[\Delta_Q])$, defined using a maximal ideal in a Hecke algebra with operators $U_{v,i}$ at places $v \in Q$. It is a perfect complex, by \cite[Theorem 5.11]{AKT}.

We write $\TT_{\chi,Q}$ for the image of the map \[\TT^{S\cup Q,\cQ_{\bar{S}_p},\bar{S}_p-\ord}_{\Gbar}\otimes_{\cO} \cO[\Delta_Q] \to \End_{\mathbf{D}(\cO[\Delta_Q])}\left(C^\bullet(K_0(Q)/K_1(Q),\cO(\chi^{-1}))_{\frakn_1^Q} \right).\] 

\begin{prop}\label{prop:levelQlgc}
There exist an integer $N \ge 1$, depending only on $[F:\Q]$, an ideal $J \subset \TT_{\chi,Q}$ such that $J^N = 0$, and a continuous surjective $\cO[\Delta_Q]$-algebra homomorphism 
\[f_{\cS_{\chi,Q}}: R_{\cS_{\chi,Q}} \to \TT_{\chi,Q}/J\] such that for each finite place $v \notin S\cup Q$ of $F$, the characteristic polynomial of $f_{\cS_{\chi,Q}}\circ \rho_{\cS_{\chi,Q}}^{\mathrm{univ}}$ equals the image of $P_v(X)$ in $\TT_{\chi,Q}/J$.
\end{prop}
\begin{proof}
   This is proved in the same way as Proposition \ref{prop:lgcnoQ}, using \cite[Theorem 3.1.1]{10author} to show that $f_{\cS_{\chi,Q}}$ is an $\cO[\Delta_Q]$-algebra homomorphism (cf.~\cite[Proposition 6.5.11]{10author} and \cite[Proposition A.13]{AKT}).
\end{proof}

It is convenient to patch complexes computing homology, so we define
\[C_{\chi,Q}:= R\Hom_{\cO[\Delta_Q]}(C^\bullet(K_0(Q)/K_1(Q),\cO(\chi^{-1}))_{\frakn_1^Q},\cO[\Delta_Q])\] and 
\[C_{\chi}:= R\Hom_{\cO}(C^\bullet(K,\cO(\chi^{-1}))_{\frakm},\cO).\]

\begin{lem}
$C_{\chi,Q}$ is a perfect complex of $\cO[\Delta_Q]$-modules, with a canonical isomorphism \[C_{\chi,Q}\otimes^{\LL}_{\cO[\Delta_Q]}\cO \cong C_{\chi}\] in $\bD(\cO)$.
\end{lem}
\begin{proof}
This follows from the fact that we can identify $C_{\chi,Q}$ and $C_{\chi}$ with the duals of perfect complexes computing localizations (at $\frakn_1^Q$ and $\m$ respectively) of equivariant homology (see the second part of \cite[Theorem 5.11]{AKT}).
\end{proof}

We now have everything we need to construct the objects required for \S\ref{ssec:patching}, using \cite[\S6.4]{10author} and \cite[Proposition A.6]{AKT} (existence of Taylor--Wiles primes) as in the proof of \cite[Theorem A.7]{AKT}. In particular, we make use of two options for the tuple of characters $\chi$: firstly, $\chi = 1$, and secondly a fixed tuple, denoted $\chi$, given by a choice of character $\chi_v: \cO_{F_v}^\times \to \cO^\times$ with $\chi_v^2 \neq 1$ which is trivial mod $\varpi$ for each $v \in R$. 

The rings $R_\infty$ and $R'_\infty$ are power series rings over $R_\loc = R_{\cS_1}^{S,\loc}$ and $R'_\loc = R_{\cS_\chi}^{S,\loc}$ respectively, which come equipped with local $R_\loc$-algebra (respectively $R'_\loc$-algebra) surjections $R_\infty \to R_{\cS_1}$ (respectively $R'_\infty \to R_{\cS_\chi}$). We can assume, extending $\cO$ if necessary, that all of the irreducible components of the local lifting rings appearing in the deformation problems $\cS_1$ and $\cS_{\chi}$ (and their special fibres) are geometrically irreducible. We do this so that we can apply \cite[Lemma 3.3]{blght}.

It follows from formal smoothness of $R_v^{\square,\epsilon_p^{-1}}$ for $v \in S - (R\cup S_p)$, \cite[Lemma A.2]{AKT}, Proposition \ref{prop:ordcomps}, Lemma \ref{lem:cegs} and \cite[Lemma 3.3]{blght} (which describes irreducible components of completed tensor products in terms of their factors) that $R_\infty$ and $R'_\infty$ satisfy Assumption \ref{assumptionsetup}. In particular, $\Spec(R'_\infty)$ has $2^{|S_p|}$ irreducible components, which biject with the irreducible components of $R_{\cS_\chi}^{S_p,\loc} = R_{\cS_1}^{S_p,\loc}$. We let $C_a$ be an irreducible component of $\Spec(R_\infty)$ containing the point $r_{\pi,\iota}$ with $C_a \subset \Spec(T_\infty)$. It exists by Proposition \ref{prop:automorphic_component}. We let $C$ be an irreducible component of $\Spec(R_\infty)$ containing the point given by $\rho$. It follows from conditions (\ref{samecomp_after_BC}) and (\ref{sameordcomp_after_BC}), Proposition \ref{prop:ordcomps}, Lemma \ref{lem:BTcomp} and \cite[Lemma 3.3]{blght} that the generic points of $C\cap \Spec(R_\infty/\varpi)$ and $C_a\cap \Spec(R_\infty/\varpi)$ lie in the same irreducible component of $\Spec(R'_\infty)$. Note that condition (\ref{sameordcomp_after_BC}) and assumption (\ref{5.2assm-ord}) ensure that both $C$ and $C_a$ lie over the component $X^{st}$ in $\Spec(R_v^\triangle)$ for $v \in S_p^{st}$. 

We will therefore be able to apply Corollary \ref{cor:char0automorphy} to deduce automorphy of $\rho$. 

This completes the proof of Proposition \ref{prop:pBT_lifting_after_BC}.
\end{proof}

\begin{proof}[The end of the proof of Theorem \ref{thm:pBT_lifting}]
We will use a special case of \cite[Lemma 4.11]{DDT}:
\begin{lem}\label{lem:DDT4.11}
Suppose that $G$ is a finite group with a representation $\rhobar: G \to \GL_2(\Fpbar)$ for an odd prime $p$.  Suppose that the character $\det\rhobar: G \to \Fpbarx$ has order $d > 1$, and for all $g$ with $\det \rhobar(g) \ne 1$ we have \begin{equation}\label{eqn:DDTtrace}(\mathrm{tr} \rhobar(g))^2 = (1+\det\rhobar(g))^2.\end{equation}
Then either $\rhobar|_{\ker(\det\rhobar)}$ is reducible, or $d = 3$ and $\rhobar$ has projective image isomorphic to $A_4$.
\end{lem}
\begin{proof}
	This is an immediate consequence of \cite[Lemma 4.11]{DDT}, applied with $\chi = \det\rhobar$. 
\end{proof}

It suffices to prove that $\rho \cong r_{\iota}(\Pi)$ for a cuspidal automorphic representation $\Pi$ of $\GL_2(\AA_F)$; then the fact that $\det(\rho) = \epsilon_p^{-1}$ implies that $\Pi$ has trivial central character. Let $L/F(\zeta_p)$ be the extension cut out by $\rhobar|_{G_{F(\zeta_p)}}$. If $F'/F$ is any finite solvable extension, we denote the base change of $\pi$ to $F'$ by $\pi_{F'}$.  

We choose a finite set $V$ of finite places of $F$ exactly as in the proof of \cite[Theorem A.14]{AKT}, so that: \begin{itemize}\item For any proper extension $L'/F$ contained in $L$, there is some $v \in V$ not splitting in $L'$. \item There is a rational prime $q \ne p$ such that $\rhobar$ is decomposed generic for $q$ and $V$ contains all $q$-adic places of $F$. 
	\item For each $v \in V$, $v \nmid 2p$ and both $\rho$ and $\pi$ are unramified at $v$. 
\end{itemize}

This ensures that if $F'/F$ is a finite Galois extension in which every place of $V$ splits, then $\rhobar|_{G_{F'}}$ remains decomposed generic, $\rhobar(G_{F'}) = \rhobar(G_{F})$ and $\rhobar(G_{F'(\zeta_p)}) = \rhobar(G_{F(\zeta_p)})$.

Now we choose a solvable, Galois, CM extension $F_0/F$ such that: \begin{itemize}
	\item Every place of $V$ splits in $F_0$.
	\item For every finite place $w$ of $F_0$, $\pi_{F_0}^{\Iw_w} \ne 0$. 
	\item For every finite place $w \nmid p$ of $F_0$, either both $\pi_{F_{0,w}}$ and $\rho|_{G_{F_{0,w}}}$ are unramified, or $\rho|_{G_{F_{0,w}}}$ is unipotently ramified, $q_w \equiv 1$ mod $p$, and $\rhobar|_{G_{F_0,w}}$ is trivial.
	\item For each $\bar{w}|p$ in $F_0^+$, $\bar{w}$ splits in $F_0$,  $\rhobar|_{G_{F_{0,w}}}$ is trivial for $w|\bar{w}$, the residue field $k_w$ is strictly bigger than $\F_p$ and there exists a place $\overline{w}' \neq \overline{w}$ of $F^+_0$ such that $\overline{w}' | p$ and 
	\[ \sum_{\overline{w}'' \neq \overline{w}, \overline{w}'} [ F^+_{0,\overline{w}''} : \bQ_p ] > \frac{1}{2} [ F^+_0 : \bQ ].
	 \]
	 \item If $w$ lies over a place $v$ of $F$ with $\rho|_{G_{F_v}}$ potentially crystalline, then $\rho|_{G_{F_{0,w}}}$ is crystalline and $\pi_{F_{0,w}}$ is unramified. Moreover, $r_{\iota}(\pi)|_{G_{F_{0,w}}}$ is crystalline and it is ordinary if and only if  $\rho|_{G_{F_{0,w}}}$ is ordinary.
	 \item If $w$ lies over a place $v$ of $F$ with $\rho|_{G_{F_v}}$ not potentially crystalline, $\rho|_{G_{F_{0,w}}}$ is a non-crystalline extension of $\epsilon_p^{-1}$ by the trivial character. 
\end{itemize}
With respect to the penultimate item, we note that it follows from Theorem \ref{thm:LGC in char 0} that when  $\pi_{F_{0,w}}$ and $\pi_{F_{0,w^c}}$ are unramified,  $r_{\iota}(\pi)|_{G_{F_{0,w}}}$ is automatically crystalline with all labelled Hodge--Tate weights equal to $(0,1)$. If $\pi_{F_{0,w}}$ is unramified but $\pi_{F_{0,w^c}}$ is not, then it follows from assumptions \eqref{assm:samecomp} and \eqref{assm:sameordcomp} that $\rho|_{G_{F_{0,w^c}}}$ is not potentially crystalline and that $\pi_{F_{0,w^c}}$ is $\iota$-ordinary of weight $0$. We can then deduce from case (cr-ord) of Theorem \ref{thm:RtoT factors thru Kisin def ring} (as in the proof of Theorem \ref{thm:LGC in char 0}) that $r_{\iota}(\pi)|_{G_{F_0,w}}$ is crystalline with all labelled Hodge--Tate weights equal to $(0,1)$. 

Making a further solvable extension $F_1/F_0$ by taking a composite with three imaginary quadratic fields (again ensuring that every place of $V$ splits in $F_1$), as in the proof of \cite[Theorem A.14]{AKT}, we furthermore satisfy:

\begin{itemize}
	\item Let $R$ be the set of finite places $w\nmid p$ of $F_1$ such that $\pi_{F_1,w}$ or $\rho|_{G_{F_{1,w}}}$ are ramified. Let $S_p$ denote the $p$-adic places of $F_1$ and set $S' = S_p \cup R$. If $l$ is a rational prime lying below an element of $S'$, or which is ramified in $F_1$, then $F_1$ contains an imaginary quadratic field in which $l$ splits. 
\end{itemize}

We can now describe the data we need to apply Proposition \ref{prop:pBT_lifting_after_BC} with $F= F_1$, $\rho = \rho|_{G_{F_1}}$ and $\pi = \pi_{F_1}$. We have already defined the set of places $R$. We let $S_p^{cr}$ be the set of places $w|p$ where $\rho|_{G_{F_{1,w}}}$ is crystalline and let $S_p^{st}$ be the set of places $w|p$ where $\rho|_{G_{F_{1,w}}}$ is non-crystalline. For $w \in S_p^{st}$, we know (by assumption) that $r_{\iota}(\pi)|_{G_{F_{1,w}}}$ is not crystalline, whilst it follows from Theorem \ref{thm:RtoT factors thru Kisin def ring} that it is ordinary.

If $\zeta_p \in F$, we set $S= S' = S_p \cup R$. If $\zeta_p\notin F$ (which entails $\zeta_p \notin F_1$), Lemma \ref{lem:DDT4.11} shows that we can find an element $g \in \rhobar(G_{F_1})$ such that $\det(g) \ne 1$, and the ratio of the eigenvalues of $g$ does not equal $\det(g)^{\pm 1}$. Using Chebotarev density, we can find infinitely many finite places $v_0$ of $F_1$ of degree $1$ over $\Q$ such that $v_0 \notin S'\cup R^c$, $q_{v_0} \not\equiv 1$ mod $p$, and the ratio of the eigenvalues of $\rhobar(\Frob_{v_0})$ does not equal $q_{v_0}^{\pm 1}$.
For such a place, $H^2(F_{1,v_0},\ad^0 \rhobar) = 0$ and the rational prime below $v_0$ splits in any quadratic subfield of $F_1$. We choose two such places $v_0, v_0'$ of distinct residue characteristic and set $S= S'\cup\{v_0,v_0'\}$. We are now in a situation where all the assumptions of Proposition \ref{prop:pBT_lifting_after_BC} are satisfied. We deduce that $\rho|_{G_{F_1}}$ is automorphic, and solvable descent \cite[Proposition 6.5.13]{10author} completes the job. \end{proof}

\section{Modularity of elliptic curves over CM fields}
In this section, our goal is to combine Theorem \ref{thm:pBT_lifting} with the results of \cite{AKT} to prove the following:
\begin{theorem}\label{thm:modularity35}
Let $F$ be an imaginary CM number field with $\zeta_5 \notin F$. Let $E/F$ be an elliptic curve satisfying one of the following two conditions:
\begin{enumerate}
\item $\rbar_{E,3}$ is decomposed generic and $\rbar_{E,3}|_{G_{F(\zeta_3)}}$ is absolutely irreducible.
\item $\rbar_{E,5}$ is decomposed generic and $\rbar_{E,5}|_{G_{F(\zeta_5)}}$ is absolutely irreducible.
\end{enumerate}
Then $E$ is modular.
\end{theorem}

By `$E$ is modular', we mean that either $E$ has CM, or there is a cuspidal automorphic representation $\pi$ of $\GL_2(\A_F)$ which is regular algebraic of weight 0, with $r_{\pi,\iota} \cong r_{E,p}^\vee$ (the dual of the Tate module) for a prime $p$ and an isomorphism $\iota: \Qpbar \to \C$. 

Before proving the theorem, we give some corollaries. These will be improved further in the next section in the special case when $F$ is imaginary quadratic.
\begin{cor}\label{cor:modularity IQF}
Let $F$ be an imaginary quadratic field. Let $E/F$ be an elliptic curve satisfying one of the following two conditions:
\begin{enumerate}
	\item $\rbar_{E,3}|_{G_{F(\zeta_3)}}$ is absolutely irreducible.
	\item $\rbar_{E,5}|_{G_{F(\zeta_5)}}$ is absolutely irreducible.
\end{enumerate}
Then $E$ is modular.
\end{cor}
\begin{proof}
	Combine Theorem \ref{thm:modularity35} and Lemma \ref{lem:quadratic irred implies generic}.
\end{proof}

\begin{cor}\label{cor:100 percent} 
Let $F$ be an imaginary CM field that is Galois over $\mathbb{Q}$ and such that $\zeta_5\notin F$. 
Then 100\% of Weierstrass equations over $F$, ordered by their height, define a modular elliptic curve. 
\end{cor}

\begin{proof} 
It follows from~\cite[Lemma 2.3]{allen-newton} that, if $F/\mathbb{Q}$ is finite Galois and $E/F$ is an elliptic curve such 
that the image of $\rbar_{E,5}$ contains $\mathrm{SL}_2(\F_5)$, 
then $\rbar_{E,5}$ is decomposed generic. Note that the decomposed generic condition used in \emph{loc.~cit.} 
is more restrictive than and therefore implies the one we are using. When 
the image of $\rbar_{E,5}$ contains $\mathrm{SL}_2(\F_5)$, we also have that 
$\rbar_{E,5}|_{G_{F(\zeta_5)}}$ is absolutely irreducible, so the hypotheses of the second part of Theorem~\ref{thm:modularity35}
are satisfied.

To conclude, we observe that a quantitative version of Hilbert irreducibility, see for example~\cite[Prop. 5.2]{zywina-maximal}, implies 
that 100\% of elliptic curves 
$E$ over a fixed number field $F$ have the property that the image of $\rbar_{E,5}$ 
contains $\mathrm{SL}_2(\F_5)$. 
\end{proof}

We recall a useful lemma from \cite{AKT}, which is proved using Varma's results on local-global compatibility \cite{ilavarma}.

\begin{lem}\label{lem:eclgc}
Let $F$ be a CM field and let $E/F$ be a modular elliptic curve without CM, with $r_{\pi,\iota} \cong r_{E,p}^\vee$ for some choice of prime $p$ and isomorphism $\iota: \Qpbar \to \C$. Then:
\begin{enumerate}
\item $\pi$ has trivial central character and weight $0$, and is uniquely determined by $E$.
\item For every prime $p$ and isomorphism $\iota:\Qpbar \to \C$, there is an isomorphism $r_{\pi,\iota} \cong r_{E,p}^\vee$.
\item For every isomorphism $\iota:\Qpbar \to \C$ and finite place $v\nmid p$ of $F$, there is an isomorphism $WD(r_{E,p}^\vee|_{G_{F_v}})^{F-ss} \cong \rec_{F_v}^T(\pi_v)$. 
\item Suppose $v|p$ is a place where $E$ has potentially multiplicative reduction. Then $\pi$ is $\iota$-ordinary of weight $0$ at $v$ for any $\iota:\Qpbar \to \C$.
\end{enumerate}
\end{lem}
\begin{proof}
The first three parts are contained in \cite[Lemma 9.1]{AKT}. The final part is proved in the same way as \cite[Corollary 9.2]{AKT}: applying the third part to $r_{E,l}^\vee$ for some $l \ne p$, we see that $\pi_v$ is a twist of the Steinberg representation by a quadratic character (quadratic since the central character of $\pi_v$ is trivial). Then the proof of \cite[Lemma 5.6]{geraghty} shows that $\pi_v$ is $\iota$-ordinary at $v$.
\end{proof}

Here is another useful lemma, taken from the proof of \cite[Corollary 9.14]{AKT}.
\begin{lem}\label{lem:condition at 5}
Let $F$ be a CM field with $\zeta_5 \notin F$ and let $E/F$ be an elliptic curve such that the projective image of $\rbar_{E,5}(G_{F(\zeta_5)})$ is conjugate to $\PSL_2(\F_5)$. Then the extension of $F$ cut out by the projective image of $\rbar_{E,5}(G_F)$ does not contain $\zeta_5$.
\end{lem}
\begin{proof}
The group $\rbar_{E,5}(G_{F(\zeta_5)})$ is a subgroup of $\SL_2(\F_5)$ surjecting onto $\PSL_2(\F_5)$, so it is equal to $\SL_2(\F_5)$. We let $G$ be the kernel of the map from $G_F$ to the projective image of $\rbar_{E,5}(G_F)$ and let $H \subset G$ be the kernel of the map from $G_{F(\zeta_5)}$ to the projective image of $\rbar_{E,5}(G_{F(\zeta_5)})$. The extension of $F$ cut out by the projective image of $\rbar_{E,5}(G_F)$ contains $\zeta_5$ if and only if $H = G$. Since $\zeta_5 \notin F$, $\det\rhobar$ has order $2$ or $4$. In the first case, the projective image of $\rbar_{E,5}(G_F)$ is again $\PSL_2(\F_5)$ so $[G:H] = [F(\zeta_5):F] > 1$. In the second case, $[F(\zeta_5):F] > [\PGL_2(\F_5):\PSL_2(\F_5)]$ which again means $H \ne G$.
\end{proof}

Now we state a variant of \cite[Proposition 9.13]{AKT}:
\begin{prop}\label{prop:potmod5}
Let $F$ be an imaginary CM number field, and let \[\rhobar: G_F \to \GL_2(\FF_5)\] be a continuous homomorphism with determinant $\bar{\epsilon}_5$. We assume $\rhobar$ is decomposed generic.

Suppose we have a decomposition $S_5 = S_5^{st}\coprod S_5^{\ord} \coprod S_5^{\ss}$ of the set of places in $F$ dividing $5$.

Let $F^\avoid/F$ be a finite Galois extension. Then we can find a solvable CM extension $L/F$ and an elliptic curve $E/L$ satisfying the following conditions:
\begin{enumerate}
\item $E$ is modular and without CM.
\item The extension $L/F$ is linearly disjoint from $F^\avoid/F$
\item For each place $v|5$ in $F$ and $w|v$ in $L$, $E_{L_w}$ has good ordinary reduction if $v \in S_5^{\ord}$, good supersingular reduction if $v \in S_5^{\ss}$ and (split) multiplicative reduction if $v \in S_5^{st}$.
\item There is an isomorphism $\rhobar|_{G_L} \cong \rbar_{E,5}$.
\item $\rhobar|_{G_L}$ is decomposed generic.
\end{enumerate}
\end{prop}
\begin{proof}
We choose $L/F$ to be a solvable CM extension such that:\begin{itemize}
\item For each place $w|2,3,5$ of $L$, $\rhobar|_{G_{L_w}}$ is trivial and $w$ is split over $L^+$.
\item For $w|5$, there are elliptic curves $E^{\ord}/L_w$, $E^{\ss}/L_w$ with good ordinary and good supersingular reduction respectively, and trivial action of $G_{L_w}$ on their $5$-torsion. 
\item For each place $w|2$ of $L$, the extension $L_w(\sqrt{-1})/L_w$ is unramified.
\item $L/F$ is linearly disjoint from $F^{\avoid}/F$.
\item There is a prime $q > 5$ which is decomposed generic for $\rhobar$ and splits in $L$.
\end{itemize}

Now we apply \cite[Lemma 9.7]{AKT} in the same way as in the proof of \cite[Proposition 9.13]{AKT} to find an $L$-rational point of the modular curve $Y_{\rhobar}$ corresponding to a modular elliptic curve $E/L$. The curve $Y_{\rhobar}$ is isomorphic to an open subset of the projective line over $F$. Combining Hilbert irreducibility, in the form of \cite[Theorem 3.5.3]{serre-topics}, and weak approximation for the projective line, we can find points in $Y_{\rhobar}(L)$ that avoid a thin subset (\cite[Definition 3.1.1]{serre-topics}) and lie in specified non-empty $w$-adically open subsets $\Omega_w \subset Y_{\rhobar}(L_w)$ for a finite set of places $w$. Compared to \cite[Proposition 9.13]{AKT}, we replace the condition that $E_{L_w}$ is a Tate curve for each place $w|5$ with the condition that for $w$ lying over a place in $S_5^{\ord}$, $E_{L_w}$ has good ordinary reduction, for $w$ lying over a place in $S_5^{\ss}$, $E_{L_w}$ has good supersingular reduction and for $w$ lying over a place in $S_5^{st}$, $E_{L_w}$ is a Tate curve. Deleting a thin subset ensures that $\rbar_{E,3}(G_{L(\zeta_3)}) = \SL_2(\F_3)$, so $E$ does not have CM (alternatively, it is also part of the construction that $E$ is a Tate curve locally at $2$- and $3$-adic places). Note that the modularity of $E$ is then proved by applying \cite[Proposition 9.12]{AKT}, which is not sensitive to the $5$-adic properties of $E$. 
\end{proof}

We have a similar statement for mod $3$ representations, which can be proved in the same way as \cite[Proposition 9.15]{AKT}:
\begin{prop}\label{prop:potmod3}
Let $F$ be an imaginary CM number field with $\zeta_5 \notin F$, and let \[\rhobar: G_F \to \GL_2(\FF_3)\] be a continuous homomorphism with determinant $\bar{\epsilon}_3$. We assume $\rhobar$ is decomposed generic.

Suppose we have a decomposition $S_3 = S_3^{st}\coprod S_3^{\ord} \coprod S_3^{\ss}$ of the set of places in $F$ dividing $3$.

Let $F^\avoid/F$ be a finite Galois extension. Then we can find a solvable CM extension $L/F$ and an elliptic curve $E/L$ satisfying the following conditions:
\begin{enumerate}
\item $E$ is modular and without CM.
\item The extension $L/F$ is linearly disjoint from $F^\avoid/F$
\item For each place $v|3$ in $F$ and $w|v$ in $L$, $E_{L_w}$ has good ordinary reduction if $v \in S_3^{\ord}$, good supersingular reduction if $v \in S_3^{\ss}$ and (split) multiplicative reduction if $v \in S_3^{st}$.
\item There is an isomorphism $\rhobar|_{G_L} \cong \rbar_{E,3}$.
\item $\rhobar|_{G_L}$ is decomposed generic.
\end{enumerate}
\end{prop}

\begin{lem}\label{lem:checklocalprops}
Let $E$ be the modular elliptic curve produced by Proposition \ref{prop:potmod5} or Proposition \ref{prop:potmod3} with $r_{E,p}^\vee \cong r_\iota(\pi)$ for $p = 5$ or $3$ respectively. Let $w|p$ be a place of $L$, lying over the place $v$ of $F$. Then $\pi$ is $\iota$-ordinary at $w$ if $v \in S_p^{st}$. The local factor $\pi_w$ is unramified if $v \in S_p^{\ord}\coprod S_p^{\ss}$.
\end{lem}
\begin{proof}
This follows from Lemma \ref{lem:eclgc}.
\end{proof}

\begin{proof}[Proof of Theorem \ref{thm:modularity35}]

Choose $p \in \{3,5\}$ so that $\rbar_{E,p}$ is decomposed generic and $\rbar_{E,p}|_{G_{F(\zeta_p)}}$ is absolutely irreducible. Now we apply
Proposition \ref{prop:potmod5} or \ref{prop:potmod3} with $\rhobar = \rbar_{E,p}$, $S_p^{\ord}$ the set of places above $p$ where $E$ has potentially good ordinary reduction,  $S_p^{\ss}$ the set of places where $E$ has potentially good supersingular reduction and $S_p^{st}$ the set of places where $E$ has potentially multiplicative reduction. We set $F^{\avoid} = F(E[p],\zeta_3,\zeta_5)$. The appropriate proposition gives us a solvable extension $L/F$ and a modular elliptic curve $A/L$ with $\rbar_{E,p}|_{G_L} \cong \rbar_{A,p}$, such that the hypotheses of Theorem \ref{thm:pBT_lifting} apply to $\rho = r_{E,p}^\vee|_{G_L}$ (we use Lemma \ref{lem:checklocalprops} here). The assumption that $\zeta_5 \notin F$ is sufficient to check the condition on the projective image of $\rhobar$ in the theorem, by Lemma \ref{lem:condition at 5}. We deduce that $E_L$ is modular, and the modularity of $E$ follows by solvable descent (cf.~\cite[Prop.~6.5.13]{10author}). 
\end{proof}

\subsection{Group theory}
We now do a little bit of group theory to optimize the statement of Theorem \ref{thm:modularity35} when $F$ is quadratic (cf.~Corollary \ref{cor:modularity IQF}). Our main tool will be the following well-known lemma: 
\begin{lem}[Goursat's lemma]
Let $G_1, G_2$ be finite groups and suppose $H \subset G_1 \times G_2$ is a subgroup with the projection map $p_i: H \to G_i$ surjective for $i = 1$ and $2$. Then we have normal subgroups $N_i = H \cap G_{i} \subset G_i$ and an isomorphism $\phi: G_1/N_1 \cong G_2/N_2$ such that 
\[H = \{(g_1,g_2):\phi(g_1N_1) = g_2N_2\}.\]
\end{lem}

\begin{lem}\label{lem:quadratic irred implies generic}
Let $F/\Q$ be a quadratic field, let $p$ be an odd prime, and let $\rhobar: G_F \to \GL_2(\F_p)$ be a homomorphism. Suppose that $\rhobar|_{G_{F(\zeta_p)}}$ is absolutely irreducible. Then $\rhobar$ is decomposed generic. Moreover, the restriction $\rhobar|_{G_{F(\zeta_p)}}$ is decomposed generic.
\end{lem}
\begin{proof}
We consider the homomorphism $P = \mathrm{Proj}(\rhobar): G_F \to \PGL_2(\F_p)$. Set $L = \overline{F}^{\ker P}$, so $P$ factors through an embedding $P: \Gal(L/F) \hookrightarrow \PGL_2(\F_p)$. We also set $L_1 = \overline{F}^{\ker P \cap G_{F(\zeta_p)}} = L(\zeta_p)$.

Let $\widetilde{L}$ be the Galois closure of $L$ over $\Q$ in $\overline{F}$. Fixing a lift $c \in \Gal(\Qbar/\Q)$ of the non-trivial element in $\Gal(F/\Q)$, we have an injective map:
\begin{align*}
\Gal(\widetilde{L}/F) &\hookrightarrow \Gal(L/F) \times \Gal(L/F)\\
\sigma &\mapsto (\sigma|_L, (c^{-1}\sigma c)|_L)
\end{align*}
whose composition with each of the two projection maps to $\Gal(L/F)$ is surjective. Injectivity follows from the fact that $\widetilde{L}$ is the composite of $L$ and $c(L)$.

The Galois closure $\widetilde{L_1}$ of $L_1$ over $\Q$ is $\widetilde{L}(\zeta_p)$, and restricting the above map to the subgroup $\Gal(\widetilde{L_1}/F(\zeta_p)) = \Gal(\widetilde{L}/\widetilde{L}\cap F(\zeta_p))$ gives an injective map \[\Gal(\widetilde{L_1}/F(\zeta_p)) \hookrightarrow \Gal(L_1/F(\zeta_p)) \times \Gal(L_1/F(\zeta_p))\]
whose composition with each of the two projection maps to $\Gal(L_1/F(\zeta_p))$ is surjective. 

We are going to show that there is an element $\tau \in \Gal(\widetilde{L_1}/F(\zeta_p))$ whose image under each projection map to $\Gal(L_1/F(\zeta_p))$ is a non-identity element of order prime to $p$. First we need to show that $\Gal(L_1/F(\zeta_p))$ itself contains a non-identity element of order prime to $p$. This group contains the image of $\rhobar(G_F(\zeta_p))$ in $\PGL_2(\Fp)$. The irreducibility of $\rhobar|_{G_{F(\zeta_p)}}$ implies that $\rhobar(G_F(\zeta_p))$ either has order prime to $p$ or contains $\SL_2(\Fp)$ (by Dickson's classification, or, more simply, \cite[Proposition 15]{serre-elliptic}). It follows that we can find a non-identity element $T \in \Gal(L_1/F(\zeta_p))$ of order prime to $p$. 

We now denote $\Gal(\widetilde{L_1}/F(\zeta_p))$ by $H$ and $\Gal({L_1}/F(\zeta_p))$ by $G$. Goursat's lemma tells us that there are normal subgroups $N_1, N_2 \triangleleft G$ and an isomorphism $\phi: G/N_1 \cong G/N_2$ such that $H = \{(g_1,g_2): \phi(g_1N_1) = g_2N_2\}$. Note that $N_1$ and $N_2$ are necessarily of the same order. We separate into two cases:
\begin{itemize}
\item The $N_i$ are $p$-groups (we include the possibility that the $N_i$ are trivial). In this case, we fix a lift $T' \in G$ of $\phi(TN_1) \in G/N_2$. The order of $\phi(TN_1)$ is equal to the order of $T$, and is equal to the order of $T'$ up to a $p$-power factor. So, replacing $T$ and $T'$ by a sufficiently large $p$th power if necessary, we have $(T,T') \in H$ with $T$ and $T'$ non-identity elements with order prime to $p$.
\item The $N_i$ are not $p$-groups. Then we can let $T_1$ be a non-identity element of $N_1$ of order prime to $p$ and $T_2$ a non-identity element of $N_2$ of order prime to $p$. The element $(T_1,T_2)$ is contained in $H$.
\end{itemize}

We have now constructed the desired element $\tau \in \Gal(\widetilde{L_1}/F(\zeta_p))$. By Chebotarev density, we can choose a rational prime $l$, unramified in $\widetilde{L_1}$, such that $\Frob_l$ is the conjugacy class of $\tau$ in $\Gal(\widetilde{L_1}/\Q)$. We can moreover ensure that $\rhobar$ is unramified at every place dividing $l$. Since $\tau$ fixes $F(\zeta_p)$, we have $l \equiv 1$ mod $p$ and $l$ splits completely in $F(\zeta_p)$. The Frobenius elements $\Frob_v$ for $v|l$ in $F$ are given by the $\Gal(\widetilde{L_1}/F)$-conjugacy classes contained in $\Frob_l$. These are the $\Gal(\widetilde{L_1}/F)$-conjugacy classes of $\tau$ and $c^{-1}\tau c$ (which could coincide). 

By construction, for each $v|l$, $\rhobar(\Frob_v)$ and $\rhobar(c^{-1}\Frob_v c)$ have image in $\PGL_2(\F_p)$ a non-identity element of order prime to $p$. It follows that they are both regular semisimple elements of $\GL_2(\Fp)$. Since $l \equiv 1$ mod $p$, we have shown that $l$ is a decomposed generic prime for $\rhobar$. Moreover, if $v|l$ in $F$, it splits completely in $F(\zeta_p)$ and for each place $w|v$ of $F(\zeta_p)$, $\rhobar(\Frob_w)$ is in the same $G_F$ conjugacy class as $\rhobar(\Frob_v)$. Hence $\rhobar|_{G_{F(\zeta_p)}}$ is also decomposed generic.  
\end{proof}

\section{Quadratic points on modular curves and modularity over quadratic fields}
In this section, inspired by the proof of modularity of elliptic curves over real quadratic fields \cite{flhs}, our goal is to extend Corollary \ref{cor:modularity IQF} to cover many of the excluded cases where $r_{E,p}|_{G_{F(\zeta_p)}}$ is absolutely reducible for $p = 3$ and $5$. We first state the main theorem of this section.

\begin{theorem}\label{thm:optimised modularity IQF}
Let $F$ be an imaginary quadratic field, and let $E/F$ be an elliptic curve such that one of the following conditions holds:
\begin{enumerate} \item The action of $G_F$ on $E[5]$ is irreducible (not necessarily absolutely irreducible).
\item The action of $G_F$ on $E[3]$ is irreducible and the image of $G_F$ in $\Aut(E[3])$ is not the normalizer of a split Cartan subgroup.
\end{enumerate} 
Then $E$ is modular.
\end{theorem}

The following lemma helps explicate the condition that $r_{E,p}|_{G_{F(\zeta_p)}}$ is absolutely reducible for $p = 3$ and $5$.
\begin{lem}\label{lem:galois image}
	Let $F$ be a number field and $\rhobar: G_F \to \GL_2(\F_p)$ an irreducible representation with determinant $\epsilonbar_p$ and $\rhobar|_{G_{F(\zeta_p)}}$ absolutely reducible. \begin{enumerate}
		\item $\rhobar(G_F)$ is a subgroup of the normalizer of a Cartan subgroup.
		\item If $p = 3$, $\rhobar(G_F)$ is conjugate to $C_{\mathrm{s}}^+(3)$ (the normalizer of a split Cartan subgroup of $\GL_2(\F_3)$) or to a subgroup of $C_{\mathrm{ns}}(3)$ (a non-split Cartan).
		\item If $p = 5$ and $[F(\zeta_5):F] = 4$, then $\rhobar(G_F)$ is conjugate to a subgroup of $C_{\mathrm{ns}}^+(5)$ (the normalizer of a non-split Cartan).
	\end{enumerate} 
\end{lem}
\begin{proof}
	Let $G = \rhobar(G_F)$. If $G$ acts absolutely reducibly on $\Fp^2$ then it is a subgroup of a non-split Cartan subgroup ($G$ acts semisimply, because it acts irreducibly). Suppose $G$ acts absolutely irreducibly. We show that $G$ is a subgroup of the normalizer of a Cartan subgroup. We have $G' = G\cap\SL_2(\F_p) = \rhobar(G_{F(\zeta_p)})$, and we apply \cite[Lemma 2.2]{flhs} to conclude the proof of the first part.
	
	Now assume that $G$ is absolutely irreducible and $\det(G) =\F_p^\times$. It follows from \cite[Lemma 2.2]{flhs} that, if $\GL_2^+(\F_p)$ is the subgroup of $\GL_2(\Fp)$ consisting of matrices with square determinant, then $G^+ = G\cap\GL_2^+(\Fp)$ is contained in a Cartan subgroup $C$ with $G$ contained in the normalizer of $C$. 
	
	For the second part, if $\zeta_3\in F$, then $G = G'$ is absolutely reducible and is therefore a subgroup of a non-split Cartan. We can now assume $G$ is absolutely irreducible and $\det(G) = \F_3^\times$. Suppose $G'$ is contained in a non-split Cartan $C_{\mathrm{ns}}$ and $G$ is contained in its normalizer $C^+_{\mathrm{ns}}$. After conjugation, we can assume that \[C_{\mathrm{ns}} = \{\begin{pmatrix}x & -y\\ y & x \end{pmatrix} : (x,y) \in \F_3^2 - \{(0,0)\}.\] Note that $C_{\mathrm{ns}}\cap\SL_2(\F_3) = \langle\left(\begin{smallmatrix}
	0 & -1\\ 1 & 0
	\end{smallmatrix}\right)\rangle$ is cyclic of order 4 and is also contained in the normalizer of the diagonal split Cartan $C_\mathrm{s}(3)$. If $G'$ is contained in the scalars, then $G$ is the pre-image of a single order $2$ element in $\PGL_2(\F_p)$ and is therefore absolutely reducible. So $G' = C_{\mathrm{ns}}\cap\SL_2(\F_3)$. Considering the possibilities for $G \subset C_{\mathrm{ns}}^+$ with $G \not\subseteq C_{\mathrm{ns}}$ which contain $G'$ with index $2$, we  can see that $G$ must also be contained in the normalizer $C^+_\mathrm{s}(3)$ of $C_\mathrm{s}(3)$. The irreducible proper subgroups of $C^+_\mathrm{s}(3)$ are also contained in $C_{\mathrm{ns}}(3)$.
	
	Finally, suppose $p = 5$ and that $G^+$ is contained in the diagonal split Cartan $C_s(5)$ with $G$ contained in $C^+_s(5)$. We assume that $G$ is irreducible (not necessarily absolutely irreducible). If $G^+$ contains an element with eigenvalues $1$ and $-1$, then we are in the situation of (the proof of) \cite[Proposition 4.1(b)]{flhs}, which shows that $G$ is a subgroup (of index 3) in the normalizer of a non-split Cartan. So we assume that $G^+$ does not contain such an element. It follows that $G^+$ is equal to the subgroup of scalar matrices in $\GL_2(\Fp)$, so $G$ is the pre-image of an order $2$ element in $\PGL_2(\F_p)$. We conclude that $G$ is absolutely reducible, hence a subgroup of a non-split Cartan.   
\end{proof}

At this point we need to introduce modular curves with special level structures at $3$ and $5$. We follow the notation of \cite[\S 2.2]{flhs}, so for a subgroup $H \subset \GL_2(\F_p)$ containing $-I$ and with $\det(H) = \Fp^\times$, we have a modular curve $X(H)/\Q$ equipped with its $j$-invariant map $X(H) \to X(1)$. If $H_1 \subset \GL_2(\F_{p_1})$ and $H_2 \subset \GL_2(\F_{p_2})$ with $p_1 \ne p_2$, $X(H_1,H_2)$ is the modular curve given by the normalization of the fibre product $X(H_1)\times_{X(1)} X(H_2)$. The cuspidal points of $X(H_1,H_2)$ are those lying over $\infty \in X(1)$.

If $F \subset \Qbar$ is a number field, the non-cuspidal $F$-rational points of $X(H_1,H_2)$ correspond to $\Qbar$-isomorphism classes of pairs $(E,[\eta])$, where $E$ is an elliptic curve over $F$ and $[\eta]$ is an $(H_1\times H_2)$-orbit of isomorphisms \[\eta: \prod_{i=1}^2 \F_{p_i}^2 \cong \prod_{i=1}^2 E[p_i](\Qbar) \] such that $\eta^{-1}\rbar_{E,p_i}(G_F)\eta$ is contained in $H_i$ for $i=1,2$.

We will be concerned with the subgroups $\rb p$, $\rs p$, $\rns p$ which are respectively the upper triangular Borel, $C_{\rs}^+(p)$ and $C_{\rns}^+(p)$ in $\GL_2(\Fp)$. We will also need the Cartan subgroups themselves: $\rs p^\circ := C_{\rs}(p)$ and $\rns p^\circ := C_{\rns}(p)$.

Before proving Theorem \ref{thm:optimised modularity IQF}, we give a corollary. Recall that $X_0(15)=X(\rb 3,\rb5)$ is an elliptic curve of rank zero over $\Q$. It is the curve with Cremona label 15A1 (see \cite[Lemma 5.6]{flhs}).
\begin{cor}\label{cor:modularity IQF improved}
	Let $F$ be an imaginary quadratic field such that $X_0(15)(F)$ is finite. Then every elliptic curve $E/F$ is modular.
\end{cor}
\begin{proof}
	By Theorem \ref{thm:optimised modularity IQF}, we only need to consider $E/F$ giving rise to an $F$-rational point $P$ of $X_0(15)= X(\rb 3,\rb5)$ or $X(\rs3,\rb5)$. We can assume that $E$ does not have CM (otherwise it would be modular), so if $(E',[\eta])$ also gives the point $P$ then $E'$ is isomorphic to $E$ or a quadratic twist of $E$. Now it suffices to show that each of the $j$-invariants coming from points of $X_0(15)(F)$ and $X(\rs3,\rb5)(F)$ is modular. The curve $X(\rs3,\rb5)$ is an elliptic curve, with Cremona label 15A3 (\cite[Lemma 5.7]{flhs}), and is isogenous to $X_0(15)$. We are assuming that $X_0(15)(F)$ is finite, so $X(\rs3,\rb5)(F)$ is also finite. It remains to check modularity for the non-cuspidal torsion points defined over imaginary quadratic fields.
	
	We have $X_0(15)(\Q) \cong \Z/2\Z\times\Z/4\Z$. A Legendre form for $X_0(15)/\Q$ is $y^2 = x(x+16)(x+25)$. It  follows from \cite[Theorem 1]{kwon-torsion} that the only quadratic fields $F$ with $X_0(15)(\Q) \subsetneq X_0(15)(F)^{\mathrm{tors}}$ are $F= \Q(\sqrt{-1})$ and $\Q(\sqrt{5})$. This information is also contained in the LMFDB. It remains to show that elliptic curves giving rise to the 8 points in $X_0(15)(\Q(\sqrt{-1}))\backslash X_0(15)(\Q)$ are modular. This reduces to checking modularity of a single elliptic curve (isogenous and conjugate curves give the other points), and its modularity can be verified using the Faltings--Serre method \cite{gdp-faltings-serre} which has been carried out as part of the LMFDB project \cite[\href{https://www.lmfdb.org/EllipticCurve/2.0.4.1/4050.1/c/3}{Elliptic curve 4050.1-c3 over number field $\Q(\sqrt{-1})$}]{lmfdb}. 
	
	We also have $X(\rs3,\rb5)(\Q) \cong \Z/2\Z\times\Z/4\Z$. A Legendre form for $X(\rs3,\rb5)/\Q$ is $y^2 = x(x+1)(x+16)$, and it follows from  \cite[Theorem 1]{kwon-torsion} that the only quadratic field $F$ with $X(s3,b5)(\Q) \subsetneq X(s3,b5)(F)^{\mathrm{tors}}$ is $F= \Q(\sqrt{5})$ (again, this information is contained in the LMFDB).
\end{proof}

Combining Theorem \ref{thm:modularity35} with Lemma \ref{lem:galois image}, to prove Theorem \ref{thm:optimised modularity IQF} we need to show modularity of the elliptic curves defined over imaginary quadratic fields giving rise to points of the following modular curves:

\begin{enumerate}
	\item $X(\rns3^\circ,\rb5)$
	\item $X(\rb3,\rns5)$
	\item $X(\rns3^\circ,\rns5)$
	\item $X(\rs3,\rns5)$
\end{enumerate}

We will do this in the following subsections. We used Magma to do the computations \cite{magma}, and there are associated Magma files available at \url{https://github.com/jjmnewton/modularity-iqf}. 

To help us find equations for these curves, we use the following well-known facts about `small' modular curves:

\begin{prop}\label{prop:small modular curves}
	\begin{enumerate}
		\item The curve $X(\rb3)$ is isomorphic to $\PP^1_\Q$ with co-ordinate $x$ and $j$-invariant $\frac{(x+27)(x+3)^3}{x}$. 
		\item The curve $X(\rb5)$ is isomorphic to $\PP^1_\Q$ with co-ordinate $x$ and $j$-invariant $\frac{(x^2+250x+5^5)^3}{x^5}$. The Fricke involution $w_5$ is given by $x\mapsto 5^3/x$.
		\item The curve $X(\rns3)$ is isomorphic to $\PP^1_\Q$ with co-ordinate $x$ and $j$-invariant $x^3$.
		\item The curve $X(\rns5)$ is isomorphic to $\PP^1_\Q$ with co-ordinate $x$ and $j$-invariant $\frac{5^3x(2x+1)^3(2x^2+7x+8)^3}{(x^2+x-1)^5}$.
		\item The curve $X(\rns3,\rns5)$ is an elliptic curve over $\Q$, isomorphic to the curve with Cremona label $225A1$, with Mordell--Weil group $X(\rns3,\rns5)(\Q) \cong \Z$.
		\item The curve $X(\rs3)$ is isomorphic to $\PP^1_\Q$ with co-ordinate $x$ and $j$-invariant $\frac{3^3(x+1)^3(x-3)^3}{x^3} $		
	\end{enumerate}
\end{prop}
\begin{proof}
For the first and last parts, see \cite[Table 1]{sutherland-zywina}.
For the second part, see \cite[Table 3]{mcmurdy}. For the third and fourth parts, see \cite[Proposition 5.1, Corollary 6.3]{chen-siegel}. For the fifth part, it follows from the third and fourth parts that a singular model for $X(\rns3,\rns5)$ is given by the equation $x^3 =   \frac{5^3y(2y+1)^3(2y^2+7y+8)^3}{(y^2+y-1)^5}$. Some simple manipulations show that this is birational to the elliptic curve $y^2 - y = x^3 + 1$, which has Cremona label $225A1$ as claimed. This is also checked in the Magma file \texttt{ns3ns5-elliptic.m}.
\end{proof}

\subsection{Quadratic points on $X(\rns3^\circ,\rb5)$}
We begin by specifying an equation for the modular curve $X(\rns3^\circ,\rb5)$, together with the properties we will need to show modularity of quadratic points.
\begin{prop}\label{prop:ns3ob5-model}
\begin{enumerate}
\item A model for $X(\rns3^\circ,\rb5)$ is given by the genus one curve \[C: y^2 = -3(x^4 +2x^3 -x^2 +10x +25),\] where (the homogenization of) this equation defines a smooth curve in the weighted projective space $\mathbb{P}(1,2,1)$. 
\item With the above equation, $w_5$ transforms $x$-co-ordinates by $x \mapsto 5/x$.
\item The Jacobian $\mathrm{Jac}_C$ is isomorphic to the elliptic curve with Cremona label $45A2$, and in particular has $\mathrm{Jac}_C(\Q) \cong \Z/2\Z \oplus \Z/2\Z$.
\end{enumerate}
\end{prop}
\begin{proof}
For the first part, we proceed in a similar way to the proof of \cite[Lemma 4.4]{zywina-possible}. It follows from Proposition \ref{prop:small modular curves} that $X(\rns3,\rb5)$ is isomorphic to $\PP^1_{\Q}$, with $j$-invariant given by $j(x) = \frac{(x^6+250x^3+5^5)^3}{x^{15}}$ and $w_5$ given by $x \mapsto 5/x$. The map $X(\rns3^\circ,\rb5) \to X(\rns3,\rb5)$ has degree two, and is ramified at 4 points with $j$-invariant 1728. The fibre of the $j$-invariant map on $X(\rns3,\rb5)$ at 1728 is cut out by the vanishing of the polynomial
\begin{multline*}(x^6+250x^3+5^5)^3-1728x^{15} \\= (x^2 - 5x-25)^2(x^2-2x+5)(x^4+2x^3-x^2+10x+25)(x^4+5x^3+50x^2-125x+625)^2.\end{multline*} Since the ramification indices for $X(\rns3^\circ,\rb5)$ over $j=1728$ divide $2$, it follows that the branch locus in $X(\rns3,\rb5)$ consists of points which are unramified over the $j$-line. It follows that this branch locus is cut out by a divisor of $(x^2-2x+5)(x^4+2x^3-x^2+10x+25)$. Since it has degree $4$ and is Galois stable, it is cut out by $(x^4+2x^3-x^2+10x+25)$. 

We conclude that $X(\rns3^\circ,\rb5)$ is defined by an equation \[y^2 = c(x^4+2x^3-x^2+10x+25)\] for some squarefree $c \in \Z$. We claim that $c = -3$. We see from the equation (since the constant term $25$ is a square) that $\Q(\sqrt{c})$ is the field of definition of the cusps of $X(\rns3^\circ,\rb5)$, or equivalently the field of definition of the cusps of $X(\rns3^\circ)$. It is well-known that this field is $\Q(\sqrt{-3})$. Here is one justification: by ramification considerations as above, the double cover $X(\rns3^\circ) \to X(\rns3)$ is given by $y^2 = d(x^2+12x+144)$ for a squarefree $d\in \Z$. The unique rational point of $X(\rns3)$ with $j$-invariant 1728 is $x = 12$. The fibre over this point in $X(\rns3^\circ)$ is given by $y^2 = 432d$, so $\Q(\sqrt{3d})$ is determined by the image of $G_{\Q(\sqrt{3d})}$ in $\Aut(E[3])$ being contained in a non-split Cartan, where $E/\Q$ has $j$-invariant 1728 (hence CM by $\Z[\sqrt{-1}]$). It follows that $\Q(\sqrt{3d}) = \Q(\sqrt{-1})$ and we deduce that $d=-3$. 

Having computed a model for our genus one curve, the Jacobian can be identified using classical invariant theory and the Mordell--Weil group determined with the help of a two-descent (see \cite{g1-jacobian} for equations for the Jacobian, we did the calculations using Magma as documented in the file \texttt{ns3ob5.m}).  
\end{proof}

We will make use of the following points in $C(\Q(\sqrt{-3}))$:
\begin{align*}
\infty^+ &= (1:\sqrt{-3}:0)\\
\infty^- &= (1:-\sqrt{-3}:0)\\
0^+ &= (0:5\sqrt{-3}:1)\\
P_1 &= (-2:-\sqrt{-3}:1)\\
P_2 &= (-5/2:5\sqrt{-3}/4:1).
\end{align*}

We have $C(\Q) = \emptyset$, as can be checked $3$-adically. Using the isomorphism between $C_{\Q(\sqrt{-3})}$ and $\mathrm{Jac}_{C,\Q(\sqrt{-3})}$ given by $P \mapsto [P]-[\infty^-]$, we can match up the order $2$ rational points of $\mathrm{Jac}_C$ with the equivalence classes of the divisors:

\begin{align*}
D_0 &= 0^+ - \infty^-\\
D_1 &= P_1 - \infty^-\\
D_2 &= P_2 - \infty^-.
\end{align*}

At this point, we want to compute which points of $\mathrm{Jac}_C(\Q)$ are represented by rational divisors, not just rational equivalence classes. We denote the group of degree 0 rational divisors modulo linear equivalence by $\Pic^0(C)$.
\begin{prop}
$\Pic^0(C) \cong \Z/2\Z$, with generator $[0^+ - \infty^-]$.
\end{prop}
\begin{proof}
Our computations here are again documented in the file \texttt{ns3ob5.m}. There is an injective map $\Pic^0(C) \hookrightarrow \mathrm{Jac}_C(\Q)$, so we need to determine which of the divisors $D_i$ (for $i = 0,1,2$) are linearly equivalent to a rational divisor. For this, we essentially follow the procedure outlined in \cite{Flynn-Bruin}. Translating by the rational divisor $\infty^+ + \infty^-$, we consider the degree $2$ divisors $E_i = D_i + \infty^+ + \infty^-$. For each divisor, we compute a basis $\{1, f_i\}$ for each of the two-dimensional Riemann--Roch spaces $L(E_i)$ defined over $\Q(\sqrt{-3})$. We get:
\begin{align*}
f_0 &= \frac{y + (\sqrt{-3}x^2 + 5\sqrt{-3})}{x}\\
f_1 &= \frac{y + (\sqrt{-3}x^2 - 5\sqrt{-3})}{x+2}\\
f_2 &= \frac{y + (\sqrt{-3}x^2 - 5\sqrt{-3})}{x + 5/2}.
\end{align*}
A divisor class $[E_i]$ contains a rational divisor if and only if there is an effective degree $2$ rational divisor $D'_i = P + \sigma(P)$ linearly equivalent to $E_i$, with $P \in C(F)$ for a quadratic field $F$ and $\sigma$ the non-trivial element of $\Gal(F/\Q)$. If this is the case, then $D'_i$ is the divisor of zeroes of a meromorphic function $f_i - \alpha$ for some $\alpha \in \Q(\sqrt{-3})$. For $i = 0,1,2$ respectively this gives us equations:
\begin{align*}
y &= \alpha x - (\sqrt{-3}x^2 + 5\sqrt{-3})\\
y &= \alpha(x+2) - (\sqrt{-3}x^2 - 5\sqrt{-3})\\
y &= \alpha(x + 5/2) - (\sqrt{-3}x^2 - 5\sqrt{-3}).
\end{align*}
Squaring both sides, using the equation for $C$, and dividing out by the linear factor in $x$ (which doesn't correspond to a zero of $f_i - \alpha$) gives us equations:
\begin{equation}
(6-2\sqrt{-3}\alpha)x^2 + (\alpha^2 - 33)x + (30-10\sqrt{-3}\alpha)  = 0\label{eqn1}\end{equation}
\begin{equation}
(6-2\sqrt{-3}\alpha)x^2 + (\alpha^2 + 15)x + (2\alpha^2 + 10\sqrt{-3}\alpha) = 0
\label{eqn2}\end{equation}
\begin{equation}
(6-2\sqrt{-3}\alpha)x^2 + (\alpha^2 + 12)x + (5\alpha^2/2 + 10\sqrt{-3}\alpha) = 0.
\label{eqn3}\end{equation}

Equation (\ref{eqn1}) has the solution $x = 1 + 2i$, $\alpha = 3 + 2\sqrt{-3}$, and $E_0$ is linearly equivalent to $P + \sigma(P)$ for $P = (1+2i:3+6i:1)$. Note that, since (\ref{eqn1}) rescales to $x^2 + \cdots + 5 = 0$, we have $x(P)x(\sigma(P)) = 5$ for any $P$ defined over a quadratic field with $E_0 \sim P + \sigma(P)$. 

For the zeroes of Equation (\ref{eqn2}) to come from a rational divisor $P+\sigma(P)$, we must have $\frac{\alpha^2+15}{6-2\sqrt{-3}\alpha} \in \Q$. Writing $\alpha = \alpha_0 + \alpha_1\sqrt{-3}$, with $\alpha_i \in \Q$, we deduce that $\alpha_0 = 0$ or $\alpha_0^2 + 3\alpha_1^2  + 6\alpha_1 + 15 = 0$. The case $\alpha_0 =0$ can be excluded by considering the equation for $y$ in terms of $x$. Indeed, this gives us an equation $y = \sqrt{-3}\cdot g(x)$ for $g(x) \in \Q[x]$ which would need to be satisfied by $P$ and $\sigma(P)$. If $x(P) \notin \Q(\sqrt{-3})$, this means $P$ is not defined over a quadratic field (points with $y=0$ are only defined over a quartic field). If $x(P) \in \Q(\sqrt{-3})$, then we have $y(\sigma(P)) = -\sqrt{-3}\cdot g(x(\sigma(P)))$, which again forces $y(P) = 0$.

The equation $\alpha_0^2 + 3\alpha_1^2  + 6\alpha_1 + 15 = 0$ has no rational solutions (check $3$-adically). The associated conic is the Brauer--Severi variety associated to the divisor class $[E_1]$, as described in \cite{Flynn-Bruin}. Similarly, from Equation (\ref{eqn3}) we get the conic equation  $\alpha_0^2 + 3\alpha_1^2 + 6\alpha_1 + 12 = 0$, which again has no rational points (check $3$-adically).
\end{proof}
\begin{remark}
The conic associated to the divisor $E_0$ has equation $\alpha_0^2+3\alpha_1^2+6\alpha_1-33 = 0$, which does have rational points (for example $(3,2)$, which corresponds to $\alpha = 3+2\sqrt{-3}$). 
\end{remark}
\begin{prop}\label{prop:ns3b5quadratic} Suppose $P \in C(F)$ for a quadratic field $F$. We assume $P$ is not one of the two points at infinity. Then either $x(P) \in \Q$ or $x(P)x(\sigma(P)) = 5$, where $\sigma$ is the non-trivial element of $\Gal(F/\Q)$. 
\end{prop}
\begin{proof}
We have $[P+\sigma(P) - \infty^+ - \infty^-] = 0$ or $[0^+ - \infty^-]$ in $\Pic^0(C)$, so we deduce that $[P+ \sigma(P)] = [\infty^+ + \infty^-]$ or $[0^+ + \infty^+]$. In the second case, the proof of the previous proposition shows that $x(P)x(\sigma(P)) = 5$. In the first case, we have $P+\sigma(P) = \infty^+ + \infty^- + \mathrm{div}(f)$, where $f \in L(\infty^+ + \infty^-)$ is a $\Q$-linear combination of $1$ and $x$. In particular, $P$ is a zero of $x -\alpha$ for some $\alpha \in \Q$, so $x(P) \in \Q$.  
\end{proof}
\begin{cor}
Let $F$ be a quadratic field and $E/F$ an elliptic curve such that $E$ is the elliptic curve underlying an $F$-point $P \in X(\mathrm{ns3}^\circ,\mathrm{b5})(F)$. Then $E$ is modular.
\end{cor}
\begin{proof}
According to Proposition \ref{prop:ns3b5quadratic} we have two cases. In the first case, $P$ maps to a point in  $X(\mathrm{ns3},\mathrm{b5})(\Q)$. In particular, $E$ has rational $j$-invariant and is therefore modular. 

In the second case, $P$ maps to a $\Q$-point of the Atkin--Lehner quotient $X(\mathrm{ns3},\mathrm{b5})/w_5$. This shows that $E$ and $E^\sigma$ are $5$-isogenous over $\Qbar$, where $\sigma$ is the non-trivial element of $\Gal(F/\Q)$. In particular, $E$ is a $\Q$-curve, and is therefore modular (cf.~\cite[\S5.1]{flhs}). 
\end{proof}

\subsection{Quadratic points on $X(\rb3,\rns5)$}
For this curve, our approach is similar to that of \cite[Lemma 5.1]{flhs}. Our computations are documented in the file \texttt{b3ns5.m}. 
\begin{prop}
\begin{enumerate}
	\item $X(\rb3,\rns5)$ is a curve of genus 2, with hyperelliptic equation \[C: y^2 = 9x^6-6x^5-35x^4+40x^2+12x-8.\]
	\item The hyperelliptic involution on $X(\rb3,\rns5)$ is equal to the Atkin--Lehner involution $w_3$.
	\item The Mordell--Weil group $\mathrm{Jac}_C(\Q)$ of the Jacobian of $C$ is isomorphic to $\Z/2\Z \oplus \Z/10\Z$. 
\end{enumerate}
\end{prop}
\begin{proof}
We used our description for $X(\rb3)$ and $X(\rns5)$ to write down (an affine patch of) the fibre product $X(\rb3)\times_{X(1)}X(\rns5)$ and then used the Magma routines \verb|IsHyperelliptic| and \verb|SimplifiedModel| to find the displayed equation for $X(\rb3,\rns5)$, together with a birational map from  $X(\rb3)\times_{X(1)}X(\rns5)$ to $C$. 
	
To show that $w_3$ is the hyperelliptic involution, it suffices to show that it has (at least) 6 fixed points. Alternatively, we can compute the automorphism group of $X(\rb3,\rns5)$ \cite{magma-quartic} and note that it has order two, so the only non-trivial automorphism is the hyperelliptic involution. The fixed points of $w_3$ come from the elliptic curves $E_0 = \C/\cO_{\Q(\sqrt{-3})}$ and $E_1 = \C/\Z[\sqrt{-3}]$ with $\rb3$ level structure coming from the kernel of the multiplication by $\sqrt{-3}$ map; using the moduli description of \cite{kohen-pacetti} we can check that there are indeed 6 fixed points for $w_3$, 2 from $E_0$ and 4 from $E_1$.

We computed the Mordell--Weil group using the Magma routine \verb|MordellWeilGroupGenus2|. For our curve, a two-descent suffices to prove that the Jacobian has rank 0. We understand that this routine was created by Michael Stoll, see also \cite{stoll-2descent}. To directly verify that the torsion subgroup has size at most 20, we can compute that $\mathrm{Jac}_C(\FF_7) \cong \Z/2\Z\oplus \Z/20\Z$ and $\mathrm{Jac}_C(\FF_{13}) \cong \Z/2\Z\oplus \Z/90\Z$.
\end{proof}

We write $D_\infty$ for the degree two (hyperelliptic) divisor at $\infty$ for $C$.
\begin{lem}
Let $P$ be a quadratic point of $C$, with conjugate $\sigma(P)$. Then we are in one of the two following cases:
\begin{enumerate}
	\item $[P+\sigma(P)] = [D_\infty]$ in $\Pic(C)$. Then $P=(x,\pm\sqrt{f(x)})$ where $f(x)$ is the polynomial defining $C$, with $x \in \Q$. 
	\item $[P+\sigma(P)-D_\infty]$ is non-zero in $\Pic(C)$. The divisor $P+\sigma(P)$ is the unique degree 2 effective divisor representing its divisor class.
\end{enumerate} 
\end{lem}
\begin{proof}
	Note that $[D_\infty]$ is the canonical divisor class. In the first case we have $L(D_\infty) = \langle 1, x\rangle$ so $P$ and $\sigma(P)$ are zeroes of $x-\alpha$ for $\alpha \in \Q$. In the second case, by Riemann--Roch we have $\dim L(P+\sigma(P)) = 1$. 
\end{proof}

\begin{prop}\label{prop:pointsonXb3ns5} Suppose $P \in C(F)$ for an imaginary quadratic field $F$. We assume $P$ is not one of the two points at infinity. Then $x(P) \in \Q$ or $F = \Q(\sqrt{-11})$ and, up to complex conjugation, $P = (\frac{-5+\sqrt{-11}}{6},\pm\frac{17-\sqrt{-11}}{6})$.
\end{prop}
\begin{proof}
	By the previous lemma, it suffices to determine the non-zero elements  $A\in\mathrm{Jac}_C(\Q)$ with $A+[D_\infty]$ represented by an effective divisor with points in its support defined over an imaginary quadratic field. Magma returns elements of $\mathrm{Jac}_C(\Q)$ using the Mumford representation, which in particular gives a minimal polynomial for the $x$ co-ordinate. So it is easy to see which points may be imaginary quadratic.  There are 9 non-zero divisor classes supported on points with rational $x$ co-ordinate, 2 supported at $\infty$, 2 supported on the imaginary quadratic points specified in the statement of this proposition (and their conjugates); the remaining 6 are supported on real quadratic points.
\end{proof}

\begin{cor}
Let $F$ be an imaginary quadratic field and $E/F$ an elliptic curve such that $E$ is the elliptic curve underlying an $F$-point $P \in X(\rb3,\rns5)(F)$. Then $E$ is modular.
\end{cor}
\begin{proof} If $P\in C(\Q)$, then $j(E)\in \mathbb{Q}$, so $E$ is modular. In particular, this applies to the case when $P$ is one of the two points at infinity, which are defined over $\mathbb{Q}$. 
The points $P \in C(F)\setminus C(\Q)$ with $x(P) \in \Q$ have $\sigma(P)= w_3(P)$. So in this case, $E$ is a $\Q$-curve, hence modular. We are now only concerned with the $\Q(\sqrt{-11})$ points identified in Proposition \ref{prop:pointsonXb3ns5}. They are related by conjugation and $w_3$, so it suffices to show modularity of one of them. We consider the elliptic curve $E_0/\Q(\sqrt{-11})$ described by \cite[\href{https://www.lmfdb.org/EllipticCurve/2.0.11.1/8100.3/a/2}{Elliptic curve 8100.3-a2 over number field $\Q(\sqrt{-11})$}]{lmfdb}. On the one hand, it is modular, because the data in loc.~cit.~shows that the image of $G_{\Q(\sqrt{-11})}$ in $\Aut(E_0[5])$ is the full normalizer of a non-split Cartan, and therefore $E_0$  satisfies the hypotheses of Corollary \ref{cor:modularity IQF} (modularity was also checked explicitly by the LMFDB using the Faltings--Serre method). On the other hand, we can write down a $\Q(\sqrt{-11})$-rational point of $X(\rb3)\times_{X(1)}X(\rns5)$ with the same $j$-invariant as $E_0$, map it to $C$, and verify that we obtain one of our points of interest. We deduce that these points are also modular. 
\end{proof}

\subsection{Quadratic points on $X(\rs3,\rns5)$ and $X(\rns3^\circ,\rns5)$}

Both curves $X(\rs3,\rns5)$ and $X(\rns3^\circ,\rns5)$ are bi-elliptic, admitting degree two maps to $X(\rns3,\rns5)$. The quadratic points with image a rational point in $X(\rns3,\rns5)$ are modular, so we need to understand quadratic points which are not pulled back from rational points. To do this, we use Siksek's relative symmetric power Chabauty method \cite{siksek-chabauty}. We can closely follow the implementation of this method by Box for some modular curves $X_0(N)$ \cite{box-quadratic}. Our computations here are documented in the Magma files \verb|ns3ons5.m| and \verb|s3ns5.m|. 

We outline Box's method, following \cite[\S2.4]{box-quadratic}. We consider a smooth geometrically irreducible projective curve $X/\Q$ with Jacobian $J$, equipped with a degree two map $\pi:X \to C$ to another smooth curve $C/\Q$. Our goal is to describe the rational points of the symmetric square $X^{(2)}$, which include pairs of conjugate quadratic points $(P,\overline{P})$ of $X$. 

We will use the following input:
\begin{enumerate}
	\item Primes $p_1, \ldots, p_r$ of good reduction for $C$ and $X$. 
	\item Divisors $D_1, \ldots, D_n$ generating a subgroup $G$ of $J(\Q)$ of finite index. We have an associated surjective homomorphism $\phi: \Z^n \to G$.
	\item A positive integer $I$ such that $I\cdot J(\Q)\subset G$.
	\item A finite non-empty set $\cL \subset X^{(2)}(\Q)$, with a fixed element $\infty \in \cL$. We have a partition $\cL= \cL_{pb}\coprod\cL_{npb}$, with $x \in \cL_{pb}$ if and only if it is pulled back from a point in $C(\Q)$.
\end{enumerate} 

For each $i = 1, \ldots, r$ we have a subset $\cL_i^{\mathrm{good}} \subset \cL$ of known points where the (relative) Chabauty $p_i$-adic criterion applies. So $x \in \cL_i^{\mathrm{good}}$ if \cite[Theorem 2.1]{box-quadratic} applies to $x$ (which entails that $x$ is the unique point of $X^{(2)}(\Q)$ in its residue class mod $p_i$) or $x = \pi^*(P)$ for $P \in C(\Q)$ and \cite[Theorem 2.4]{box-quadratic} applies to $x$ (which entails that every point in the residue class of $x$ mod $p_i$ is pulled back from a point in $C(\Q)$). The following proposition is immediate from the definition.

\begin{prop}\label{prop:box-chabauty}
	Suppose $x \in X^{(2)}(\Q)$ and $\mathrm{red}_{p_i}(x) \in \mathrm{red}_{p_i}(\cL_i^{\mathrm{good}})$ for at least one $i = 1, \ldots, r$. Then $x \in \cL \cup \pi^*C(\Q)$.
\end{prop}

In practice, we will use a version of this proposition which lends itself to explicit computation. Consider the commutative diagram
\[
\xymatrix{ X^{(2)}(\Q) \ar[r]^{\iota}\ar[d]_{\mathrm{red}_p} & G\ar[d]_{\mathrm{red}_p} \\ 
X^{(2)}(\Fp) \ar[r]^{\iota_p} &J(\Fp)}
\]
where $\iota(x) = I([x]-\infty)$ and $\iota_p$ is the same map on the reduction mod $p$. For $i=1,\ldots,r$, let $\cM^{\mathrm{bad}}_i$ be the subset \[\cM^{\mathrm{bad}}_i:=\iota_{p_i}^{-1}(\mathrm{red}_{p_i}(G)) - \mathrm{red}_{p_i}(\cL_i^{\mathrm{good}}) \subset X^{(2)}(\F_{p_i}).\] It follows from Proposition \ref{prop:box-chabauty} that if $x$ is a rational point in $X^{(2)}(\Q)$ which is not in $\cL \cup \pi^*C(\Q)$, then $\mathrm{red}_{p_i}(x) \in \cM^{\mathrm{bad}}_i$ for each $i= 1 \ldots r$.

\begin{thm}{\cite[Theorem 2.6]{box-quadratic}}
If the set \[\bigcap_{i=1}^r \mathrm{red}_{p_i}^{-1}(\iota_{p_i}(\cM_i^{\mathrm{bad}})) \subset G\] is empty, then $X^{(2)}(\Q) = \cL \cup\pi^*C(\Q)$.
\end{thm}
\begin{proof}
	Suppose $x \in X^{(2)}(\Q) - (\cL \cup\pi^*C(\Q))$. We have already observed that $\mathrm{red}_{p_i}(x) \in \cM^{\mathrm{bad}}_i$ for each $i$. Then $\iota(x) \in G$  and we have $\mathrm{red}_{p_i}(\iota(x)) \in \iota_{p_i}(\cM_i^{\mathrm{bad}})$ for each $i$.
\end{proof}

For each $i$, $\mathrm{red}_{p_i}^{-1}(\iota_{p_i}(\cM_i^{\mathrm{bad}}))$ will give us a union of cosets of $\ker(\mathrm{red}_{p_i})$ in $G$ (possibly an empty union, if $\cM_i^{\mathrm{bad}} = \emptyset$). 

To get started applying this method in practice, we need equations for the curves and a formula for the bi-elliptic involution with quotient $X = X(\rns3,\rns5)$. It turns out that both curves have automorphism groups of order two, so the non-trivial automorphism must be the bi-elliptic involution.

\begin{prop}\label{prop:facts plane quartics}
	\begin{enumerate}
		\item $X(\rns3^\circ,\rns5)$ is isomorphic to the plane quartic $C_1$ with equation \[9x^4+19x^2y^2+y^4+9x^3+19x^2y+22xy^2+2y^3+10x^2+22xy+13y^2+7x+12y+11 = 0.\] It has a unique automorphism $w_1$ of order $2$ defined over $\Q$, given by $(x,y) \mapsto (x,-y-1)$.
		
		\item $X(\rs3,\rns5)$ is isomorphic to the plane quartic $C_2$ with equation \[-x^4+2x^3y+x^2y^2+8x^3+2x^2y-2xy^2-y^3-3x^2-3xy+3y^2+2x-3y+1 = 0.\] It has a unique automorphism $w_2$ of order $2$ defined over $\Q$, given by $(x,y) \mapsto (\frac{3x+y+2}{4x-2y+1},\frac{8x+y-8}{4x-2y+1})$.
		
		\item The $\Q(\sqrt{-55})$ points \[P_1 = (\frac{1+\sqrt{-55}}{28},\frac{27-\sqrt{-55}}{56}), P_2 = (\frac{3-\sqrt{-55}}{4},\frac{3+3\sqrt{-55}}{4})\] of $C_2$ have $j$-invariant $-32768$.
	\end{enumerate}
\end{prop}
\begin{proof}
	We determined a model $y^2=-3(x^2+12x+144)$ for the conic $X(\rns3^\circ)$ in Proposition \ref{prop:ns3ob5-model}. Together with our model for $X(\rns5)$ this gives us a singular model for $X(\rns3^\circ,\rns5)$. Using Magma, we determined a model for the desingularization and checked it was isomorphic to the plane quartic $C_1$ with the given equation. The automorphisms of $C_1$ can be computed using the function  \verb|AutomorphismGroupOfPlaneQuartic| written by Lercier, Sijsling and Ritzenthaler \cite{magma-quartic}. The same procedure was used for $X(\rs3,\rns5)$. We also wrote down the points in the singular model (where the $j$-invariant is easy to compute) which map to $P_1$ and $P_2$.
\end{proof}

We write $\pi_i: C_i \to X := X(\rns3,\rns5)$ for the bi-elliptic quotient map.
	\begin{prop}\begin{enumerate}
\item We have $\mathrm{rk}(\mathrm{Jac}_{C_i}(\Q)) = \mathrm{rk}(\mathrm{Jac}_{X}(\Q)) = 1$ for $i=1$ and $2$.\label{item:samerank}
\item The torsion subgroups $\mathrm{Jac}_{C_i}(\Q)^{\mathrm{tors}}$ satisfy
\begin{itemize}
\item $\mathrm{Jac}_{C_1}(\Q)^{\mathrm{tors}}$ is isomorphic to a subgroup of $\Z/2\Z \oplus \Z/2\Z$
\item $\mathrm{Jac}_{C_2}(\Q)^{\mathrm{tors}} \cong \Z/2\Z \oplus \Z/10\Z$ (explicit generators will be identified in the proof).
\end{itemize}
\item Pick a generator $D \in \mathrm{Jac}_{X}(\Q) \cong \ZZ$. Set $G_i = \langle \pi_i^* D \rangle \subset \mathrm{Jac}_{C_i}(\Q)$. Then $4(\mathrm{Jac}_{C_1}(\Q)) \subset 2G_1$ and $10(\mathrm{Jac}_{C_2}(\Q)) \subset \langle 5G_2,\mathrm{Jac}_{C_2}(\Q)[2] \rangle$.
\end{enumerate}
	\end{prop}
\begin{proof}
Using Chen's isogeny (in the form of \cite[Th\'{e}or\`{e}me 2]{deS-edix-chen}) and the fact that $X(\rns3^\circ)$ and $X(\rns5)$ have genus 0, we see that $\mathrm{Jac}_{C_1}$ is isogenous to the new part of the Jacobian of $X_0(225)/w_{25}$. The relevant space of cuspforms $S_2(\Gamma_0(225))^{w_{25}=1,\mathrm{new}}$ has dimension $3$. It has a basis of Hecke eigenforms with rational $q$-expansions, with LMFDB labels $225.2.a.b$, $225.2.a.c$ and $225.2.a.e$. The associated isogeny classes of elliptic curves are those with Cremona labels $225c$, $225a$ and $225d$ respectively. They have Mordell--Weil rank $0$, $1$ and $0$ respectively, so we deduce that $\mathrm{rk}(\mathrm{Jac}_{C_1}(\Q)) = 1$.

Similarly, $\mathrm{Jac}_{C_2}$ is isogenous to the $5$-new part of the Jacobian of $X_0(225)/\langle w_9 ,w_{25}\rangle$. The relevant space of cuspforms is now  $S_2(\Gamma_0(225))^{w_9=1,w_{25}=1,5\mathrm{-new}}$. The newform $225.2.a.c$, with associated rank $1$ elliptic curve $225a$, still contributes (as it should, since $C_2$ maps to $X(\rns3,\rns5)$). We also get contributions from the $2$ rational eigenforms in $S_2(\Gamma_0(75))^{w_{25}=1,5\mathrm{-new}}$, labels $75.2.a.a$ and $75.2.a.b$ in the LMFDB ($w_9$ has characteristic polynomial $X^2-1$ on the oldspaces generated by each eigenform). The associated isogeny classes of elliptic curves (Cremona labels $75c$ and $75b$ respectively) both have rank $0$.

It now follows from the equality of ranks in part \eqref{item:samerank} that \[2(\mathrm{Jac}_{C_i}(\Q)) \subset \langle G_i, \mathrm{Jac}_{C_i}(\Q)^{\mathrm{tors}} \rangle\] for $i = 1,2$ (cf.~\cite[Proposition 3.1]{box-quadratic}, the running assumption there that the curves have a rational point is not used in the proof). We claim that
\begin{enumerate}
	\item $2(\mathrm{Jac}_{C_1}(\Q)^{\mathrm{tors}}) = 0$
	\item $10(\mathrm{Jac}_{C_2}(\Q)^{\mathrm{tors}}) = 0$.
\end{enumerate}
In the first case, the orders of $\mathrm{Jac}_{C_1}(\F_p)$ for $p = 7,11$ and $13$ have $\gcd$ $4$, so $\mathrm{Jac}_{C_1}(\Q)^{\mathrm{tors}}$ has order dividing $4$. Moreover, $\mathrm{Jac}_{C_1}(\F_{13}) \cong \Z/2\Z \oplus \Z/1710\Z$, so $\mathrm{Jac}_{C_1}(\Q)^{\mathrm{tors}}$ must be isomorphic to a subgroup of $\Z/2\Z \oplus \Z/2\Z$.

In the second case, the orders of $\mathrm{Jac}_{C_2}(\F_p)$ for $p = 7,11$ and $13$ have $\gcd$ $20$, so $\mathrm{Jac}_{C_2}(\Q)^{\mathrm{tors}}$ has order dividing $20$. On the other hand, we can write down a subgroup of $\mathrm{Jac}_{C_2}(\Q)^{\mathrm{tors}}$ isomorphic to $\Z/2\Z\oplus\Z/10\Z$, so we deduce that this is the full torsion subgroup. Generators for the torsion subgroup are given by $D_1$ of order 10 and $D_2$ of order 2, for \begin{align*}
D_1 &= [(0,1,1)]-[(-3,7,1)],\\
D_2 &= 5([(0,1,0)]+[(-1/2,-1/2,1)] + [P_2] + [\overline{P}_2] - 2[P_1] - 2[\overline{P}_1]),
\end{align*}
where $P_1 = (\frac{5+\sqrt{21}}{2},-\sqrt{21}-4)$, $P_2 = (\frac{-1-\sqrt{5}}{2},\frac{3\sqrt{5}+9}{2})$ and $\overline{P}_i$ denotes the Galois conjugate of $P_i$.

We deduce from this that $4(\mathrm{Jac}_{C_1}(\Q)) \subset  2G_1 $ and $10(\mathrm{Jac}_{C_2}(\Q)) \subset \langle 5G_2,\mathrm{Jac}_{C_2}(\Q)[2] \rangle$.
\end{proof}

This gives us all the inputs we need to apply the relative symmetric power Chabauty method. The description of vanishing differentials and their reduction goes through exactly as in \cite[\S3.4]{box-quadratic}, replacing the Atkin--Lehner involution of $X_0(N)$ which appears there with our bi-elliptic involutions. We used Box's code, available at \url{https://github.com/joshabox/quadraticpoints/}, to carry out the computations. This also includes code written by Ozman and Siksek to search for rational points in $X^{(2)}(\Q)$. 

\begin{prop}
\begin{enumerate}
	\item Let $X = X(\rns3^\circ,\rns5)$, $C = X(\rns3,\rns5)$ and let $\pi:X \to C$ be the natural quotient map. We have $X^{(2)}(\Q) = \pi^*C(\Q)$.
	\item Let $X = X(\rs3,\rns5)$, $C = X(\rns3,\rns5)$ and let $\pi:X \to C$ be the natural quotient map. There are eight conjugate pairs of quadratic points of $X$ which do not have image in $C(\Q)$. Two are defined over an imaginary quadratic field, corresponding to the points $P_1, P_2$ identified in Proposition \ref{prop:facts plane quartics}.
\end{enumerate}
\end{prop}
\begin{proof}
For the first part, it turns out that we just need to apply the relative Chabauty criterion of \cite[Theorem 2.4]{box-quadratic} for $p = 43$ and with $\cL$ consisting of 8 degree two divisors (all pulled back from $C(\Q)$). Digging into what is happening, it turns out that $\mathrm{red}_{43}(G)$ is cyclic of order $7$ and $\iota_{43}^{-1}(\mathrm{red}_{43}(G))$ also has size $7$. The divisors in $\cL$ cover all these possibilities, and they all satisfy the hypotheses of \cite[Theorem 2.4]{box-quadratic}.
	
For the second part, we run Box's Mordell--Weil sieve for the primes $11,43$, $G = \langle 5G_2,\mathrm{Jac}_{C_2}(\Q)[2] \rangle$ and $I = 10$. The set $\cL$ includes 16 rational degree two divisors which are not pulled back from $C(\Q)$, 8 of which are sums of two points in $X(\Q)$ (not interchanged by the bi-elliptic involution).
\end{proof}
This proposition has the immediate corollary:
\begin{cor}
	\begin{enumerate}
\item Let $F$ be a quadratic field and $E/F$ an elliptic curve such that $E$ is the elliptic curve underlying an $F$-point $P \in X(\rns3^\circ,\rns5)(F)$. Then $E$ is modular. 

\item Let $F$ be an imaginary quadratic field and $E/F$ an elliptic curve such that $E$ is the elliptic curve underlying an $F$-point $P \in X(\rs3,\rns5)(F)$. Then $E$ is modular. 
\end{enumerate}
\end{cor}
\begin{proof}
	We have shown that all relevant quadratic points map to rational points of $X(\rns3,\rns5)$, with the exception of the points $P_1, P_2 \in X(\rs3,\rns5)(\Q(\sqrt{-55}))$ and their conjugates. The latter points have rational $j$-invariant so we are done in all cases. 
\end{proof}

We have now completed the proof of Theorem \ref{thm:optimised modularity IQF}!

\newcommand{\etalchar}[1]{$^{#1}$}
\renewcommand{\MR}[1]{}
\providecommand{\bysame}{\leavevmode\hbox to3em{\hrulefill}\thinspace}
\providecommand{\MR}{\relax\ifhmode\unskip\space\fi MR }
\providecommand{\MRhref}[2]{%
	\href{http://www.ams.org/mathscinet-getitem?mr=#1}{#2}
}
\providecommand{\href}[2]{#2}

\end{document}